\documentclass{article} 
\usepackage{latexsym,amsmath,amscd,amssymb,graphics}
\usepackage[backend=bibtex,sorting=none]{biblatex}
\usepackage{enumerate}
\usepackage[]{mcode}
\usepackage[toc,page]{appendix}
\usepackage[a4paper,margin=.75in]{geometry}
\usepackage{array}
\usepackage[section]{placeins}
\usepackage{afterpage}
\usepackage{verbatim}

\usepackage{graphicx}
\usepackage{subfig}
\usepackage{framed}
\usepackage[colorlinks]{hyperref}
\usepackage{url}
\usepackage{mathrsfs}
\usepackage{epstopdf}

\usepackage[all]{xy}
\usepackage{indentfirst}
\usepackage[mathscr]{euscript}
\usepackage{algorithm}
\usepackage{algpseudocode}
\usepackage{authblk}

\makeatletter
\renewcommand*{\@fnsymbol}[1]{\ensuremath{\ifcase#1\or \dagger\or *\or \ddagger\or
		\mathsection\or \mathparagraph\or \|\or **\or \dagger\dagger
		\or \ddagger\ddagger \else\@ctrerr\fi}}
\makeatother

\makeatletter

\@addtoreset{figure}{section}
\def\thefigure{\thesection.\@arabic\c@figure}
\def\fps@figure{h, t}
\@addtoreset{table}{bsection}
\def\thetable{\thesection.\@arabic\c@table}
\def\fps@table{h, t}
\@addtoreset{equation}{section}

\makeatother


\makeatletter
\AtBeginDocument{%
	\expandafter\renewcommand\expandafter\subsection\expandafter{%
		\expandafter\@fb@secFB\subsection
	}%
}
\makeatother

\newcommand{\unaryminus}{\scalebox{0.75}[1.0]{\( - \)}}

\newtheorem{theorem}{Theorem}

\newtheorem{definition}[theorem]{Definition}

\numberwithin{theorem}{subsection}

\def\bea{\begin{eqnarray}}
\def\eea{\end{eqnarray}}
\def\ba{\begin{array}}
\def\ea{\end{array}}

\def\bOm{\boldsymbol{\Omega}}

\def\diag{\mathrm{ \textbf{diag}}}

\def\bx{{\boldsymbol {x} }}

\let\<\langle
\let \>\rangle

\newcommand{\sgn}{\mathop{\mathrm{sgn}}}

\newcommand{\rem}[1]{}

\newcommand{\de}{\delta}

\newcommand{\bu}{\boldsymbol{u}}

\newcommand{\bv}{\boldsymbol{v}}
\newcommand{\bz}{\boldsymbol{z}}

\newcommand{\bGam}{\boldsymbol{\Gamma}}
\newcommand{\bGamma}{\boldsymbol{\Gamma}}
\newcommand{\bzeta}{\boldsymbol{\zeta}}

\newcommand{\bsigma}{\boldsymbol{\sigma}}

\newcommand{\bpsi}{\boldsymbol{\psi}}

\newcommand{\bxi}{\boldsymbol{\xi}}

\newcommand{\bY}{\mathbf{Y}}

\newcommand{\bPi}{\boldsymbol{\Pi}}

\newcommand{\bPhi}{\boldsymbol{\Phi}}
\newcommand{\bUpsilon}{\boldsymbol{\Upsilon}}

\newcommand{\bkappa}{\boldsymbol{\kappa}}
\newcommand{\bpi}{\boldsymbol{\pi}}
\newcommand{\blam}{\boldsymbol{\lambda}}
\newcommand{\bnu}{\boldsymbol{\nu}}

\newcommand{\btheta}{\boldsymbol{\theta}}
\newcommand{\inertia}{\mathbb{I}}

\newcommand{\pp}[2]{\frac{\partial #1}{\partial #2}}
\newcommand{\dd}[2]{\frac{\mathrm{d} #1}{\mathrm{d} #2}}

\newcommand{\dt}{\mathrm{d}t}
\newcommand{\mrmd}{\mathrm{d}}
\newcommand{\ACP}{\lambda}
\newcommand{\fricudirbase}{\sigma}
\newcommand{\fricudir}{\boldsymbol{\fricudirbase}}

\newcommand{\dprime}{\prime \prime}

\DeclareMathOperator{\ReLU}{ReLU}

\newcommand{\revision}[2]{#2} 

\newcommand{\revisionNACO}[2]{#2}

\usepackage[normalem]{ulem}

\graphicspath{ {./gen/} {./dsim8/} {./dsim8_tp/} {./bsim10/} {./bsim12/}  }

\title{Numerical Simulations of a Rolling Ball Robot Actuated by Internal Point Masses}
\author[1,2]{Vakhtang Putkaradze\thanks{Email address: \texttt{putkarad@ualberta.ca}}}
\author[3]{Stuart Rogers\thanks{Email address: \texttt{srogers@umn.edu}}}
\affil[1]{Department of Mathematical and Statistical Sciences, Faculty of Science, University of Alberta, CAB 632, Edmonton, AB T6G 2G1, Canada}
\affil[2]{ATCO SpaceLab, 5302 Forand ST SW, Calgary, AB T3E 8B4, Canada}
\affil[3]{Institute for Mathematics and its Applications, College of Science and Engineering, University of Minnesota, 207 Church ST SE, 306 Lind Hall, Minneapolis, MN 55455, USA}
\date{\today}

\providecommand{\keywords}[1]{\textbf{\textit{Keywords:}} #1}
\begin{document}

\maketitle

\abstract{\noindent 
The controlled motion of a rolling ball actuated by internal point masses that move along arbitrarily-shaped rails fixed within the ball is considered. The controlled equations of motion are solved numerically using a predictor-corrector continuation method, starting from an initial solution obtained via a direct method, to realize trajectory tracking and obstacle avoidance maneuvers.} 
\\
\\
\tableofcontents 

\section{Introduction} \label{sec_introduction}

\subsection{Overview}
This paper is a continuation of \cite{putkaradze2020optimal}, providing numerical solutions of the controlled equations of motion for several special cases of the rolling disk and ball actuated by moving internal point masses. The paper \cite{putkaradze2020optimal} invokes Pontryagin's minimum principle to derive the theoretical background for the optimal control of the rolling disk and ball having general performance indexes. This paper implements the theory derived in \cite{putkaradze2020optimal} to solve several practical examples, such as trajectory tracking for the rolling disk and obstacle avoidance for the rolling ball. The key contributions of this paper are listed below.
\begin{itemize}
	\item The controlled equations of motion, for a rolling ball actuated by internal point masses that move along arbitrarily-shaped rails fixed within the ball, are solved numerically by a predictor-corrector continuation method, starting from an initial solution provided by a direct method.
	\item Jacobians of the ordinary differential equations (ODEs) and boundary conditions (BCs) which constitute the controlled equations of motion (i.e. an ordinary differential equation two-point boundary value problem (ODE TPBVP)) corresponding to the optimal control of a dynamical system are derived. These Jacobians are useful for numerically solving the controlled equations of motion.
	\item Algorithms for solving an ODE TPBVP by predictor-corrector continuation are developed and were implemented in \mcode{MATLAB} to numerically solve the controlled equations of motion for the rolling ball. There are not very many predictor-corrector continuation methods publicly available for solving dynamical systems. The idea of using a monotonic continuation ODE TPBVP solver in conjunction with a predictor-corrector continuation method to advance (or ``sweep'') as far along the tangent as possible is new,  and this novel technique was used to obtain all the numerical results in this paper.
\end{itemize}

The paper is organized as follows.  Subsections~\ref{ssec_ball_uncontrolled} and \ref{ssec_disk_uncontrolled} review the specific types of rolling disk and ball considered, define coordinate systems and notation used to describe this rolling disk and ball, and present the uncontrolled equations of motion for this rolling disk and ball derived earlier in \cite{Putkaradze2018dynamicsP,putkaradze2018normalpub}. In Sections~\ref{sec_disk_sim} and \ref{sec_ball_sim}, the controlled equations of motion for the rolling disk and ball are formulated and solved numerically via a predictor-corrector continuation method, starting from an initial solution provided by a direct method. Subsection~\ref{ssec_numerical_methods} provides details of the numerical methods used to solve the controlled equations of motion. Section~\ref{sec_conclusions} summarizes the results of the paper and discusses topics for future work.  The background material for this paper is contained in several appendices. In particular, Appendix~\ref{sec_optimal_control} reviews the theory of optimal control needed to derive the controlled equations of motion for a generic dynamical system given initial and final conditions, given a performance index to be minimized, and in the absence of path inequality constraints. Appendix~\ref{app_imp_details} provides details for numerically solving the controlled equations of motion. Appendices~\ref{app_predictor_corrector} and \ref{app_sweep_predictor_corrector} develop  two predictor-corrector continuation algorithms which numerically solve an ODE TPBVP, the latter of which is utilized to numerically solve the controlled equations of motion for the rolling disk and ball. 

\subsection{Rolling Ball} \label{ssec_ball_uncontrolled}
\begin{figure}[h]
	\centering
	\includegraphics[width=0.5\linewidth]{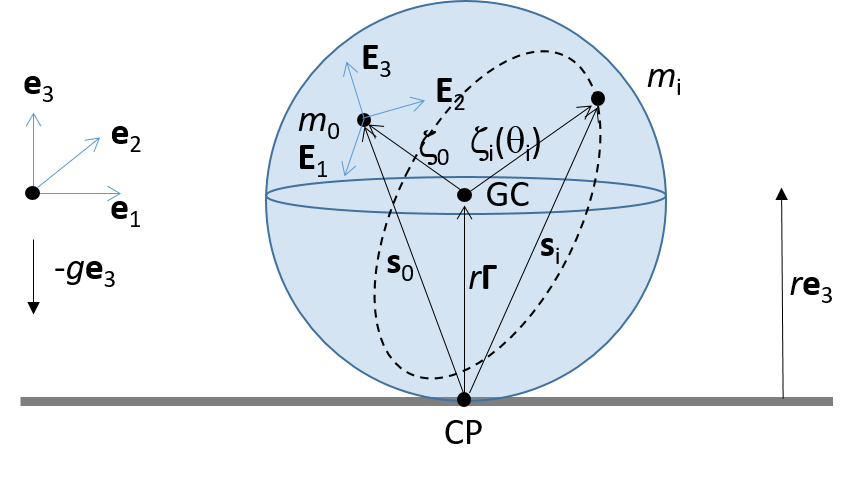}
	\caption{A ball of radius $r$ and mass $m_0$ rolls without slipping on a horizontal surface in the presence of a uniform gravitational field of magnitude $g$. The ball's geometric center, center of mass, and contact point with the horizontal surface are denoted by GC, $m_0$, and CP, respectively. The spatial frame has origin located at height $r$ above the horizontal surface and orthonormal axes $\mathbf{e}_1$, $\mathbf{e}_2$, and $\mathbf{e}_3$. The body frame has origin located at the ball's center of mass (denoted by $m_0$) and orthonormal axes $\mathbf{E}_1$, $\mathbf{E}_2$, and $\mathbf{E}_3$. The ball's motion is actuated by $n$ point masses, each of mass $m_i$, $1 \le i \le n$, and each moving along its own rail fixed inside the ball. The $i^\mathrm{th}$ rail is depicted here by the dashed hoop. The trajectory of the $i^\mathrm{th}$ rail, with respect to the body frame translated to the GC, is denoted by $\bzeta_i$ and is parameterized by $\theta_i$. All vectors inside the ball are expressed with respect to the body frame, while all vectors outside the ball are expressed with respect to the spatial frame. }
	\label{fig:detailed_1dparam_rolling_ball}
\end{figure}
Consider a rigid ball of radius $r$ containing some static internal structure as well as $n$ point masses, \revision{R1Q7}{where either $n$ is a positive integer denoting the number of moving masses or $n=0$ if no moving masses are used and the structure of the ball is static}. 
This ball rolls without slipping on a horizontal  surface in the presence of a uniform gravitational field of magnitude $g$, as illustrated in Figure~\ref{fig:detailed_1dparam_rolling_ball}. The ball with its static internal structure has mass $m_0$ and the $i^\mathrm{th}$ point mass has mass $m_i$ for $1 \le i \le n$. Let $M = \sum_{i=0}^n m_i$ denote the mass of the total system. The total mechanical system consisting of the ball with its static internal structure and the $n$ point masses is referred to as the ball or the rolling ball, the ball with its static internal structure but without the $n$  point masses may also be referred to as $m_0$, and the $i^\mathrm{th}$  point mass may also be referred to as $m_i$ for $1 \le i \le n$. Note that this system is the Chaplygin ball \cite{Ho2011_pII} equipped with point masses.

Two coordinate systems, or frames of reference, will be used to describe the motion of the rolling ball, an inertial spatial coordinate system and a body coordinate system in which each particle within the ball is always fixed. For brevity, the spatial coordinate system will be referred to as the spatial frame and the body coordinate system will be referred to as the body frame. These two frames are depicted in Figure~\ref{fig:detailed_1dparam_rolling_ball}. The spatial frame has orthonormal axes $\mathbf{e}_1$, $\mathbf{e}_2$, $\mathbf{e}_3$, such that the $\mathbf{e}_1$-$\mathbf{e}_2$ plane is parallel to the horizontal surface and passes through the ball's geometric center (i.e. the $\mathbf{e}_1$-$\mathbf{e}_2$ plane is a height $r$ above the horizontal surface), such that $\mathbf{e}_3$ is vertical (i.e. $\mathbf{e}_3$ is perpendicular to the horizontal surface) and points ``upward" and away from the horizontal surface, and such that $\left(\mathbf{e}_1, \mathbf{e}_2, \mathbf{e}_3 \right)$ forms a right-handed coordinate system. For simplicity, the spatial frame axes are chosen to be
\begin{equation}
\mathbf{e}_1 = \begin{bmatrix} 1 & 0 & 0 \end{bmatrix}^\mathsf{T}, \quad \mathbf{e}_2 = \begin{bmatrix} 0 & 1 & 0 \end{bmatrix}^\mathsf{T}, \quad \mathrm{and} \quad \mathbf{e}_3 = \begin{bmatrix} 0 & 0 & 1 \end{bmatrix}^\mathsf{T}.
\end{equation}
The acceleration due to gravity in the uniform gravitational field is $\mathfrak{g} = -g \mathbf{e}_3  = \begin{bmatrix} 0 & 0 & -g  \end{bmatrix}^\mathsf{T}$ in the spatial frame.

The body frame's origin is chosen to coincide with the position of $m_0$'s center of mass. The body frame has orthonormal axes $\mathbf{E}_1$, $\mathbf{E}_2$, and $\mathbf{E}_3$, chosen to coincide with $m_0$'s principal axes, in which $m_0$'s inertia tensor $\inertia$ is diagonal, with corresponding \revision{R1Q9}{principal} moments of inertia $d_1$, $d_2$, and $d_3$. That is, in this body frame the inertia tensor is the diagonal matrix $\inertia = \diag \left( \begin{bmatrix} d_1 & d_2 & d_3 \end{bmatrix} \right)$.
Moreover, $\mathbf{E}_1$, $\mathbf{E}_2$, and $\mathbf{E}_3$ are chosen so that $\left(\mathbf{E}_1, \mathbf{E}_2, \mathbf{E}_3 \right)$ forms a right-handed coordinate system. For simplicity, the body frame axes are chosen to be
\begin{equation}
\mathbf{E}_1 = \begin{bmatrix} 1 & 0 & 0 \end{bmatrix}^\mathsf{T}, \quad \mathbf{E}_2 = \begin{bmatrix} 0 & 1 & 0 \end{bmatrix}^\mathsf{T}, \quad \mathrm{and} \quad \mathbf{E}_3 = \begin{bmatrix} 0 & 0 & 1 \end{bmatrix}^\mathsf{T}.
\end{equation}
In the spatial frame, the body frame is the moving frame $\left(\Lambda \left(t\right) \mathbf{E}_1, \Lambda \left(t\right) \mathbf{E}_2, \Lambda \left(t\right) \mathbf{E}_3  \right)$, where $\Lambda \left(t\right) \in SO(3)$ defines the orientation (or attitude) of the ball at time $t$ relative to its reference configuration, for example at some initial time.   

For $1 \le i \le n$, it is assumed that $m_i$ moves along its own 1-d rail. It is further assumed that the $i^\mathrm{th}$  rail is parameterized by a 1-d parameter $\theta_i$, so that the \revision{R1Q8}{trajectory} $\bzeta_i$ of the $i^\mathrm{th}$  rail, in the body frame translated to the ball's geometric center, as a function of $\theta_i$ is $\bzeta_i(\theta_i)$. Refer to Figure~\ref{fig:detailed_1dparam_rolling_ball} for an illustration. Therefore, the body frame vector from the ball's geometric center to $m_i$'s center of mass is denoted by $\bzeta_i(\theta_i (t))$. Since $m_0$ is stationary in the body frame and to be consistent with the positional notation for $m_i$ for $1 \le i \le n$, $\bzeta_0 \equiv \bzeta_0(\theta_0) \equiv \bzeta_0(\theta_0 (t))$ is the constant \revision{R1Q10}{(time-independent)} vector from the ball's geometric center to $m_0$'s center of mass for any scalar-valued, time-varying function $\theta_0(t)$. In addition, suppose a time-varying external force $\mathbf{F}_\mathrm{e}(t)$ acts at the ball's geometric center.

Let $\mathbf{z}_i(t)$ denote the position of $m_i$'s center of mass in the spatial frame so that the position of $m_i$'s center of mass in the spatial frame is $\mathbf{z}_i(t)=\mathbf{z}_0(t)+\Lambda(t) \left[\bzeta_i(\theta(t))-\bzeta_0\right]$. In general, a particle with position $\mathbf{w}(t)$ in the body frame has position $\mathbf{z}(t) = \mathbf{z}_0(t)+\Lambda(t) \mathbf{w}(t)$ in the spatial frame and has position $\mathbf{w}(t)+\bzeta_0$ in the body frame translated to the ball's geometric center.

For conciseness, the ball's geometric center is often denoted GC, $m_0$'s center of mass is often denoted CM, and the ball's contact point with the surface is often denoted CP. The GC is located at $\mathbf{z}_\mathrm{GC}(t) = \mathbf{z}_0(t)-\Lambda(t) \bzeta_0$ in the spatial frame, at $-\bzeta_0$ in the body frame, and at $\mathbf{0} = \begin{bmatrix} 0 & 0 & 0 \end{bmatrix}^\mathsf{T}$ in the body frame translated to the GC. The CM is located at $\mathbf{z}_0(t)$ in the spatial frame, at $\mathbf{0}$ in the body frame, and at $\bzeta_0$ in the body frame translated to the GC. The CP is located at $\mathbf{z}_\mathrm{CP}(t) = \mathbf{z}_0(t)-\Lambda(t) \left[r\bGam(t)+\bzeta_0 \right]$ in the spatial frame, at $-\left[r\bGam(t)+\bzeta_0 \right]$ in the body frame, and at $-r\bGam(t)$ in the body frame translated to the GC, where $\bGamma(t) \equiv \Lambda^{-1}(t) \mathbf{e}_3$. Since the third spatial coordinate of the ball's GC is always $0$ and of the ball's CP is always $-r$, only the first two spatial coordinates of the ball's GC and CP, denoted by $\bz(t)$, are needed to determine the spatial location of the ball's GC and CP. 

For succintness, the explicit time dependence of variables is often dropped. That is, the orientation of the ball at time $t$ is denoted simply $\Lambda$ rather than $\Lambda(t)$, the position of $m_i$'s center of mass in the spatial frame at time $t$ is denoted $\mathbf{z}_i$ rather than $\mathbf{z}_i(t)$,  the position of $m_i$'s center of mass in the body frame translated to the GC at time $t$ is denoted $\bzeta_i$ or $\bzeta_i(\theta_i)$ rather than $\bzeta_i(\theta_i(t))$, the spatial $\mathbf{e}_1$-$\mathbf{e}_2$ position of the ball's GC at time $t$ is denoted $\bz$ rather than $\bz(t)$, and the external force is denoted $\mathbf{F}_\mathrm{e}$ rather than $\mathbf{F}_\mathrm{e}(t)$. 

As shown in \cite{Putkaradze2018dynamicsP,putkaradze2018normalpub}, the uncontrolled equations of motion for this rolling ball are
\begin{equation}
\begin{split} \label{uncon_ball_eqns_explicit_1d}
\dot \bOm &= \left[\sum_{i=0}^n m_i \widehat{\mathbf{s}_i}^2  -\inertia \right]^{-1}  \Bigg[\bOm \times \inertia \bOm+r \tilde \bGamma \times \bGamma 
\\ &\hphantom{=} + \sum_{i=0}^n m_i \mathbf{s}_i \times  \left\{ g \bGamma+ \bOm \times \left(\bOm \times \bzeta_i +2 \dot \theta_i \bzeta_i^{\prime} \right) + \dot \theta_i^2 \bzeta_i^{\dprime} + \ddot \theta_i  \bzeta_i^{\prime} \right\}  \Bigg], \\
\dot \Lambda &= \Lambda \widehat{\bOm}, \\
\dot \bz &= \left( \Lambda \bOm \times r \mathbf{e}_3 \right)_{12},
\end{split}
\end{equation}
where $\mathbf{s}_i \equiv r \bGamma +\bzeta_i$ \revisionNACO{R1Q3}{is the body frame vector from the CP to $m_i$} for $0\le i\le n$, $\bOm \equiv \left( \Lambda^{-1} \dot \Lambda \right)^\vee$ \revisionNACO{R1Q3}{is the ball's body angular velocity}, $\bGamma \equiv \Lambda^{-1} \mathbf{e}_3$ \revisionNACO{R1Q3}{is the spatial unit normal expressed in the body frame}, and $\tilde \bGamma \equiv \Lambda^{-1} \mathbf{F}_\mathrm{e}$ \revisionNACO{R1Q3}{is the external force expressed in the body frame}. For $\mathbf{v} = \begin{bmatrix} v_1 & v_2 & v_3  \end{bmatrix}^\mathsf{T} \in \mathbb{R}^3$, $\widehat{\mathbf{v}}^2=\widehat{\mathbf{v}}\widehat{\mathbf{v}}$ is the symmetric matrix given by
\begin{equation}
\widehat{\mathbf{v}}^2 = \begin{bmatrix}
-(v_2^2+v_3^2) & v_1 v_2 & v_1 v_3 \\
v_1 v_2 & -(v_1^2+v_3^2)  & v_2 v_3 \\
v_1 v_3 & v_2 v_3 & -(v_1^2+v_2^2) 
\end{bmatrix}
\end{equation}
and $\mathbf{v}_{12}$ is the projected vector consisting of the first two components of $\mathbf{v}$ so that
\begin{equation}
\mathbf{v}_{12} = \begin{bmatrix}v_1 & v_2 \end{bmatrix}^\mathsf{T} \in \mathbb{R}^2.
\end{equation}
\revisionNACO{R1Q5}{Let $N$ denote the magnitude of the normal force acting at the ball's CP. Let $f_\mathrm{s}$ and $\fricudir$ denote the magnitude of and unit-length direction antiparallel to the static friction acting at the ball's CP, respectively. As shown in \cite{putkaradze2018normalpub}, the magnitude of the normal force is
\begin{equation} \label{eq_normal_1d}
N = Mg+\left<\sum_{i=0}^n m_i \left[\dot \bOm \times \mathbf{s}_i + \bOm \times \left(\bOm \times \bzeta_i +2 \dot \theta_i \bzeta_i^{\prime} \right) + \dot \theta_i^2 \bzeta_i^{\dprime} + \ddot \theta_i \bzeta_i^{\prime} \right],\bGamma\right>-F_{\mathrm{e},3}
\end{equation}
and the static friction is  
\begin{equation} \label{eq_static_friction_1d}
-f_\mathrm{s} \fricudir = \begin{bmatrix}
\left(\Lambda \sum_{i=0}^n m_i \left[\dot \bOm \times \mathbf{s}_i +  \bOm \times \left(\bOm \times \bzeta_i +2 \dot \theta_i \bzeta_i^{\prime} \right) + \dot \theta_i^2 \bzeta_i^{\dprime} + \ddot \theta_i \bzeta_i^{\prime} \right] - \mathbf{F}_\mathrm{e} \right)_{12}\\
0
\end{bmatrix}.
\end{equation}	
The dynamics encapsulated by \eqref{uncon_ball_eqns_explicit_1d} are valid only if the ball does not detach from the surface ($N>0$) and rolls without slipping ($\mu_\mathrm{s} N \ge f_\mathrm{s}$), where $\mu_\mathrm{s}$ denotes the coefficient of static friction between the ball and the surface.}

\subsection{Rolling Disk} \label{ssec_disk_uncontrolled}
\begin{figure}[h]
	\centering
	\includegraphics[width=0.6\linewidth]{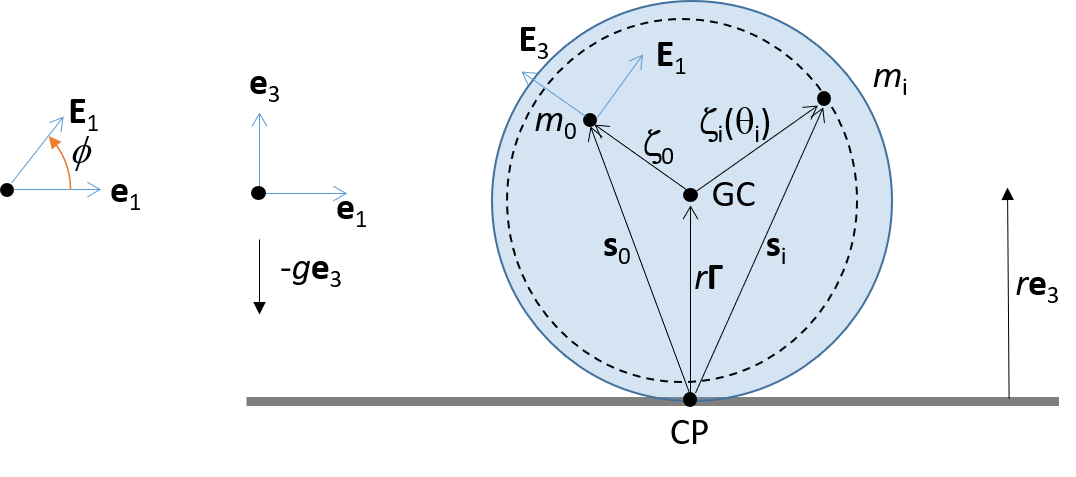}
	\caption{A disk of radius $r$ and mass $m_0$ rolls without slipping in the $\mathbf{e}_1$-$\mathbf{e}_3$ plane. $\mathbf{e}_2$ and $\mathbf{E}_2$ are directed into the page and are omitted from the figure. The disk's center of mass is denoted by $m_0$. The disk's motion is actuated by $n$ point masses, each of mass $m_i$, $1 \le i \le n$, and each moving along its own rail fixed inside the disk. The point mass depicted here by $m_i$ moves along a circular hoop in the disk that is not centered on the disk's geometric center (GC). The disk's orientation is determined by $\phi$, the angle measured counterclockwise from $\mathbf{e}_1$ to $\mathbf{E}_1$.  }
	\label{fig:detailed_rolling_disk}
\end{figure}
Now suppose that  the ball's inertia is such that one of  the ball's principal axes, say the one labeled $\mathbf{E}_2$, is orthogonal to the plane containing the GC and CM. 
Also assume that all the  point masses move along 1-d  rails which lie in the plane containing the GC and CM. Moreover, suppose that the ball is oriented initially so that the plane containing the GC and CM coincides with the $\mathbf{e}_1$-$\mathbf{e}_3$ plane and that the external force $\mathbf{F}_\mathrm{e}$ acts in the $\mathbf{e}_1$-$\mathbf{e}_3$ plane. Then for all time, the ball will remain oriented so that the plane containing the GC and CM coincides with the $\mathbf{e}_1$-$\mathbf{e}_3$ plane and the ball will only move in the $\mathbf{e}_1$-$\mathbf{e}_3$ plane, with the ball's rotation axis always parallel to $\mathbf{e}_2$. Note that the dynamics of this system are equivalent to that of the Chaplygin disk \cite{Ho2011_pII}, equipped with point masses, rolling in the $\mathbf{e}_1$-$\mathbf{e}_3$ plane, and where the Chaplygin disk (minus the  point masses) has polar moment of inertia $d_2$. Therefore, henceforth, this particular ball with this special inertia, orientation,  and placement of the rails and point masses, may be referred to as the disk or the rolling disk. Figure~\ref{fig:detailed_rolling_disk} depicts the rolling disk. Let $\phi$ denote the angle between $\mathbf{e}_1$ and $\mathbf{E}_1$, measured counterclockwise from $\mathbf{e}_1$ to $\mathbf{E}_1$. Thus, if $\dot \phi > 0$, the disk rolls in the $-\mathbf{e}_1$ direction and $\bOm$ has the same direction as $-\mathbf{e}_2$, and  if $\dot \phi < 0$, the disk rolls in the $\mathbf{e}_1$ direction and $\bOm$ has the same direction as $\mathbf{e}_2$. 

As shown in \cite{Putkaradze2018dynamicsP,putkaradze2018normalpub}, the uncontrolled equation of motion for this rolling disk is
\begin{equation} \label{eqmo_chap_disk_4}
\ddot \phi = \frac{ -r F_{\mathrm{e},1}+ \sum_{i=0}^n m_i K_i }{d_2+\sum_{i=0}^n m_i \left[\left( r \sin \phi + \zeta_{i,1} \right)^2+\left( r \cos \phi+ \zeta_{i,3} \right)^2 \right]} \equiv \kappa\left(t,\btheta,\dot \btheta,\phi,\dot \phi,\ddot \btheta\right),
\end{equation}
where 
\begin{equation} \label{eq_K_i}
\begin{split}
K_i &\equiv  \left(g+ r {\dot \phi}^2 \right) \left(\zeta_{i,3} \sin \phi - \zeta_{i,1} \cos \phi  \right)+
\left(r \cos \phi + \zeta_{i,3} \right) \left(- 2 \dot \phi {\dot \theta}_i \zeta_{i,3}^{\prime} + {\dot \theta}_i^2 \zeta_{i,1}^{\dprime} + {\ddot \theta}_i \zeta_{i,1}^{\prime} \right)\\
&\hphantom{\equiv}   - \left(r \sin \phi + \zeta_{i,1} \right) \left( 2 \dot \phi {\dot \theta}_i \zeta_{i,1}^{\prime}+ {\dot \theta}_i^2 \zeta_{i,3}^{\dprime} + {\ddot \theta}_i \zeta_{i,3}^{\prime} \right). 
\end{split}
\end{equation}
In \eqref{eqmo_chap_disk_4}, $\kappa$ is a function that depends on time ($t$) through the possibly time-varying external force $F_{\mathrm{e},1}(t)$, on the  point mass parameterized positions ($\btheta$), velocities ($\dot \btheta$), and accelerations ($\ddot \btheta$), and on the disk's orientation angle ($\phi$) and its time-derivative ($\dot \phi$). The spatial $\mathbf{e}_1$ position $z$ of the disk's GC is given by 
\begin{equation} \label{eq_disk_GC}
z = z_a-r \left(\phi-\phi_a \right),
\end{equation}
where $z_a$ is the spatial $\mathbf{e}_1$ position of the disk's GC at time $t=a$ and $\phi_a$ is the disk's angle at time $t=a$. 
\revisionNACO{R1Q5}{Let $N$ denote the magnitude of the normal force acting at the disk's CP. Let $f_\mathrm{s}$ and $\fricudir$ denote the magnitude of and unit-length direction antiparallel to the static friction acting at the disk's CP, respectively. As shown in \cite{putkaradze2018normalpub}, the magnitude of the normal force is
\begin{equation} \label{eq_normal_disk}
\begin{split}
N&=Mg+\sum_{i=0}^n m_i \Big[
\left(-\ddot \phi \zeta_{i,3} - {\dot \phi}^2 \zeta_{i,1} - 2 \dot \phi {\dot \theta}_i \zeta_{i,3}^{\prime} + {\dot \theta}_i^2 \zeta_{i,1}^{\dprime} + {\ddot \theta}_i \zeta_{i,1}^{\prime} \right) \sin \phi \\
&\hphantom{=Mg+\sum_{i=0}^n m_i \Big[}   + \left(\ddot \phi \zeta_{i,1} - {\dot \phi}^2 \zeta_{i,3}+ 2 \dot \phi {\dot \theta}_i \zeta_{i,1}^{\prime}+ {\dot \theta}_i^2 \zeta_{i,3}^{\dprime} + {\ddot \theta}_i \zeta_{i,3}^{\prime} \right) \cos \phi \Big]-F_{\mathrm{e},3}
\end{split}
\end{equation}
and the static friction is
\begin{equation} \label{eq_static_friction_disk}
\begin{split}
-f_\mathrm{s} \fricudir &= -\Big\{ Mr \ddot \phi + \sum_{i=0}^n m_i \Big[\left( \ddot \phi \zeta_{i,3} + {\dot \phi}^2 \zeta_{i,1} + 2 \dot \phi {\dot \theta}_i \zeta_{i,3}^{\prime} - {\dot \theta}_i^2 \zeta_{i,1}^{\dprime} - {\ddot \theta}_i \zeta_{i,1}^{\prime} \right) \cos \phi \\
&\hphantom{= -\Big\{ Mr \ddot \phi + \sum_{i=0}^n m_i \Big[} +\left(\ddot \phi \zeta_{i,1} - {\dot \phi}^2 \zeta_{i,3}+ 2 \dot \phi {\dot \theta}_i \zeta_{i,1}^{\prime}+ {\dot \theta}_i^2 \zeta_{i,3}^{\dprime} + {\ddot \theta}_i \zeta_{i,3}^{\prime} \right) \sin \phi \Big] + F_{\mathrm{e},1} \Big\} \mathbf{e}_1.
\end{split}
\end{equation}
The dynamics encapsulated by \eqref{eqmo_chap_disk_4} are valid only if the disk does not detach from the surface ($N>0$) and rolls without slipping ($\mu_\mathrm{s} N \ge f_\mathrm{s}$), where $\mu_\mathrm{s}$ denotes the coefficient of static friction between the disk and the surface.}

\subsection{Numerical Methods} \label{ssec_numerical_methods}

In Sections~\ref{sec_disk_sim} and \ref{sec_ball_sim}, the motions of the rolling disk and ball are simulated in \mcode{MATLAB} R2019b by numerically solving the controlled equations of motion \eqref{eq_pmp_bvp_disk} and \eqref{eq_pmp_bvp_ball} corresponding to the optimal control problems \eqref{dyn_opt_problem_disk} and \eqref{dyn_opt_problem_ball} for the rolling disk and ball, respectively. Subsection~\ref{ssec_disk_sim} simulates the rolling disk, while Subsections~\ref{ssec_ball_sim} and \ref{ssec_ball_sim_redux} simulate the rolling ball. Because the controlled equations of motion have a very small radius of convergence \cite{BrHo1975applied,bryson1999dynamic,betts2010practical}, a direct method, namely the \mcode{MATLAB} toolbox \mcode{GPOPS-II} \cite{patterson2014gpops} version 2.5, is first used to construct a good initial guess. In these simulations, \mcode{GPOPS-II} is configured to use the NLP solver SNOPT \cite{snopt76,GilMS05} version 7.6.0, though \mcode{GPOPS-II} can also be configured to use the NLP solver IPOPT \cite{wachter2006implementation}. For the rolling disk, the direct method is used to solve the rolling disk optimal control problem \eqref{dyn_opt_problem_disk}. When using the direct method to solve the rolling ball optimal control problem, the differential-algebraic equation (DAE) formulation \eqref{dyn_opt_problem_dae2_ball} is solved first. The direct method solution to the DAE formulation is then used as an initial guess to solve the ODE formulation \eqref{dyn_opt_problem_ball}, which is consistent with the controlled equations of motion \eqref{eq_pmp_bvp_ball} for the rolling ball, by the direct method. The \mcode{MATLAB} automatic differentiation toolbox \mcode{ADiGator} \cite{weinstein2017algorithm,weinstein2015utilizing} version 1.5 is used to supply vectorized first derivatives (i.e. Jacobians) to the direct method solver \mcode{GPOPS-II}, since SNOPT accepts first, but not second, derivatives.

Starting from the initial guess provided by the direct method, the controlled equations of motion \eqref{eq_pmp_bvp_disk} and \eqref{eq_pmp_bvp_ball} are solved by predictor-corrector continuation in the parameter $\mu$, utilizing the algorithm described in Appendix~\ref{app_sweep_predictor_corrector}. The predictor-corrector continuation method uses the \mcode{MATLAB} global method ODE TPBVP solvers \mcode{sbvp} \cite{auzinger2003collocation} version 1.0 or \mcode{bvptwp} \cite{cash2013algorithm} version 1.0. By vectorized automatic differentiation of $H_{\bx}$, $\bpi$, and $\hat{\mathbf{f}}$, \mcode{ADiGator} is used to numerically construct the Jacobians of the normalized ODE velocity function \eqref{eq_DPhi_z} and \eqref{eq_DPhi_mu}. By non-vectorized automatic differentiation of the Hamiltonian $H$, the initial condition function $\bsigma$, the final condition function $\bpsi$, and the endpoint function $G$, \mcode{ADiGator} is used to numerically construct the normalized BC function \eqref{eq_tUps} and the Jacobians of the normalized BC function \eqref{eq_Jac_tUps_tz0}, \eqref{eq_Jac_tUps_tz1}, and \eqref{eq_Jac_tUps_mu}. These functions are needed by the ODE TPBVP solvers \mcode{sbvp} and \mcode{bvptwp} to solve the controlled equations of motion \eqref{eq_pmp_bvp_disk} and \eqref{eq_pmp_bvp_ball} by predictor-corrector continuation in the parameter $\mu$. 

In contrast to the direct method, the controlled equations of motion obtained via the indirect method have a very small radius of convergence \cite{BrHo1975applied,bryson1999dynamic,betts2010practical}. Therefore, the direct method is needed to initialize the predictor-corrector continuation of the controlled equations of motion. Predictor-corrector continuation is used in conjunction with the indirect, rather than direct, method, because a predictor-corrector continuation direct method requires a predictor-corrector continuation NLP solver. Even though predictor-corrector continuation NLP solver algorithms are provided in \cite{zangwill1981pathways,kungurtsev2017predictor}, there do not seem to be any publicly available predictor-corrector continuation NLP solvers.

\section{Trajectory Tracking for the Rolling Disk} \label{sec_disk_sim}

\subsection{Optimal Control Problem and Controlled Equations of Motion} \label{ssec_disk_controlled}
\begin{figure}[h] 
	\centering
	\includegraphics[scale=.6]{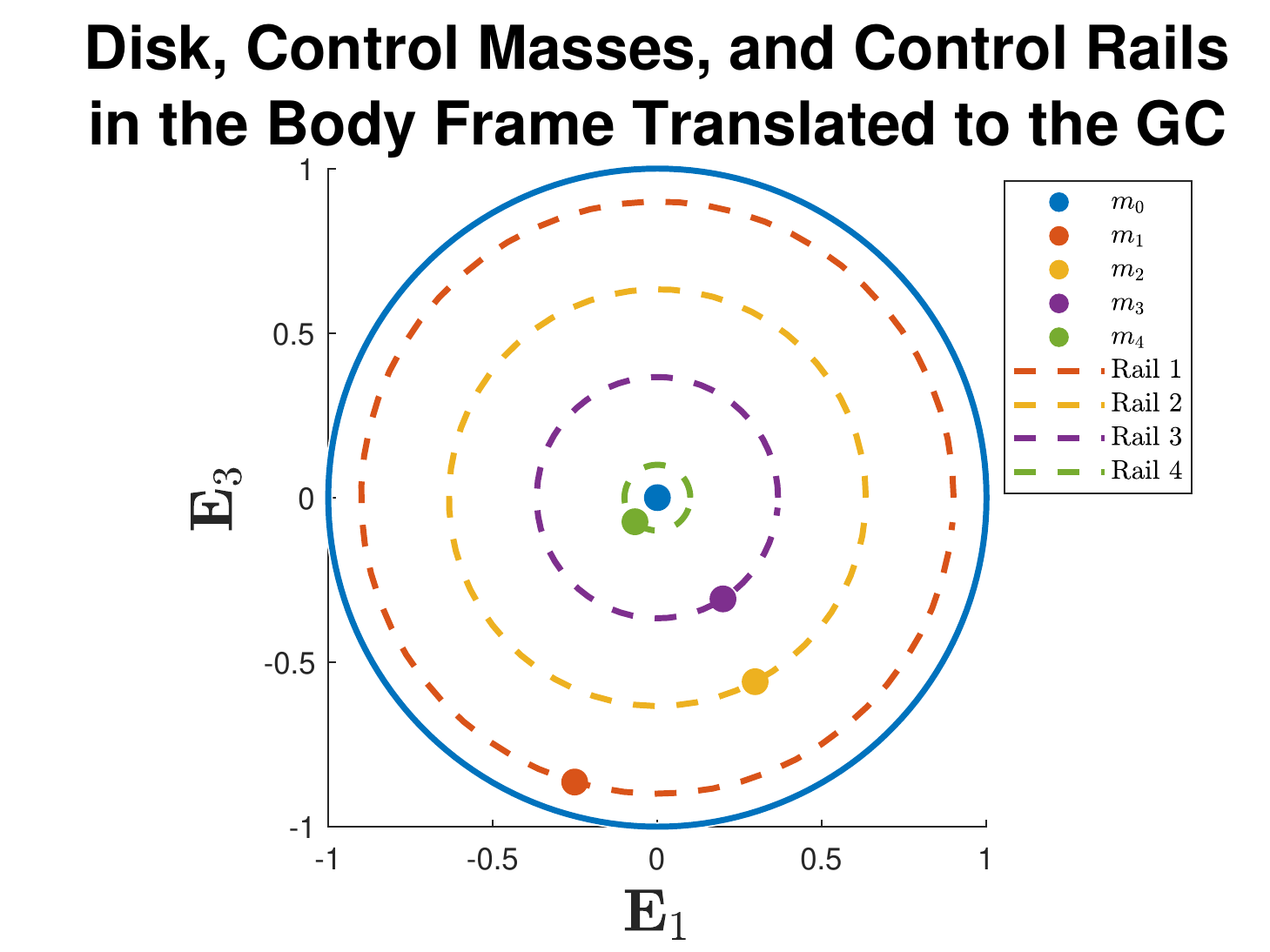}
	\caption{The disk of radius $r=1$ actuated by $4$ control masses, $m_1$, $m_2$, $m_3$, and $m_4$, each on its own circular control rail. The control rail radii are $r_1=.9$, $r_2=.6\overline{3}$, $r_3=.3\overline{6}$, and $r_4=.1$. The location of the disk's CM is denoted by $m_0$.}
	\label{fig_dsim8_control_masses_rails}
\end{figure}

In the next subsection, numerical solutions of the controlled equations of motion for the rolling disk are presented, where the goal is to move the disk between a pair of points while the disk's GC tracks a prescribed trajectory. \revisionNACO{R1Q5}{A rolling disk of mass $m_0=4$, radius $r=1$, polar moment of inertia $d_2=1$, and with the CM coinciding with the GC (i.e. $\bzeta_0=\mathbf{0}$) is simulated. These physical parameters for the disk are consistent with the necessary and sufficient conditions stipulated by Inequality 3.2 in \cite{rozenblat2016choice}; in that inequality, note that $d_1+d_3=d_2=1$ since the disk is a planar distribution of mass.} There are $n=4$ control masses, each of mass $1$ so that $m_1=m_2=m_3=m_4=1$, located on concentric circular control rails centered on the GC of radii $r_1=.9$, $r_2=.6\overline{3}$, $r_3=.3\overline{6}$, and $r_4=.1$, as shown in Figure~\ref{fig_dsim8_control_masses_rails}. \revision{R1Q18}{For $1 \le i \le 4$, the position of $m_i$ in the body frame centered on the GC is
\begin{equation}
\bzeta_i\left(\theta_i\right) = r_i \begin{bmatrix} \cos \theta_i \\ 0 \\ \sin \theta_i  \end{bmatrix}.
\end{equation}}

The total mass of the system  is \revisionNACO{R1Q5}{$M=8$}, and gravity is \revisionNACO{R1Q5}{$g=9.81$}. There is no external force acting on the disk's GC, so that $F_{\mathrm{e},1}=0$ in \eqref{eqmo_chap_disk_4}. The initial time is fixed to $a=0$ and the final time is fixed to \revisionNACO{R1Q5}{$b=12$}. The disk's GC starts at rest at $z_a=0$ at time $a=0$ and stops at rest at $z_b=1$ at time \revisionNACO{R1Q5}{$b=12$}. Table~\ref{table_disk_ICs} shows parameter values used in the rolling disk's initial conditions \eqref{eq_disk_initial_conds} and final conditions \eqref{eq_disk_final_conds}. Since the initial orientation of the disk is $\phi_a=0$ and since the initial configurations of the control masses are given by $\btheta_a=\begin{bmatrix} \unaryminus \frac{\pi}{2} & \unaryminus \frac{\pi}{2} & \unaryminus \frac{\pi}{2} & \unaryminus \frac{\pi}{2} \end{bmatrix}^\mathsf{T}$, all the control masses are initially located directly below the GC. In order for the disk to start and stop at rest, $\dot \btheta_a = \dot \btheta_b =\begin{bmatrix} 0 & 0 & 0 & 0 \end{bmatrix}^\mathsf{T}$ and $\dot z_a = \dot z_b = 0$. Table~\ref{table_disk_FCs} shows parameter values used in the rolling disk's final conditions \eqref{eq_disk_final_conds}. 

\begin{table}[h!]
	\centering 
	{ 
		\setlength{\extrarowheight}{1.5pt}
		\begin{tabular}{| c | c |} 
			\hline
			\textbf{Parameter} & \textbf{Value} \\ 
			\hline\hline 
			$\btheta_a$ & $\begin{bmatrix} \unaryminus \frac{\pi}{2} & \unaryminus \frac{\pi}{2} & \unaryminus \frac{\pi}{2} & \unaryminus \frac{\pi}{2} \end{bmatrix}^\mathsf{T}$  \\  
			\hline
			$\dot \btheta_a$ & $\begin{bmatrix} 0 & 0 & 0 & 0 \end{bmatrix}^\mathsf{T}$ \\ 
			\hline
			$\phi_a$ & $0$ \\ 
			\hline
			$z_a$ & $0$ \\
			\hline
			$\dot z_a$ & $0$ \\ 
			\hline
		\end{tabular} 
	}
	\caption{Initial condition parameter values for the rolling disk. Refer to \eqref{eq_disk_initial_conds} and \eqref{eq_disk_final_conds}.}
	\label{table_disk_ICs}
\end{table}

\begin{table}[h!]
	\centering 
	{ 
		\setlength{\extrarowheight}{1.5pt}
		\begin{tabular}{| c | c |} 
			\hline
			\textbf{Parameter} & \textbf{Value} \\ 
			\hline\hline 
			$\dot \btheta_b$ & $\begin{bmatrix} 0 & 0 & 0 & 0 \end{bmatrix}^\mathsf{T}$ \\ 
			\hline
			$z_b$ & $1$ \\
			\hline
			$\dot z_b$ & $0$ \\ 
			\hline
		\end{tabular} 
	}
	\caption{Final condition parameter values for the rolling disk. Refer to \eqref{eq_disk_final_conds}.}
	\label{table_disk_FCs}
\end{table}

\begin{figure}[!ht]
	\centering
	\subfloat[The GC tracks the desired path crudely  when $\alpha=.1$.]{\includegraphics[scale=.5]{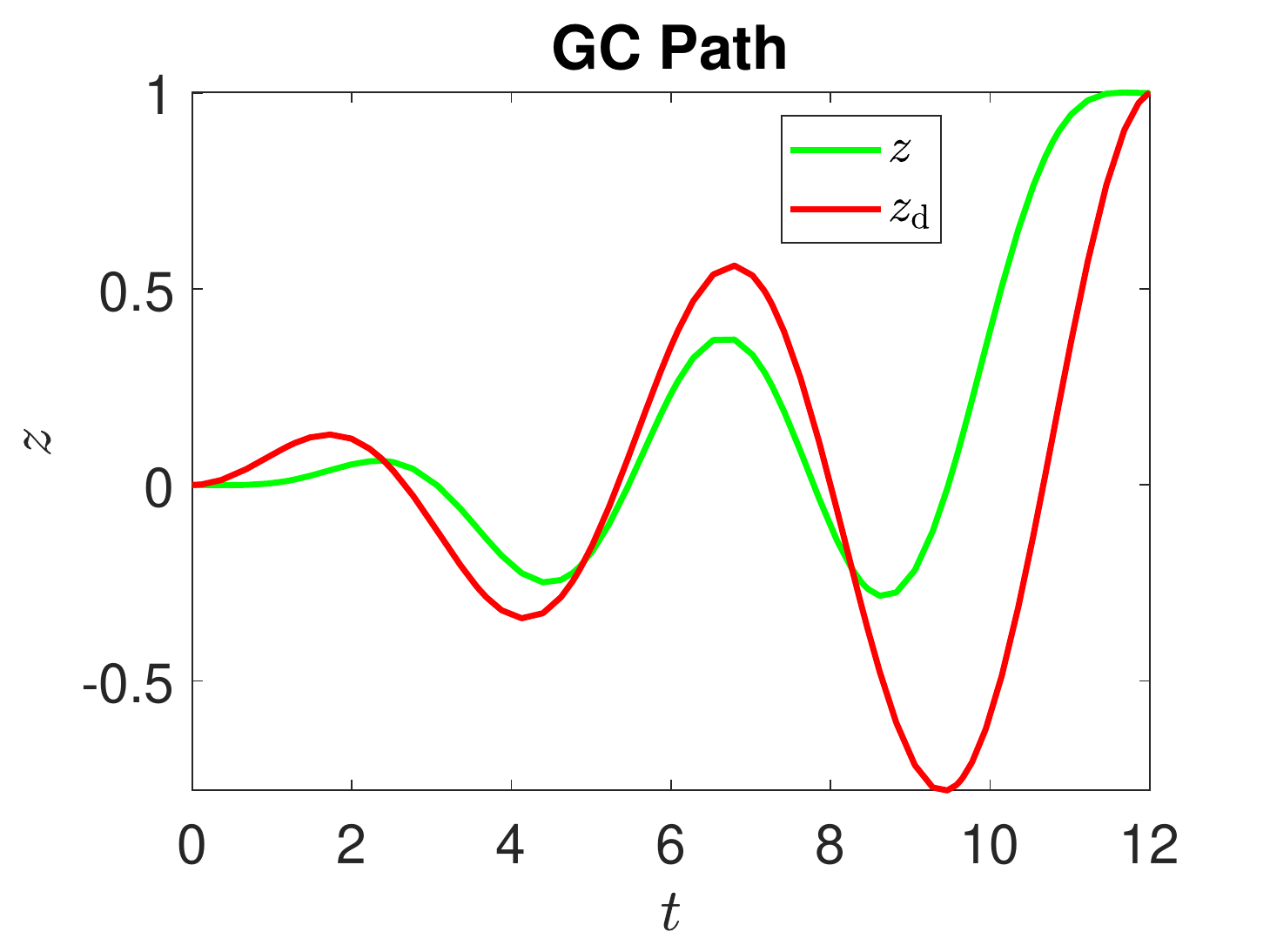}\label{fig_dsim8_dm_gc_path}}
	\hspace{5mm}
	\subfloat[The GC tracks the desired path very accurately when $\alpha \approx 272$.]{\includegraphics[scale=.5]{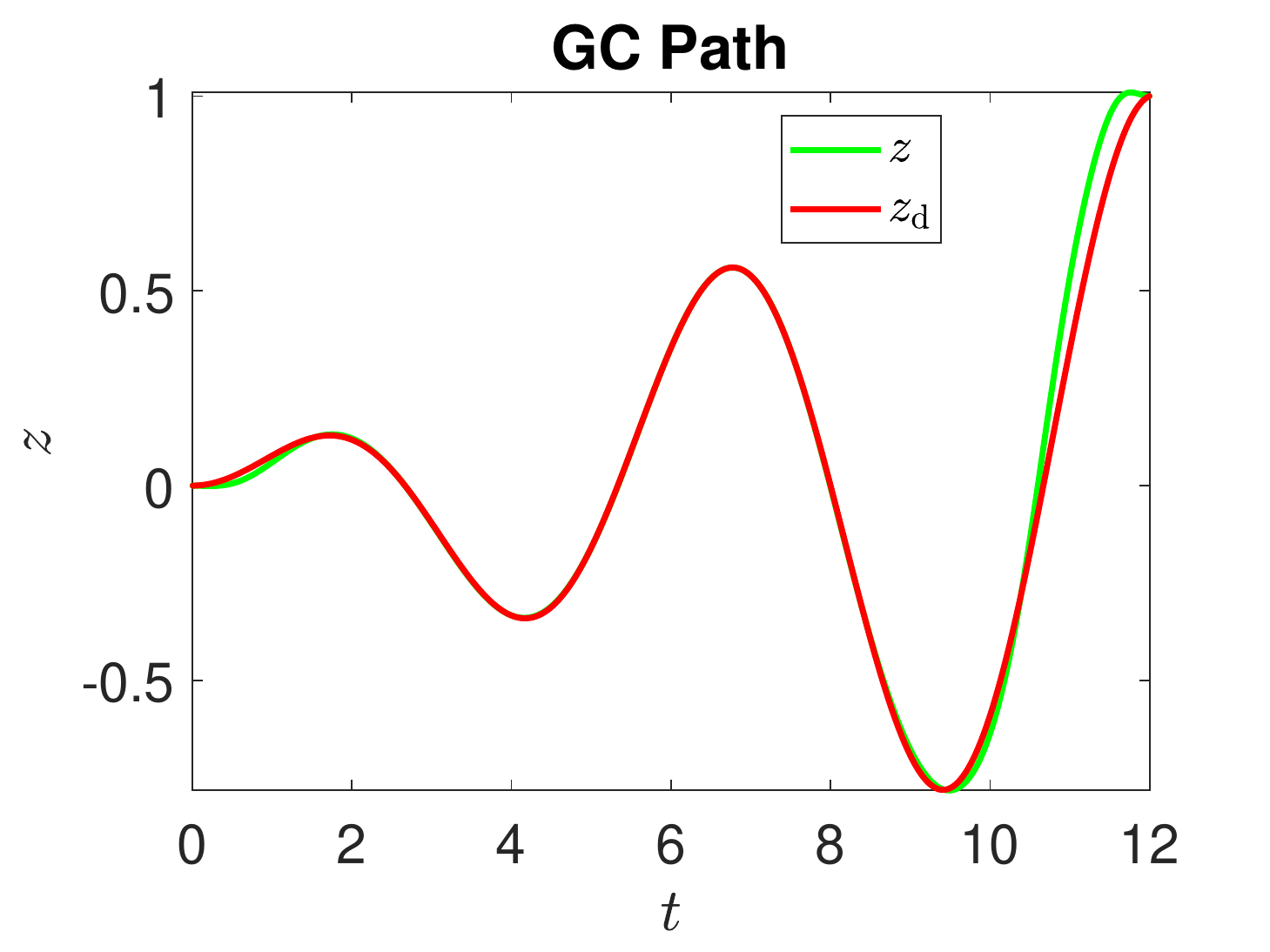}\label{fig_dsim8_pc_gc_path}}
	\\
	\subfloat[The motions of the center of masses are modest when $\alpha=.1$.]{\includegraphics[scale=.5]{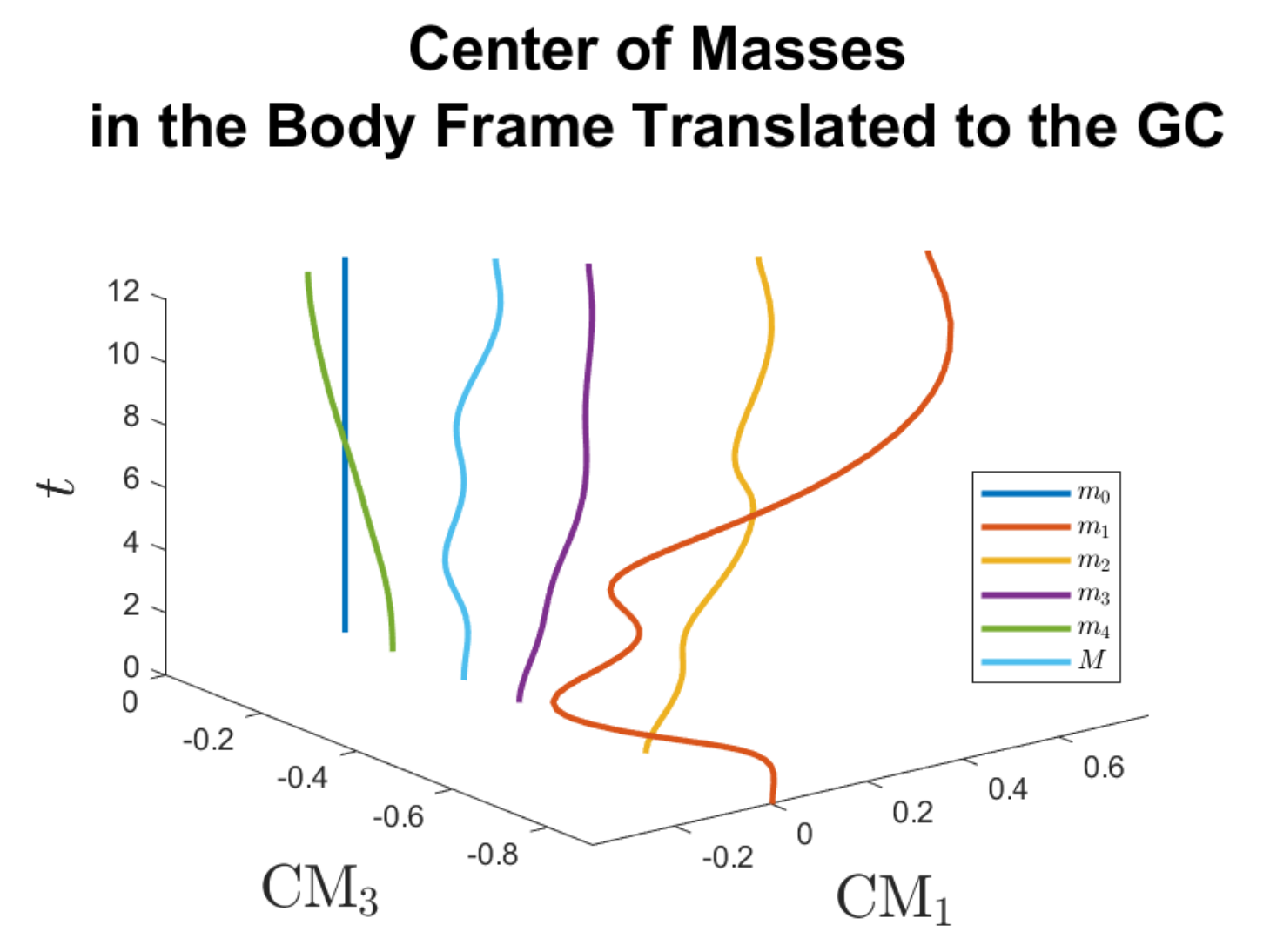}\label{fig_dsim8_dm_cm_bf}}
	\hspace{5mm}
	\subfloat[The motions of the center of masses are more serpentine when $\alpha \approx 272$.]{\includegraphics[scale=.5]{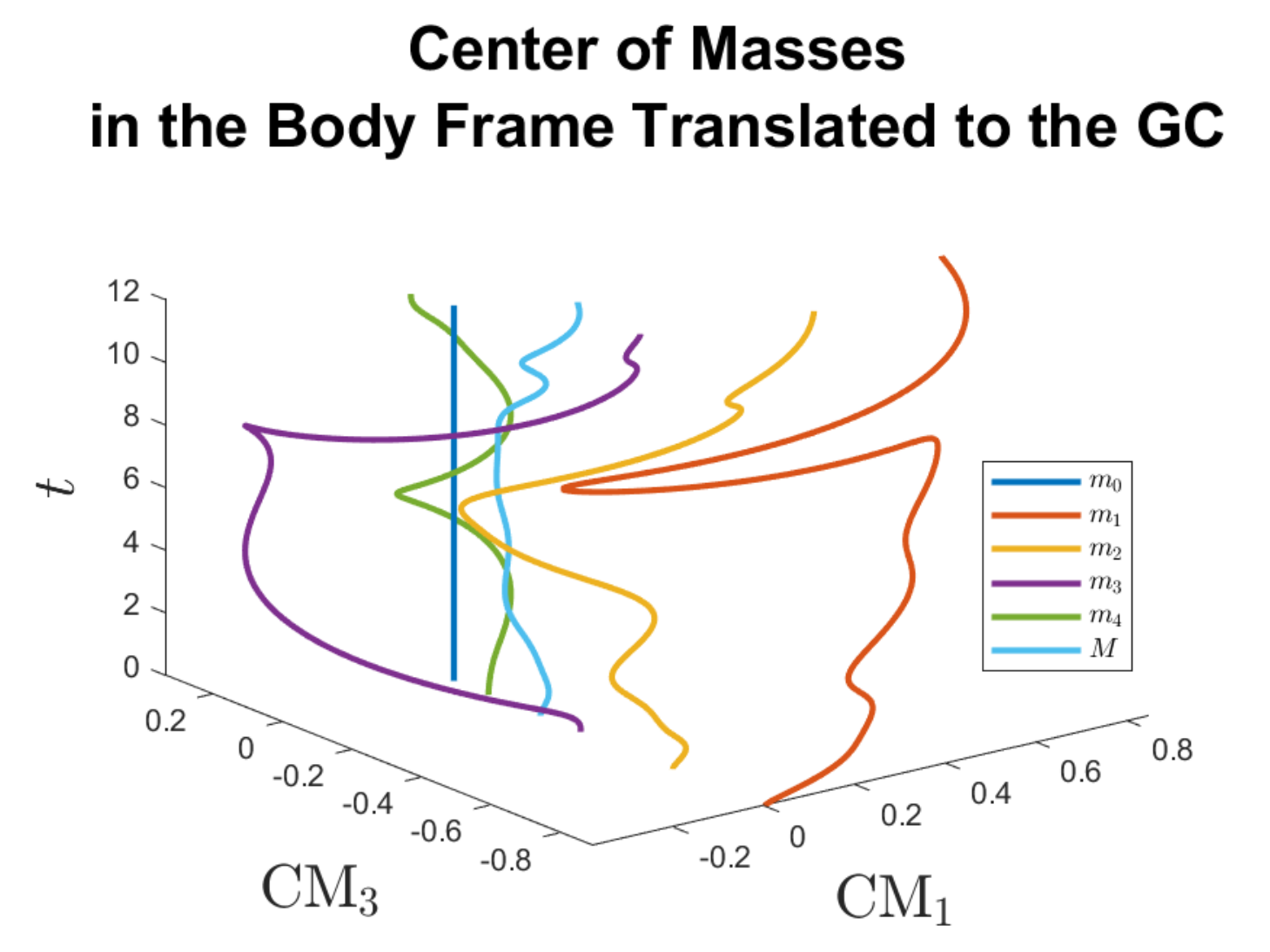}\label{fig_dsim8_pc_cm_bf}}
	\\
	\subfloat[The controls have relatively small magnitudes when $\alpha=.1$.]{\includegraphics[scale=.5]{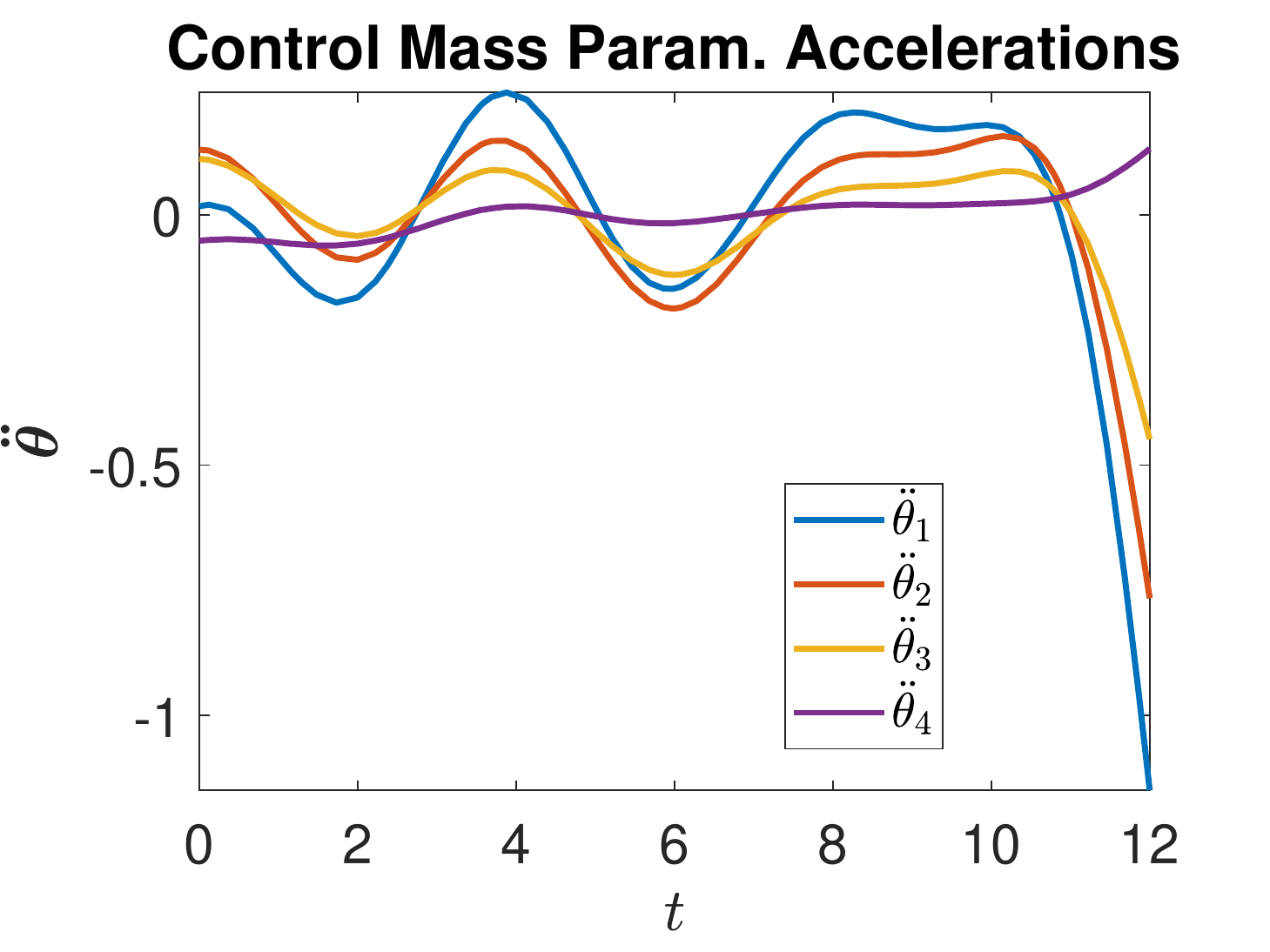}\label{fig_dsim8_dm_controls}}
	\hspace{5mm}
	\subfloat[The controls have relatively larger magnitudes when $\alpha \approx 272$.]{\includegraphics[scale=.5]{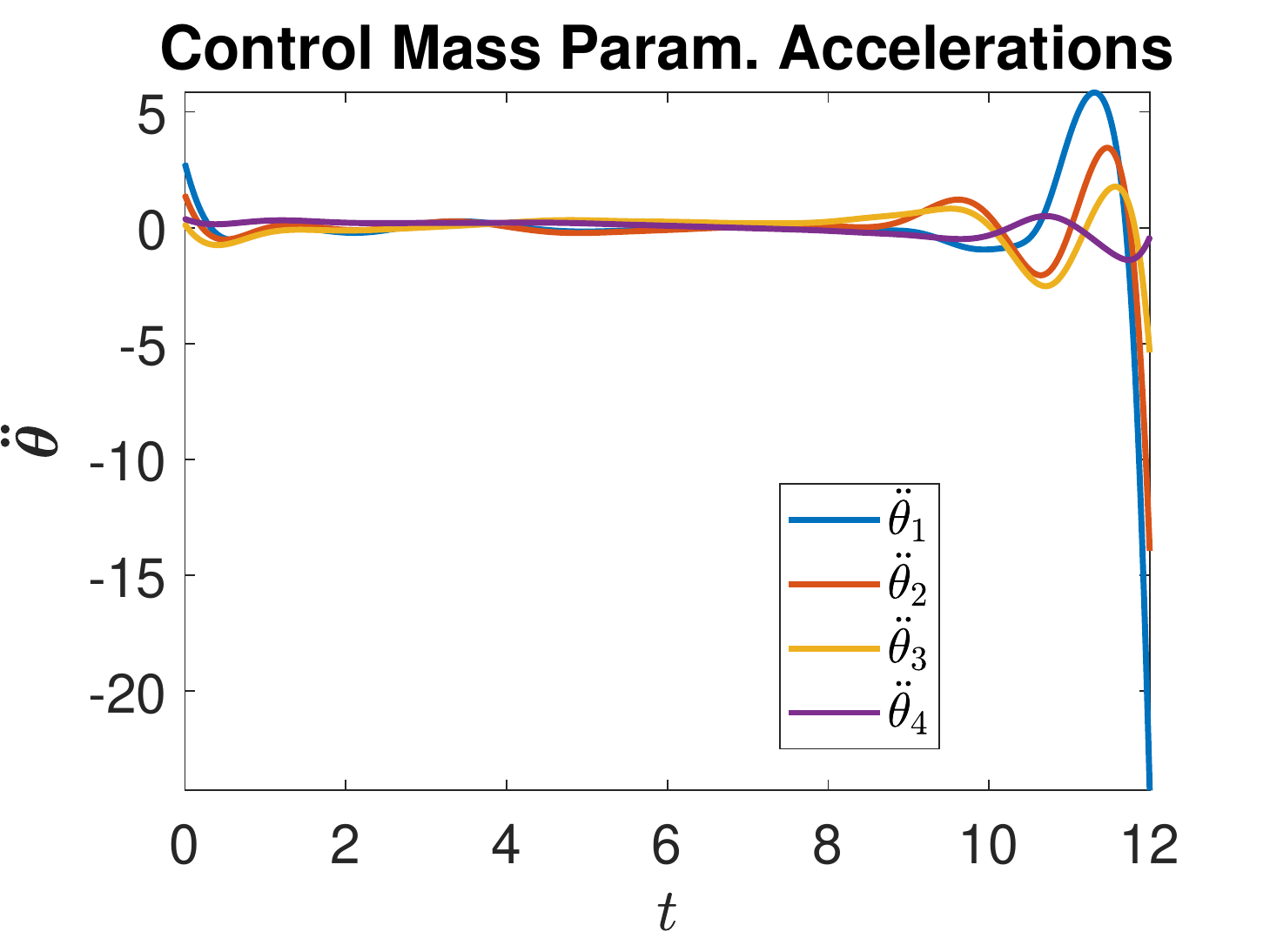}\label{fig_dsim8_pc_controls}}
	\caption{Numerical solutions of the rolling disk optimal control problem \eqref{dyn_opt_problem_disk} using $4$ control masses for $ \gamma_1=\gamma_2=\gamma_3=\gamma_4=.1$ and for fixed initial and final times. The direct method results for $\alpha=.1$ are shown in the left column, while the predictor-corrector continuation indirect method results for $\alpha \approx 272$ are shown in the right column. The direct method solution tracks the desired GC path crudely, whereas the indirect method solution tracks the desired GC path very accurately at the expense of larger magnitude controls.}
	\label{fig_dsim8}
\end{figure}

The desired GC path $z_\mathrm{d}$ is depicted by the red curve in Figures~\ref{fig_dsim8_dm_gc_path} and \ref{fig_dsim8_pc_gc_path}. $z_\mathrm{d}$ encourages the disk's GC to track a sinusoidally-modulated linear trajectory \revision{R1Q19}{connecting $z(0)=0$ with $z(12)=1$}. That is, the disk is encouraged to roll right, then left, then right, then left, and finally to the right, with the amplitude of each successive roll increasing from the previous one. Specifically, $z_\mathrm{d}$ is given by
\begin{equation} \label{eq_z_d}
z_\mathrm{d}(t) \equiv \left[z_a w(t)+\tilde{z}_\mathrm{d}(t) \left(1-w(t)\right) \right] \left(1-y(t)\right)+z_b y(t),
\end{equation}
where
\begin{equation} \label{eq_sigmoid}
S(t) \equiv \frac{1}{2} \left[1+\tanh{\left( \frac{-t}{\epsilon}\right)} \right],
\end{equation}
\begin{equation} \label{eq_tran_sigmoidl}
w(t) \equiv S\left(t-a\right),
\end{equation}
\begin{equation} \label{eq_tran_sigmoidr}
y(t) \equiv S\left(-t+b\right),
\end{equation}
and
\begin{equation}
\tilde{z}_\mathrm{d}(t) \equiv \left[ z_a+\left(z_b-z_a\right) \frac{t-a}{b-a} \right] \sin\left(\frac{9\pi}{2}\frac{t-a}{b-a}\right),
\end{equation}
with $\epsilon=.01$ in \eqref{eq_sigmoid}. The reader is referred to \cite{putkaradze2020optimal} for further details about the properties and construction of \eqref{eq_z_d}.

The optimal control problem for the rolling disk is
\begin{equation}
\min_{\bu} J 
\mbox{\, s.t. \,}
\left\{
\begin{array}{ll}
\dot {\bx} = \mathbf{f}\left(\bx,\bu\right), \\
\bsigma\left(\bx(a)\right) = \mathbf{0},\\
\bpsi\left(\bx(b)\right) = \mathbf{0}.
\end{array}
\right.
\label{dyn_opt_problem_disk}
\end{equation}
In \eqref{dyn_opt_problem_disk}, the system state $\bx$ and control $\bu$ are
\begin{equation}
\bx \equiv \begin{bmatrix} \btheta \\ \dot \btheta \\ \phi \\ \dot \phi  \end{bmatrix} \quad \mathrm{and} \quad \bu \equiv \ddot \btheta,
\end{equation}
where $\btheta, \dot \btheta, \ddot \btheta \in \mathbb{R}^4$ and $\phi, \dot \phi  \in \mathbb{R}$. In \eqref{dyn_opt_problem_disk}, the system dynamics defined for $a \le t \le b$ are
\begin{equation} \label{eq_disk_dynamics}
\dot {\bx} = \begin{bmatrix} \dot \btheta \\ \ddot \btheta \\ \dot \phi \\ \ddot \phi  \end{bmatrix}  = \mathbf{f}\left(\bx,\bu\right) \equiv \begin{bmatrix} \dot \btheta \\ \bu  \\ \dot \phi \\ \kappa\left(\bx,\bu\right)  \end{bmatrix},
\end{equation}
where $\kappa\left(\bx,\bu \right)$ is given by the right-hand side of \eqref{eqmo_chap_disk_4}. In \eqref{eq_disk_dynamics}, the time-dependence of $\kappa$ is dropped since $F_{\mathrm{e},1}=0$ in \eqref{eqmo_chap_disk_4} for these simulations. In \eqref{dyn_opt_problem_disk}, the prescribed initial conditions at time $t=a$ are
\begin{equation} \label{eq_disk_initial_conds}
\bsigma\left(\bx(a)\right) \equiv \begin{bmatrix} \btheta(a) - \btheta_a \\ \dot \btheta(a) - {\dot \btheta}_a \\ \phi(a)-\phi_a \\ -r \dot \phi(a) - {\dot z}_a  \end{bmatrix} = \mathbf{0},
\end{equation}
and the prescribed final conditions at time $t=b$ are
\begin{equation} \label{eq_disk_final_conds}
\bpsi\left(\bx(b)\right) \equiv \begin{bmatrix} \Pi \left(\tilde \Lambda\left(\phi(b)\right) \left[ \frac{1}{M} \sum_{i=0}^4 m_i \bzeta_i\left(\theta_i(b)\right) \right] \right) \\ \dot \btheta(b) - {\dot \btheta}_b  \\ z_a - r \left(\phi(b)-\phi_a\right) - z_b \\ -r \dot \phi(b) - {\dot z}_b  \end{bmatrix} = \mathbf{0}.
\end{equation} 
In \eqref{eq_disk_final_conds},
\begin{equation}
\frac{1}{M} \sum_{i=0}^4 m_i \bzeta_i\left(\theta_i(t)\right)
\end{equation}
is the total system CM expressed in the body frame translated to the disk's GC at time $t$,
\begin{equation}
\tilde \Lambda (\phi(t)) = \Lambda(t) = \begin{bmatrix} \cos \phi(t) & 0 & - \sin \phi(t) \\ 0 & 1 & 0 \\ \sin \phi(t) & 0 & \cos \phi(t) \end{bmatrix}
\end{equation} 
is the rotation matrix that maps the body to spatial frame at time $t$, and $\Pi$ is \revision{R1Q13}{the} projection onto the first component. Therefore, the first constraint in \eqref{eq_disk_final_conds} ensures that the total system CM is above or below the disk's GC in the spatial frame at the final time $t=b$, so that, in conjunction with the final condition parameter values given in Table~\ref{table_disk_FCs}, the disk stops at rest. In \eqref{dyn_opt_problem_disk}, the performance index is
\begin{equation} \label{eq_disk_J}
J \equiv \int_a^b L\left(t,\bx,\bu,\mu\right) \dt 
= \int_a^b \left[ \frac{\alpha(\mu)}{2} \left(z_a - r \left(\phi-\phi_a \right) - z_\mathrm{d} \right)^2  + \sum_{i=1}^4 \frac{\gamma_i}{2} {\ddot \theta}_i^2 \right] \dt,
\end{equation} 
where the integrand cost function is
\begin{equation} \label{eq_disk_integrand_cost}
L\left(t,\bx,\bu,\mu\right) \equiv \frac{\alpha(\mu)}{2} \left(z_a - r \left(\phi-\phi_a\right) - z_\mathrm{d} \right)^2 + \sum_{i=1}^4 \frac{\gamma_i}{2} {\ddot \theta}_i^2,
\end{equation}
for positive coefficients $\alpha(\mu)$ and $\gamma_i$, $1 \le i \le 4$.
The first summand $\frac{\alpha(\mu)}{2} \left(z_a - r (\phi-\phi_a) - z_\mathrm{d} \right)^2$ in $L$ encourages the disk's GC to track the desired spatial $\mathbf{e}_1$ path $z_\mathrm{d}$, and \revision{R1Q4}{$\mu$ is a scalar continuation parameter used to construct a sequence of optimal control problems.} The next $4$ summands $\frac{\gamma_i}{2} {\ddot \theta}_i^2$, $1 \le i \le 4$, in $L$ limit the magnitude of the acceleration of the $i^\mathrm{th}$ control mass parameterization. Table~\ref{table_disk_integrand} shows the values set for the integrand cost function coefficients in \eqref{eq_disk_integrand_cost}.  

\begin{table}[h!]
	\centering 
	{ 
		\setlength{\extrarowheight}{1.5pt}
		\begin{tabular}{| c | c |} 
			\hline
			\textbf{Parameter} & \textbf{Value} \\ 
			\hline\hline 
			$\alpha(\mu)$ & $.1+\frac{.95-\mu}{.95-.00001}\left(5000-.1\right)$ \\
			\hline
			$\gamma_1=\gamma_2=\gamma_3=\gamma_4$ & $.1$ \\ 
			\hline
		\end{tabular} 
	}
	\caption{Integrand cost function coefficient values for the rolling disk when predictor-corrector continuation is performed in $\alpha$. Refer to \eqref{eq_disk_integrand_cost}.}
	\label{table_disk_integrand}
\end{table} 

As explained in Appendix~\ref{sec_optimal_control}, the controlled equations of motion for the rolling disk's optimal control problem \eqref{dyn_opt_problem_disk} are encapsulated by the ODE TPBVP:
\begin{equation} \label{eq_pmp_bvp_disk}
\begin{split}
\dot {\bx} &= \hat{H}_{\blam}^\mathsf{T} \left(t,\bx,\blam,\mu\right) = \hat{\mathbf{f}}\left(\bx,\blam\right)\equiv \mathbf{f}\left(\bx,\bpi\left(\bx,\blam\right)\right), \\
\dot {\blam} &= - \hat{H}_{\bx}^\mathsf{T} \left(t,\bx,\blam,\mu\right)=-H_{\bx}^\mathsf{T}\left(t,\bx,\blam,\bpi\left(\bx,\blam\right),\mu\right), \\
\left. \blam \right|_{t=a} &= -G_{\bx(a)}^\mathsf{T}, \quad G_{\bxi}^\mathsf{T} = \bsigma\left(\bx(a)\right) = \mathbf{0}, \\
\left. \blam \right|_{t=b} &= G_{\bx(b)}^\mathsf{T}, \quad G_{\bnu}^\mathsf{T} = \bpsi\left(\bx(b)\right) = \mathbf{0}.
\end{split}
\end{equation}
Subappendix~A.1 of \cite{putkaradze2020optimal} derives the formulas for constructing $H_{\bx}^\mathsf{T}$.
In \eqref{eq_pmp_bvp_disk}, $G$ is the endpoint function 
\begin{equation} \label{eq_disk_endpoint_fcn}
\begin{split}
G\left(\bx(a),\bxi,\bx(b),\bnu\right) &\equiv \bxi^\mathsf{T} \bsigma\left(\bx(a)\right)+\bnu^\mathsf{T} \bpsi\left(\bx(b)\right) \\
&=\bxi^\mathsf{T}  \begin{bmatrix} \btheta(a) - \btheta_a \\ \dot \btheta(a) - {\dot \btheta}_a \\ \phi(a)-\phi_a \\ -r \dot \phi(a) - {\dot z}_a  \end{bmatrix}+\bnu^\mathsf{T} \begin{bmatrix} \Pi \left(\tilde \Lambda\left(\phi(b)\right) \left[ \frac{1}{M} \sum_{i=0}^4 m_i \bzeta_i\left(\theta_i(b)\right) \right] \right) \\ \dot \btheta(b) - {\dot \btheta}_b  \\ z_a - r \left(\phi(b)-\phi_a\right) - z_b \\ -r \dot \phi(b) - {\dot z}_b  \end{bmatrix},
\end{split}
\end{equation}
\revision{R1Q4}{where $\bxi \in \mathbb{R}^{10}$ and $\bnu \in \mathbb{R}^{7}$ are constant Lagrange multiplier vectors enforcing the initial and final conditions, \eqref{eq_disk_initial_conds} and \eqref{eq_disk_final_conds}, respectively.} In \eqref{eq_pmp_bvp_disk}, $H$ is the Hamiltonian
\begin{equation} \label{eq_disk_ham}
\begin{split}
H\left(t,\bx,\blam,\bu,\mu\right) &\equiv L\left(t,\bx,\bu,\mu\right) + \blam^\mathsf{T} \mathbf{f}\left(\bx,\bu\right) \\
&= \frac{\alpha(\mu)}{2} \left(z_a - r \left(\phi-\phi_a \right) - z_\mathrm{d} \right)^2  + \sum_{i=1}^4 \frac{\gamma_i}{2} {\ddot \theta}_i^2 + \blam^\mathsf{T} \begin{bmatrix} \dot \btheta \\ \bu  \\ \dot \phi \\ \kappa\left(\bx,\bu \right)  \end{bmatrix},
\end{split}
\end{equation}
\revision{R1Q4}{where $\blam \in \mathbb{R}^{10}$ is a time-varying Lagrange multiplier vector enforcing the dynamics \eqref{eq_disk_dynamics}.}  In \eqref{eq_pmp_bvp_disk}, $\bpi$ is an analytical formula expressing the control $\bu$ as a function of the state $\bx$ and the costate $\blam$. The components of $\bpi$ are given by
\begin{equation}  \label{eq_ddtheta_exp}
{\ddot \theta}_i = \pi_i \left(\bx,\blam\right) \equiv -\gamma_i^{-1} \left\{ \lambda_{4+i} + \lambda_{10} \frac{m_i \left[ \left(r \cos \phi + \zeta_{i,3} \right)\zeta_{i,1}^{\prime}- \left(r \sin \phi + \zeta_{i,1} \right) \zeta_{i,3}^{\prime} \right]}{d_2+\sum_{k=0}^4 m_k \left[\left( r \sin \phi + \zeta_{k,1} \right)^2+\left( r \cos \phi+ \zeta_{k,3} \right)^2 \right]} \right\},
\end{equation}
for $1 \le i \le 4$. In \eqref{eq_pmp_bvp_disk}, $\hat H$ is the regular Hamiltonian
\begin{equation} \label{eq_disk_regular_Hamiltonian}
\begin{split}
\hat H\left(t,\bx,\blam,\mu\right) &\equiv H\left(t,\bx,\blam,\bpi \left(\bx,\blam\right),\mu\right) \\
&= \frac{\alpha(\mu)}{2} \left(z_a - r \left(\phi-\phi_a \right) - z_\mathrm{d} \right)^2 + \sum_{i=1}^4 \frac{\gamma_i}{2} {\pi}_i^2 \left(\bx,\blam\right)  + \blam^\mathsf{T} \begin{bmatrix} \dot \btheta \\ \bpi \left(\bx,\blam\right)  \\ \dot \phi \\ \kappa\left(\bx,\bpi \left(\bx,\blam\right) \right)  \end{bmatrix}.
\end{split}
\end{equation}
The reader is referred to \cite{putkaradze2020optimal} for a more general description of the rolling disk's optimal control problem \eqref{dyn_opt_problem_disk} and the associated controlled equations of motion \eqref{eq_pmp_bvp_disk}. 
 
\revisionNACO{R1Q4 \\ R1Q5}{\subsection{Numerical Solutions: Trajectory Tracking}} \label{ssec_disk_sim}
The direct method solver \mcode{GPOPS-II} is used to solve the optimal control problem \eqref{dyn_opt_problem_disk} when the integrand cost function coefficient is $\alpha=.1$. Predictor-corrector continuation  is then used to solve the controlled equations of motion \eqref{eq_pmp_bvp_disk}, starting from the direct method solution. The continuation parameter is $\mu$, which is used to adjust $\alpha$ according to the linear homotopy given in Table~\ref{table_disk_integrand}, so that $\alpha=.1$ when $\mu=.95$ and $\alpha \approx 272$ when $\mu \approx .8983$. The predictor-corrector continuation begins at $\mu=.95$, which is consistent with the direct method solution obtained at $\alpha=.1$.

\begin{figure}[!ht]
	\centering
	\subfloat[The magnitude of the normal force  when $\alpha=.1$.]{\includegraphics[scale=.5]{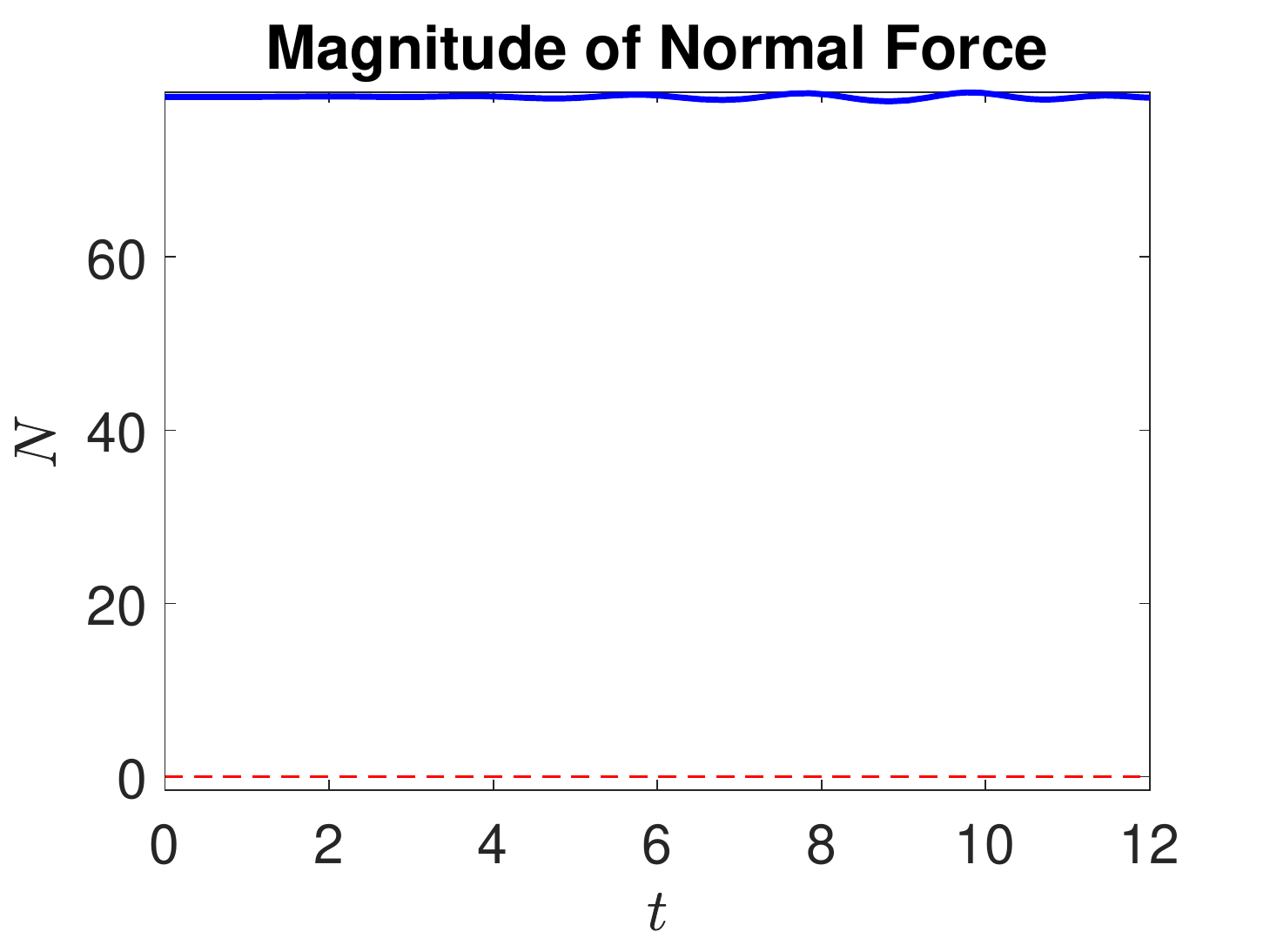}\label{fig_dsim8_dm_normal}}
	\hspace{5mm}
	\subfloat[The magnitude of the normal force when $\alpha \approx 272$.]{\includegraphics[scale=.5]{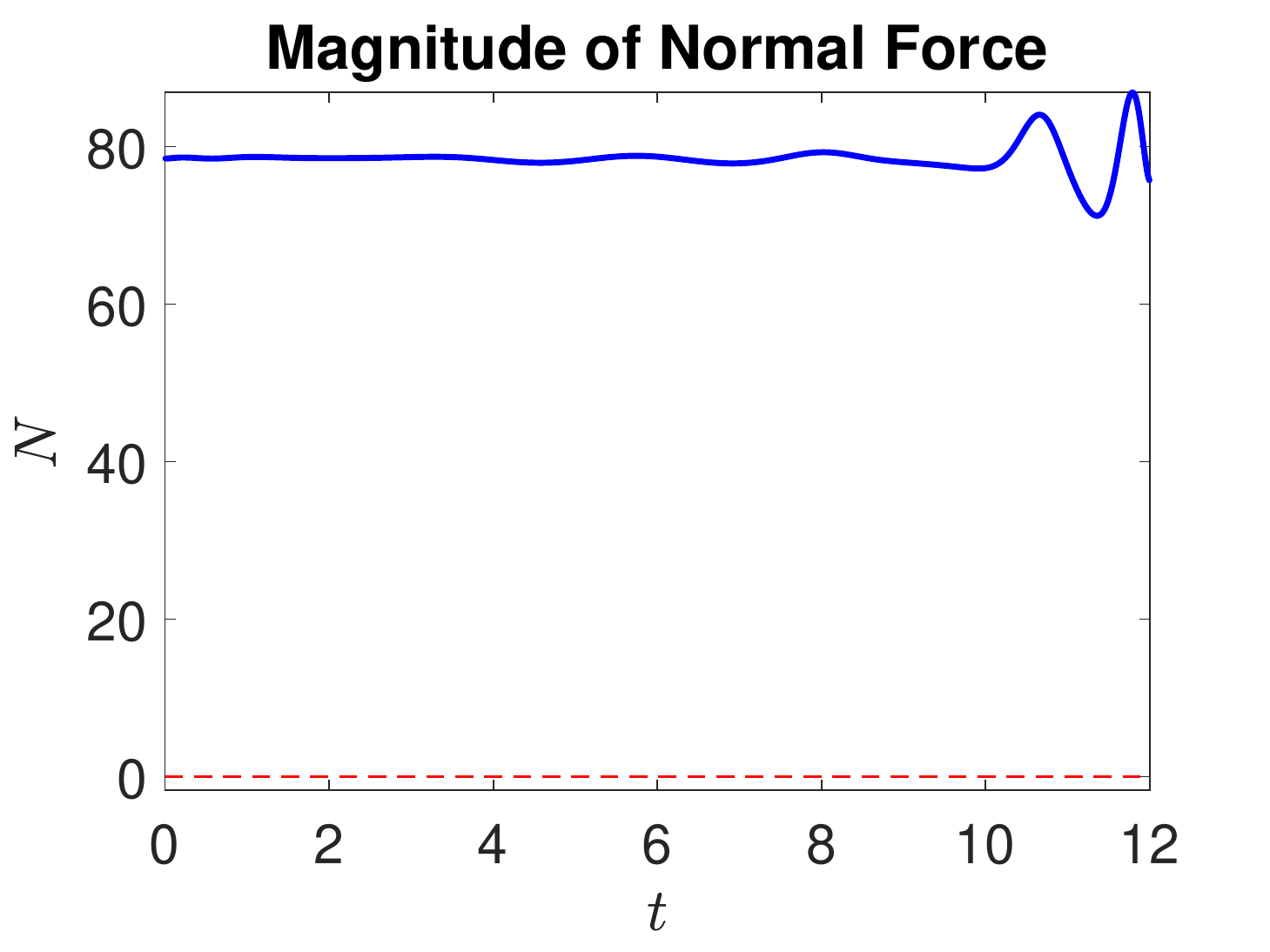}\label{fig_dsim8_pc_normal}}
	\\
	\subfloat[The minimum coefficient of static friction to prevent slipping is $.07799$ when $\alpha=.1$.]{\includegraphics[scale=.5]{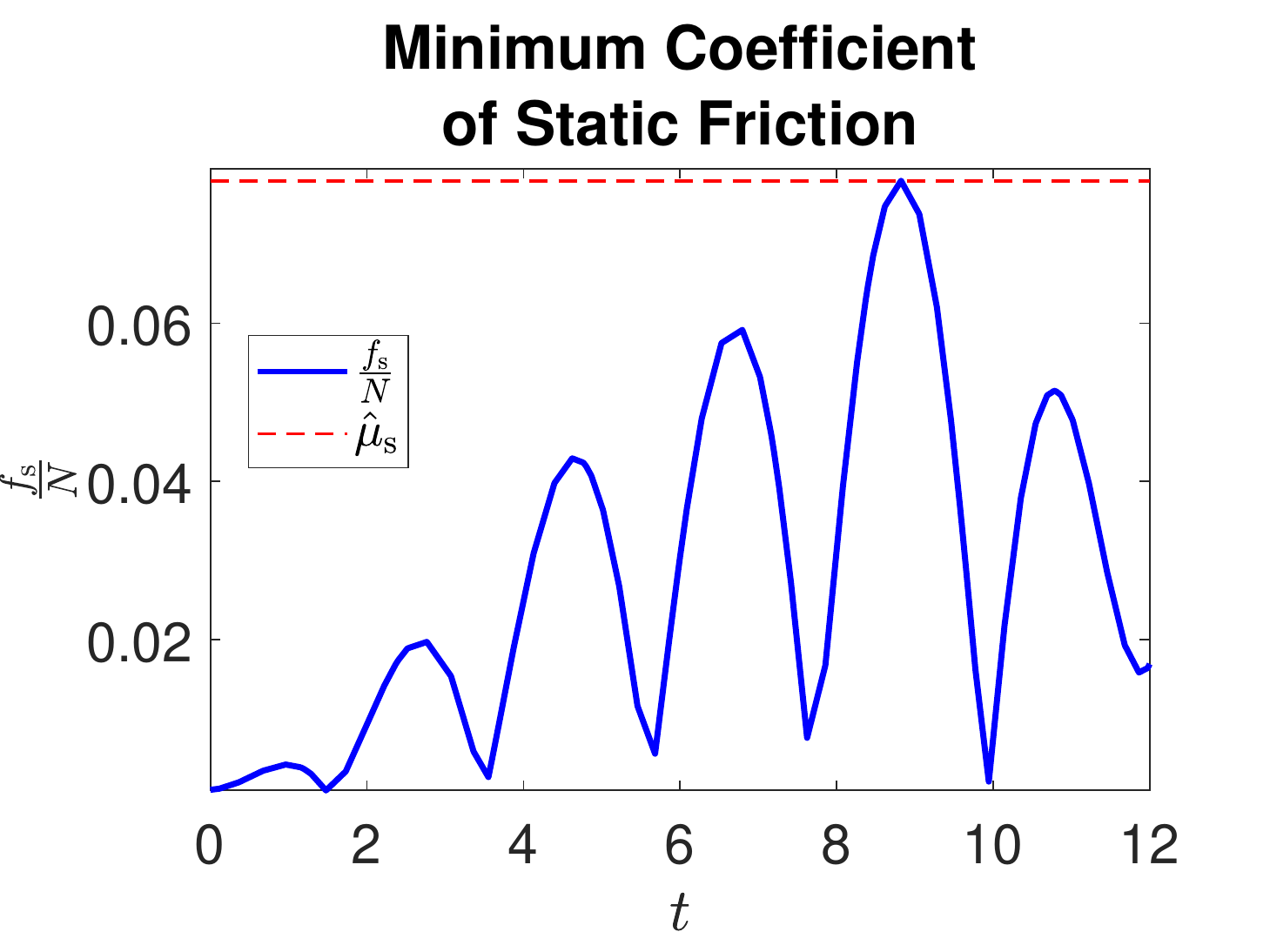}\label{fig_dsim8_dm_mu_s}}
	\hspace{5mm}
	\subfloat[The minimum coefficient of static friction to prevent slipping is $.3502$ when $\alpha \approx 272$.]{\includegraphics[scale=.5]{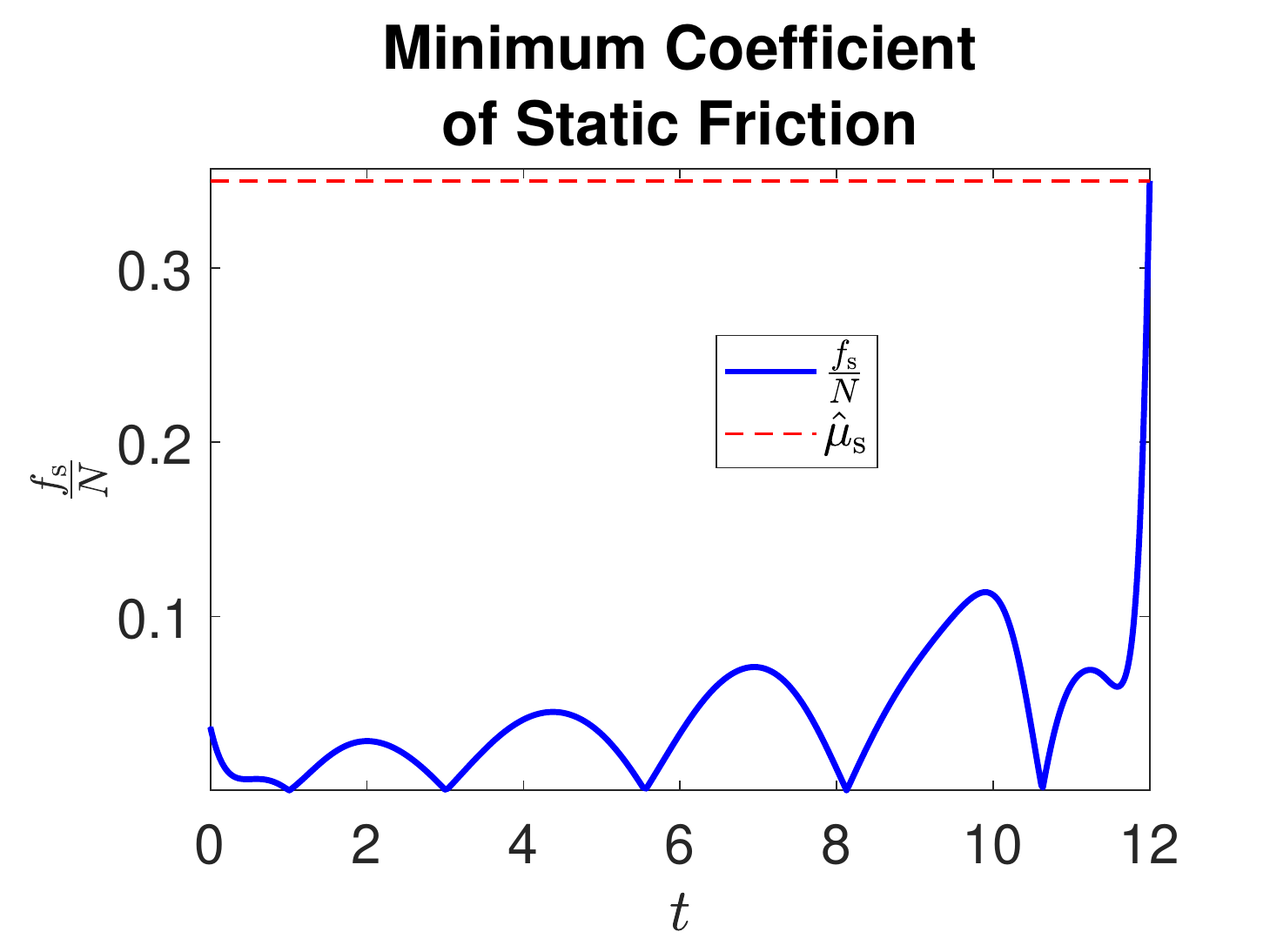}\label{fig_dsim8_pc_mu_s}}
	\caption{Numerical solutions of the rolling disk optimal control problem \eqref{dyn_opt_problem_disk} using $4$ control masses for $ \gamma_1=\gamma_2=\gamma_3=\gamma_4=.1$ and for fixed initial and final times. The direct method results for $\alpha=.1$ are shown in the left column, while the predictor-corrector continuation indirect method results for $\alpha \approx 272$ are shown in the right column. The disk does not detach from the surface since the magnitude of the normal force is always positive. The disk rolls without slipping if $\mu_\mathrm{s} \ge .07799$ for the direct method solution and if $\mu_\mathrm{s} \ge .3502$ for the indirect method solution. That is, the indirect method solution requires a much larger coefficient of static friction.}
	\label{fig_dsim8_normal}
\end{figure}

\begin{figure}[!ht] 
	\centering 
	\subfloat[Evolution of the continuation parameter $\mu$ from $.95$ down to $.8983$.]{\includegraphics[scale=.5]{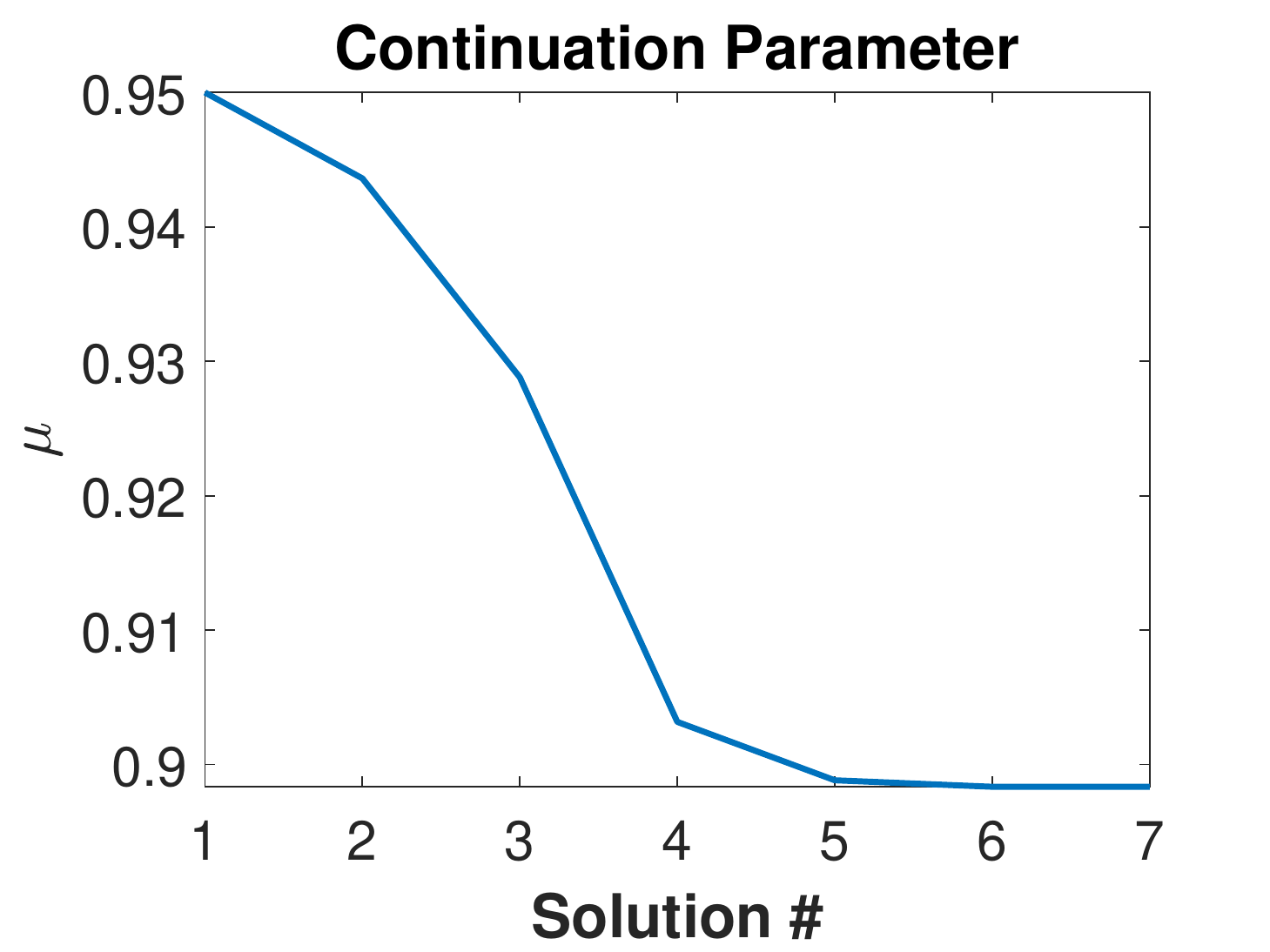}\label{fig_dsim8_pc_cont_param}}
	\hspace{5mm}
	\subfloat[Evolution of the GC path weighting factor $\alpha$ from $.1$ up to $272$.]{\includegraphics[scale=.5]{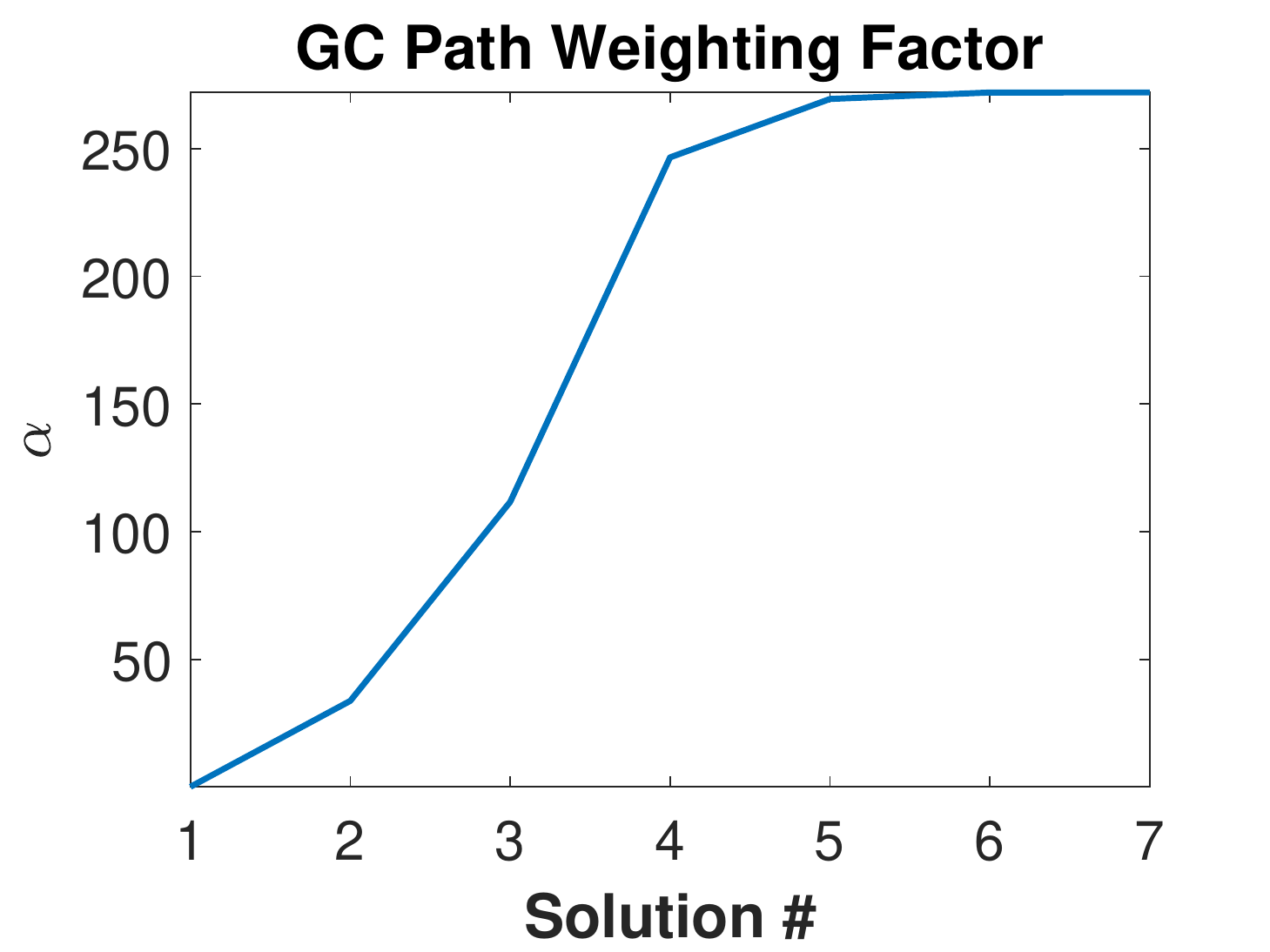}\label{fig_dsim8_pc_alpha}}
	\\
	\subfloat[Evolution of the performance index $J$.]{\includegraphics[scale=.5]{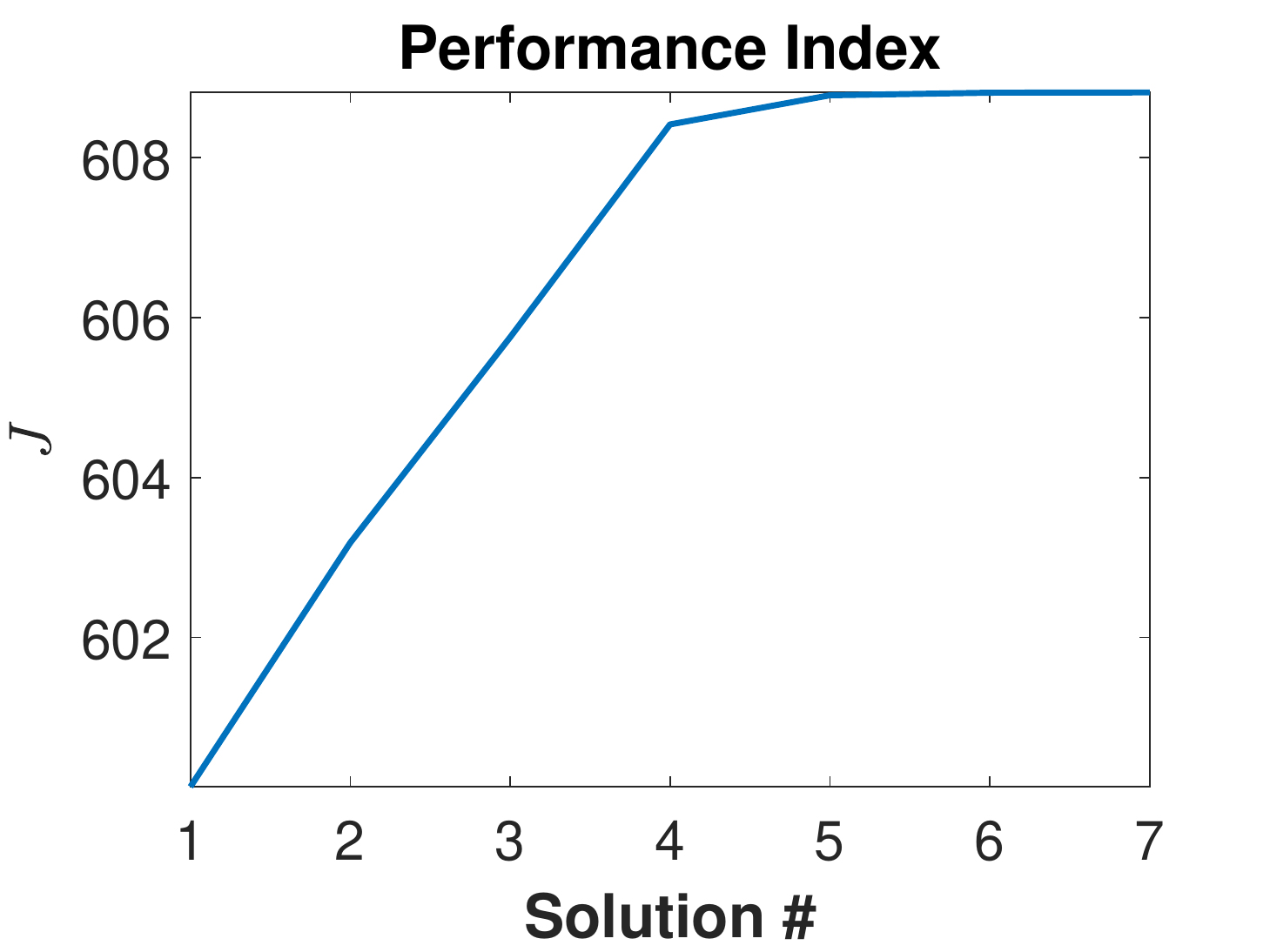}\label{fig_dsim8_pc_J}}
	\hspace{5mm}
	\subfloat[Evolution of the tangent steplength $\sigma$.]{\includegraphics[scale=.5]{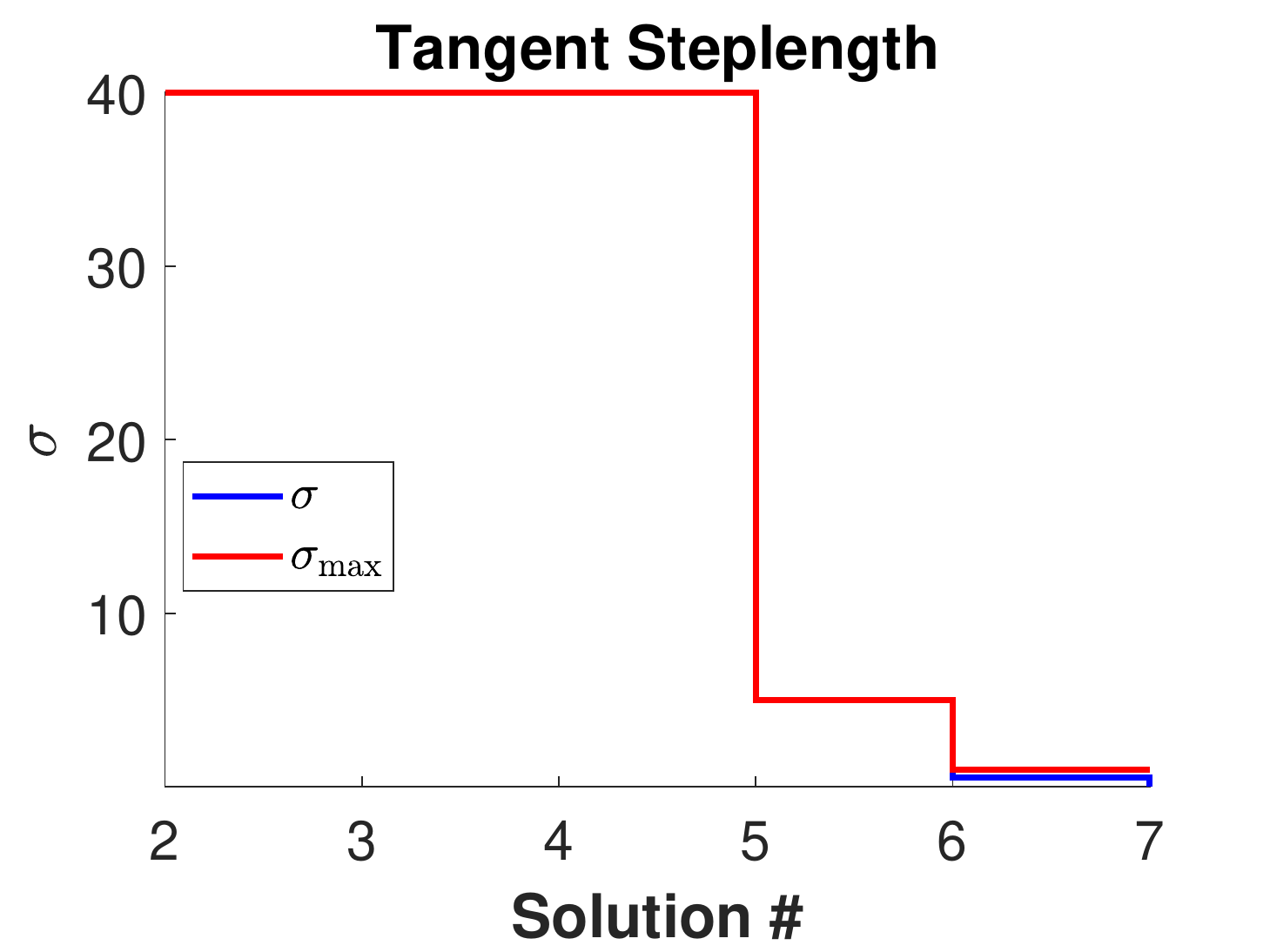}\label{fig_dsim8_pc_sigma}}
	\caption{Evolution of various parameters and variables during the predictor-corrector continuation indirect method, which starts from the direct method solution, used to solve the rolling disk optimal control problem \eqref{dyn_opt_problem_disk}. $\mu$ decreases monotonically, while $\alpha$ and $J$ increase monotonically.}
	\label{fig_dsim8_pc}
\end{figure}

For the direct method, \mcode{GPOPS-II} is run using the NLP solver SNOPT. The \mcode{GPOPS-II} mesh error tolerance is $1\mathrm{e}\unaryminus 6$ and the SNOPT error tolerance is $1\mathrm{e}\unaryminus 7$. In order to encourage convergence of SNOPT, a constant $C=50$ is added to the integrand cost function $L$ in \eqref{eq_disk_integrand_cost}. The sweep predictor-corrector continuation method discussed in Appendix~\ref{app_sweep_predictor_corrector} is used by the indirect method. For the sweep predictor-corrector continuation method, the maximum tangent steplength $\sigma_\mathrm{max}$ is adjusted according to Figure~\ref{fig_dsim8_pc_sigma} over the course of $6$ predictor-corrector steps, the maximum tangent steplength in each step is $\sigma_\mathrm{max}=\begin{bmatrix} 40 & 40 & 40 & 5 & 1 & 1 \end{bmatrix}$, the direction of the initial unit tangent is determined by setting $d=\unaryminus 2$ to force the continuation parameter $\mu$ to initially decrease, the relative error tolerance is $1\mathrm{e}\unaryminus 8$, the unit tangent solver is \mcode{twpbvpc_m}, and the monotonic ``sweep'' continuation solver is \mcode{acdcc}. The numerical results are shown in Figures~\ref{fig_dsim8}, \ref{fig_dsim8_normal}, and \ref{fig_dsim8_pc}. As $\mu$ decreases from $.95$ down to $.8983$ during continuation (see Figure~\ref{fig_dsim8_pc_cont_param}), $\alpha$ increases from $.1$ up to $272$ (see Figure~\ref{fig_dsim8_pc_alpha}). Since $\alpha$ is ratcheted up during continuation, thereby increasing the penalty in the integrand cost function \eqref{eq_disk_integrand_cost} for deviation between the disk's GC and $z_\mathrm{d}$, by the end of continuation, the disk's GC tracks $z_\mathrm{d}$ very accurately (compare Figures~\ref{fig_dsim8_dm_gc_path} vs \ref{fig_dsim8_pc_gc_path}), at the expense of more serpentine control mass trajectories (compare Figures~\ref{fig_dsim8_dm_cm_bf} vs \ref{fig_dsim8_pc_cm_bf}) and larger magnitude controls (compare Figures~\ref{fig_dsim8_dm_controls} vs \ref{fig_dsim8_pc_controls}). The disk does not detach from the surface since the magnitude of the normal force is always positive (see Figures~\ref{fig_dsim8_dm_normal} and \ref{fig_dsim8_pc_normal}). The disk rolls without slipping if the coefficient of static friction $\mu_\mathrm{s}$ is at least $\hat \mu_\mathrm{s} \approx .07799$ for the direct method solution (see Figure~\ref{fig_dsim8_dm_mu_s}) and if $\mu_\mathrm{s}$ is at least $\hat \mu_\mathrm{s} \approx .3502$ for the indirect method solution (see Figure~\ref{fig_dsim8_pc_mu_s}). That is, the indirect method solution requires a much larger coefficient of static friction. 

\subsubsection{Turning Points} \label{sssec_disk_sim_tp}
\begin{figure}[!ht] 
	\centering 
	\subfloat[Evolution of the continuation parameter $\mu$.]{\includegraphics[scale=.5]{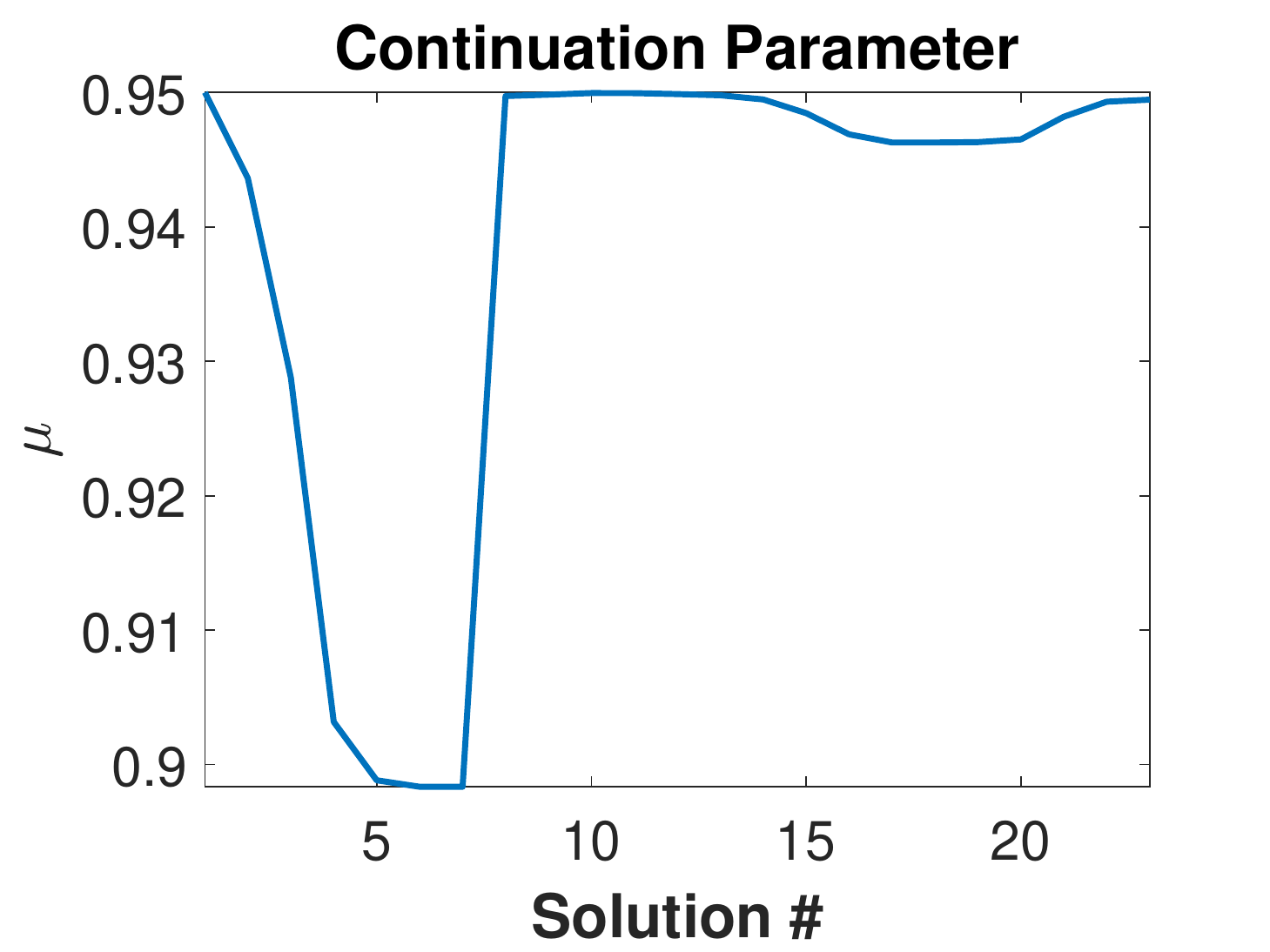}\label{fig_dsim8_tp_pc_cont_param}}
	\hspace{5mm}
	\subfloat[Evolution of the GC path weighting factor $\alpha$.]{\includegraphics[scale=.5]{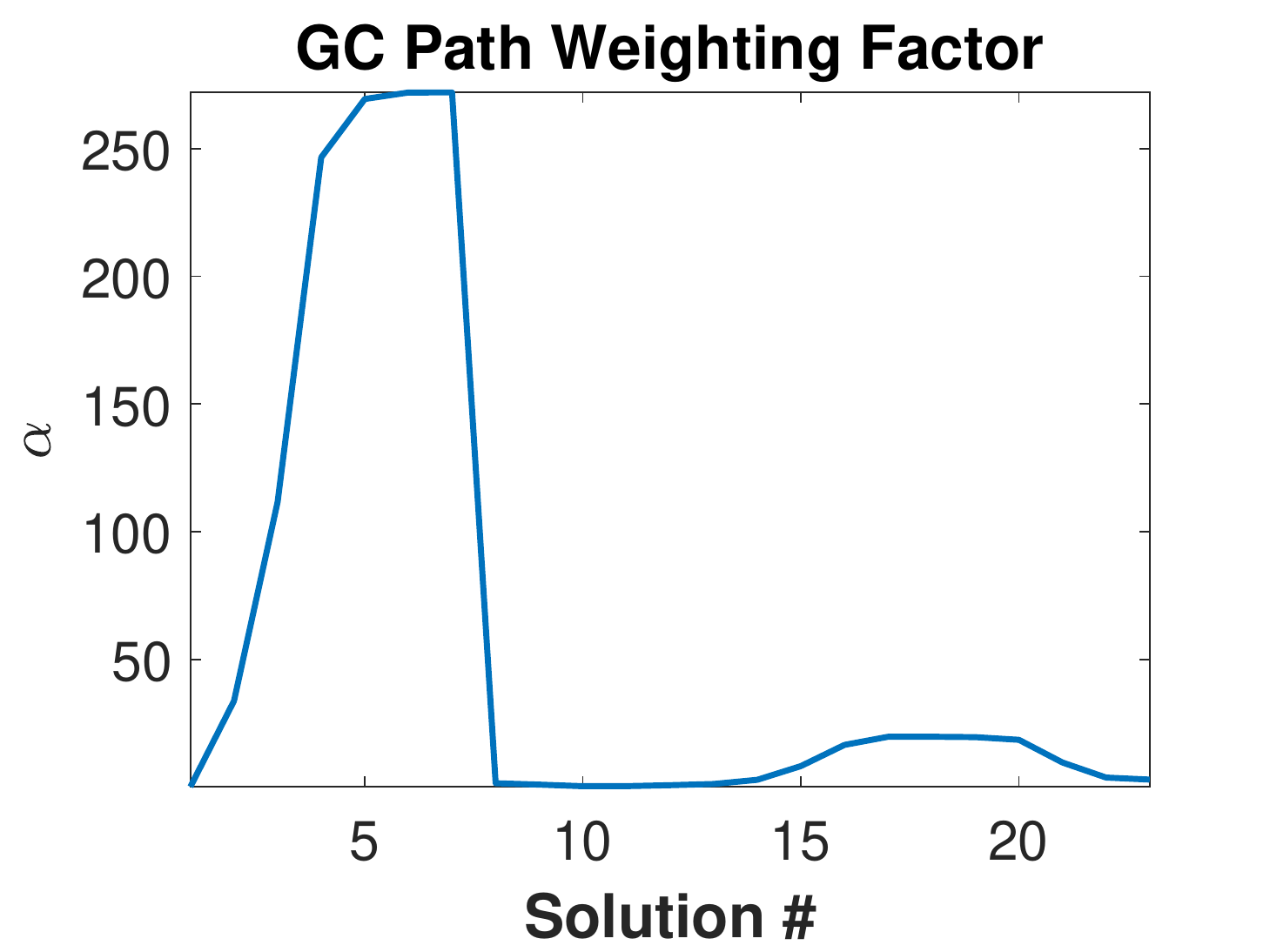}\label{fig_dsim8_tp_pc_alpha}}
	\\
	\subfloat[Evolution of the performance index $J$.]{\includegraphics[scale=.5]{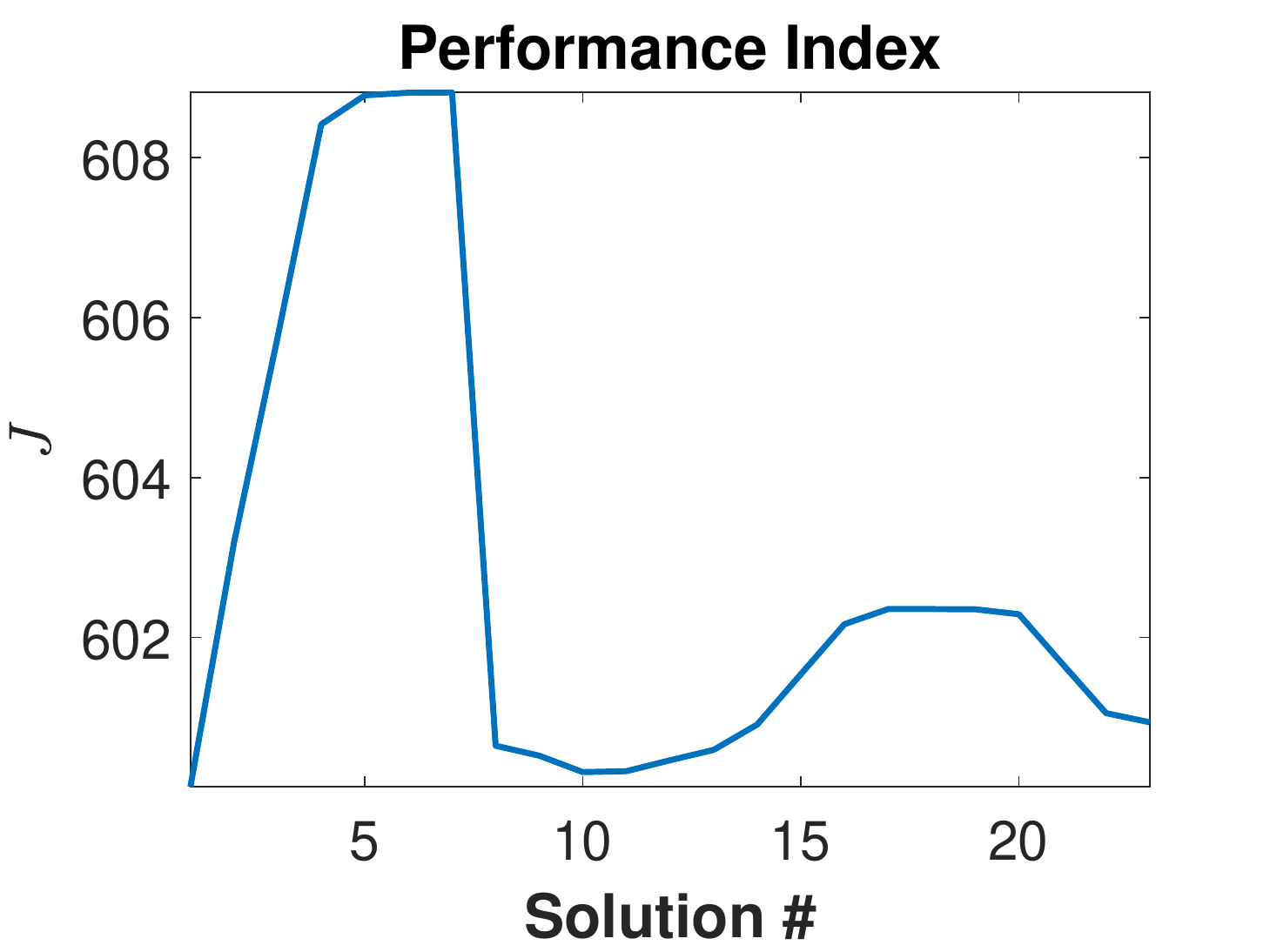}\label{fig_dsim8_tp_pc_J}}
	\hspace{5mm}
	\subfloat[Evolution of the GC tracking error, measured by the square of the $L^2$ norm.]{\includegraphics[scale=.5]{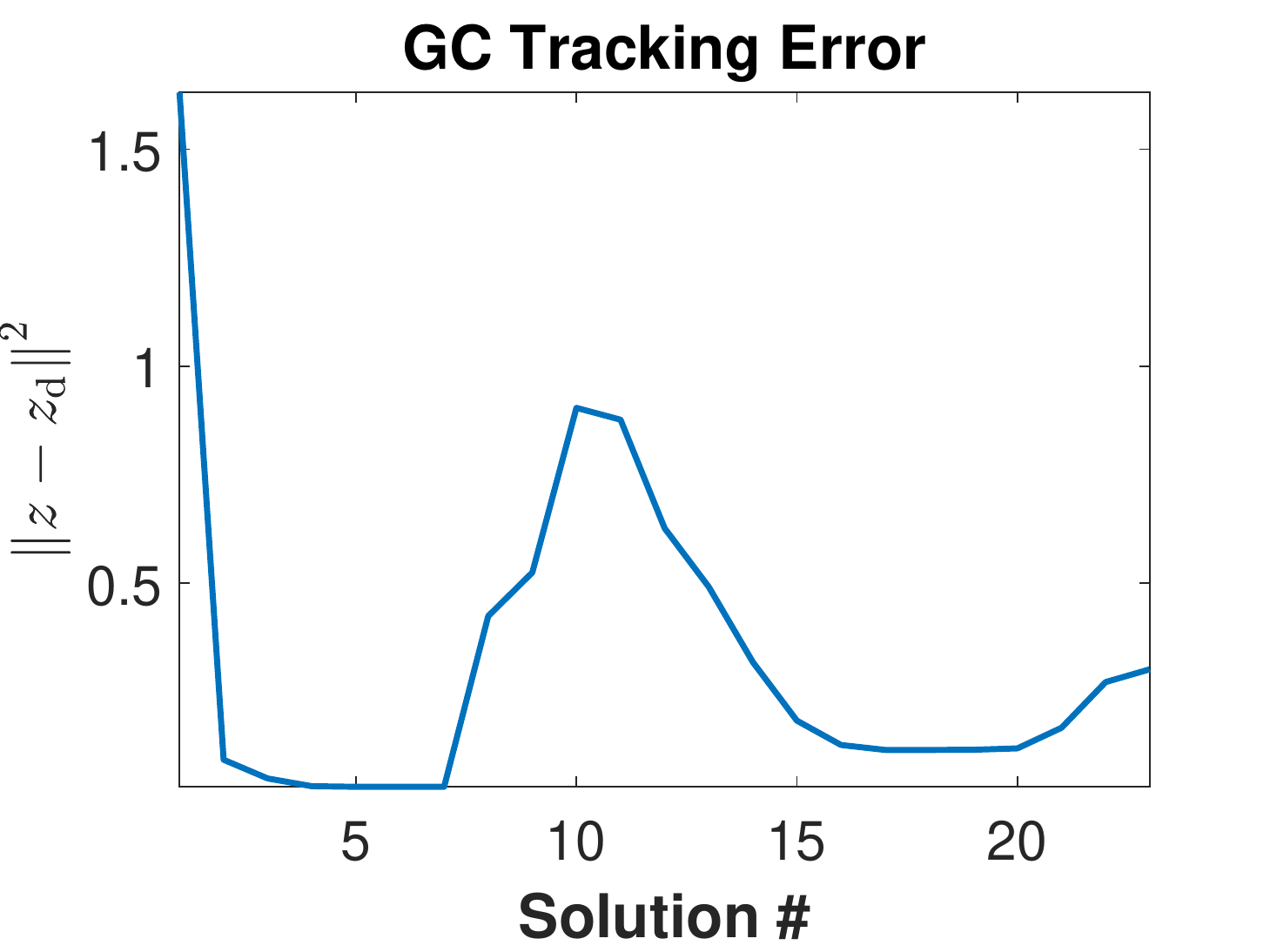}\label{fig_dsim8_tp_pc_trk_err}}
	\\
	\subfloat[Evolution of the tangent steplength $\sigma$.]{\includegraphics[scale=.5]{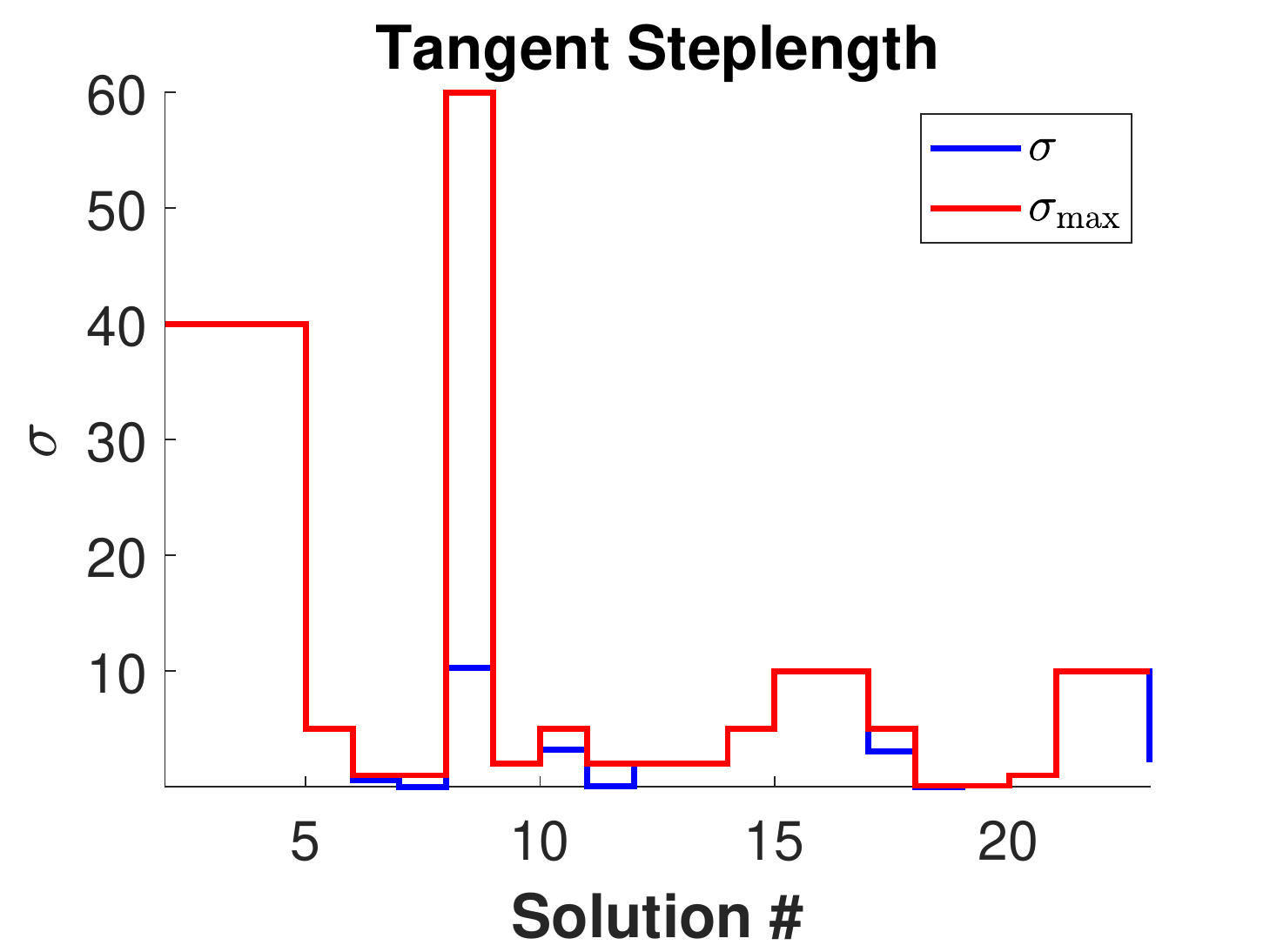}\label{fig_dsim8_tp_pc_sigma}}
	\caption{Evolution of various parameters and variables during an extended run of the predictor-corrector continuation indirect method, which starts from the direct method solution, used to solve the rolling disk optimal control problem \eqref{dyn_opt_problem_disk}. Note the turning points at solutions 7, 10, and 18. The minimum of the GC tracking error occurs at solution 7.}
	\label{fig_dsim8_tp_pc}
\end{figure}

The previous sweep predictor-corrector continuation indirect method is repeated, but this time with 22 steps, where the maximum tangent steplength in each step is 
\begin{equation} \label{eq_tp_sigma_max}
\sigma_\mathrm{max}=\begin{bmatrix} 40 & 40 & 40 & 5 & 1 & 1 & 60 & 2 & 5 & 2 & 2 & 2 & 5 & 10 & 10 & 5 & .1 & .1 & 1 & 10 & 10 & 10 \end{bmatrix}.
\end{equation}
Note that the first 6 maximum tangent steplengths in \eqref{eq_tp_sigma_max} agree with those used in the previous simulation, so that the two simulations agree for the first 7 solutions. Figure~\ref{fig_dsim8_tp_pc} shows the evolution of the continuation parameter $\mu$, GC path weighting factor $\alpha$, performance index $J$, GC tracking error $\left \Vert z- z_\mathrm{d} \right \Vert ^2$, and tangent steplength $\sigma$ over the course of the 23 solutions (the first solution is initialized by the direct method) constructed by the sweep predictor-corrector continuation indirect method. Note the turning points (local maxima or minima) at solutions 7, 10, and 18 in Figures~\ref{fig_dsim8_tp_pc_cont_param}-\ref{fig_dsim8_tp_pc_trk_err}. Figure~\ref{fig_dsim8_tp_pc_trk_err} shows that the GC tracking error realizes a minimum at solution 7, which explains how the stopping point for the previous simulation was selected. 

\section{Obstacle Avoidance for the Rolling Ball} \label{sec_ball_sim}

\subsection{Optimal Control Problem and Controlled Equations of Motion} \label{ssec_ball_controlled}

\begin{figure}[h] 
	\centering
	\includegraphics[scale=.5]{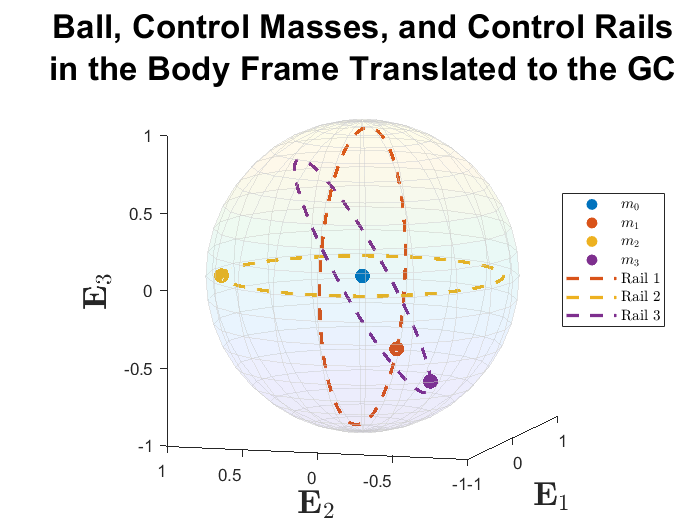}
	\caption{The ball of radius $r=1$ actuated by $3$ control masses, $m_1$, $m_2$, and $m_3$, each on its own circular control rail. The control rail radii are $r_1=.95$, $r_2=.9$, and $r_3=.85$. The location of the ball's CM is denoted by $m_0$.}
	\label{fig_bsim10_control_masses_rails}
\end{figure}

In the next two subsections, numerical solutions of the controlled equations of motion for the rolling ball are presented, where the goal is to move the ball between a pair of points while the ball's GC avoids a pair of obstacles. \revisionNACO{R1Q5}{A rolling ball of mass $m_0=4$, radius $r=1$, principal moments of inertia $d_1=d_2=d_3=1$, and with the CM coinciding with the GC (i.e. $\bzeta_0=\mathbf{0}$) is simulated. These physical parameters for the ball are consistent with the necessary and sufficient conditions stipulated by Inequalities 3.1 and 3.2 in \cite{rozenblat2016choice}.} There are $n=3$ control masses, each of mass $1$ so that $m_1=m_2=m_3=1$, located on circular control rails centered on the GC of radii $r_1=.95$, $r_2=.9$, and $r_3=.85$, oriented as \revision{R1Q18}{shown in Figure~\ref{fig_bsim10_control_masses_rails}. For $1 \le i \le 3$, the position of $m_i$ in the body frame centered on the GC is
\begin{equation}
\bzeta_i\left(\theta_i\right) = r_i \mathcal{B}_i \left(\boldsymbol{\varsigma} \left(\mathbf{v}_i\right) \right) \begin{bmatrix} \cos \theta_i \\ 0 \\ \sin \theta_i  \end{bmatrix},
\end{equation}
where $\mathcal{B}_i \left(\mathbf{n} \right) \in SO(3)$ is a rotation matrix whose columns are the right-handed orthonormal basis constructed from the unit vector $\mathbf{n} \in \mathbb{R}^3$ based on the algorithm given in Section 4 and Listing 2 of \cite{frisvad2012building}, $\boldsymbol{\varsigma} \colon \mathbb{R}^3 \to \mathbb{R}^3$ maps spherical coordinates to Cartesian coordinates:
\begin{equation}
\boldsymbol{\varsigma}\left(\begin{bmatrix} \phi \\ \theta \\ \rho \end{bmatrix} \right) = \begin{bmatrix} \rho \cos \theta \cos \phi  \\ \rho \cos \theta \sin \phi \\ \rho \sin \theta \end{bmatrix},
\end{equation}
and
\begin{equation}
\mathbf{v}_1 = \begin{bmatrix} 0 & 0 & 1 \end{bmatrix}^\mathsf{T}, \quad \mathbf{v}_2 = \begin{bmatrix} \frac{\pi}{2} & 0 & 1 \end{bmatrix}^\mathsf{T}, \quad \mathrm{and} \quad \mathbf{v}_3 = \begin{bmatrix} \frac{\pi}{4} & \frac{\pi}{4} & 1 \end{bmatrix}^\mathsf{T}
\end{equation}
are spherical coordinates of unit vectors in $\mathbb{R}^3$.} 

\begin{figure}[!ht] 
	\centering
	\subfloat[The GC plows through the obstacles when the obstacle heights $h_1=h_2$ are $0$.]{\includegraphics[scale=.5]{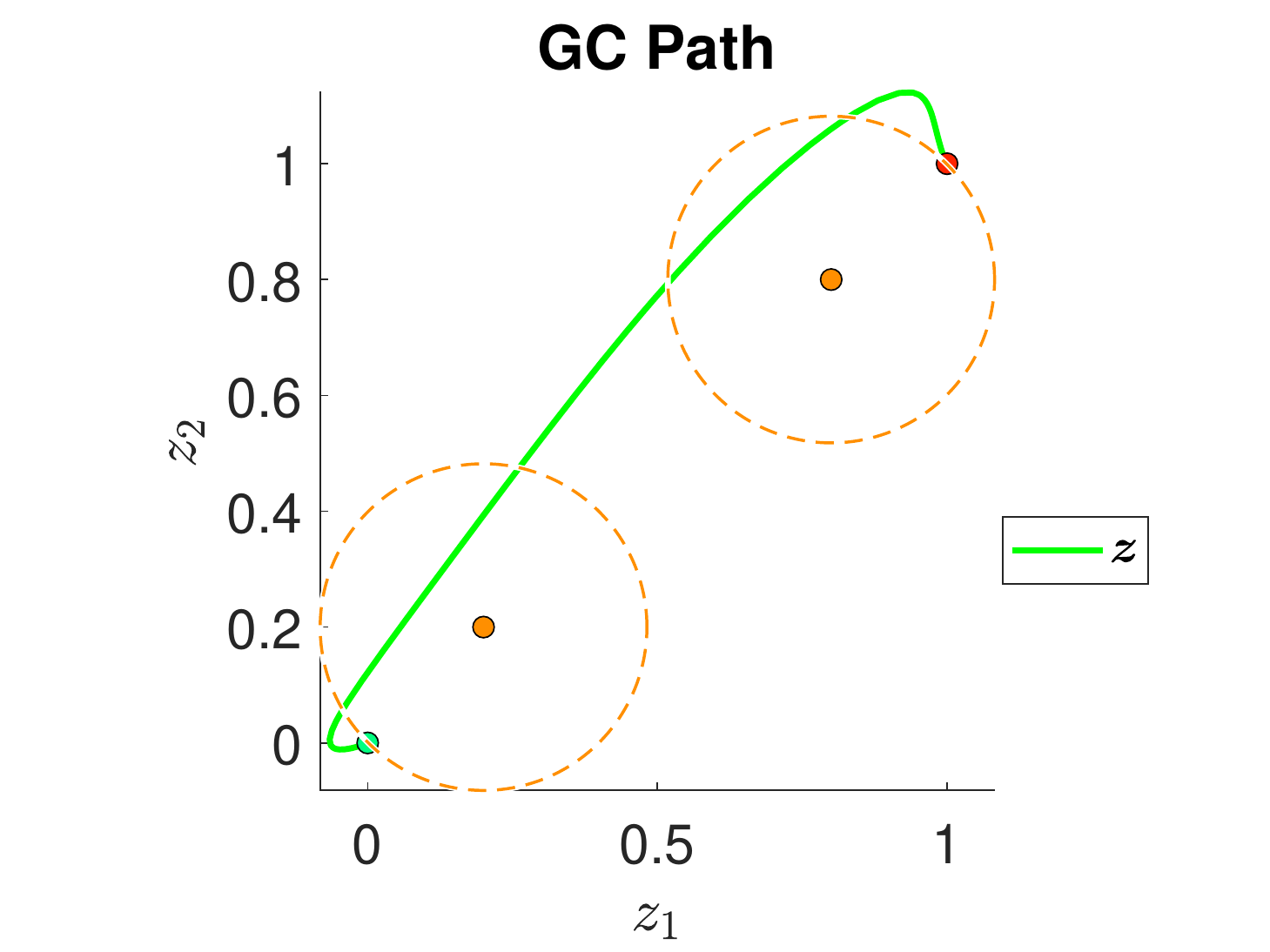}\label{fig_bsim10_dm_gc_path}}
	\hspace{5mm}
	\subfloat[The GC veers around the obstacles when the obstacle heights $h_1=h_2$ are $9.93$.]{\includegraphics[scale=.5]{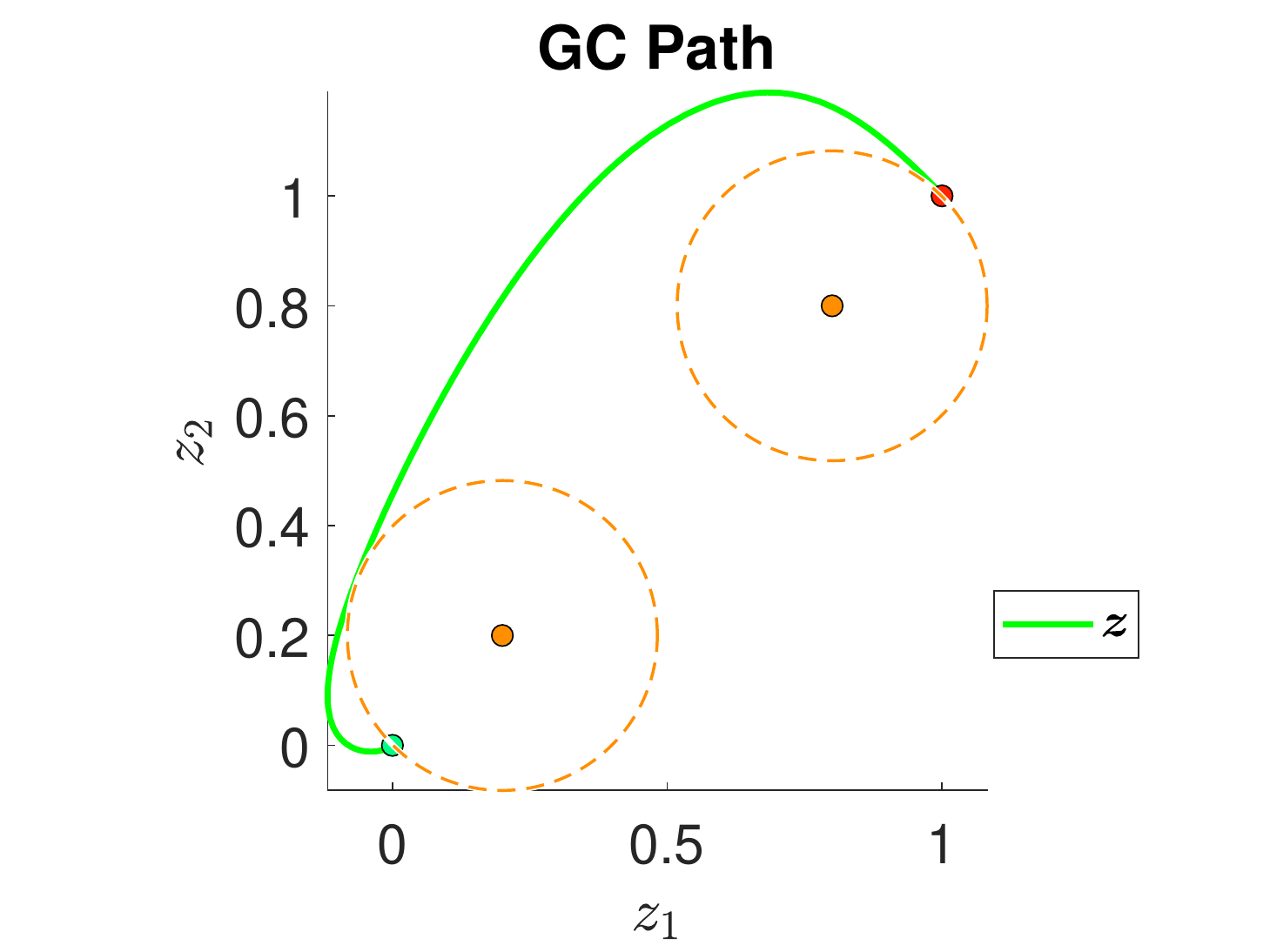}\label{fig_bsim10_pc_gc_path}}
	\\
	\subfloat[Motion of the center of masses when the obstacle heights $h_1=h_2$ are $0$.]{\includegraphics[scale=.5]{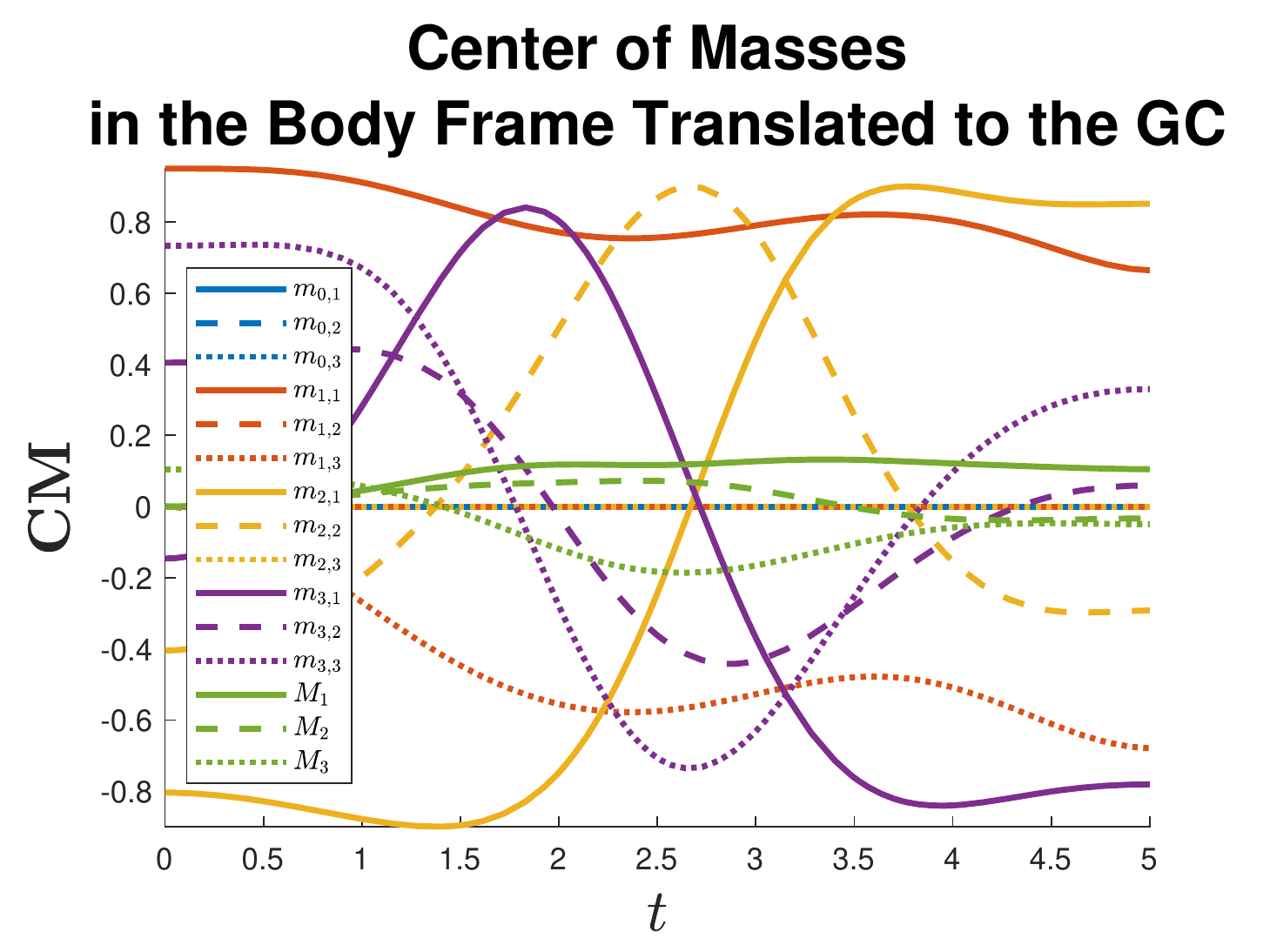}\label{fig_bsim10_dm_cm_bf}}
	\hspace{5mm}
	\subfloat[Motion of the center of masses when the obstacle heights $h_1=h_2$ are $9.93$.]{\includegraphics[scale=.5]{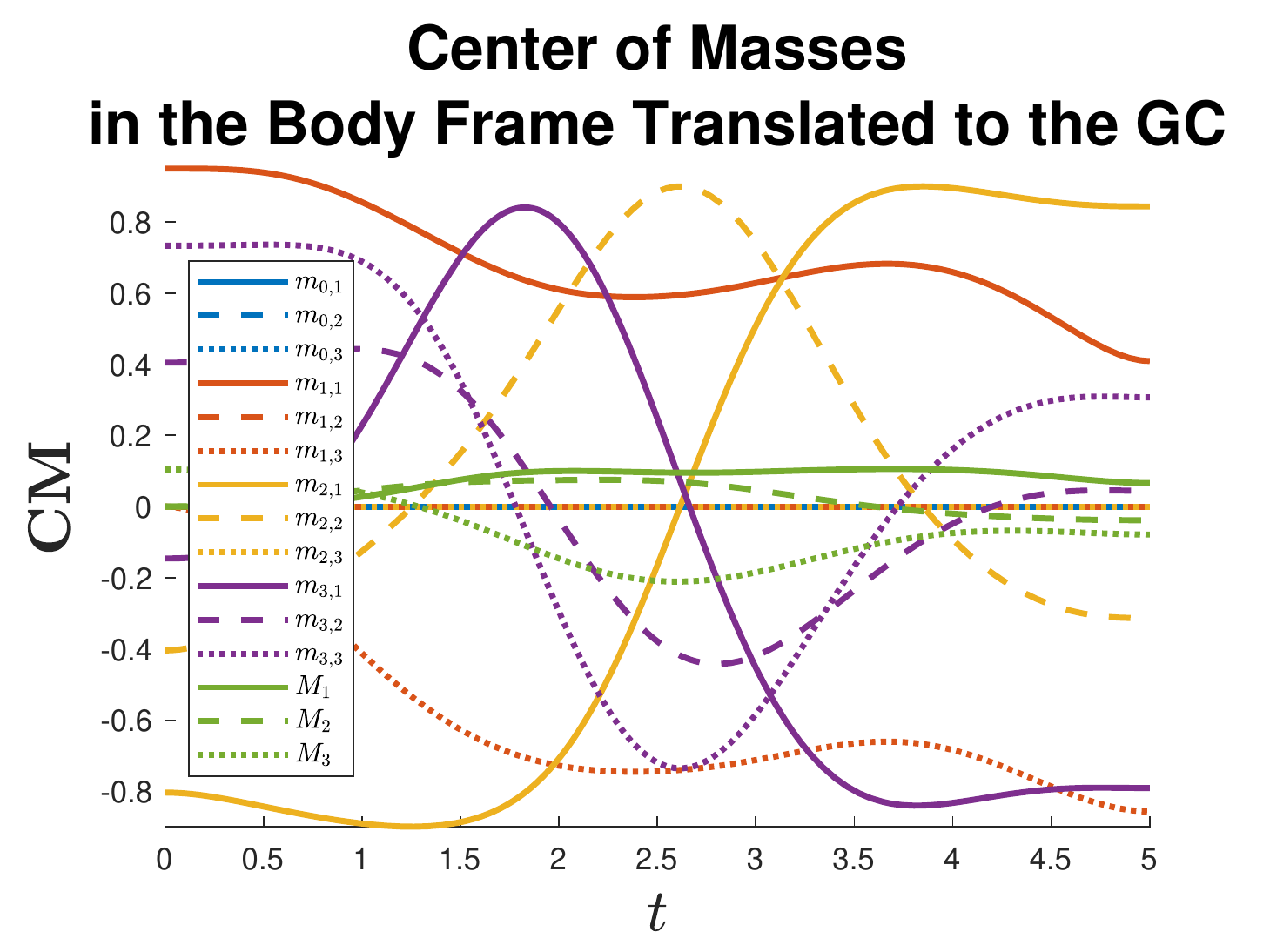}\label{fig_bsim10_pc_cm_bf}}
	\\
	\subfloat[The controls when the obstacle heights $h_1=h_2$ are $0$.]{\includegraphics[scale=.5]{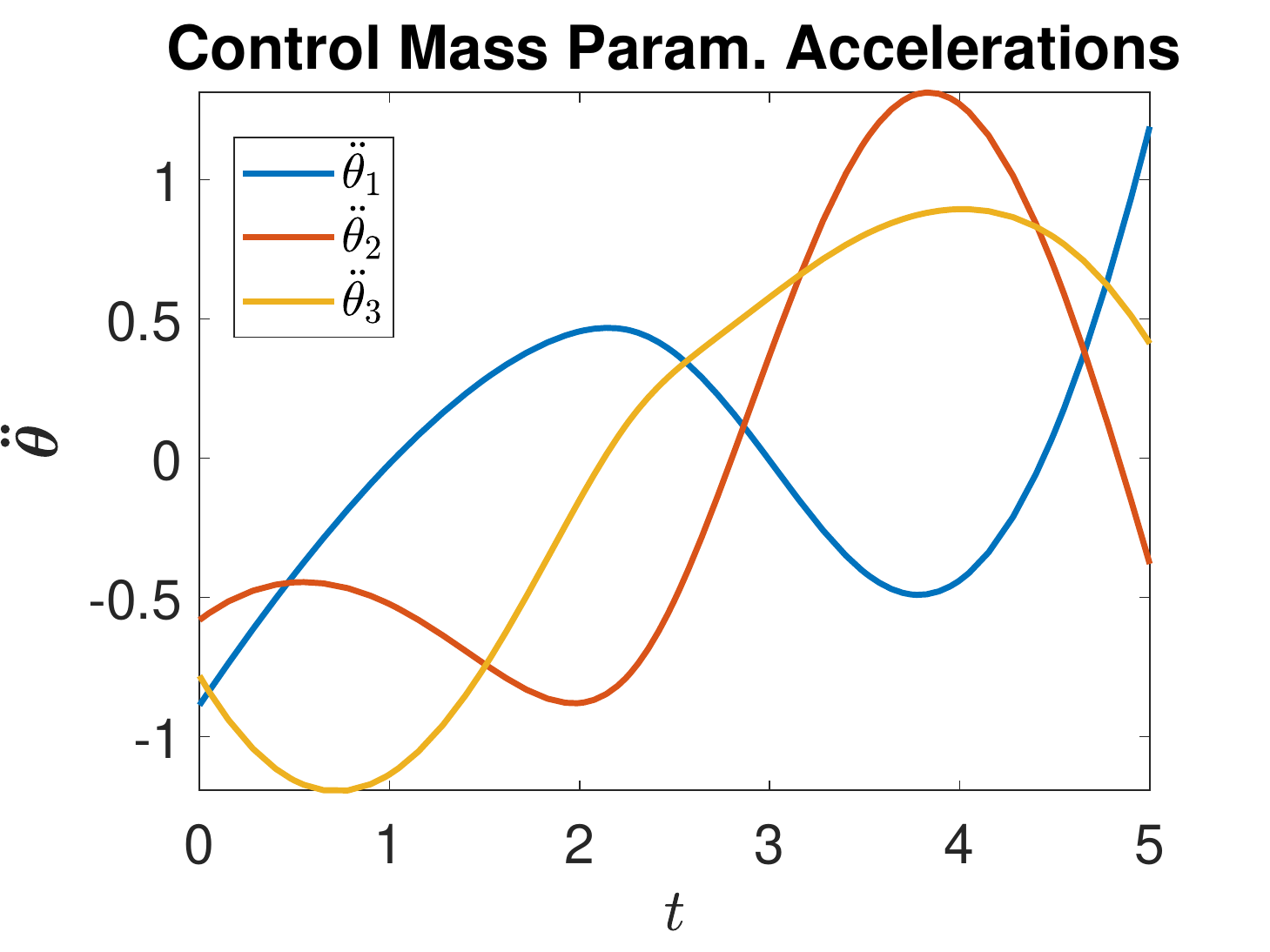}\label{fig_bsim10_dm_controls}}
	\hspace{5mm}
	\subfloat[The controls increase slightly in magnitude when the obstacle heights $h_1=h_2$ are  $9.93$.]{\includegraphics[scale=.5]{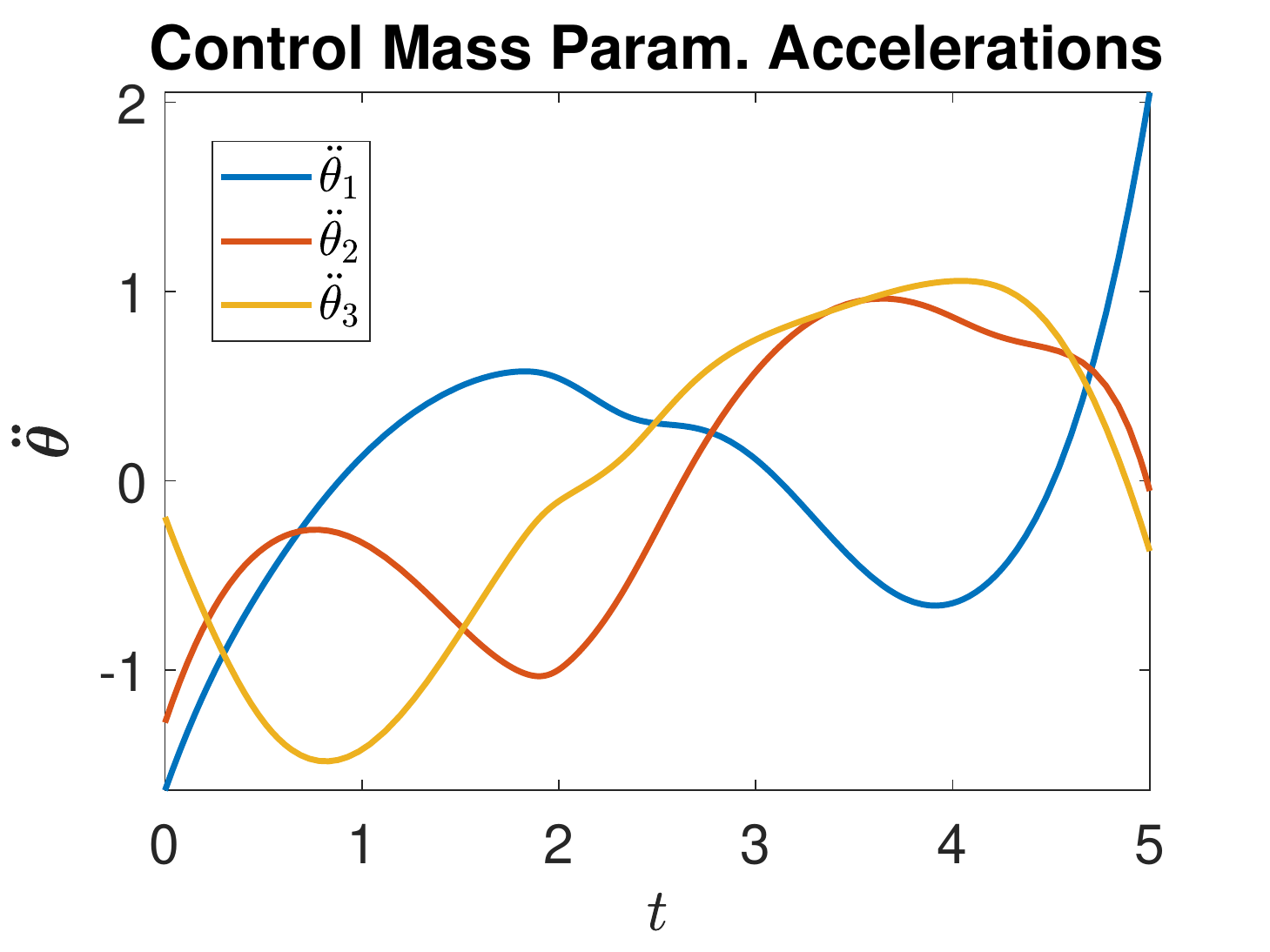}\label{fig_bsim10_pc_controls}}
	\caption{Numerical solutions of the rolling ball optimal control problem \eqref{dyn_opt_problem_ball} for sigmoid obstacle avoidance using $3$ control masses for $\gamma_1=\gamma_2=\gamma_3=10$ and for fixed initial and final times. The obstacle centers are located at $\bv_1=\protect\begin{bmatrix}  v_{1,1} & v_{1,2} \protect\end{bmatrix}^\mathsf{T} = \protect\begin{bmatrix}  .2 & .2 \protect\end{bmatrix}^\mathsf{T}$ and $\bv_2 = \protect\begin{bmatrix}  v_{2,1} & v_{2,2} \protect\end{bmatrix}^\mathsf{T} = \protect\begin{bmatrix}  .8 & .8 \protect\end{bmatrix}^\mathsf{T}$ and the obstacle radii are $\rho_1=\rho_2=.282$. The direct method results for obstacle heights at $h_1=h_2=0$ are shown in the left column, while the predictor-corrector continuation indirect method results for obstacle heights at $h_1=h_2 \approx 9.93$ are shown in the right column.}
	\label{fig_bsim10}
\end{figure}

The total mass of the system  is \revisionNACO{R1Q5}{$M=7$}, and gravity is \revisionNACO{R1Q5}{$g=9.81$}. There is no external force acting on the \revision{R1Q17}{ball}'s GC, so that $\tilde \bGamma \equiv \Lambda^{-1} \mathbf{F}_\mathrm{e}=\mathbf{0}$ in \eqref{uncon_ball_eqns_explicit_1d}. The initial time is fixed to $a=0$ and the final time is fixed to \revisionNACO{R1Q5}{$b=5$}. The ball's GC starts at rest at $\bz_a=\begin{bmatrix}  0 & 0 \end{bmatrix}^\mathsf{T}$ at time $a=0$ and stops at rest at $\bz_b=\begin{bmatrix}  1 & 1 \end{bmatrix}^\mathsf{T}$ at time \revisionNACO{R1Q5}{$b=5$}. The ball's GC should avoid $2$ circular obstacles, depicted in Figures~\ref{fig_bsim10_dm_gc_path} and \ref{fig_bsim10_pc_gc_path}. The obstacles each have radius $\rho_1 = \rho_2 = .282$ and are centered at $\bv_1 = \begin{bmatrix}  .2 & .2 \end{bmatrix}^\mathsf{T}$ and $\bv_2 = \begin{bmatrix}  .8 & .8 \end{bmatrix}^\mathsf{T}$.

The ODE formulation of the optimal control problem for the rolling ball is
\begin{equation}
\min_{\bu} J 
\mbox{\, s.t. \,}
\left\{
\begin{array}{ll}
\dot {\bx} = \mathbf{f}\left(\bx,\bu\right), \\
\bsigma\left(\bx(a)\right) = \mathbf{0},\\
\bpsi\left(\bx(b)\right) = \mathbf{0}.
\end{array}
\right.
\label{dyn_opt_problem_ball}
\end{equation}
In \eqref{dyn_opt_problem_ball}, the system state $\bx$ and control $\bu$ are
\begin{equation} \label{eq_ball_state_control}
\bx \equiv \begin{bmatrix} \btheta \\ \dot \btheta \\ \mathfrak{q} \\ \bOm \\ \bz  \end{bmatrix} \quad \mathrm{and} \quad \bu \equiv \ddot \btheta,
\end{equation}
where $\btheta, \, \dot \btheta, \, \ddot \btheta \in \mathbb{R}^3$, $\mathfrak{q} \in \mathscr{S} \cong \mathbb{S}^3 \subset \mathbb{R}^4$, $\bOm \in \mathbb{R}^3$, and $\bz \in \mathbb{R}^2$. In the system state defined in \eqref{eq_ball_state_control}, the rolling ball's orientation matrix $\Lambda \in SO(3)$ is parameterized by $\mathfrak{q} \in \mathscr{S}$, where $\mathscr{S}$ denotes the set of versors (i.e. unit quaternions) \cite{Ho2011_pII,graf2008quaternions,stevens2015aircraft,baraff2001physically}. The properties of versors and the notation used to manipulate versors are explained in Appendix D of \cite{Putkaradze2018dynamicsP}. \revision{R1Q14}{Recall from \cite{Putkaradze2018dynamicsP} that given a column vector $\bv \in \mathbb{R}^3$, $\bv^\sharp$ is the quaternion 
	\begin{equation}
	\bv^\sharp = \begin{bmatrix}0 \\ \bv \end{bmatrix},
	\end{equation}
	and given a quaternion $\mathfrak{p} \in \mathbb{H}$, $\mathfrak{p}^\flat \in \mathbb{R}^3$ is the column vector such that 
	\begin{equation}
	\mathfrak{p} = \begin{bmatrix} p_0 \\ \mathfrak{p}^\flat \end{bmatrix}.
	\end{equation} } 
\revisionNACO{R1Q3}{As explained in \cite{Putkaradze2018dynamicsP},  the transformation of a body frame vector $\bY \in \mathbb{R}^3$ into the spatial frame by the ball's orientation matrix $\Lambda \in SO(3)$ can be realized using the versor $\mathfrak{q} \in \mathscr{S}$ via the Euler-Rodrigues formula
\begin{equation} \label{eq_euler_rod}
\Lambda \bY = \left[\mathfrak{q} \bY^\sharp \mathfrak{q}^{-1} \right]^\flat.
\end{equation}}
In \eqref{dyn_opt_problem_ball}, the system dynamics defined for $a \le t \le b$ are
\begin{equation} \label{rolling_ball_opt_con_f}
\dot {\bx} = \begin{bmatrix} \dot \btheta \\ \ddot \btheta \\ \dot {\mathfrak{q}} \\ \dot \bOm \\ \dot \bz  \end{bmatrix}  = \mathbf{f}\left(\bx,\bu\right) \equiv \begin{bmatrix} \dot \btheta \\ \bu  \\ \frac{1}{2} \mathfrak{q} \bOm^\sharp \\ \bkappa\left(\bx,\bu \right) \\ \left( \left[\mathfrak{q} \bOm^\sharp \mathfrak{q}^{-1} \right]^\flat \times r \mathbf{e}_3  \right)_{12}  \end{bmatrix},
\end{equation}
where $\bkappa\left(\bx,\bu \right)$ is given by the right-hand side of the formula for $\dot \bOm$ in \eqref{uncon_ball_eqns_explicit_1d}. In \eqref{rolling_ball_opt_con_f}, the time-dependence of $\bkappa$ is dropped since  $\tilde \bGamma \equiv \Lambda^{-1} \mathbf{F}_\mathrm{e}=\mathbf{0}$ in \eqref{uncon_ball_eqns_explicit_1d} for these simulations.  In \eqref{dyn_opt_problem_ball}, the prescribed initial conditions at time $t=a$ are
\begin{equation} \label{eq_ball_initial_conds}
\bsigma\left(\bx(a)\right) \equiv \begin{bmatrix} \btheta(a) - \btheta_a \\ \dot \btheta(a) - {\dot \btheta}_a \\ \mathfrak{q}(a)-\mathfrak{q}_a \\ \bOm(a) - \bOm_a \\ \bz(a)- \bz_a \end{bmatrix} = \mathbf{0},
\end{equation}
and the prescribed final conditions at time $t=b$ are
\begin{equation} \label{eq_ball_final_conds}
\bpsi\left(\bx(b)\right) \equiv \begin{bmatrix} \bPi \left(\left[\mathfrak{q}(b) \left[ \frac{1}{M} \sum_{i=0}^3 m_i \bzeta_i\left(\theta_i(b)\right) \right]^\sharp \mathfrak{q}(b)^{-1} \right]^\flat \right) \\ \dot \btheta(b) - {\dot \btheta}_b \\ \bOm(b) - \bOm_b \\ \bz(b)- \bz_b  \end{bmatrix} = \mathbf{0}.
\end{equation}
Table~\ref{table_ball_ICs} shows the parameter values used in the rolling ball's initial conditions \eqref{eq_ball_initial_conds}. The initial configurations of the control masses are selected so that the total system CM in the spatial frame is initially located above or below the ball's GC. Hence, in conjunction with the other initial condition parameter values given in Table~\ref{table_ball_ICs}, the ball starts at rest. In \eqref{eq_ball_final_conds}, 
\begin{equation} \label{eq_ball_stab_fin_prep}
\frac{1}{M} \sum_{i=0}^3 m_i \bzeta_i\left(\theta_i(t)\right)
\end{equation}
is the total system CM in the body frame translated to the ball's GC at time $t$,
\revisionNACO{R1Q3}{
\begin{equation} \label{eq_ball_stab_fin}
\left[\mathfrak{q}(t) \left[ \frac{1}{M} \sum_{i=0}^3 m_i \bzeta_i\left(\theta_i(t)\right) \right]^\sharp \mathfrak{q}(t)^{-1} \right]^\flat = \Lambda(t) \left[ \frac{1}{M} \sum_{i=0}^3 m_i \bzeta_i\left(\theta_i(t)\right) \right]
\end{equation}} 
is the total system CM in the spatial frame translated to the ball's GC at time $t$, and $\bPi$ is the projection onto the first two components. Therefore, the first two constraints in \eqref{eq_ball_final_conds} ensure that the total system CM in the spatial frame is above or below the ball's GC  at the final time $t=b$. Hence, in conjunction with the final condition parameter values given in Table~\ref{table_ball_FCs}, the ball stops at rest. In \eqref{dyn_opt_problem_ball},  the performance index is
\begin{equation} \label{eq_ball_J}
J \equiv \int_a^b L\left(\bx,\bu,\mu\right) \dt = \int_a^b \left[ \sum_{i=1}^3 \frac{\gamma_i(\mu)}{2} {\ddot \theta}_i^2+\sum_{j=1}^2 h_j(\mu) \, S \left(\left| \bz - \bv_j \right| - \rho_j \right) \right] \dt,
\end{equation}
where the integrand cost function is
\begin{equation} \label{eq_ball_integrand_cost}
L\left(\bx,\bu,\mu\right) \equiv \sum_{i=1}^3 \frac{\gamma_i(\mu)}{2} {\ddot \theta}_i^2+\sum_{j=1}^2 h_j(\mu) \, S \left(\left| \bz - \bv_j \right| - \rho_j \right),
\end{equation}
for positive coefficients $\gamma_i(\mu)$, $1 \le i \le 3$, and $ h_j(\mu)$, $1 \le j \le 2$. 
The first $3$ summands $\frac{\gamma_i(\mu)}{2} {\ddot \theta}_i^2$, $1 \le i \le 3$, in $L$ limit the magnitude of the acceleration of the $i^\mathrm{th}$ control mass parameterization and  the final $2$ summands $h_j(\mu) \, S \left(\left| \bz - \bv_j \right| - \rho_j \right)$, $1 \le j \le 2$, in $L$ encourage the ball's GC to avoid the pair of obstacles. For the obstacle avoidance function in $L$, $S$ is either the time-reversed sigmoid function \eqref{eq_sigmoid} or the $C^2$ cutoff function
\begin{equation} \label{eq_obst_cutoff}
S(y)\equiv \left[\max \left\{0,-y \right\} \right]^4  =\ReLU^4\left(-y\right). 
\end{equation}   
In \eqref{eq_obst_cutoff}, $\ReLU$ is the rectified linear unit  function frequently used in the machine learning literature \cite{lecun2015deep}. \revision{R1Q4}{In \eqref{eq_ball_integrand_cost}, the coefficients $\gamma_i(\mu)$, $1 \le i \le 3$, and $h_j(\mu)$, $1 \le j \le 2$, depend on the scalar continuation parameter $\mu$ so that a sequence of optimal control problems may be constructed.} Note that a solution  obtained by the optimal control procedure, that minimizes \eqref{eq_disk_J} for the rolling disk or \eqref{eq_ball_J} for the rolling ball, is a ``compromise" between several, often conflicting, components, where some components of the performance index can be made more prominent by making their coefficients appropriately larger. The minimization of the  performance index  does not guarantee the minimization of each  component individually.

\begin{table}[h!]
	\centering 
	{ 
		\setlength{\extrarowheight}{1.5pt}
		\begin{tabular}{| c | c |} 
			\hline
			\textbf{Parameter} & \textbf{Value} \\ 
			\hline\hline 
			$\btheta_a$ & $\begin{bmatrix}  0 & 2.0369 & .7044 \end{bmatrix}^\mathsf{T}$  \\  
			\hline
			$\dot \btheta_a$ & $\begin{bmatrix} 0 & 0 & 0 \end{bmatrix}^\mathsf{T}$ \\ 
			\hline
			$\mathfrak{q}_a$ & $\begin{bmatrix} 1 &  0 & 0 & 0 \end{bmatrix}^\mathsf{T}$ \\
			\hline
			$\bOm_a$ & $\begin{bmatrix}  0 & 0 & 0 \end{bmatrix}^\mathsf{T}$ \\ 
			\hline
			$\bz_a$ & $\begin{bmatrix}  0 & 0 \end{bmatrix}^\mathsf{T}$ \\ 
			\hline
		\end{tabular} 
	}
	\caption{Initial condition parameter values for the rolling ball. Refer to \eqref{eq_ball_initial_conds}.}
	\label{table_ball_ICs}
\end{table}

\begin{table}[h!]
	\centering 
	{ 
		\setlength{\extrarowheight}{1.5pt}
		\begin{tabular}{| c | c |} 
			\hline
			\textbf{Parameter} & \textbf{Value} \\ 
			\hline\hline 
			$\dot \btheta_b$ & $\begin{bmatrix} 0 & 0 & 0 \end{bmatrix}^\mathsf{T}$ \\ 
			\hline
			$\bOm_b$ & $\begin{bmatrix}  0 & 0 & 0 \end{bmatrix}^\mathsf{T}$ \\ 
			\hline
			$\bz_b$ & $\begin{bmatrix}  1 & 1 \end{bmatrix}^\mathsf{T}$ \\ 
			\hline
		\end{tabular} 
	}
	\caption{Final condition parameter values for the rolling ball. Refer to \eqref{eq_ball_final_conds}.}
	\label{table_ball_FCs}
\end{table}

There is also a DAE formulation of the optimal control problem for the rolling ball which explicitly enforces the algebraic versor constraint on $\mathfrak{q}$ and which is mathematically equivalent to \eqref{dyn_opt_problem_ball}. In the DAE formulation, the first component, $q_0$, of the versor $\mathfrak{q}$ is moved from the state $\bx$ to the control $\bu$ and an imitator state, $\tilde q_0$, is used to replace $q_0$ in $\bx$. $\tilde q_{a,0} = q_{a,0}$, so that with perfect integration (i.e. no numerical integration errors), $\tilde q_0(t) = q_0(t)$ for $a \le t \le b$. Defining 
\begin{equation}
\tilde {\mathfrak{q}} \equiv \begin{bmatrix} \tilde q_0 \\  \mathfrak{q}^\flat \end{bmatrix},
\end{equation}
then with perfect integration, 
\begin{equation}
\tilde {\mathfrak{q}}(t) \equiv \begin{bmatrix} \tilde q_0(t) \\  \mathfrak{q}^\flat(t) \end{bmatrix}=\begin{bmatrix} q_0(t) \\  \mathfrak{q}^\flat(t) \end{bmatrix}=\mathfrak{q}(t)
\end{equation}
for $a \le t \le b$. $\tilde q_0$ is added to the state since the final conditions require knowledge of $q_0$, which is unavailable if it has been moved to the control since the final conditions are not a function of the control. The DAE formulation of the rolling ball's optimal control problem is
\begin{equation}
\min_{\bu_2} J 
\mbox{\, s.t. \,}
\left\{
\begin{array}{ll}
\dot {\bx}_2 = \mathbf{f}_2\left(\bx_2,\bu_2\right), \\
h_2\left(\bx_2,\bu_2\right)=1, \\
\bsigma_2\left(\bx_2(a)\right) = \mathbf{0},\\
\bpsi_2\left(\bx_2(b)\right) = \mathbf{0},
\end{array}
\right.
\label{dyn_opt_problem_dae2_ball}
\end{equation}
where
\revision{R1Q15}{
	\begin{equation}
	\bx_2 \equiv \begin{bmatrix} \btheta \\ \dot \btheta \\ \tilde {\mathfrak{q}} \\ \bOm \\ \bz  \end{bmatrix}, \quad
	\bu_2 \equiv \begin{bmatrix} \ddot \btheta \\ q_0 \end{bmatrix}, \quad \mathbf{f}_2\left(\bx_2,\bu_2\right) \equiv \begin{bmatrix} \dot \btheta \\ \ddot \btheta  \\ \frac{1}{2} \mathfrak{q} \bOm^\sharp \\ \bkappa\left(\bx,\ddot \btheta \right) \\ \left( \left[\mathfrak{q} \bOm^\sharp \mathfrak{q}^{-1} \right]^\flat \times r \mathbf{e}_3  \right)_{12}  \end{bmatrix}, \quad h_2\left(\bx_2,\bu_2\right) \equiv \left| \mathfrak{q} \right|^2,
	\end{equation}}
\begin{equation}
\bsigma_2\left(\bx_2(a)\right) \equiv \begin{bmatrix} \btheta(a) - \btheta_a \\ \dot \btheta(a) - {\dot \btheta}_a \\ \tilde {\mathfrak{q}}(a)-\mathfrak{q}_a \\ \bOm(a) - \bOm_a \\ \bz(a)- \bz_a \end{bmatrix},
\end{equation}
and
\begin{equation}
\bpsi_2\left(\bx_2(b)\right) \equiv \begin{bmatrix} \bPi \left(\left[\tilde {\mathfrak{q}}(b) \left[ \frac{1}{M} \sum_{i=0}^3 m_i \bzeta_i\left(\theta_i(b)\right) \right]^\sharp \tilde {\mathfrak{q}}(b)^{-1} \right]^\flat \right) \\ \dot \btheta(b) - {\dot \btheta}_b \\ \bOm(b) - \bOm_b \\ \bz(b)- \bz_b  \end{bmatrix}.
\end{equation}
Even though the DAE formulation \eqref{dyn_opt_problem_dae2_ball} is mathematically equivalent to the ODE formulation \eqref{dyn_opt_problem_ball}, the DAE formulation \eqref{dyn_opt_problem_dae2_ball} tends to be numerically more stable to solve than the ODE formulation \eqref{dyn_opt_problem_ball}, as explained in Example 6.12 ``Reorientation of an Asymmetric Rigid Body'' of \cite{betts2010practical}.

As explained in Appendix~\ref{sec_optimal_control}, the controlled equations of motion for the ODE formulaton \eqref{dyn_opt_problem_ball} of the rolling ball's optimal control problem are encapsulated by the ODE TPBVP:
\begin{equation} \label{eq_pmp_bvp_ball}
\begin{split}
\dot {\bx} &= \hat{H}_{\blam}^\mathsf{T} \left(\bx,\blam,\mu\right) = \hat{\mathbf{f}}\left(\bx,\blam\right) \equiv \mathbf{f}\left(\bx,\bpi\left(\bx,\blam\right)\right), \\
\dot {\blam} &= - \hat{H}_{\bx}^\mathsf{T} \left(\bx,\blam,\mu\right) =-H_{\bx}^\mathsf{T}\left(\bx,\blam,\bpi\left(\bx,\blam\right),\mu\right), \\
\left. \blam \right|_{t=a} &= -G_{\bx(a)}^\mathsf{T}, \quad G_{\bxi}^\mathsf{T} = \bsigma\left(\bx(a)\right) = \mathbf{0}, \\
\left. \blam \right|_{t=b} &= G_{\bx(b)}^\mathsf{T}, \quad G_{\bnu}^\mathsf{T} = \bpsi\left(\bx(b)\right) = \mathbf{0}.
\end{split}
\end{equation}
Subappendix~A.2 of \cite{putkaradze2020optimal} derives the formulas for constructing $H_{\bx}^\mathsf{T}$. In \eqref{eq_pmp_bvp_ball}, $G$ is the endpoint function 
\begin{equation} \label{eq_ball_endpoint_fcn}
\begin{split}
G\left(\bx(a),\bxi,\bx(b),\bnu\right) &\equiv \bxi^\mathsf{T} \bsigma\left(\bx(a)\right)+\bnu^\mathsf{T} \bpsi\left(\bx(b)\right) \\
&= \bxi^\mathsf{T} \begin{bmatrix} \btheta(a) - \btheta_a \\ \dot \btheta(a) - {\dot \btheta}_a \\ \mathfrak{q}(a)-\mathfrak{q}_a \\ \bOm(a) - \bOm_a \\ \bz(a)- \bz_a \end{bmatrix}+\bnu^\mathsf{T} \begin{bmatrix} \bPi \left(\left[\mathfrak{q}(b) \left[ \frac{1}{M} \sum_{i=0}^3 m_i \bzeta_i\left(\theta_i(b)\right) \right]^\sharp \mathfrak{q}(b)^{-1} \right]^\flat \right) \\ \dot \btheta(b) - {\dot \btheta}_b \\ \bOm(b) - \bOm_b \\ \bz(b)- \bz_b  \end{bmatrix},
\end{split}
\end{equation}
\revision{R1Q4}{where $\bxi \in \mathbb{R}^{15}$ and $\bnu \in \mathbb{R}^{10}$ are constant Lagrange multiplier vectors enforcing the initial and final conditions, \eqref{eq_ball_initial_conds} and \eqref{eq_ball_final_conds}, respectively.} In \eqref{eq_pmp_bvp_ball}, $H$ is the Hamiltonian
\begin{equation} \label{eq_ball_ham}
\begin{split}
H\left(\bx,\blam,\bu,\mu\right) &\equiv L\left(\bx,\bu,\mu\right) + \blam^\mathsf{T} \mathbf{f}\left(\bx,\bu\right) \\
&= \sum_{i=1}^3 \frac{\gamma_i(\mu)}{2} {\ddot \theta}_i^2+\sum_{j=1}^2 h_j(\mu) \, S \left(\left| \bz - \bv_j \right| - \rho_j \right)   + \blam^\mathsf{T} \begin{bmatrix} \dot \btheta \\ \bu  \\ \frac{1}{2} \mathfrak{q} \bOm^\sharp \\ \bkappa\left(\bx,\bu \right) \\ \left( \left[\mathfrak{q} \bOm^\sharp \mathfrak{q}^{-1} \right]^\flat \times r \mathbf{e}_3  \right)_{12}  \end{bmatrix},
\end{split}
\end{equation}
\revision{R1Q4}{where $\blam \in \mathbb{R}^{15}$ is a time-varying Lagrange multiplier vector enforcing the dynamics \eqref{rolling_ball_opt_con_f}.} In \eqref{eq_pmp_bvp_ball}, $\bpi$ is an analytical formula expressing the control $\bu$ as a function of the state $\bx$ and the costate $\blam$. The components of $\bpi$ are given by
\begin{equation}  \label{eq_ddtheta_exp_rball}
{\ddot \theta}_i = \pi_i \left(\bx,\blam\right) \equiv -\gamma_i^{-1} \left\{ \lambda_{3+i} + \blam_{\bOm}^\mathsf{T} \left[\sum_{k=0}^3 m_k \widehat{\mathbf{s}_k}^2  -\inertia \right]^{-1}  \left[ m_i \mathbf{s}_i \times  \bzeta_i^{\prime}  \right]  \right\},
\end{equation}
for $1 \le i \le 3$ and where $\blam_{\bOm} \equiv \begin{bmatrix} \lambda_{11} & \lambda_{12} & \lambda_{13} \end{bmatrix}^\mathsf{T}$. In \eqref{eq_pmp_bvp_ball}, $\hat H$ is the regular Hamiltonian
\begin{equation} \label{eq_ball_regular_Hamiltonian}
\begin{split}
\hat H\left(\bx,\blam,\mu\right) &\equiv H\left(\bx,\blam,\bpi \left(\bx,\blam\right),\mu\right) \\
&= \sum_{i=1}^3 \frac{\gamma_i(\mu)}{2} {\pi}_i^2 \left(\bx,\blam\right)+\sum_{j=1}^2 h_j(\mu) \, S \left(\left| \bz - \bv_j \right| - \rho_j \right)  + \blam^\mathsf{T} \begin{bmatrix} \dot \btheta \\ \bpi \left(\bx,\blam\right)  \\ \frac{1}{2} \mathfrak{q} \bOm^\sharp \\ \bkappa\left(\bx,\bpi \left(\bx,\blam\right) \right) \\ \left( \left[\mathfrak{q} \bOm^\sharp \mathfrak{q}^{-1} \right]^\flat \times r \mathbf{e}_3  \right)_{12}  \end{bmatrix}.
\end{split}
\end{equation}
The reader is referred to \cite{putkaradze2020optimal} for a more general description of the ODE and DAE formulations, \eqref{dyn_opt_problem_ball} and \eqref{dyn_opt_problem_dae2_ball}, of the rolling ball's optimal control problem and the controlled equations of motion \eqref{eq_pmp_bvp_ball} correpsonding to \eqref{dyn_opt_problem_ball}. The DAE TPBVP encapsulating the controlled equations of motion corresponding to the  DAE formulation \eqref{dyn_opt_problem_dae2_ball} of the rolling ball's optimal control problem were not investigated since a robust DAE TPBVP solver is not readily available in \mcode{MATLAB}. \revisionNACO{R1Q6}{COLDAE is a robust DAE TPBVP solver that uses collocation \cite{ascher1994collocation}; however, COLDAE is only available in Fortran and thus was not used in our calculations. }
 
\revisionNACO{R1Q5}{\subsection{Numerical Solutions: Sigmoid Obstacle Avoidance}} \label{ssec_ball_sim}

The controlled equations of motion \eqref{eq_pmp_bvp_ball} for the rolling ball are solved numerically to move the ball between the pair of points while avoiding the pair of obstacles, where the obstacle avoidance function $S$ in \eqref{eq_ball_integrand_cost} is realized via the time-reversed sigmoid function \eqref{eq_sigmoid} with $\epsilon=.01$. The direct method solver \mcode{GPOPS-II} is used to solve the DAE formulation \eqref{dyn_opt_problem_dae2_ball} of the optimal control problem, where the obstacle heights appearing in the integrand cost function \eqref{eq_ball_integrand_cost} are $h_1=h_2=0$ and where the values for the other integrand cost function coefficients are given in Table~\ref{table_ball_integrand}. Using this direct method solution as an initial guess, \mcode{GPOPS-II} is used again to solve the ODE formulation \eqref{dyn_opt_problem_ball} of the same optimal control problem.  Predictor-corrector continuation is then used to solve the controlled equations of motion \eqref{eq_pmp_bvp_ball}, starting from the second direct method solution. The continuation parameter is $\mu$, which is used to adjust $h_1=h_2$ according to the linear homotopy shown in Table~\ref{table_ball_integrand}, so that $h_1=h_2=0$ when $\mu=.95$ and $h_1=h_2=1{,}000$ when $\mu=.00001$. The predictor-corrector continuation begins at $\mu=.95$, which is consistent with the direct method solution obtained at $h_1=h_2=0$.

\begin{figure}[!ht]
	\centering
	\subfloat[The magnitude of the normal force  when the obstacle heights  $h_1=h_2$ are $0$.]{\includegraphics[scale=.5]{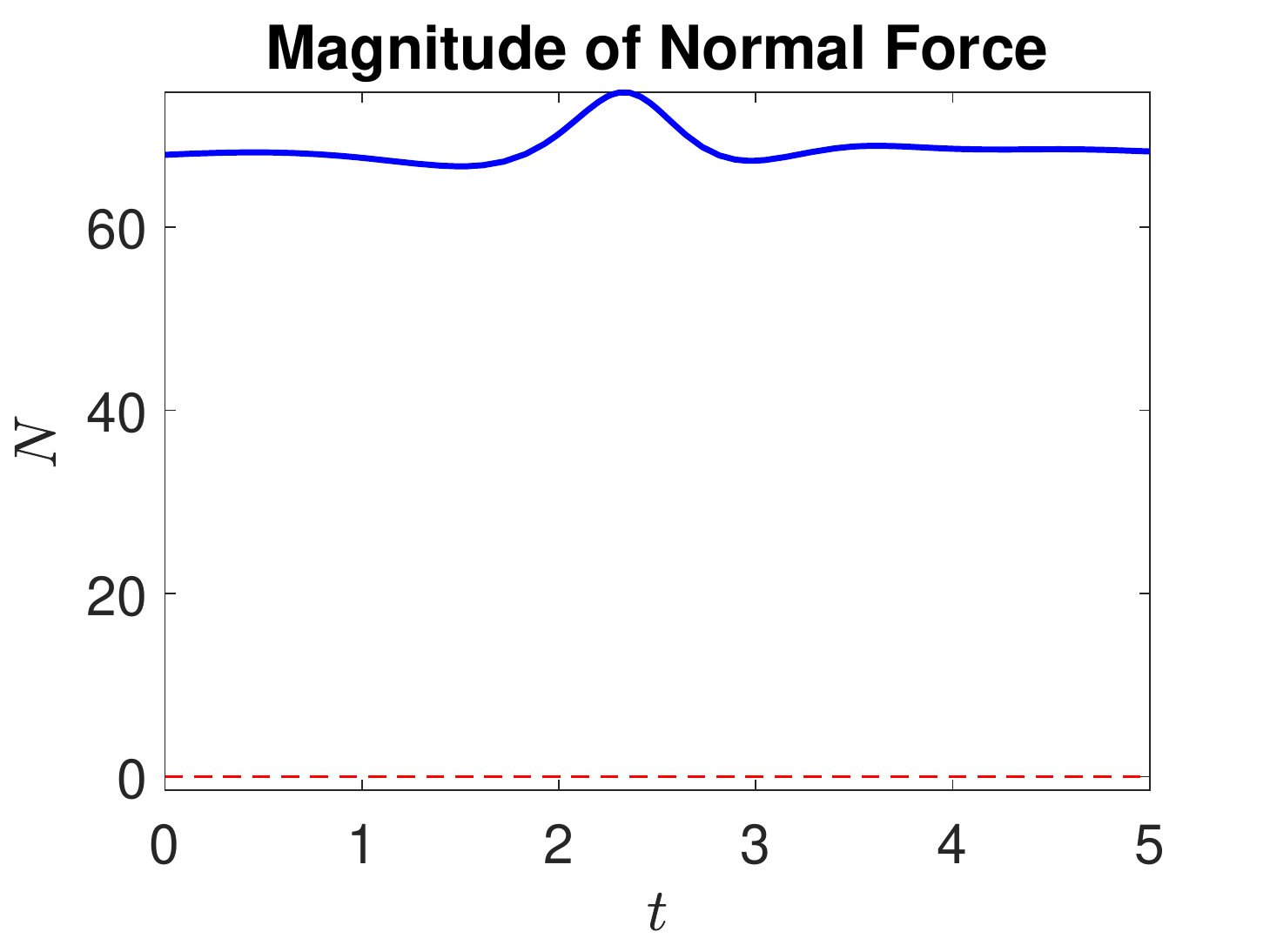}\label{fig_bsim10_dm_normal}}
	\hspace{5mm}
	\subfloat[The magnitude of the normal force when the obstacle heights $h_1=h_2$ are $9.93$.]{\includegraphics[scale=.5]{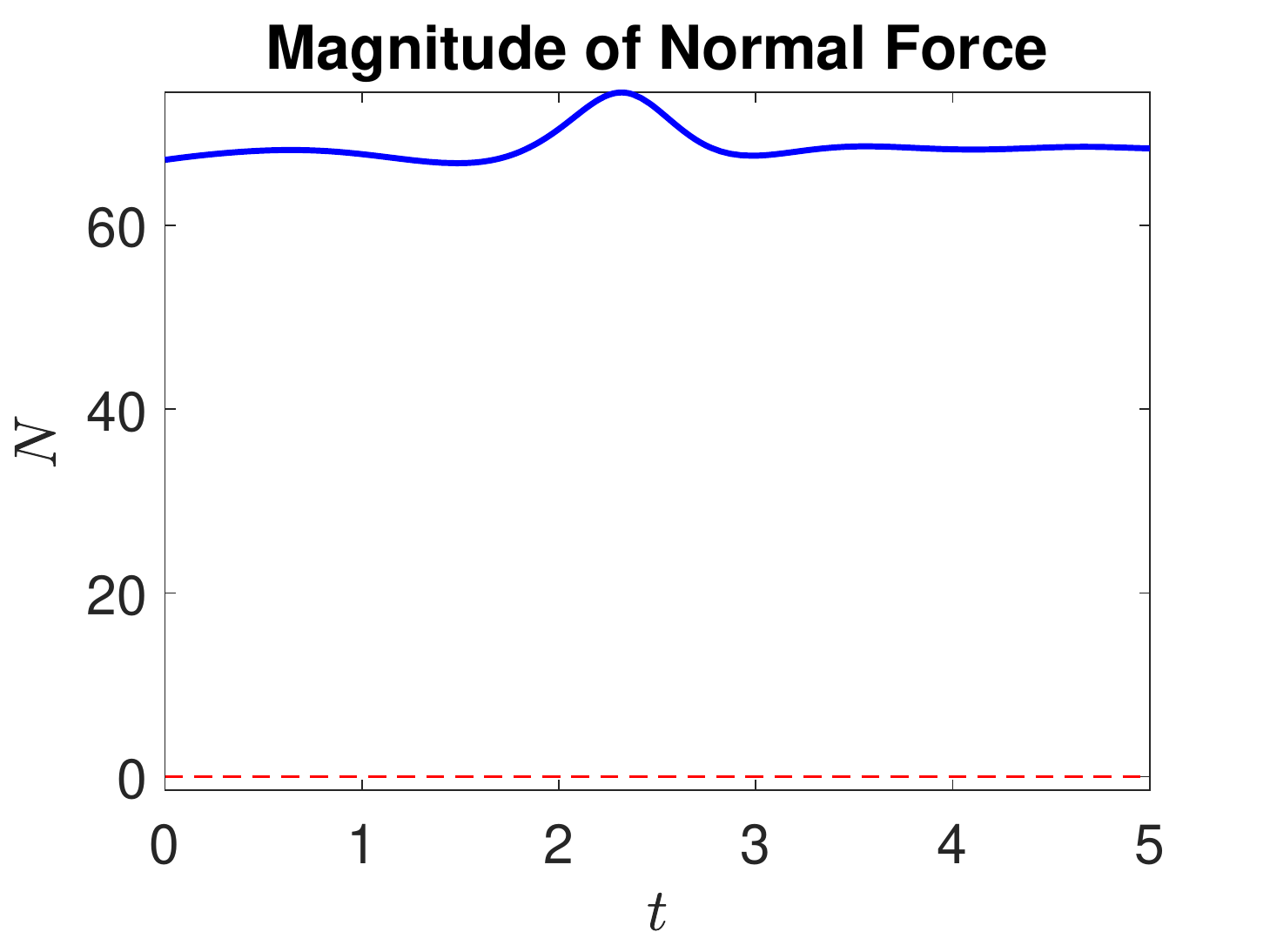}\label{fig_bsim10_pc_normal}}
	\\
	\subfloat[The minimum coefficient of static friction to prevent slipping is $.1055$ when the obstacle heights $h_1=h_2$ are $0$.]{\includegraphics[scale=.5]{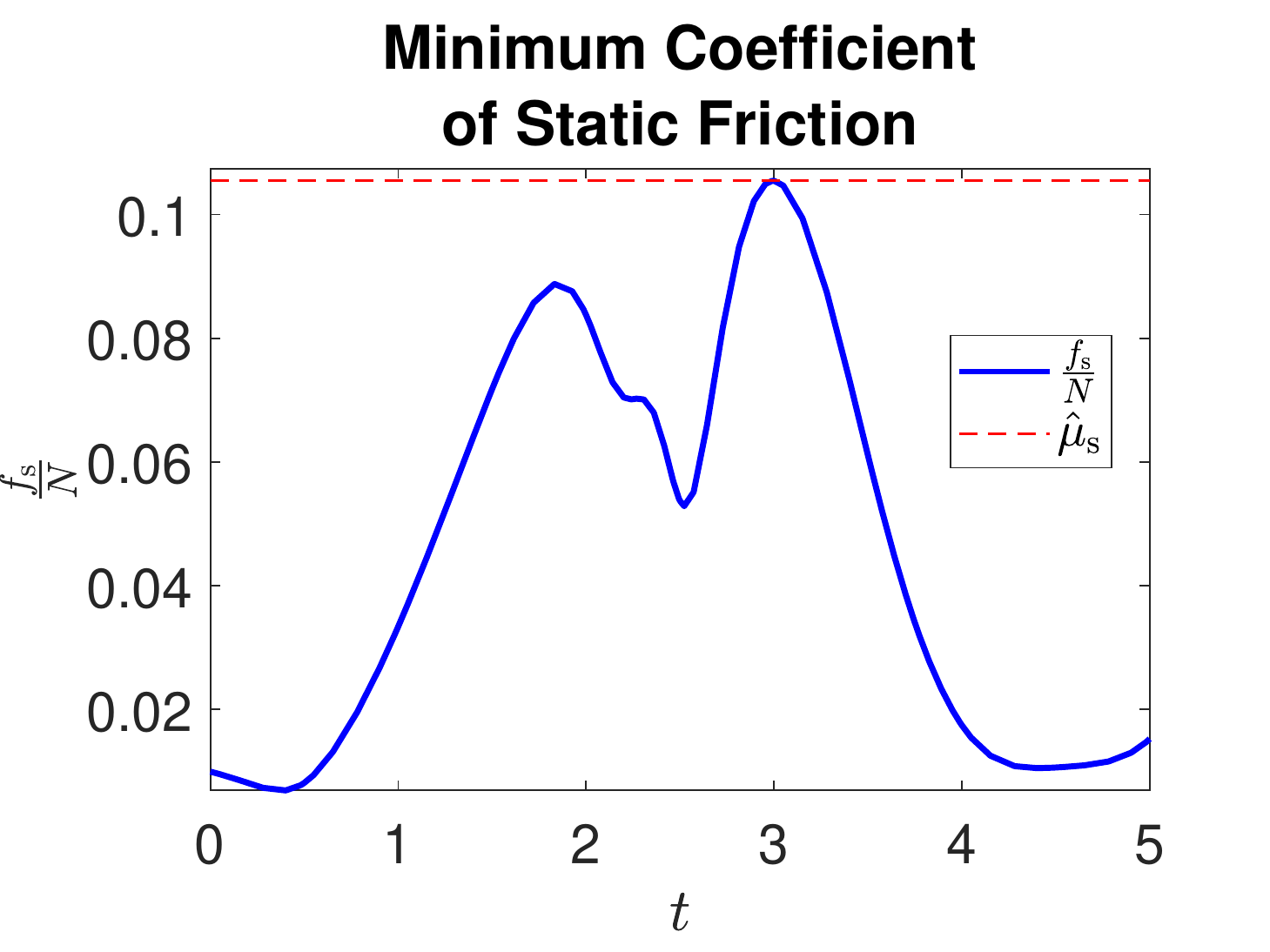}\label{fig_bsim10_dm_mu_s}}
	\hspace{5mm}
	\subfloat[The minimum coefficient of static friction to prevent slipping is $.0988$ when the obstacle heights $h_1=h_2$ are  $9.93$.]{\includegraphics[scale=.5]{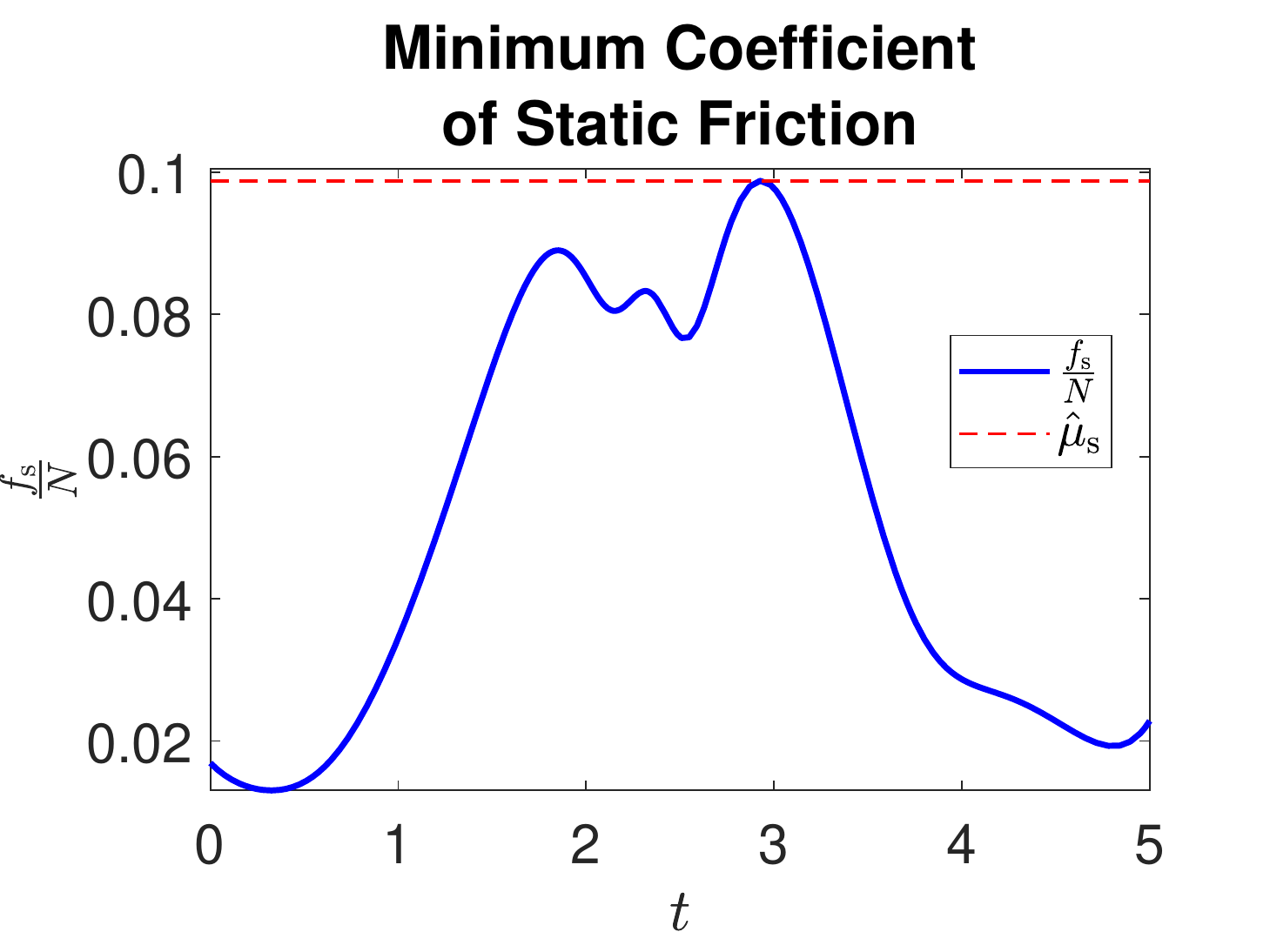}\label{fig_bsim10_pc_mu_s}}
	\caption{Numerical solutions of the rolling ball optimal control problem \eqref{dyn_opt_problem_ball} for sigmoid obstacle avoidance using $3$ control masses for $ \gamma_1=\gamma_2=\gamma_3=10$ and for fixed initial and final times. The direct method results for obstacle heights at $h_1=h_2=0$ are shown in the left column, while the predictor-corrector continuation indirect method results for obstacle heights at $h_1=h_2 \approx 9.93$ are shown in the right column. The ball does not detach from the surface since the magnitude of the normal force is always positive. The ball rolls without slipping if $\mu_\mathrm{s} \ge .1055$ for the direct method solution and if $\mu_\mathrm{s} \ge .0988$ for the indirect method solution. }
	\label{fig_bsim10_normal}
\end{figure}

\begin{figure}[!ht] 
	\centering 
	\subfloat[Evolution of the continuation parameter $\mu$. ]{\includegraphics[scale=.5]{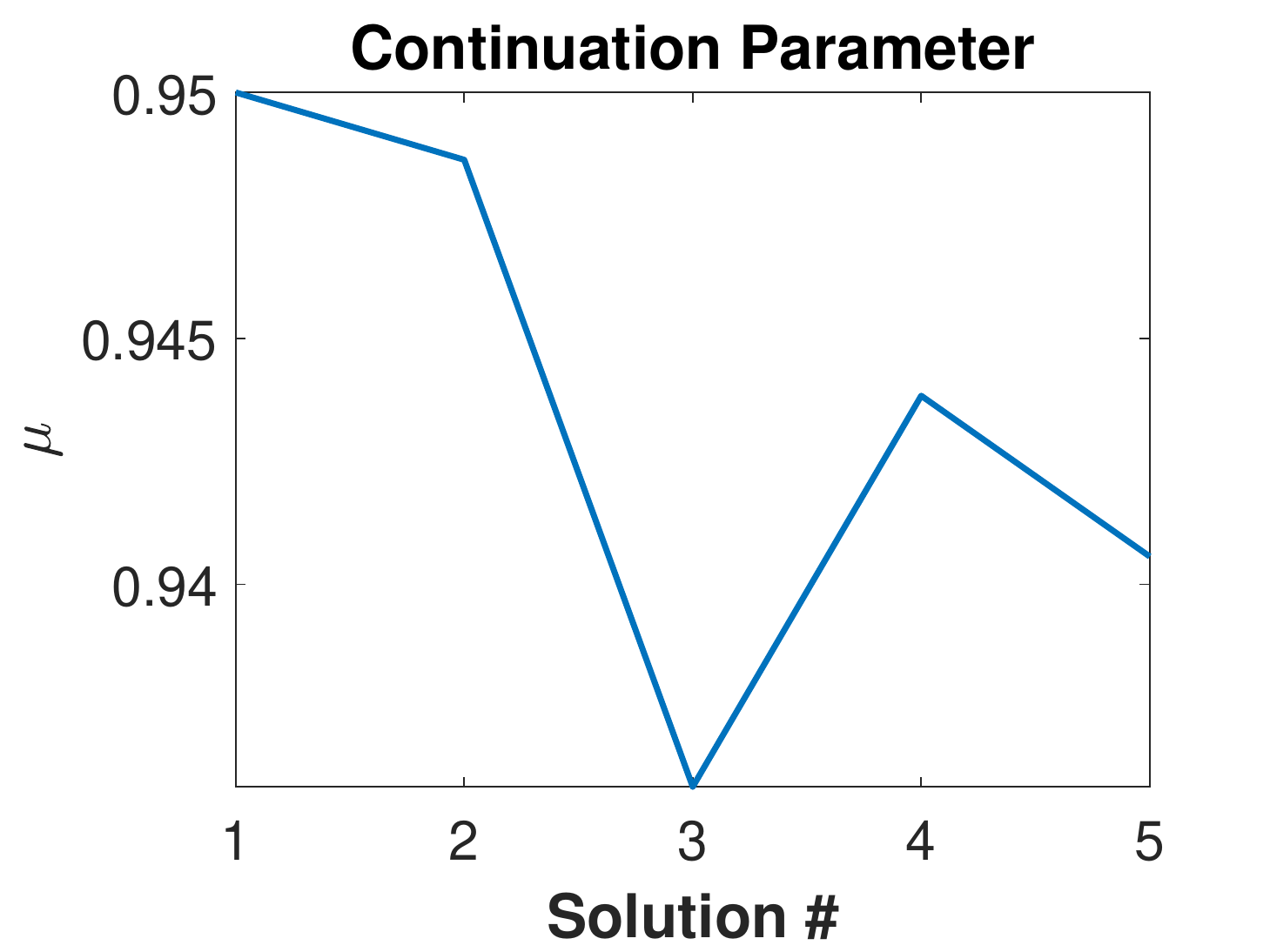}\label{fig_bsim10_pc_cont_param}}
	\hspace{5mm}
	\subfloat[Evolution of the obstacle heights $h_1=h_2$.]{\includegraphics[scale=.5]{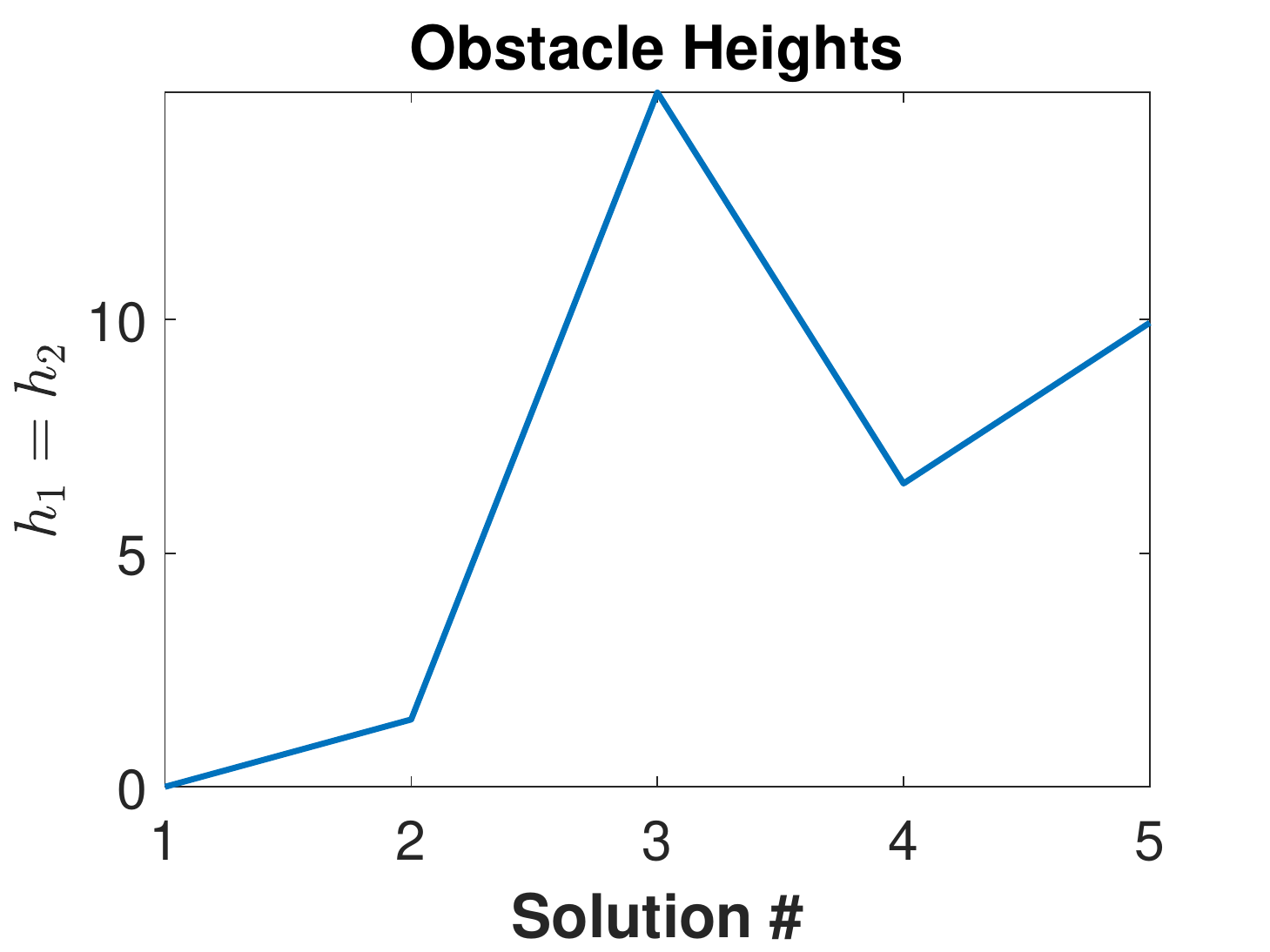}\label{fig_bsim10_pc_obst_height}}
	\\
	\subfloat[Evolution of the performance index $J$.]{\includegraphics[scale=.5]{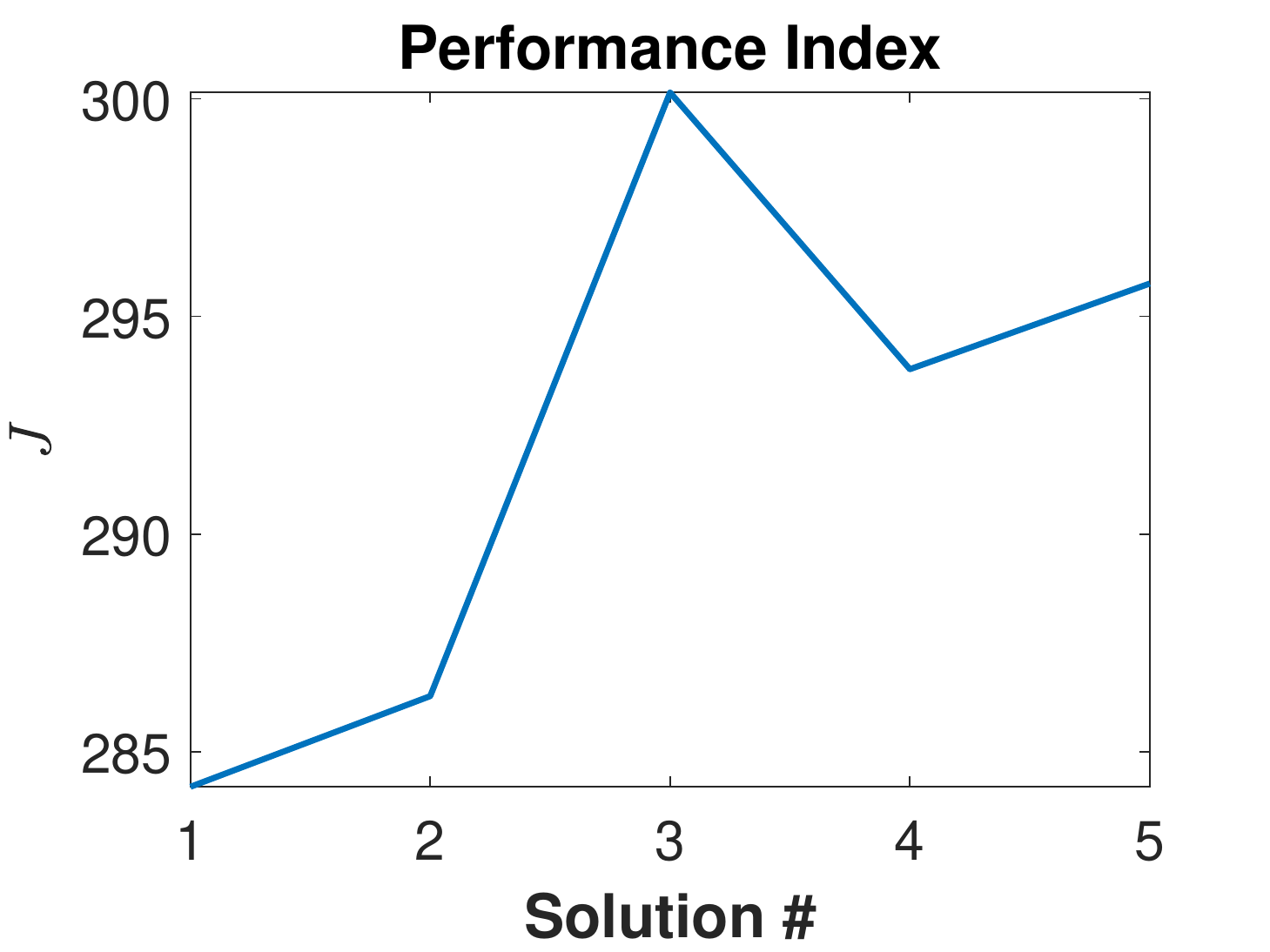}\label{fig_bsim10_pc_J}}
	\hspace{5mm}
	\subfloat[Evolution of the tangent steplength $\sigma$.]{\includegraphics[scale=.5]{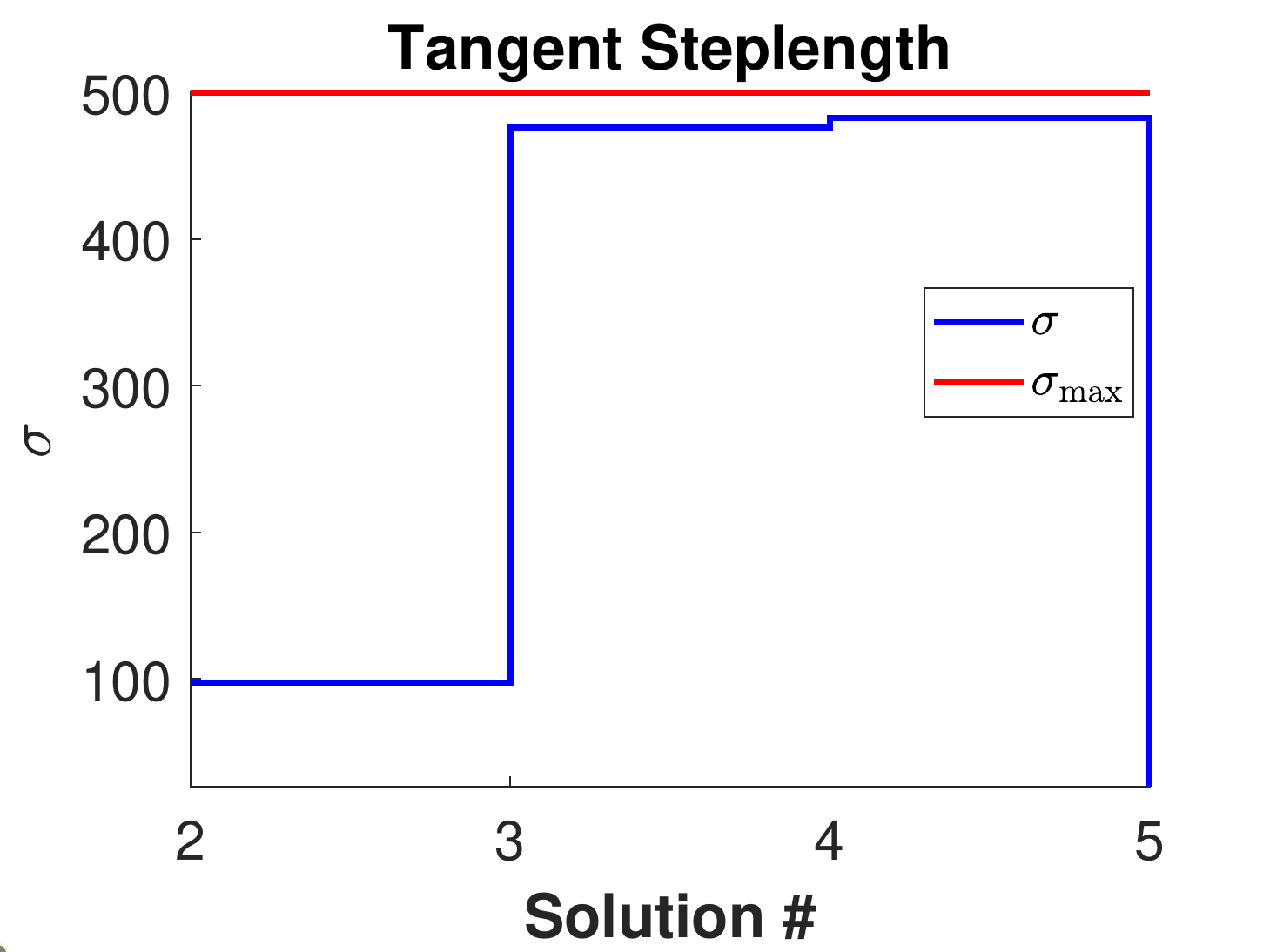}\label{fig_bsim10_pc_sigma}}
	\caption{Evolution of various parameters and variables during the predictor-corrector continuation indirect method, which starts from the direct method solution, used to solve the rolling ball optimal control problem \eqref{dyn_opt_problem_ball} for sigmoid obstacle avoidance. Note the turning points at solutions 3 and 4.}
	\label{fig_bsim10_pc}
\end{figure}

For the direct method, \mcode{GPOPS-II} is run using the NLP solver SNOPT. The \mcode{GPOPS-II} mesh error tolerance is $1\mathrm{e}\unaryminus 6$ and the SNOPT error tolerance is $1\mathrm{e}\unaryminus 7$. In order to encourage convergence of SNOPT, a constant $C=50$ is added to the integrand cost function $L$ in \eqref{eq_ball_integrand_cost}. The sweep predictor-corrector continuation method discussed in Appendix~\ref{app_sweep_predictor_corrector} is used by the indirect method. For the sweep predictor-corrector continuation method, there are $4$ predictor-corrector steps, the maximum tangent steplength in each step is $\sigma_\mathrm{max}= 500$, the direction of the initial unit tangent is determined by setting $d=\unaryminus 2$ to force the continuation parameter $\mu$ to initially decrease, the relative error tolerance is $1\mathrm{e}\unaryminus 6$, the unit tangent solver is \mcode{twpbvpc_m}, and the monotonic ``sweep'' continuation solver is \mcode{acdcc}. The numerical results are shown in Figures~\ref{fig_bsim10}, \ref{fig_bsim10_normal}, and \ref{fig_bsim10_pc}. As $\mu$ decreases from $.95$ down to $.9406$ during continuation (see Figure~\ref{fig_bsim10_pc_cont_param}), $h_1=h_2$ increases from $0$ up to $9.93$ (see Figure~\ref{fig_bsim10_pc_obst_height}). Since $h_1=h_2$ is ratcheted up during continuation, thereby increasing the penalty in the integrand cost function \eqref{eq_ball_integrand_cost} when the GC intrudes into the obstacles, by the end of continuation, the ball's GC avoids both obstacles while veering smartly around the first obstacle (compare Figures~\ref{fig_bsim10_dm_gc_path} vs \ref{fig_bsim10_pc_gc_path}), at the expense of slightly larger magnitude controls (compare Figures~\ref{fig_bsim10_dm_controls} vs \ref{fig_bsim10_pc_controls}). The ball does not detach from the surface since the magnitude of the normal force is always positive (see Figures~\ref{fig_bsim10_dm_normal} and \ref{fig_bsim10_pc_normal}). The ball rolls without slipping if the coefficient of static friction $\mu_\mathrm{s}$ is at least $\hat \mu_\mathrm{s} \approx .1055$ for the direct method solution (see Figure~\ref{fig_bsim10_dm_mu_s}) and if $\mu_\mathrm{s}$ is at least $\hat \mu_\mathrm{s} \approx .0988$ for the indirect method solution (see Figure~\ref{fig_bsim10_pc_mu_s}). As shown in Figures~\ref{fig_bsim10_pc_cont_param}-\ref{fig_bsim10_pc_J}, the sweep predictor-corrector continuation indirect method encounters turning points at solutions 3 and 4. 

\begin{table}[h!]
	\centering 
	{ 
		\setlength{\extrarowheight}{1.5pt}
		\begin{tabular}{| c | c |} 
			\hline
			\textbf{Parameter} & \textbf{Value} \\ 
			\hline\hline 
			$\gamma_1=\gamma_2=\gamma_3$ & $10$ \\ 
			\hline
			$h_1(\mu)=h_2(\mu)$ & $\frac{.95-\mu}{.95-.00001}\left(1000\right)$ \\
			\hline
			$\bv_1$ & $\begin{bmatrix}  .2 & .2 \end{bmatrix}^\mathsf{T}$ \\
			\hline
			$\bv_2$ & $\begin{bmatrix}  .8 & .8 \end{bmatrix}^\mathsf{T}$ \\
			\hline
			$\rho_1=\rho_2$ & $.282$ \\
			\hline
		\end{tabular} 
	}
	\caption{Integrand cost function coefficient values for the rolling ball when predictor-corrector continuation is performed in the obstacle heights. Refer to \eqref{eq_ball_integrand_cost}.}
	\label{table_ball_integrand}
\end{table} 

\revisionNACO{R1Q5}{\subsection{Numerical Solutions: $\ReLU$ Obstacle Avoidance}} \label{ssec_ball_sim_redux}
The controlled equations of motion \eqref{eq_pmp_bvp_ball} for the rolling ball are solved numerically again to move the ball between the pair of points while avoiding the pair of obstacles, but this time the obstacle avoidance function $S$ in \eqref{eq_ball_integrand_cost} is realized via the $C^2$ cutoff function \eqref{eq_obst_cutoff}. With the obstacle heights appearing in the integrand cost function \eqref{eq_ball_integrand_cost} set to $h_1=h_2=0$ and with the values for the other integrand cost function coefficients set according to Table~\ref{table_ball_integrand}, the same double pass (DAE formulation followed by ODE formulation) direct method, with all the same settings, parameters, and initial conditions, is used as in the previous subsection to generate the same initial solution. Starting from the direct method solution, two rounds of the same sweep predictor-corrector continuation method that was used in the previous subsection are again used here to solve the controlled equations of motion \eqref{eq_pmp_bvp_ball}, with all the same settings and parameters, except for the number of predictor-corrector steps and maximum tangent steplengths used. In the first round, the continuation parameter $\mu$ is used to adjust $h_1=h_2$ according to the linear homotopy shown in Table~\ref{table_ball_integrand}, so that $h_1=h_2=0$ when $\mu=.95$ and $h_1=h_2=1{,}000$ when $\mu=.00001$. This predictor-corrector continuation begins at $\mu=.95$, which is consistent with the direct method solution obtained at $h_1=h_2=0$. In the first round, two predictor-corrector steps are made with the maximum tangent steplengths $\sigma_\mathrm{max}=\begin{bmatrix} 20000 & 100000 \end{bmatrix}$. Starting from the predictor-corrector continuation solution obtained in the first round, a second predictor-corrector continuation is used to solve the controlled equations of motion \eqref{eq_pmp_bvp_ball}, where the continuation parameter is $\mu$, which is now used to adjust $\gamma_1=\gamma_2=\gamma_3$ according to the linear homotopy shown in Table~\ref{table_ball_integrand_pc2}, so that $\gamma_1=\gamma_2=\gamma_3=10$ when $\mu=.95$ and $\gamma_1=\gamma_2=\gamma_3=\unaryminus1{,}000$ when $\mu=.00001$. The second predictor-corrector continuation begins at $\mu=.95$, which is consistent with the first predictor-corrector continuation solution obtained at $\gamma_1=\gamma_2=\gamma_3=10$. Moreover, during the second predictor-corrector continuation, the obstacle heights $h_1=h_2$ are fixed at $7.846\mathrm{e}{8}$, which is consistent with the final obstacle heights obtained by the first predictor-corrector continuation solution. In the second continuation, four predictor-corrector steps are made with the maximum tangent steplengths $\sigma_\mathrm{max}=\begin{bmatrix} 5000 & 200 & 1 & 1 \end{bmatrix}$.

\begin{figure}[!ht] 
	\centering
	\subfloat[The GC plows through the obstacles when the obstacle heights $h_1=h_2$ are $0$.]{\includegraphics[scale=.35]{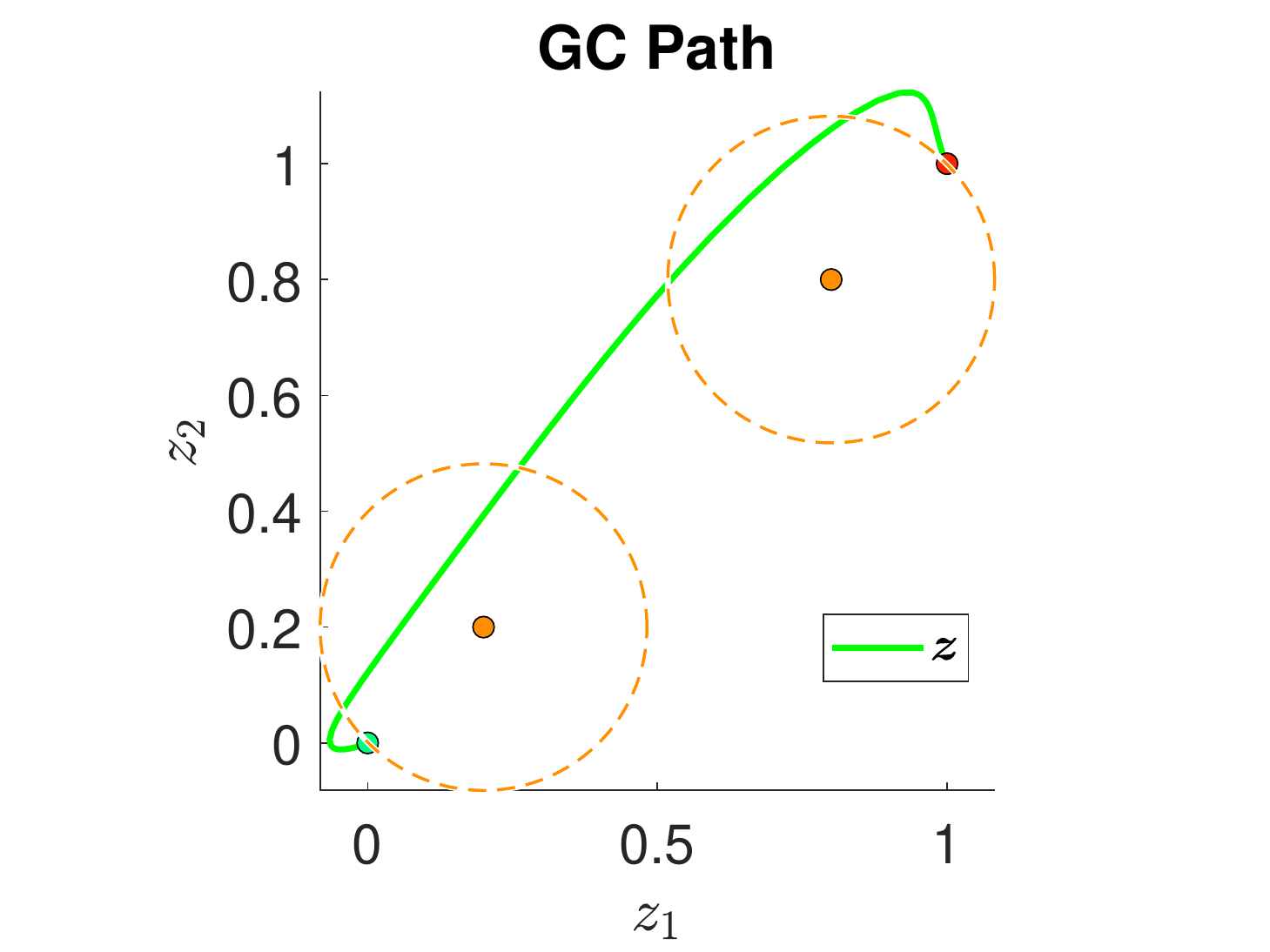}\label{fig_bsim12_dm_gc_path}}
	\hspace{5mm}
	\subfloat[The GC nearly clears the obstacles after ramping the obstacle heights $h_1=h_2$ from $0$ up to $7.846\mathrm{e}{8}$.]{\includegraphics[scale=.35]{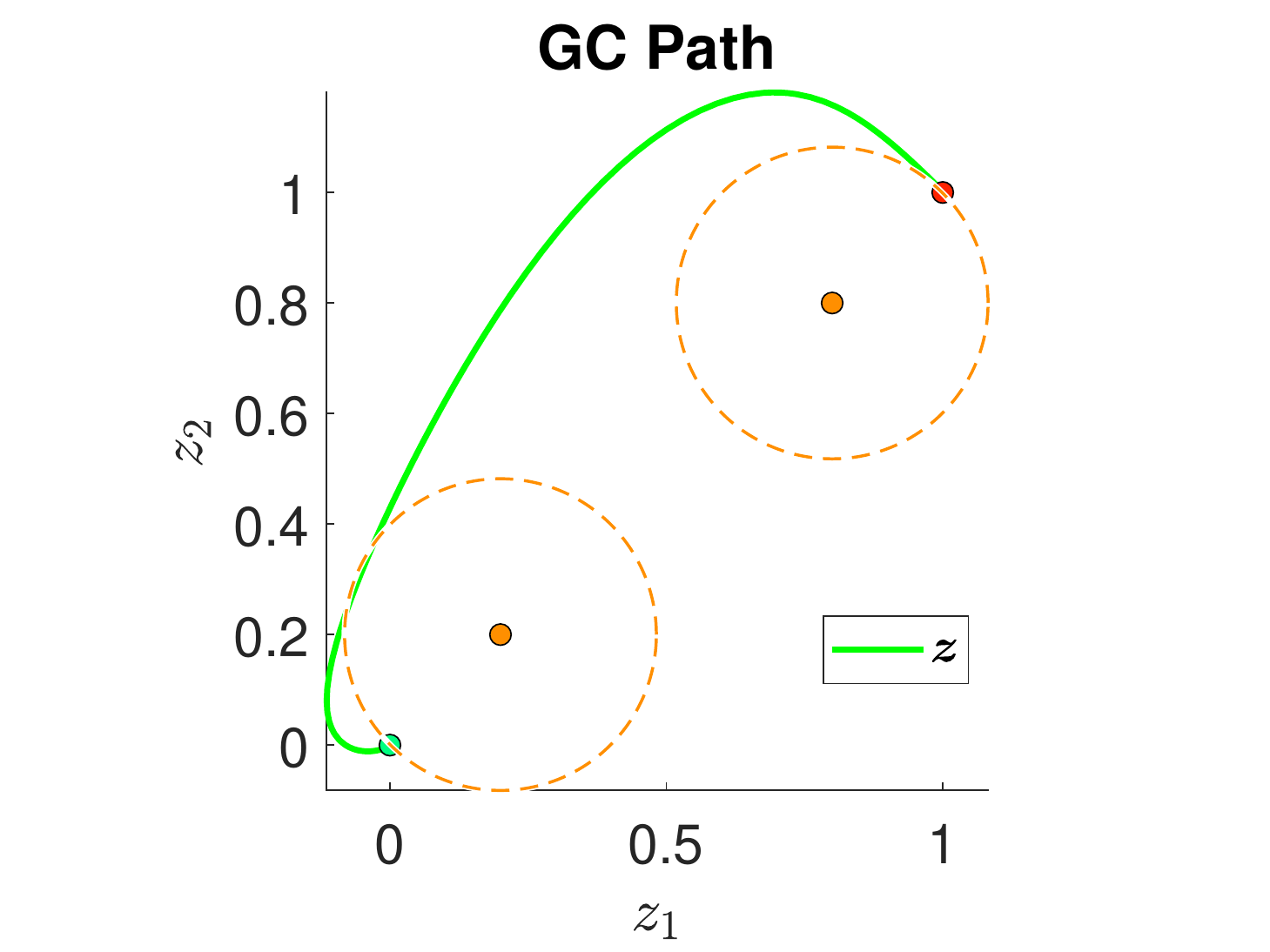}\label{fig_bsim12_pc_h_gc_path}}
	\hspace{5mm}
	\subfloat[The GC veers tightly around the obstacle boundaries after relaxing the control coefficients $\gamma_1=\gamma_2=\gamma_3$ from $10$ down to $3.602\mathrm{e}{-5}$.]{\includegraphics[scale=.35]{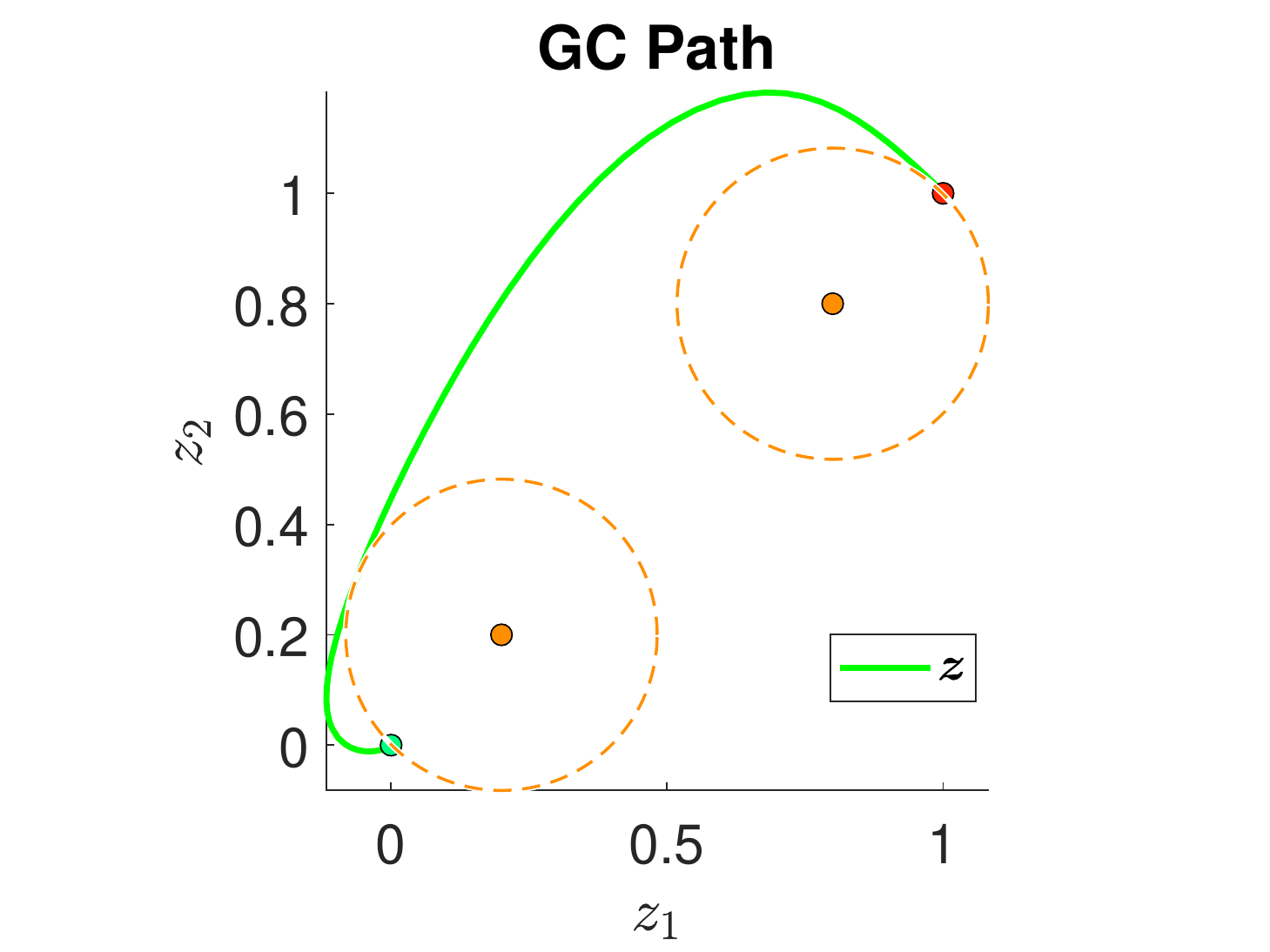}\label{fig_bsim12_pc_gamma_gc_path}}
	\\
	\subfloat[Motion of the center of masses when the obstacle heights $h_1=h_2$ are $0$.]{\includegraphics[scale=.35]{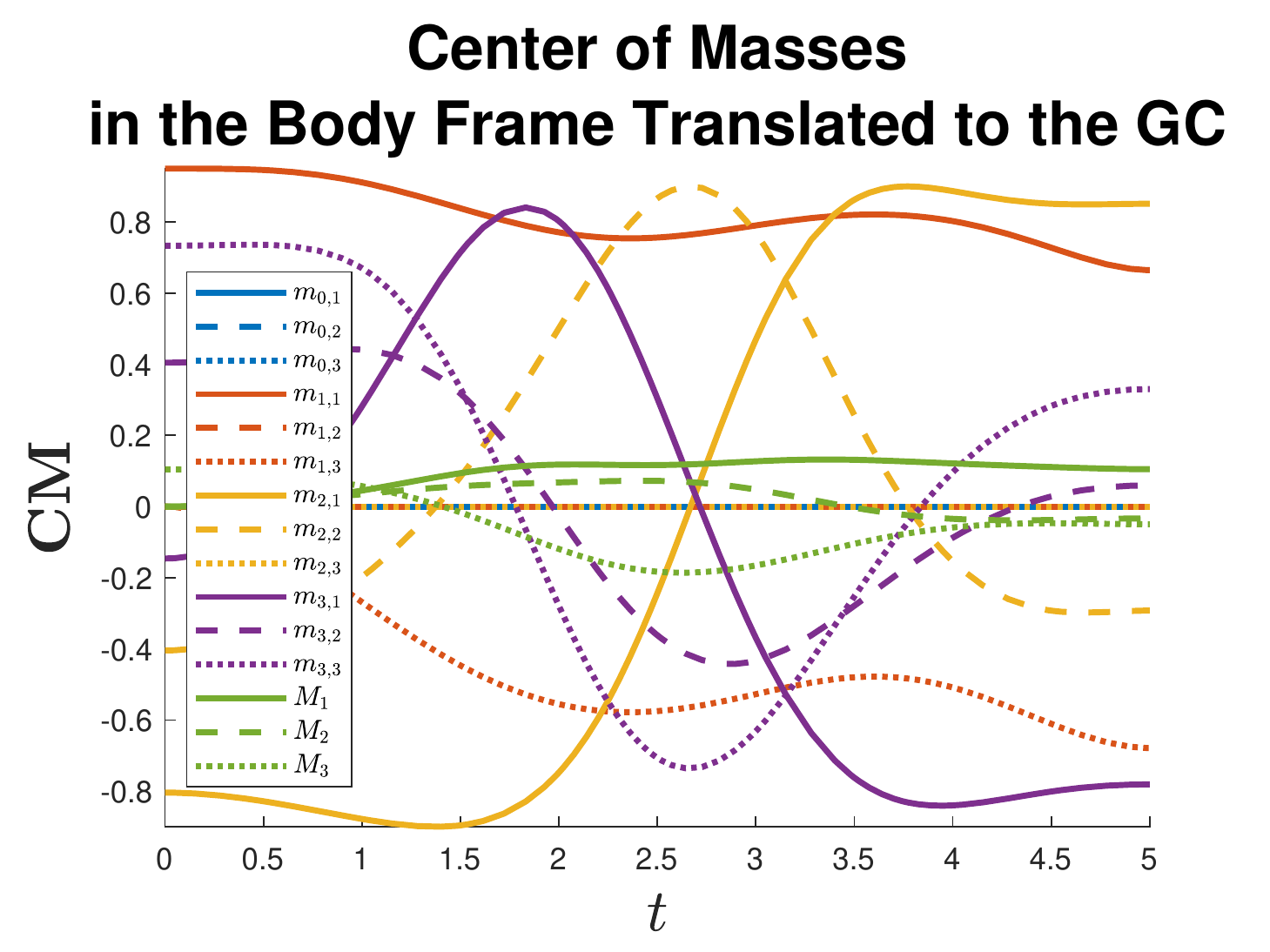}\label{fig_bsim12_dm_cm_bf}}
	\hspace{5mm}
	\subfloat[Motion of the center of masses after ramping the obstacle heights $h_1=h_2$ from $0$ up to $7.846\mathrm{e}{8}$.]{\includegraphics[scale=.35]{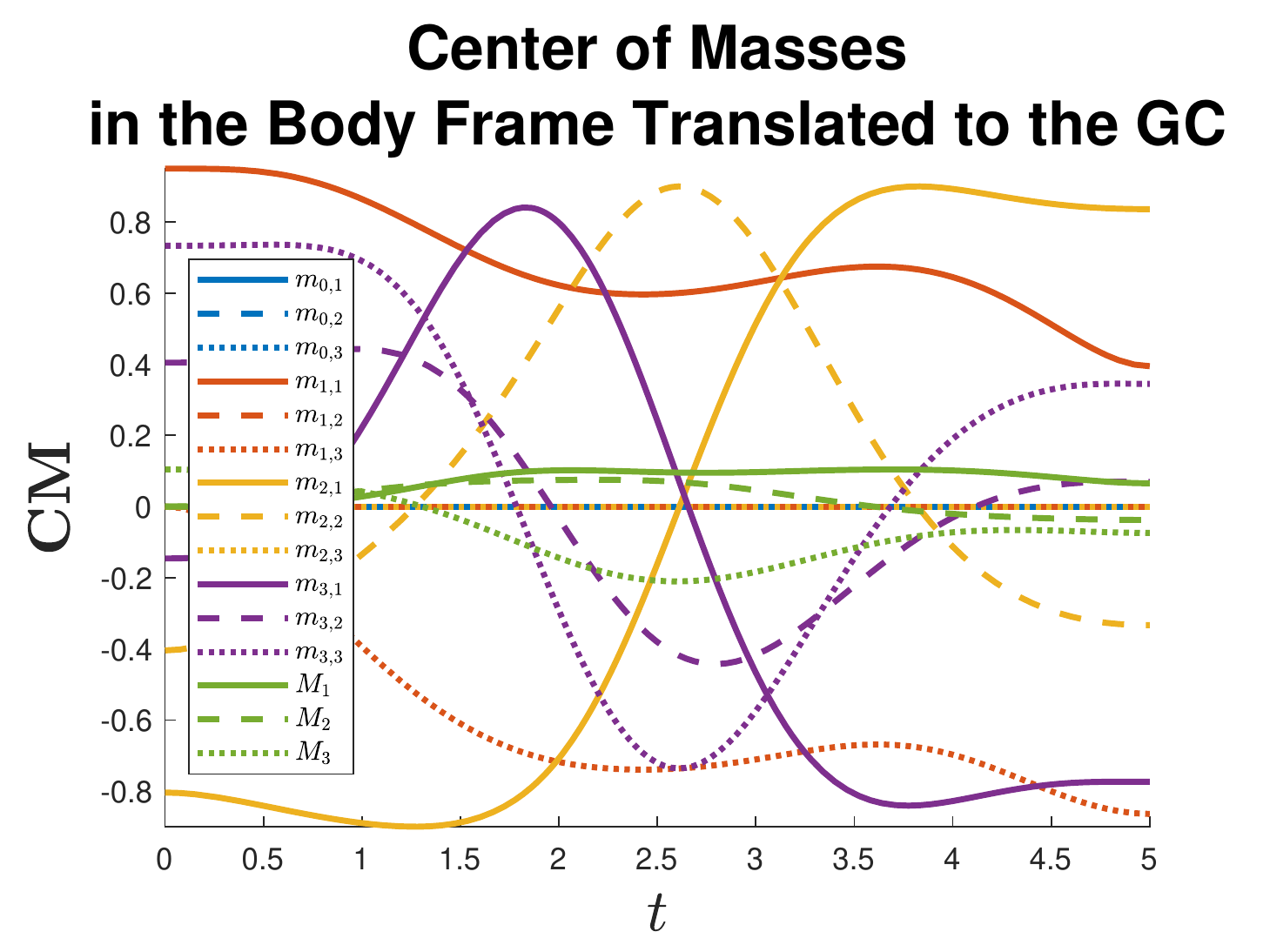}\label{fig_bsim12_pc_h_cm_bf}}
	\hspace{5mm}
	\subfloat[Motion of the center of masses after relaxing the control coefficients $\gamma_1=\gamma_2=\gamma_3$ from $10$ down to $3.602\mathrm{e}{-5}$.]{\includegraphics[scale=.35]{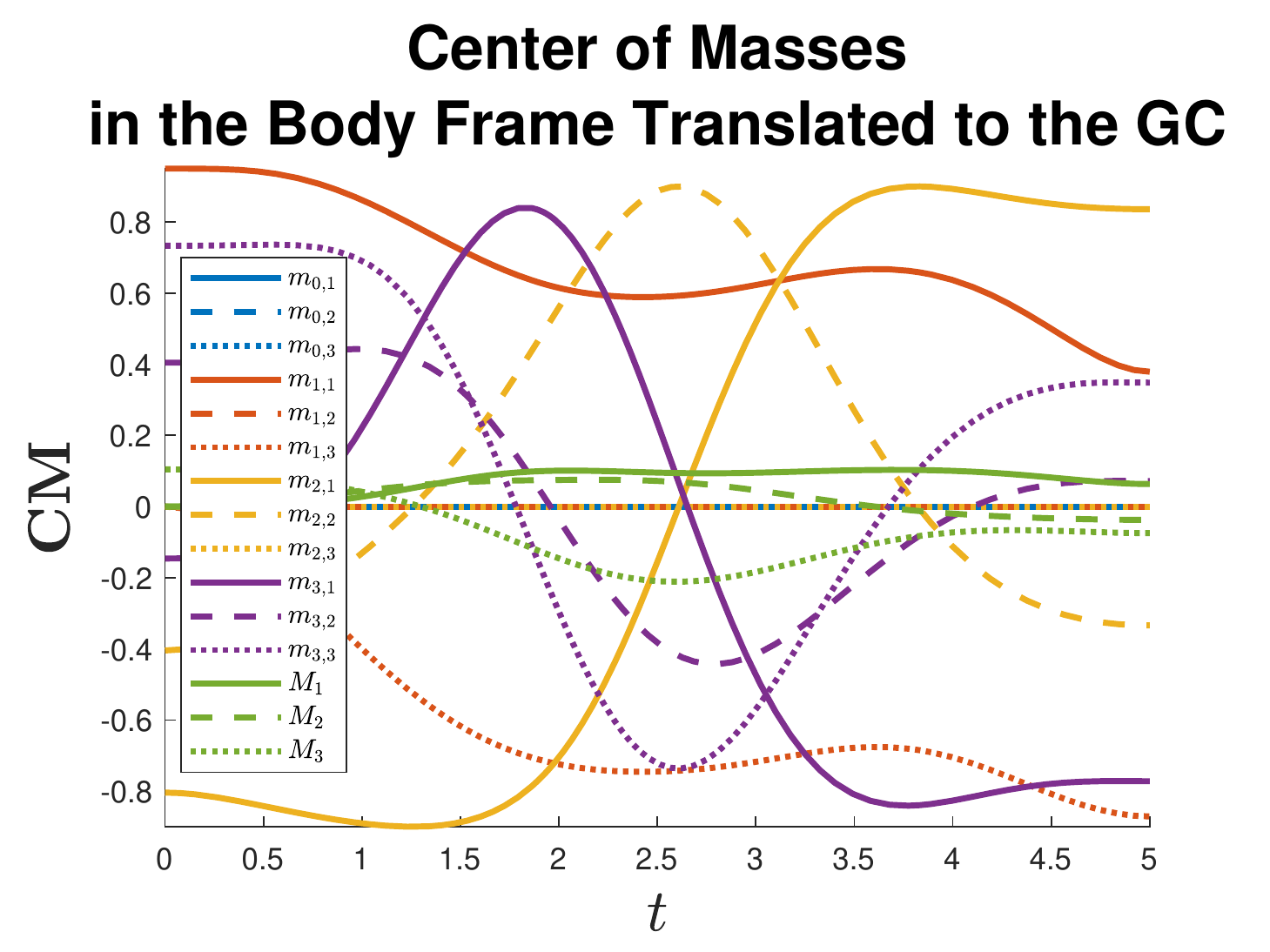}\label{fig_bsim12_pc_gamma_cm_bf}}
	\\
	\subfloat[The controls when the obstacle heights $h_1=h_2$ are $0$.]{\includegraphics[scale=.35]{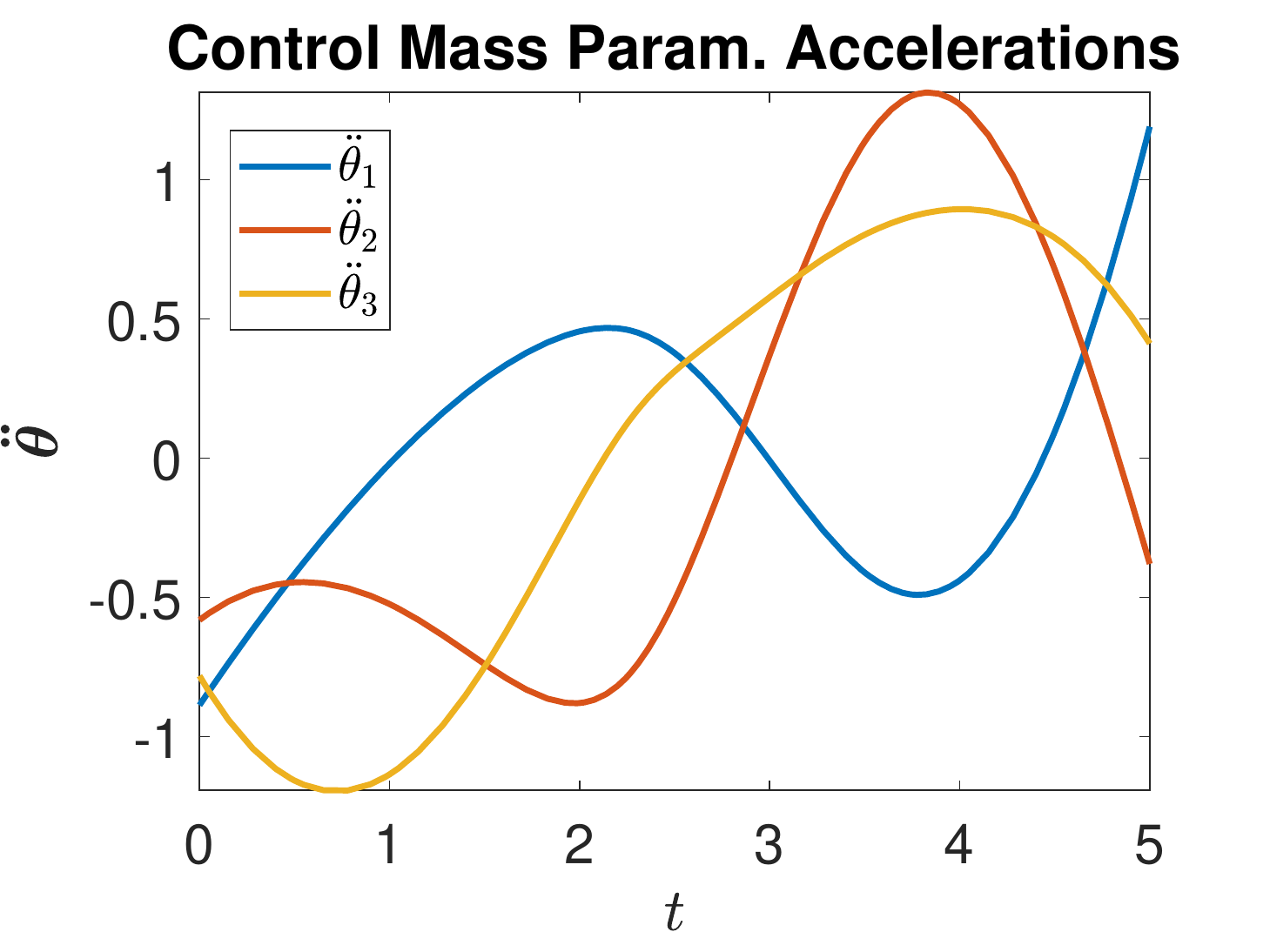}\label{fig_bsim12_dm_controls}}
	\hspace{5mm}
	\subfloat[The controls after ramping the obstacle heights $h_1=h_2$ from $0$ up to $7.846\mathrm{e}{8}$.]{\includegraphics[scale=.35]{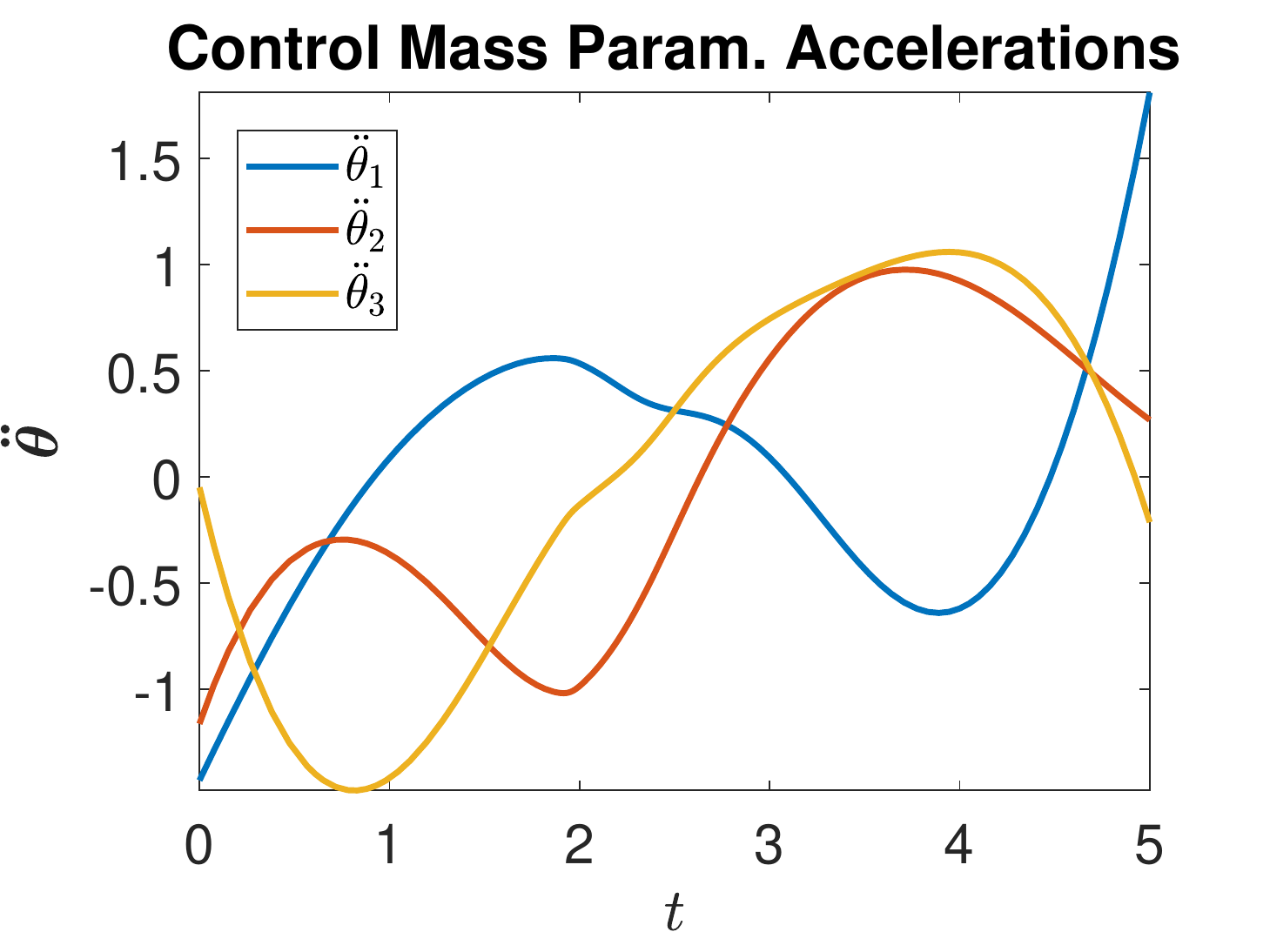}\label{fig_bsim12_pc_h_controls}}
	\hspace{5mm}
	\subfloat[The controls after relaxing the control coefficients $\gamma_1=\gamma_2=\gamma_3$ from $10$ down to $3.602\mathrm{e}{-5}$.]{\includegraphics[scale=.35]{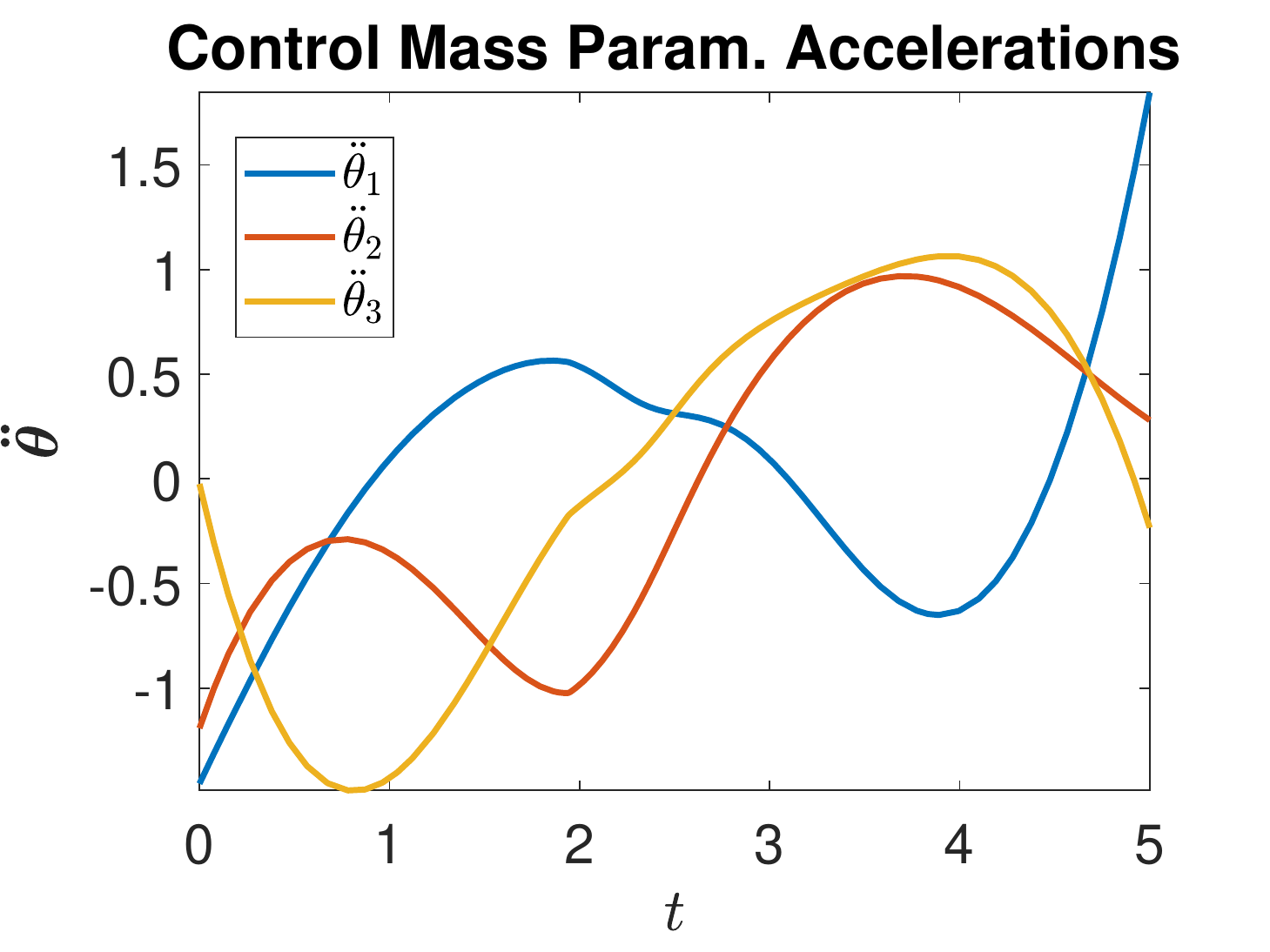}\label{fig_bsim12_pc_gamma_controls}}
	\caption{Direct and predictor-corrector continuation indirect methods realize a $\ReLU$ obstacle avoidance maneuver for the rolling ball. First, a direct method (left column) solves to provide an initial solution that plows through the obstacles. Then, predictor-corrector continuation in the obstacle heights $h_1=h_2$ (middle column) solves to transform the ball's trajectory to the obstacle boundaries. Finally, predictor-corrector continuation in the control coefficients $\gamma_1=\gamma_2=\gamma_3$ (right column) slightly perturbs the control masses so that the ball's trajectory is just outside the obstacle boundaries. }
	
	\label{fig_bsim12}
\end{figure}

\begin{figure}[!ht] 
	\centering
	\subfloat[The magnitude of the normal force  when the obstacle heights $h_1=h_2$ are $0$.]{\includegraphics[scale=.35]{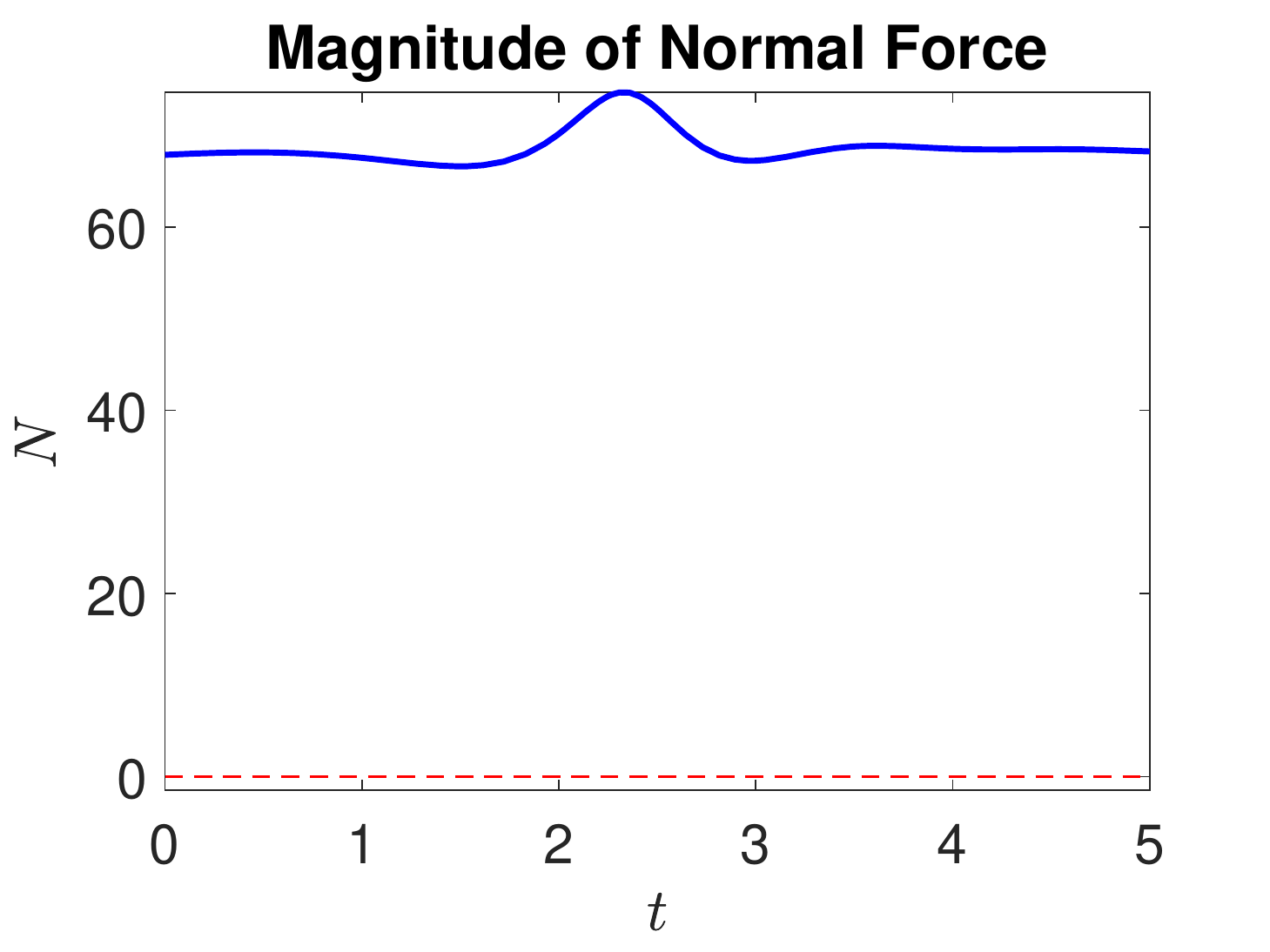}\label{fig_bsim12_dm_normal}}
	\hspace{5mm}
	\subfloat[The magnitude of the normal force after ramping the obstacle heights $h_1=h_2$ from $0$ up to $7.846\mathrm{e}{8}$.]{\includegraphics[scale=.35]{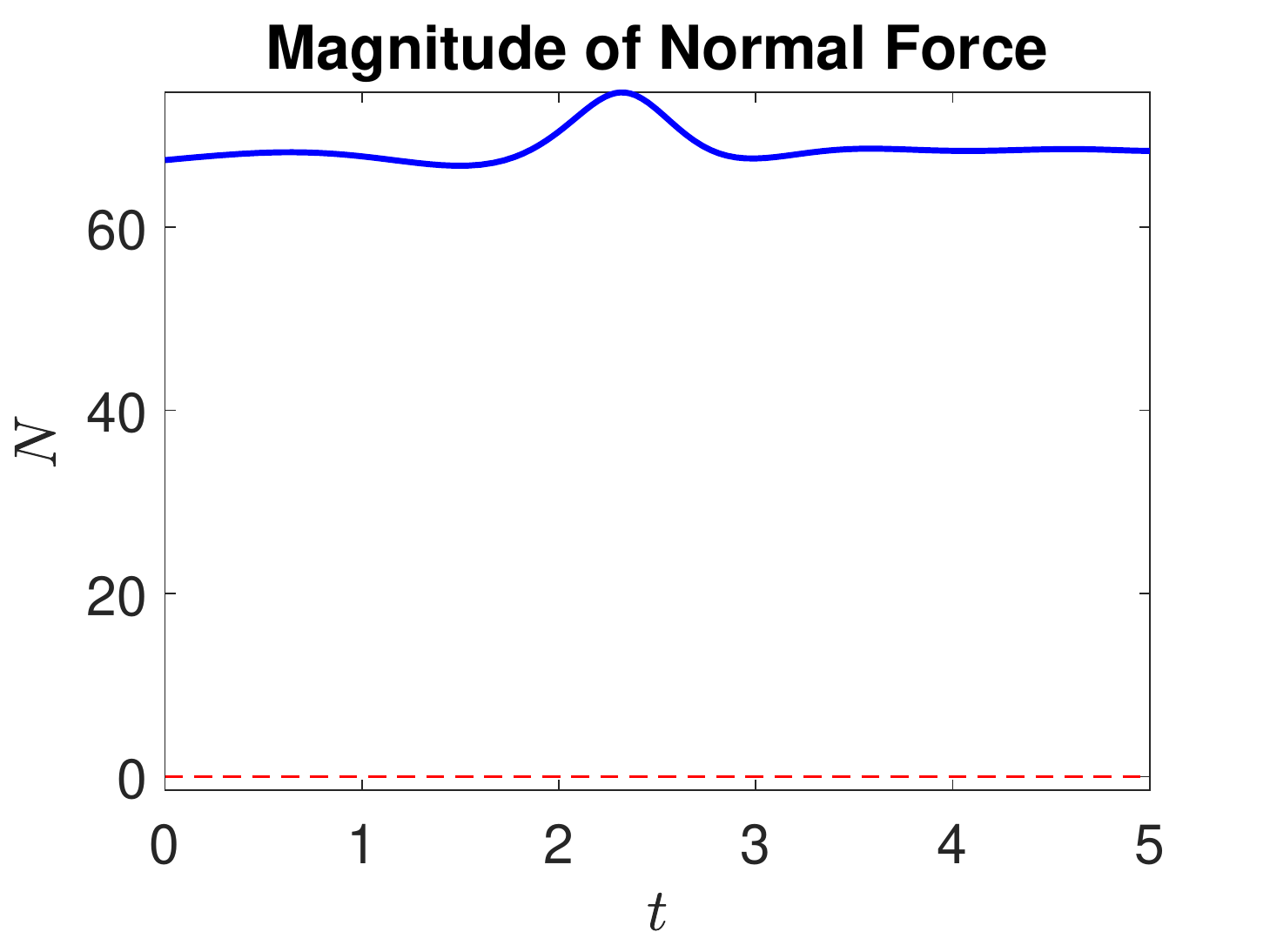}\label{fig_bsim12_pc_h_normal}}
	\hspace{5mm}
	\subfloat[The magnitude of the normal force  after relaxing the control coefficients $\gamma_1=\gamma_2=\gamma_3$ from $10$ down to $3.602\mathrm{e}{-5}$.]{\includegraphics[scale=.35]{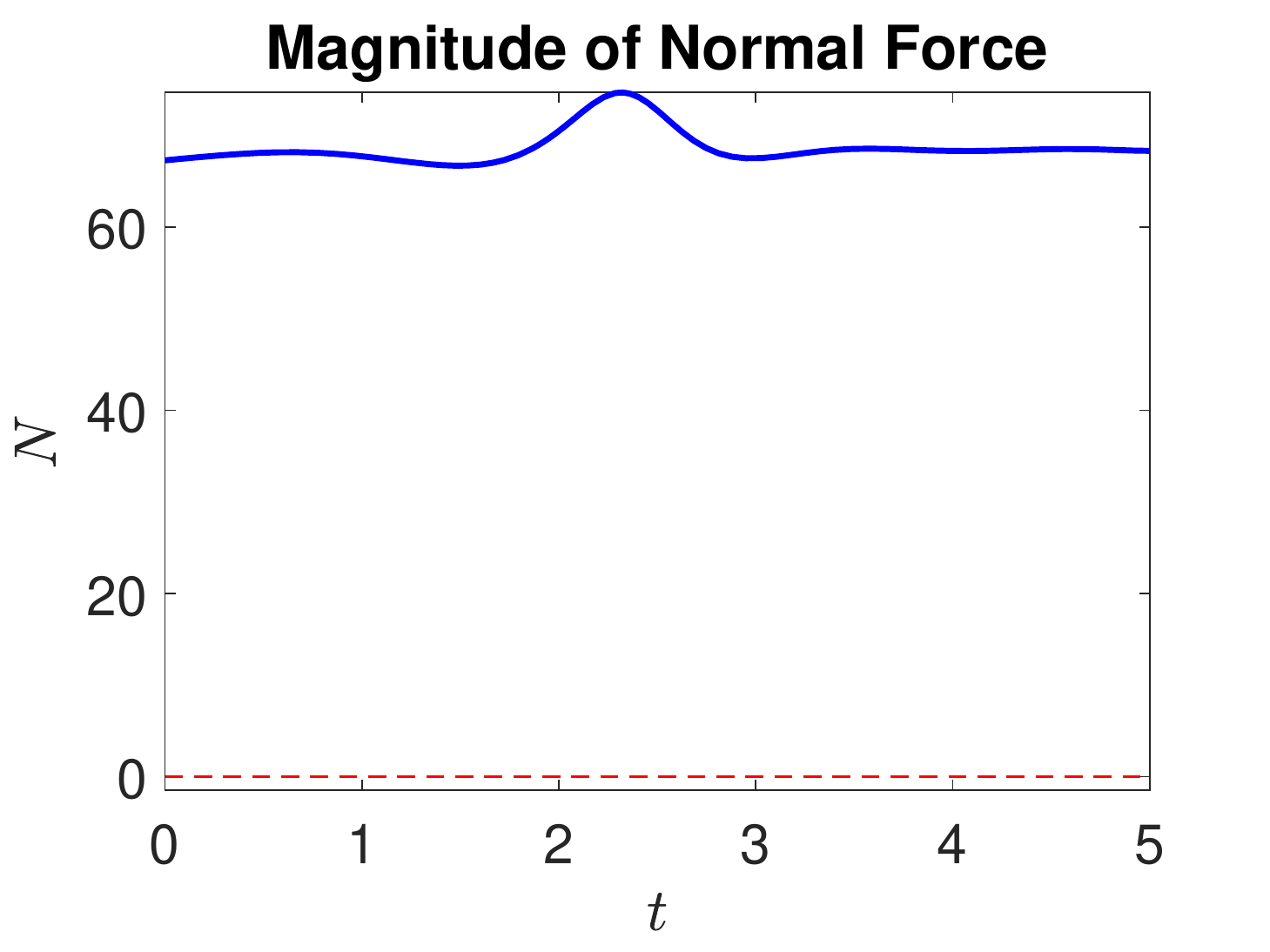}\label{fig_bsim12_pc_gamma_normal}}
	\\
	\subfloat[The minimum coefficient of static friction to prevent slipping is $.1055$ when the obstacle heights $h_1=h_2$ are $0$.]{\includegraphics[scale=.35]{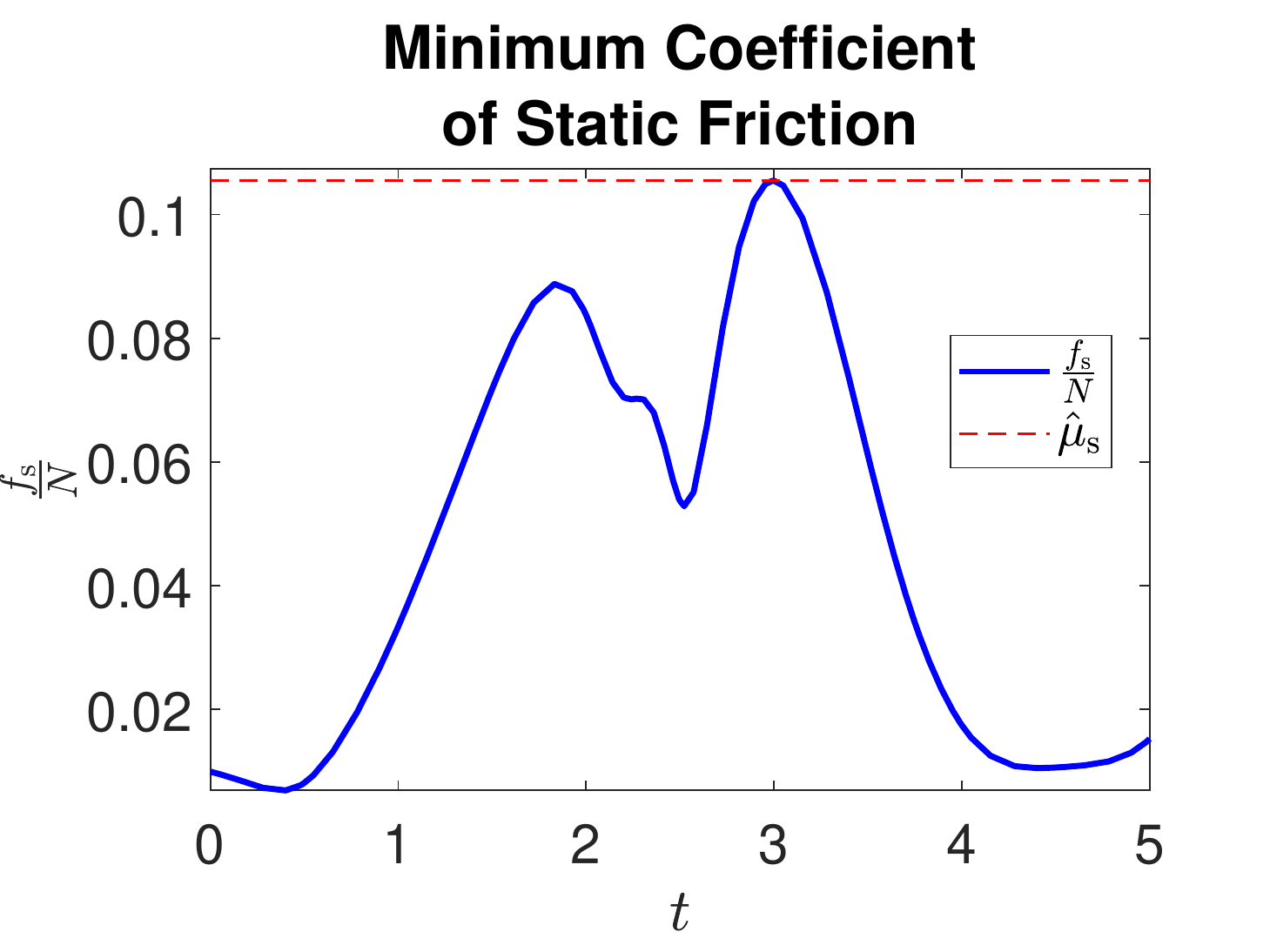}\label{fig_bsim12_dm_mu_s}}
	\hspace{5mm}
	\subfloat[The minimum coefficient of static friction to prevent slipping is $.09942$ after ramping the obstacle heights $h_1=h_2$ from $0$ up to $7.846\mathrm{e}{8}$.]{\includegraphics[scale=.35]{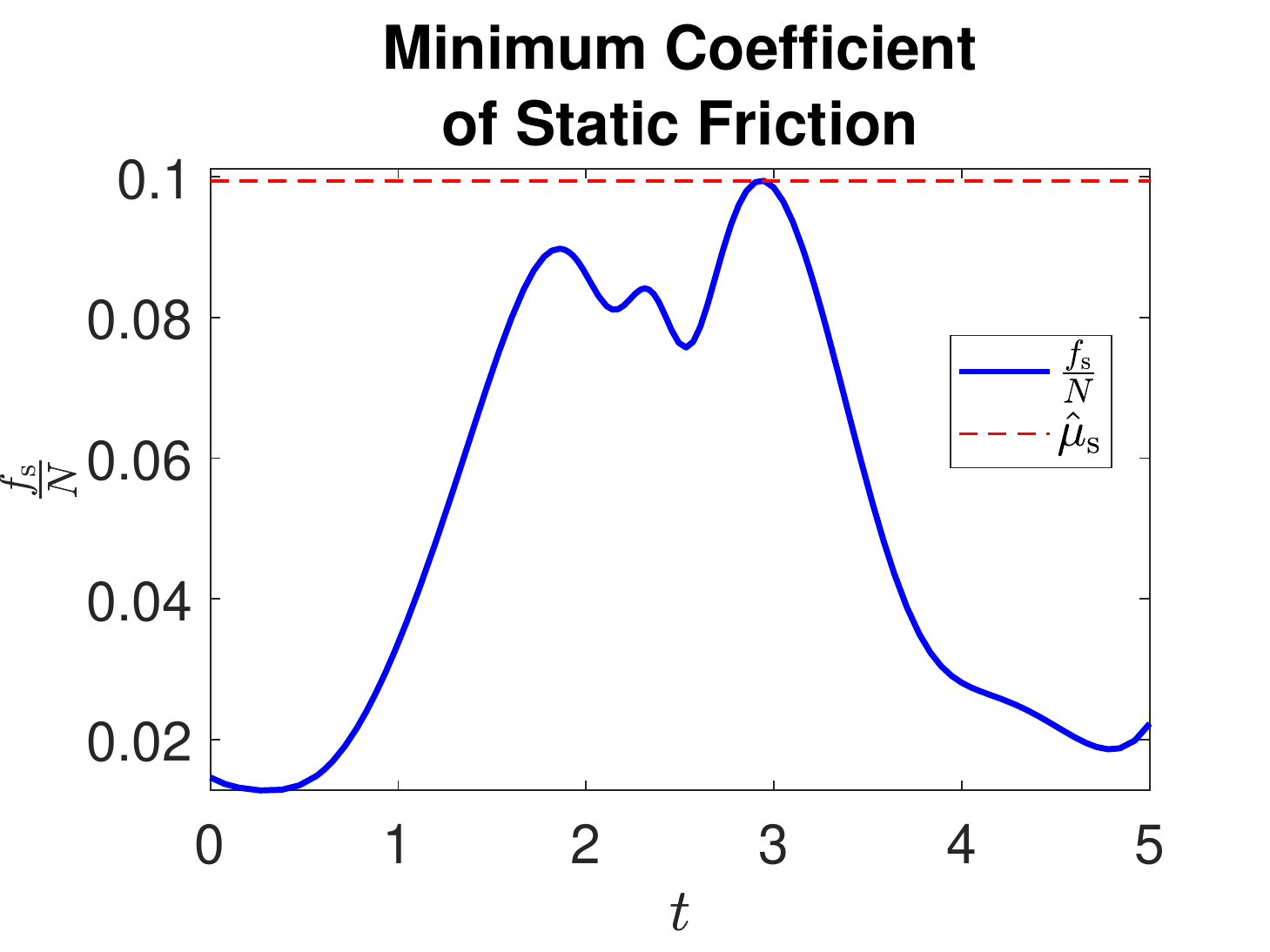}\label{fig_bsim12_pc_h_mu_s}}
	\hspace{5mm}
	\subfloat[The minimum coefficient of static friction to prevent slipping is $.09917$ after relaxing the control coefficients $\gamma_1=\gamma_2=\gamma_3$ from $10$ down to $3.602\mathrm{e}{-5}$.]{\includegraphics[scale=.35]{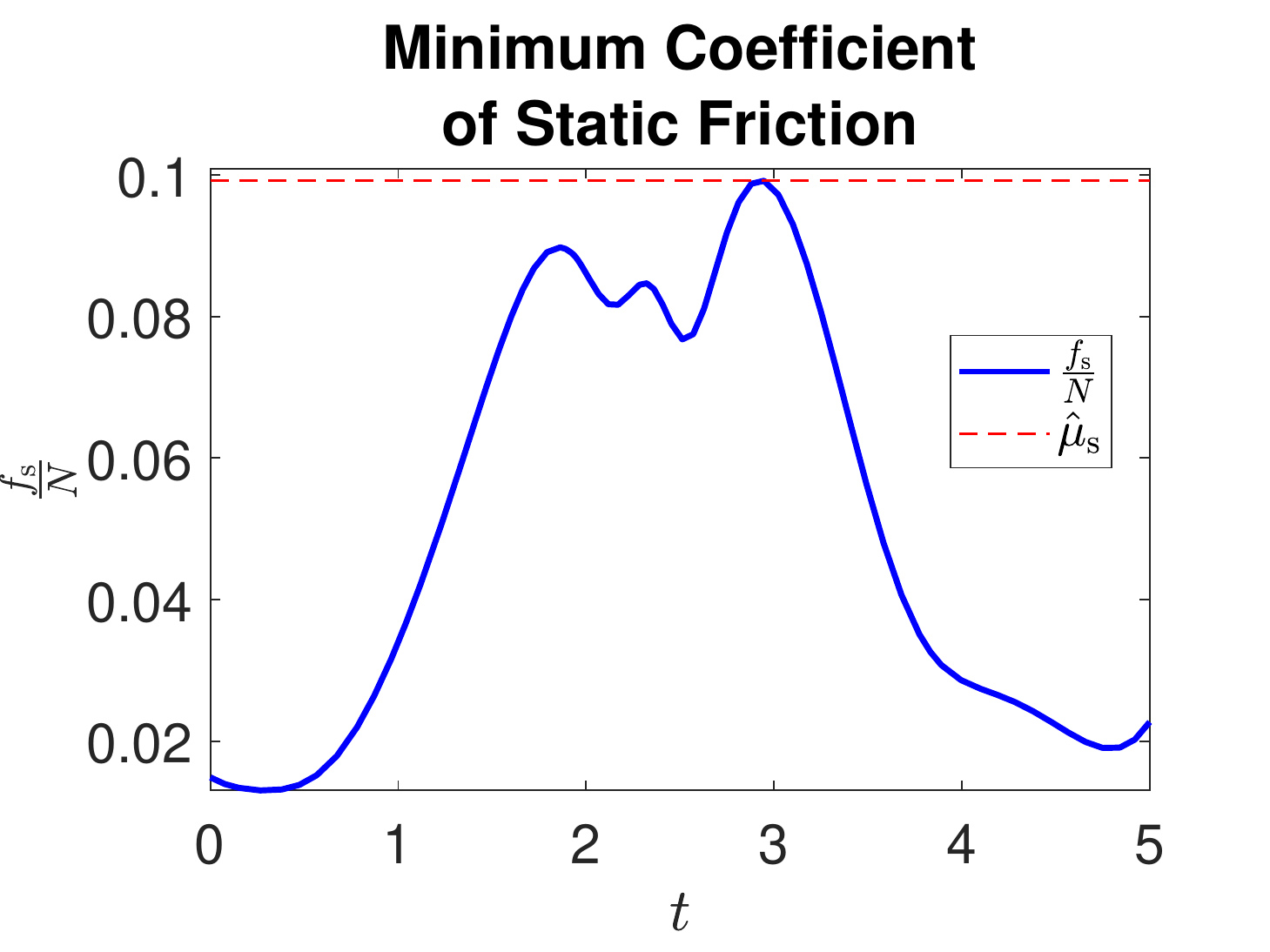}\label{fig_bsim12_pc_gamma_mu_s}}
	\caption{A direct method (left column) is followed by two rounds of a predictor-corrector continuation indirect method to realize a $\ReLU$ obstacle avoidance maneuver for the rolling ball. The first round (middle column) of predictor-corrector continuation increases the obstacle heights $h_1=h_2$, and the second round (right column) of predictor-corrector continuation decreases the control coefficients $\gamma_1=\gamma_2=\gamma_3$. The ball does not detach from the surface since the magnitude of the normal force is always positive. The ball rolls without slipping if $\mu_\mathrm{s} \ge .1055$ for the direct method solution, if $\mu_\mathrm{s} \ge .09942$ for the first indirect method solution, and if $\mu_\mathrm{s} \ge .09917$ for the second indirect method solution. 
 }
	
	\label{fig_bsim12_normal}
\end{figure}

\begin{figure}[!ht] 
	\centering
	\subfloat[Evolution of the obstacle heights $h_1=h_2$ as they increase from $0$ up to $7.846\mathrm{e}{8}$.]{\includegraphics[scale=.5]{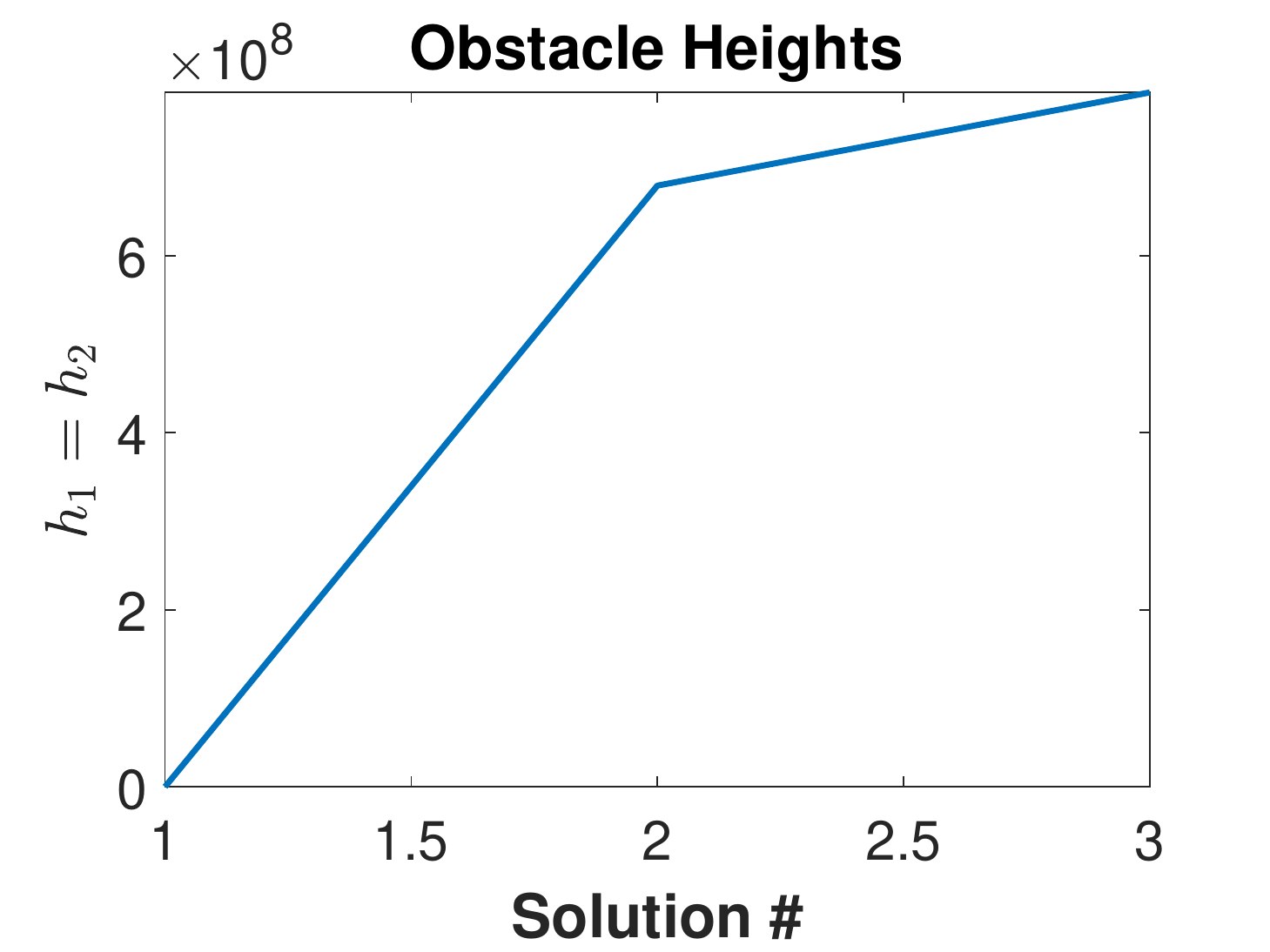}\label{fig_bsim12_pc_h_obst_height}}
	\hspace{5mm}
	\subfloat[Evolution of the control coefficients $\gamma_1=\gamma_2=\gamma_3$ as they decrease from $10$ down to $3.602\mathrm{e}{-5}$.]{\includegraphics[scale=.5]{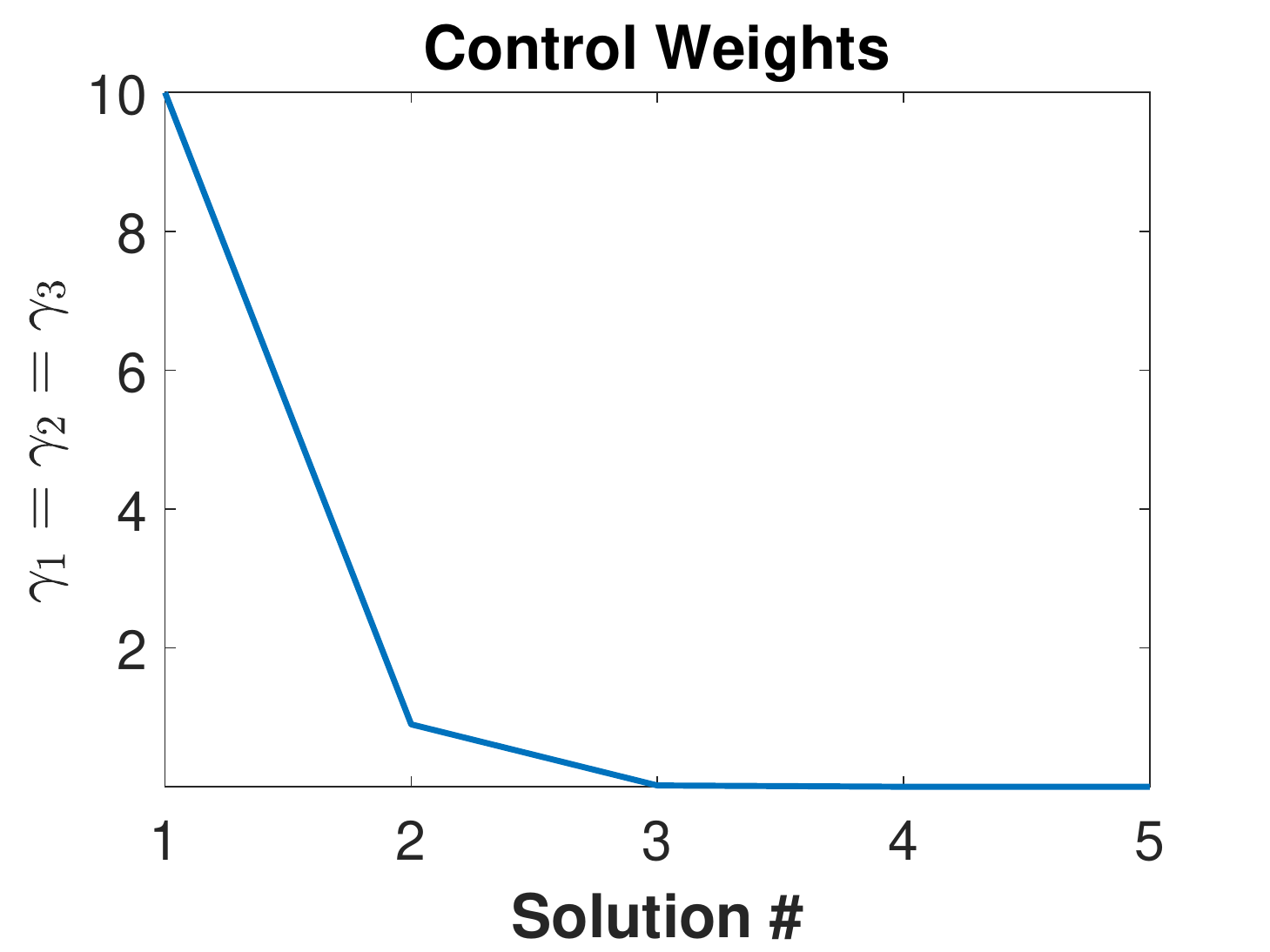}\label{fig_bsim12_pc_gamma_weights}}
	\\
	\subfloat[Evolution of the performance index $J$ from $284.2$ up to $289.7$ as the obstacle heights $h_1=h_2$ increase from $0$ up to $7.846\mathrm{e}{8}$.]{\includegraphics[scale=.5]{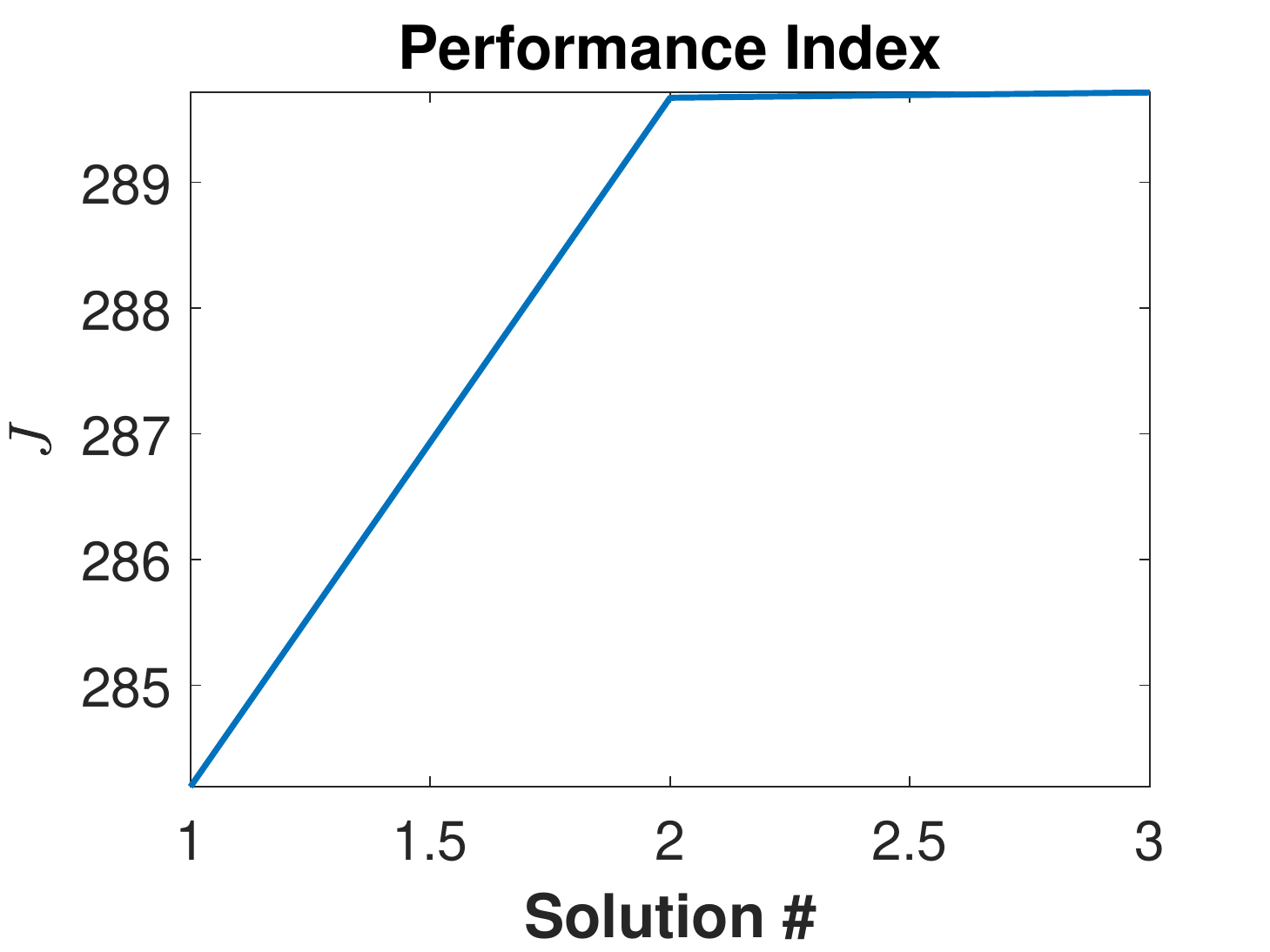}\label{fig_bsim12_pc_h_J}}
	\hspace{5mm}
	\subfloat[Evolution of the performance index $J$ from $289.7$ down to $250$ as the control coefficients $\gamma_1=\gamma_2=\gamma_3$  decrease from $10$ down to $3.602\mathrm{e}{-5}$.]{\includegraphics[scale=.5]{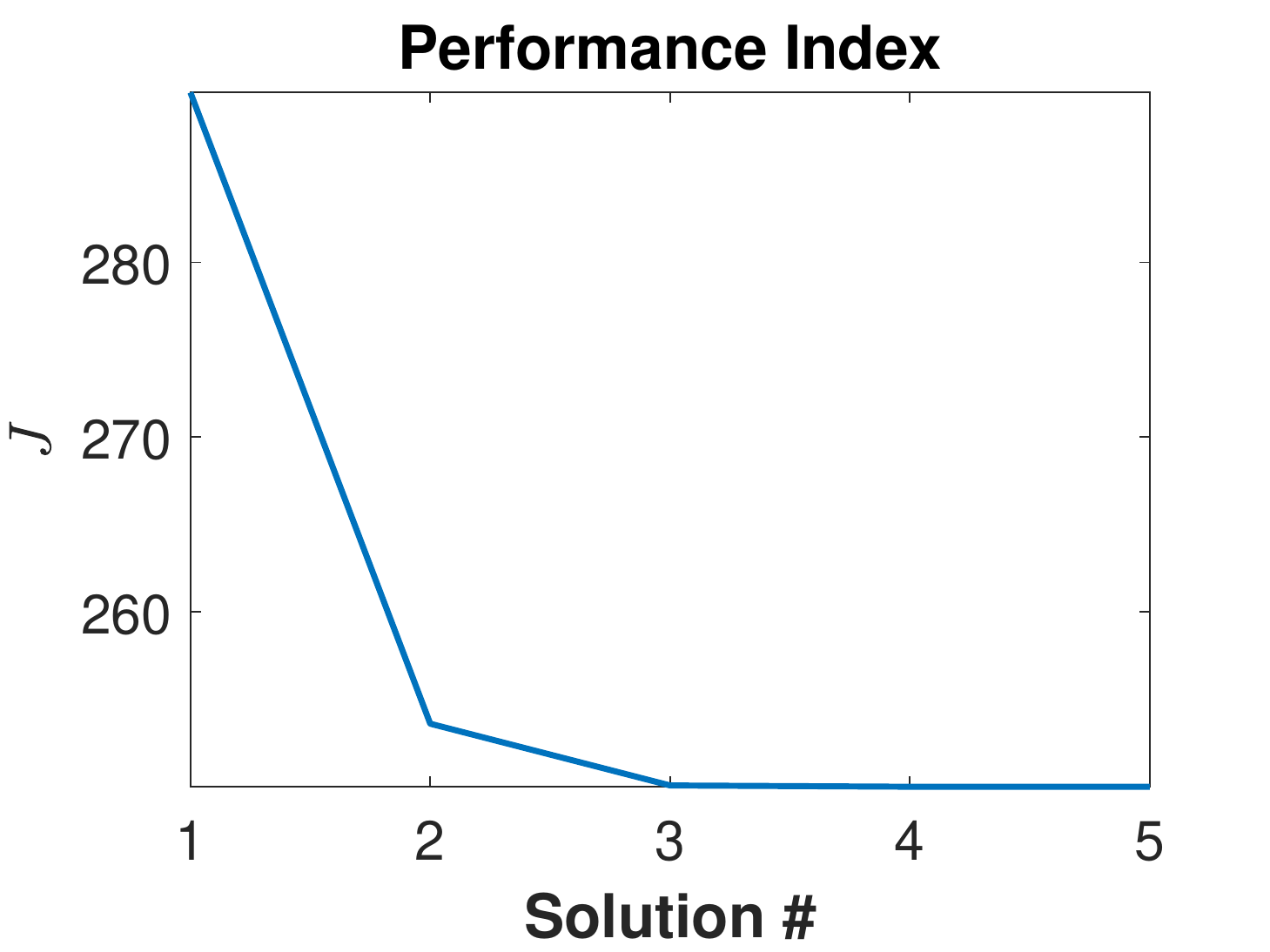}\label{fig_bsim12_pc_gamma_J}}
		\\
	\subfloat[Evolution of the tangent steplength $\sigma$ as the obstacle heights $h_1=h_2$ increase from $0$ up to $7.846\mathrm{e}{8}$.]{\includegraphics[scale=.5]{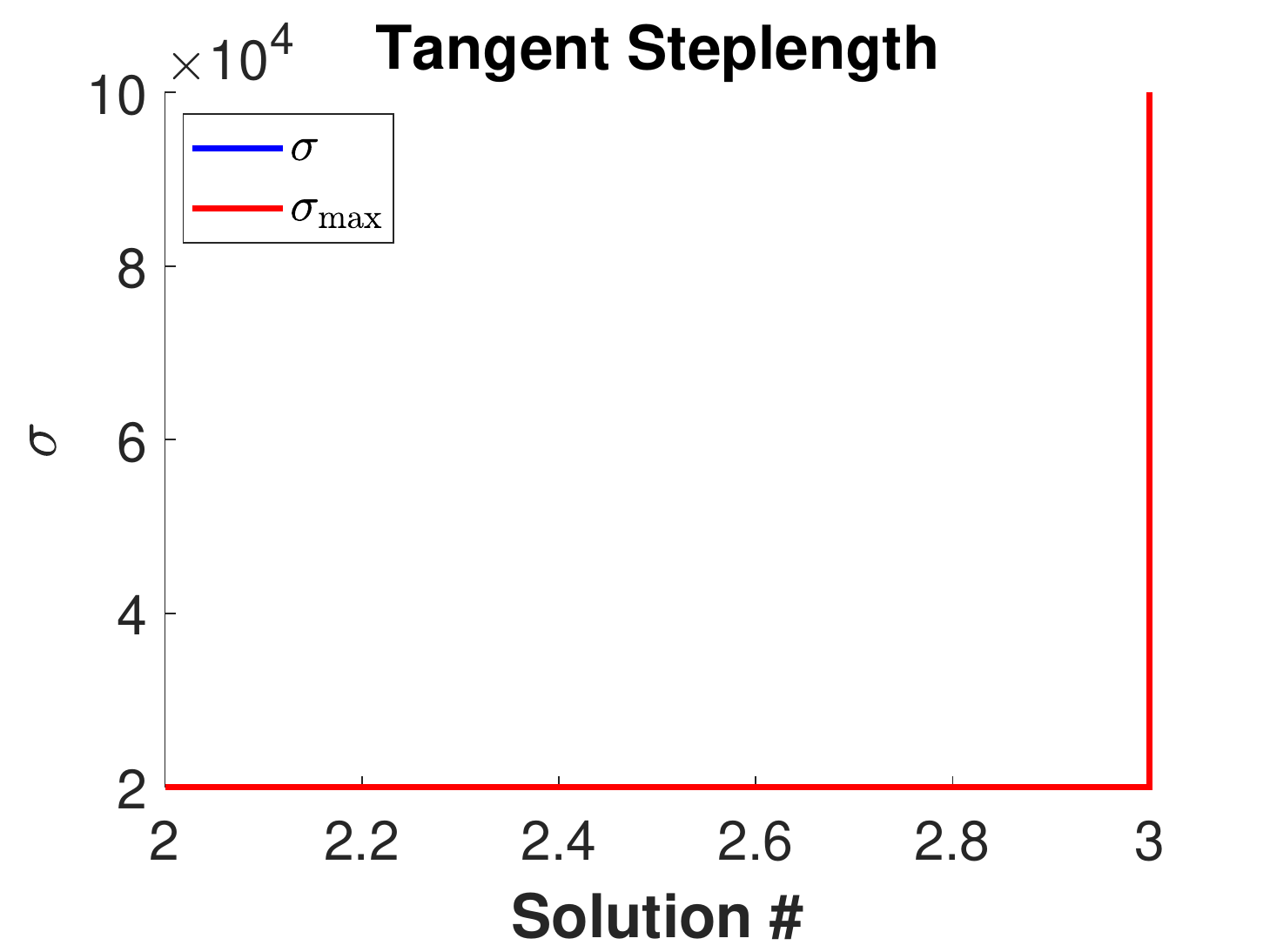}\label{fig_bsim12_pc_h_sigma}}
	\hspace{5mm}
	\subfloat[Evolution of the tangent steplength $\sigma$ as the control coefficients $\gamma_1=\gamma_2=\gamma_3$ decrease from $10$ down to $3.602\mathrm{e}{-5}$.]{\includegraphics[scale=.5]{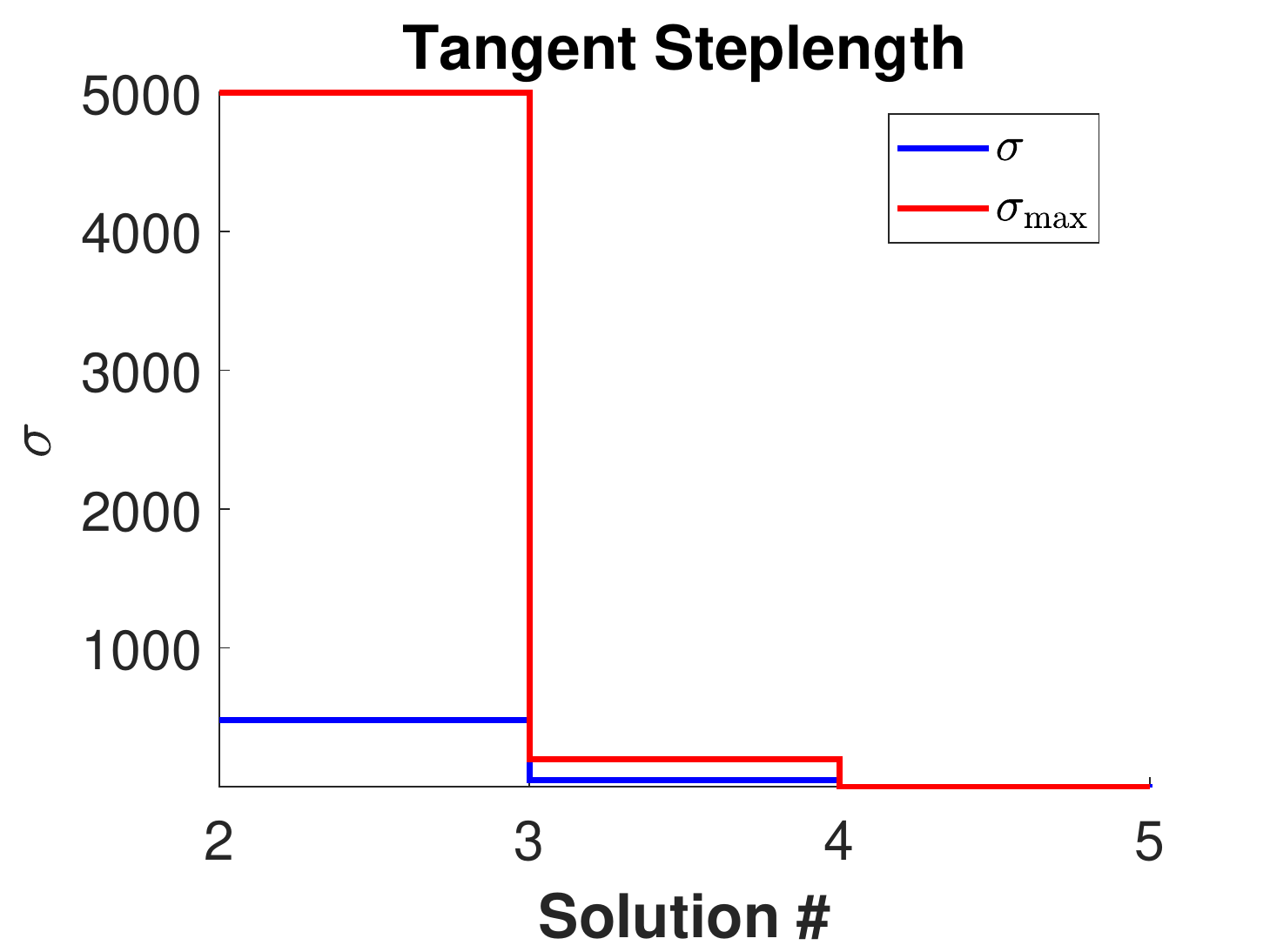}\label{fig_bsim12_pc_gamma_sigma}}
	\caption{Predictor-corrector continuation in the obstacle heights $h_1=h_2$ (left column) is followed by predictor-corrector continuation in the control coefficients $\gamma_1=\gamma_2=\gamma_3$ (right column) to realize a $\ReLU$ obstacle avoidance maneuver for the rolling ball.}
	\label{fig_bsim12_pc}
\end{figure}

The numerical results are shown in Figures~\ref{fig_bsim12}, \ref{fig_bsim12_normal}, and \ref{fig_bsim12_pc}. As $\mu$ decreases from $.95$ down to $\unaryminus 7.453\mathrm{e}{5}$ during the first round of predictor-corrector continuation, $h_1=h_2$ increases from $0$ up to $7.846\mathrm{e}{8}$ (see Figure~\ref{fig_bsim12_pc_h_obst_height}). Since $h_1=h_2$ is ratcheted up during continuation, thereby increasing the penalty in the integrand cost function \eqref{eq_ball_integrand_cost} when the GC intrudes into the obstacles, by the end of continuation, the ball's GC completely exits the second obstacle and approaches the boundary of the first obstacle (compare Figures~\ref{fig_bsim12_dm_gc_path} vs \ref{fig_bsim12_pc_h_gc_path}). As $\mu$ decreases from $.95$ down to $.9406$ during the second round of predictor-corrector continuation, $\gamma_1=\gamma_2=\gamma_3$ decreases from $10$ down to $3.602\mathrm{e}{-5}$ (see Figure~\ref{fig_bsim12_pc_gamma_weights}). Since $\gamma_1=\gamma_2=\gamma_3$ is ratcheted down during continuation, thereby decreasing the penalty in the integrand cost function \eqref{eq_ball_integrand_cost} for large magnitude accelerations of the control mass parameterizations, and since the obstacle heights are held fixed at $7.846\mathrm{e}{8}$, by the end of continuation, the ball's GC avoids both obstacles while veering smartly around the first obstacle (compare Figures~\ref{fig_bsim12_pc_h_gc_path} vs \ref{fig_bsim12_pc_gamma_gc_path}). Figure~\ref{fig_bsim12_pc_h_J}  shows that the performance index $J$ increases from $284.2$ up to $289.7$ as the obstacle heights are ramped up in the first round of predictor-corrector continuation; Figure~\ref{fig_bsim12_pc_gamma_J} shows that the performance index $J$ then decreases down to $250$ as the control coefficients are ramped down and the ball's GC fully departs the tall obstacles in the second round of predictor-corrector continuation.

\begin{table}[h!]
	\centering 
	{ 
		\setlength{\extrarowheight}{1.5pt}
		\begin{tabular}{| c | c |} 
			\hline
			\textbf{Parameter} & \textbf{Value} \\ 
			\hline\hline 
			$\gamma_1(\mu)=\gamma_2(\mu)=\gamma_3(\mu)$ & $10+\frac{.95-\mu}{.95-.00001}\left(-1000-10\right)$ \\
			\hline
			$h_1=h_2$ & $7.846\mathrm{e}{8}$ \\ 
			\hline
			$\bv_1$ & $\begin{bmatrix}  .2 & .2 \end{bmatrix}^\mathsf{T}$ \\
			\hline
			$\bv_2$ & $\begin{bmatrix}  .8 & .8 \end{bmatrix}^\mathsf{T}$ \\
			\hline
			$\rho_1=\rho_2$ & $.282$ \\
			\hline
		\end{tabular} 
	}
	\caption{Integrand cost function coefficient values for the rolling ball when a second round of predictor-corrector continuation is performed in the control coefficients. Refer to \eqref{eq_ball_integrand_cost}.}
	\label{table_ball_integrand_pc2}
\end{table}

\section{Summary, Discussion, and Future Work} \label{sec_conclusions}

The controlled equations of motion for the rolling disk and ball were solved numerically using predictor-corrector continuation, starting from an initial solution obtained via a direct method, to solve trajectory tracking problems for the rolling disk and obstacle avoidance problems for the rolling ball. These optimal control maneuvers were achieved by performing predictor-corrector continuation in weighting factors that scale penalty functions in the integrand cost function of the performance index. 
 
This paper focused on the indirect, rather than direct, method to numerically solve the optimal control problems. Because the indirect and direct methods only converge to a local minimum solution near the initial guess, a robust continuation algorithm capable of handling turning points is needed to obtain indirect and direct method solutions of complicated, nonconvex optimal control problems. A continuation indirect method requires a continuation ODE or DAE TPBVP solver, while a continuation direct method requires a continuation NLP solver. Predictor-corrector continuation ODE TPBVP algorithms were presented in Appendices~\ref{app_predictor_corrector} and \ref{app_sweep_predictor_corrector} and implemented in \mcode{MATLAB} to realize the continuation indirect method used to solve the rolling disk and ball optimal control problems. Even though predictor-corrector continuation NLP solver algorithms are provided in the literature (e.g. see \cite{zangwill1981pathways,kungurtsev2017predictor}), there do not seem to be any publicly available predictor-corrector continuation NLP solvers, which inhibited the use of a continuation direct method in this paper. When compared against the direct method, the indirect method suffers from two major deficiencies:
\begin{enumerate}
	\item Unlike the direct method, the indirect method has a very small radius of convergence and therefore requires a very accurate initial solution guess \cite{betts2010practical,bryson1999dynamic,BrHo1975applied}. Moreover, unlike the direct method, the indirect method requires a guess of the costates, which are unphysical.
	\item Unlike the direct method, the indirect method is unable to construct the switching structure (i.e. the times when the states and/or controls enter and exit the boundary) of an optimal control problem having path inequality constraints.
\end{enumerate}
Since predictor-corrector continuation was used in this paper, the first deficiency in the indirect method only applied when constructing the solution of the initial ODE TPBVP, and this deficiency was circumvented by using a direct method to solve the optimal control problem corresponding to that initial ODE TPBVP. To circumvent the second deficiency in the indirect method, path inequality constraints were incorporated into the optimal control problems as soft constraints through penalty functions in the integrand cost functions.

The predictor-corrector continuation methods presented in Appendices~\ref{app_predictor_corrector} and \ref{app_sweep_predictor_corrector} work if the control can be expressed analytically as a function of the state and costate, e.g. if the Hamiltonian is quadratic in the control. That is, those methods perform continuation only in the state and costate after the original optimal control DAE TPBVP is transformed into an ODE TPBVP through elimination of the control. For more complicated Hamiltonians (such as when penalty functions are added to the integrand cost function to softly enforce bounded body frame accelerations of the control masses, the no-detachment constraint, and the no-slip constraint), numerical methods (such as Newton's method) must be used to construct the control numerically from the state, costate, and a good initial guess of the control. In these cases, the predictor-corrector continuation method of Appendix~\ref{app_predictor_corrector} must be extended to perform continuation in the state, costate, and control by  solving the optimal control DAE TPBVP. This will be investigated in subsequent work.

In future work, instead of using \mcode{MATLAB}, the simulation code could be reimplemented in the higher performance programming languages Julia or C++, while relying on Fortran routines like COLNEW \cite{bader1987new}, COLMOD \cite{cash2001automatic}, TWPBVP(C) \cite{cash1991deferred,cash2005new}, TWPBVPL(C) \cite{bashir1998lobatto,cash2006hybrid}, and ACDC(C) \cite{cash2001automatic,cash2013algorithm} to solve the underlying ODE TPBVPs, to obtain faster numerical results. Julia and C++ feature several mature and efficient automatic differentiation libraries \cite{autodiff} capable of constructing the Jacobians and Hessians needed by the ODE TPBVP solvers. In addition, a more efficient and robust predictor-corrector adaptive tangent steplength algorithm, such as described in \cite{bader1989continuation,deuflhard2011newton}, could be implemented.   

Another avenue for future investigation is to use a neighboring extremal optimal control (NEOC) method \cite{gupta2017combined}, which constructs a homotopy between the controlled equations of motion and their linearization about a nominal solution; however, the NEOC method in \cite{gupta2017combined} could be made more robust by using predictor-corrector, rather than monotonic, continuation in the homotopy parameter. Yet another avenue for future investigation is to perform  predictor-corrector continuation in a weighting factor that scales a term in the endpoint cost function measuring the deviation between the actual and prescribed final conditions. 

\revision{R1Q20}{Additionally, throughout the paper we have kept the initial and final times fixed. It would be interesting to perform additional studies when the time duration is free, as outlined in Appendix~\ref{sec_optimal_control}, especially regarding problems of navigation over complex terrains and slippery surfaces.}
Complex terrains will affect both the  uncontrolled equations of motion, due to gravity, and the performance index, for example, through potential energy penalty functions that discourage the ball from  ascending steep slopes. If the substrate is slippery, for example, due to the presence of moisture, one can imagine situations where the valleys are wet and the slopes are dry and thus less slippery. Then, one can introduce an additional term in the performance index penalizing motion through the valleys where the possibility of a slip is high. The interplay between the terms in the performance index  discouraging and encouraging motion up the slopes will give a very interesting control system.  These and other interesting questions will be considered in future work.

\section*{Acknowledgements}

We are indebted to our colleagues  A.M. Bloch, D.M. de Diego, F. Gay-Balmaz, D.D. Holm, M. Leok, A. Lewis, T. Ohsawa, and D.V. Zenkov for useful and fruitful discussions. 
M.J. Weinsten provided copious advice on using the \mcode{MATLAB} automatic differentiation toolbox \mcode{ADiGator} and fixed numerous bugs in \mcode{ADiGator} that were revealed in the course of this research. A.V. Rao provided a free license to use the \mcode{MATLAB} direct method optimal control solver \mcode{GPOPS-II} for some of this research. \revisionNACO{R1Q5}{P. Tallapragada observed that the reaction forces exerted on the ball by the accelerating internal point masses  may cause the ball to detach from the surface, which prompted the inclusion of plots depicting the normal force and minimum coefficient of static friction. G.M. Rozenblat pointed out the necessary and sufficient conditions that must be satisfied by a ball's physical parameters in \cite{rozenblat2016choice}.}

This research was partially  supported by  the NSERC Discovery Grant, the University of Alberta Centennial Fund, and the Alberta Innovates Technology Funding (AITF) which came through the Alberta Centre for Earth Observation Sciences (CEOS). S.M. Rogers also received support from the University of Alberta Doctoral Recruitment Scholarship, the FGSR Graduate Travel Award, the IGR Travel Award, the GSA Academic Travel Award, the AMS Fall Sectional Graduate Student Travel Grant, Target Corporation, and the Institute for Mathematics and its Applications at the University of Minnesota, Twin Cities.

\phantomsection
\addcontentsline{toc}{section}{References}
\printbibliography
\hypertarget{References}{}

\appendix

\section{Optimal Control: Variational Pontryagin's Minimum Principle} \label{sec_optimal_control}

\rem{ 
\subsection{Calculus of Variations} \label{ssec_calc_var}

Before proceeding with the variational Pontryagin's minimum principle, some terminology from the calculus of variations is briefly reviewed. Suppose that $y$ is a time-dependent function, $w$ is a time-independent variable, and $Q$ is a scalar-valued function or functional that depends on $y$ and $w$. The variation of $y$ is $\de y  \equiv \left. \pp {y}{\epsilon} \right|_{\epsilon=0}$, the differential of $y$ is $\mathrm{d} y \equiv \de y + \dot y \mathrm{d}t = \left. \pp {y}{\epsilon} \right|_{\epsilon=0}+\left. \pp{y}{t} \right|_{\epsilon=0} \mathrm{d}t$, and the differential of $w$ is $\mathrm{d} w \equiv \left. \dd {w}{\epsilon} \right|_{\epsilon=0}$, where $\epsilon$ represents an independent ``variational'' variable. The variation of $Q$ with respect to $y$ is $\de_y Q \equiv \pp{Q}{y} \de y$, while the differential of $Q$ with respect to $w$ is $\mathrm{d}_{w} Q \equiv \pp{Q}{w} \mathrm{d} w$. The total differential (or for brevity ``the differential'') of $Q$ is $\mathrm{d} Q \equiv \de_y Q + \mathrm{d}_{w} Q =\pp{Q}{y} \de y + \pp{Q}{w} \mathrm{d} w$. Colloquially, the variation of $Q$ with respect to $y$ means the change in $Q$ due to a small change in $y$, the differential of $Q$ with respect to $w$ means the change in $Q$ due to a small change in $w$, and the total differential of $Q$ means the change in $Q$ due to small changes in $y$ and $w$. The extension to vectors of time-dependent functions and time-independent variables is staightforward. If $\mathbf{y}$ is a vector of time-dependent functions, $\mathbf{w}$ is a vector of time-independent variables, and $Q$ is a scalar-valued function or functional depending on $\mathbf{y}$ and $\mathbf{w}$, then the variation of $\mathbf{y}$ is $\de \mathbf{y} \equiv \left.  \pp {\mathbf{y}}{\epsilon} \right|_{\epsilon=0}$, the differential of $\mathbf{y}$ is $\mathrm{d} \mathbf{y} \equiv \de \mathbf{y} + \dot {\mathbf{y}} \mathrm{d}t = \left. \pp {\mathbf{y}}{\epsilon} \right|_{\epsilon=0}+\left. \pp{\mathbf{y}}{t} \right|_{\epsilon=0} \mathrm{d}t$, the differential of $\mathbf{w}$ is $\mathrm{d} \mathbf{w} \equiv \left. \dd {\mathbf{w}}{\epsilon} \right|_{\epsilon=0}$,   the variation of $Q$ with respect to $\mathbf{y}$ is $\de_\mathbf{y} Q \equiv \pp{Q}{\mathbf{y}} \de \mathbf{y}$, the differential of $Q$ with respect to $\mathbf{w}$ is $\mathrm{d}_{\mathbf{w}} Q \equiv \pp{Q}{\mathbf{w}} \mathrm{d} \mathbf{w}$, and the total differential (or for brevity ``the differential'') of $Q$ is $\mathrm{d} Q \equiv \de_\mathbf{y} Q + \mathrm{d}_{\mathbf{w}} Q =\pp{Q}{\mathbf{y}} \de \mathbf{y} + \pp{Q}{\mathbf{w}} \mathrm{d} \mathbf{w}$.

To illustrate these definitions, consider the integral $I$ of the function $F\left(t,\mathbf{y}(t)\right)$ with respect to $t$ with free upper limit of integration $b$ and free lower limit of integration $a$
\begin{equation}
I = \int_a^b F\left(t,\mathbf{y}(t)\right) \dt.
\end{equation}
Applying the above definitions and using the Fundamental Theorem of Calculus, the differential of $I$ is
\begin{equation} \label{eq_leibnitz}
\begin{split}
\mathrm{d} I &= \de_\mathbf{y} \int_a^b F \left(t,\mathbf{y}(t)\right) \dt + \pp{}{b} \left[ \int_a^b F\left(t,\mathbf{y}(t)\right) \dt \right] \mathrm{d}b + \pp{}{a} \left[ \int_a^b F\left(t,\mathbf{y}(t)\right) \dt \right] \mathrm{d}a \\
&= \int_a^b \de_\mathbf{y} F \left(t,\mathbf{y}(t)\right) \dt + F\left(b,\mathbf{y}(b)\right) \mathrm{d}b - F\left(a,\mathbf{y}(a)\right) \mathrm{d}a = \int_a^b \pp{F}{\mathbf{y}} \de \mathbf{y} \dt + \left[ F\left(t,\mathbf{y}(t)\right) \dt \right]_a^b,
\end{split}
\end{equation}
which is Leibnitz's Rule.
} 

This appendix presents necessary conditions, called the variational Pontryagin's minimum principle, which a solution to an optimal control problem lacking path inequality constraints must satisfy; there is a more general version of Pontryagin's minimum principle that applies to optimal control problems possessing path inequality constraints. In this paper, these necessary conditions, in the context of describing the optimal control of the rolling ball, are referred to as the controlled equations of motion. In the literature, application of Pontryagin's minimum principle to solve an optimal control problem is called the indirect method. Let $n,m \in \mathbb{N}$. Let $a$ be a prescribed or free initial time and let $k_1 \in \mathbb{N}^0$ be such that $0 \le k_1 \le n$ if $a$ is prescribed and $1 \le k_1 \le n+1$ if $a$ is free. Let $b$ be a prescribed or free final time and let $k_2 \in \mathbb{N}^0$ be such that $0 \le k_2 \le n$ if $b$ is prescribed and $1 \le k_2 \le n+1$ if $b$ is free. Suppose a dynamical system has state $\bx \in \mathbb{R}^n$ and control $\bu \in \mathbb{R}^m$ and the control $\bu$ is sought that minimizes the performance index
\begin{equation} \label{eq_PI}
J \equiv p\left(a,\bx(a),b,\bx(b),\mu\right)+ \int_a^b L\left(t,\bx,\bu,\mu\right) \dt
\end{equation}
subject to the system dynamics defined for $a \le t \le b$
\begin{equation} \label{eq_pmp_sys_dyn}
\dot {\bx} = \mathbf{f}\left(t,\bx,\bu,\mu\right),
\end{equation}
the prescribed initial conditions at time $t=a$
\begin{equation}
\bsigma\left(a,\bx(a),\mu\right) = \mathbf{0}_{k_1 \times 1},
\end{equation}
and the prescribed final conditions at time $t=b$
\begin{equation}
\bpsi\left(b,\bx(b),\mu\right) = \mathbf{0}_{k_2 \times 1}.
\end{equation}
$p$ is a scalar-valued function called the endpoint cost function, $L$ is a scalar-valued function called the integrand cost function, $\bx$ and $\mathbf{f}$ are $n \times 1$ vector-valued functions, $\bu$ is an $m \times 1$ vector-valued function, $\bsigma$ is a $k_1 \times 1$ vector-valued function, and $\bpsi$ is a $k_2 \times 1$ vector-valued function. $\mu$ is a prescribed scalar parameter which may be exploited to numerically solve this problem via continuation. More concisely, this optimal control problem may be stated as \rem{
\begin{equation}
\min_{a,\bx(a),b,\bu} \left[ p\left(a,\bx(a),b,\bx(b),\mu\right)+  \int_a^b L\left(t,\bx,\bu,\mu\right) \dt \right] 
\mbox{\, s.t. \,}
\left\{
\begin{array}{ll}
\dot {\bx} = \mathbf{f}\left(t,\bx,\bu,\mu\right), \\
\bsigma\left(a,\bx(a),\mu\right) = \mathbf{0}_{k_1 \times 1},\\
\bpsi\left(b,\bx(b),\mu\right) = \mathbf{0}_{k_2 \times 1},
\end{array}
\right.
\label{dyn_opt_problem_pmp1}
\end{equation}
or even more concisely as}
\begin{equation}
\min_{a,\bx(a),b,\bu} J 
\mbox{\, s.t. \,}
\left\{
\begin{array}{ll}
\dot {\bx} = \mathbf{f}\left(t,\bx,\bu,\mu\right), \\
\bsigma\left(a,\bx(a),\mu\right) = \mathbf{0}_{k_1 \times 1},\\
\bpsi\left(b,\bx(b),\mu\right) = \mathbf{0}_{k_2 \times 1}.
\end{array}
\right.
\label{dyn_opt_problem_pmp}
\end{equation}
Observe that the optimal control problem encapsulated by \eqref{dyn_opt_problem_pmp} ignores path inequality constraints such as $\mathbf{D}\left(t,\bx,\bu,\mu\right) \le \mathbf{0}$, where $\mathbf{D}$ is an $r \times 1$ vector-valued function for $r \in \mathbb{N}^0$. Path inequality constraints can be incorporated into \eqref{dyn_opt_problem_pmp} as soft constraints through penalty functions in the integrand cost function $L$ or the endpoint cost function $p$. By omitting hard path inequality constraints from \eqref{dyn_opt_problem_pmp}, a solution of \eqref{dyn_opt_problem_pmp} does not lie on the boundary of a compact set and the calculus of variations may be applied to derive necessary conditions, called the variational Pontryagin's minimum principle, which a solution of \eqref{dyn_opt_problem_pmp} must satisfy.

\rem{To be concrete, suppose that for a particular optimal control problem the magnitude of the $j^{\mathrm{th}}$ control is bounded by a positive constant $B_j \in \mathbb{R}^{+}$, so that $\left|u_j \right| \le B_j \iff u_j^2 \le B_j^2 \iff u_j^2 - B_j^2 \le 0$, for $1 \le j \le m$. In this case, the path inequality constraint $\mathbf{D}$ is a function of $\bu$ and has the form
\begin{equation} \label{eq_D_example}
\mathbf{D}\left( \bu \right) \equiv \begin{bmatrix} u_1^2 - B_1^2 \\ u_2^2 - B_2^2 \\ \vdots \\ u_m^2 - B_m^2 \end{bmatrix} \le \mathbf{0}_{m \times 1}.
\end{equation}
Emulating the technique utilized in \cite{gupta2017combined}, to incorporate the path inequality constraint \eqref{eq_D_example} into \eqref{dyn_opt_problem_pmp}, the penalty function 
\begin{equation} \label{eq_penalty}
\sum_{j=1}^m \upsilon_j \left[\max\left\{0,u_j^2 - B_j^2\right\}\right]^4
\end{equation} 
could be added to the integrand cost function $L$, where $\upsilon_j \in \mathbb{R}^{+}$ is a  large positive penalty weighting factor which deters the magnitude of $u_j$ from exceeding $B_j$. Note that each summand in \eqref{eq_penalty} is $C^2$ so that \eqref{eq_penalty} is twice differentiable.
}

\rem{To begin the derivation of the variational Pontryagin's minimum principle, the augmented performance index for the optimal control problem \eqref{dyn_opt_problem_pmp} is obtained by adjoining the dynamic, initial, and final constraints to the original performance index \eqref{eq_PI} via Lagrange multipliers:
\begin{equation} \label{eq_aug_PI}
\begin{split}
\tilde{J} &\equiv J+\bxi^\mathsf{T} \bsigma\left(a,\bx(a),\mu\right)+\bnu^\mathsf{T} \bpsi\left(b,\bx(b),\mu\right)+\int_a^b \blam^\mathsf{T} \left(\mathbf{f}\left(t,\bx,\bu,\mu\right)-\dot {\bx} \right)  \dt \\
&= p\left(a,\bx(a),b,\bx(b),\mu\right)+\bxi^\mathsf{T} \bsigma\left(a,\bx(a),\mu\right)+\bnu^\mathsf{T} \bpsi\left(b,\bx(b),\mu\right)\\
&\hphantom{=}+\int_a^b \left[ L\left(t,\bx,\bu,\mu\right) + \blam^\mathsf{T} \left(\mathbf{f}\left(t,\bx,\bu,\mu\right)-\dot {\bx} \right) \right] \dt \\
&= G\left(a,\bx(a),\bxi,b,\bx(b),\bnu,\mu\right)+\int_a^b \left[ H\left(t,\bx,\blam,\bu,\mu\right) -  \blam^\mathsf{T} \dot {\bx}  \right] \dt,
\end{split}
\end{equation}
where the endpoint function $G$ and the Hamiltonian $H$ are defined by
\begin{equation}
\begin{split}
G\left(a,\bx(a),\bxi,b,\bx(b),\bnu,\mu\right) &\equiv p\left(a,\bx(a),b,\bx(b),\mu\right)+\bxi^\mathsf{T} \bsigma\left(a,\bx(a),\mu\right)+\bnu^\mathsf{T} \bpsi\left(b,\bx(b),\mu\right) \\
H\left(t,\bx,\blam,\bu,\mu\right) &\equiv L\left(t,\bx,\bu,\mu\right) + \blam^\mathsf{T} \mathbf{f}\left(t,\bx,\bu,\mu\right),
\end{split}
\end{equation}
and where $\bxi$ is a $k_1 \times 1$ constant Lagrange multiplier vector, $\bnu$ is a $k_2 \times 1$ constant Lagrange multiplier vector, and $\blam$ is an $n \times 1$ time-varying Lagrange multiplier vector. In the literature, the time-varying Lagrange multiplier vector used to adjoin the system dynamics to the integrand cost function is often called the adjoint variable or the costate. Henceforth, the time-varying Lagrange multiplier vector is referred to as the costate and the elements in this vector are referred to as the costates.} 

Define the endpoint function $G$ and the Hamiltonian $H$ by
\begin{equation}
\begin{split}
G\left(a,\bx(a),\bxi,b,\bx(b),\bnu,\mu\right) &\equiv p\left(a,\bx(a),b,\bx(b),\mu\right)+\bxi^\mathsf{T} \bsigma\left(a,\bx(a),\mu\right)+\bnu^\mathsf{T} \bpsi\left(b,\bx(b),\mu\right), \\
H\left(t,\bx,\blam,\bu,\mu\right) &\equiv L\left(t,\bx,\bu,\mu\right) + \blam^\mathsf{T} \mathbf{f}\left(t,\bx,\bu,\mu\right),
\end{split}
\end{equation}
where $\bxi$ is a $k_1 \times 1$ constant Lagrange multiplier vector, $\bnu$ is a $k_2 \times 1$ constant Lagrange multiplier vector, and $\blam$ is an $n \times 1$ time-varying Lagrange multiplier vector. In the literature, the time-varying Lagrange multiplier vector used to adjoin the system dynamics to the integrand cost function is often called the adjoint variable or the costate. Henceforth, the time-varying Lagrange multiplier vector is referred to as the costate and the elements in this vector are referred to as the costates. The necessary conditions \cite{hull2013optimal} on $\bx$, $\blam$, and $\bu$ which a solution of \eqref{dyn_opt_problem_pmp} must satisfy are the DAEs defined for $a \le t \le b$
\begin{equation} \label{eq_pmp_dae}
\begin{split}
\dot {\bx} &= H_{\blam}^\mathsf{T} \left(t,\bx,\blam,\bu,\mu\right) = \mathbf{f}\left(t,\bx,\bu,\mu\right), \\
\dot {\blam} &= - H_{\bx}^\mathsf{T} \left(t,\bx,\blam,\bu,\mu\right), \\
\mathbf{0}_{m \times 1} &= H_{\bu}^\mathsf{T} \left(t,\bx,\blam,\bu,\mu\right),
\end{split}
\end{equation}
the left boundary conditions defined at time $t=a$
\begin{equation} \label{eq_pmp_lbc}
\left. H \right|_{t=a} = G_{a}, \quad \left. \blam \right|_{t=a} = -G_{\bx(a)}^\mathsf{T}, \quad G_{\bxi}^\mathsf{T} = \bsigma\left(a,\bx(a),\mu\right) = \mathbf{0}_{k_1 \times 1},
\end{equation}
and the right boundary conditions defined at time $t=b$
\begin{equation} \label{eq_pmp_rbc}
\left. H \right|_{t=b} = -G_{b}, \quad \left. \blam \right|_{t=b} = G_{\bx(b)}^\mathsf{T}, \quad G_{\bnu}^\mathsf{T} = \bpsi\left(b,\bx(b),\mu\right) = \mathbf{0}_{k_2 \times 1}.
\end{equation}
If the initial time $a$ is prescribed, then the left boundary condition $\left. H \right|_{t=a} = G_{a}$ is dropped. If the final time $b$ is prescribed, then the right boundary condition $\left. H \right|_{t=b} = -G_{b}$ is dropped. The necessary conditions \eqref{eq_pmp_dae}, \eqref{eq_pmp_lbc}, and \eqref{eq_pmp_rbc} constitute a differential-algebraic equation two-point boundary value problem (DAE TPBVP).

\rem {\paragraph{Note on Normal and Abnormal Extremals} There is a slightly more general formulation of the endpoint function and the Hamiltonian where $G\left(a,\bx(a),\xi_0,\bxi,b,\bx(b),\bnu,\mu\right) \equiv \xi_0 p\left(a,\bx(a),b,\bx(b),\mu\right)+\bxi^\mathsf{T} \bsigma\left(a,\bx(a),\mu\right)+\bnu^\mathsf{T} \bpsi\left(b,\bx(b),\mu\right)$ and $H\left(t,\bx,\lambda_0,\blam,\bu,\mu\right) \equiv \lambda_0 L\left(t,\bx,\bu,\mu\right) + \blam^\mathsf{T} \mathbf{f}\left(t,\bx,\bu,\mu\right)$, where $\xi_0$ is a constant Lagrange multiplier scalar and $\lambda_0$ is a time-varying Lagrange multiplier scalar. But it can be shown (Pontryagin showed this in \cite{boltyanskiy1962mathematical}) that $\lambda_0$ must be a nonnegative constant for an optimal solution! If $\lambda_0 > 0$, the extremal is normal, and if $\lambda_0 = 0$, the extremal is abnormal. If $\lambda_0 > 0$, then the Hamiltonian can be normalized so that $\lambda_0 = 1$.
\rem{ \cite{BrHo1975applied,hull2013optimal} do not consider this more general approach in their derivations, and I don't know how hard it is to extend their derivations to include $\xi_0$. The derivation is given in \cite{boltyanskiy1962mathematical} for the Lagrange formulation of the optimal control problem. }

\vspace{4mm}
Taking the differential of $\tilde{J}$ yields
\begin{equation}
\begin{split}
\mrmd \tilde{J} &= G_{a} \mrmd a + G_{\bx(a)} \mrmd \bx(a) + G_{\bxi} \mrmd \bxi +G_{b} \mrmd b +G_{\bx(b)} \mrmd \bx(b) + G_{\bnu} \mrmd \bnu \\
&\hphantom{=}+\left[ \left( H -  \blam^\mathsf{T} \dot {\bx} \right) \dt \right]_a^b+ \int_a^b \left[H_{\bx} \de \bx - \blam^\mathsf{T} \de \dot \bx+ \left(H_{\blam} - {\dot \bx}^\mathsf{T} \right) \de \blam + H_{\bu} \de \bu   \right] \dt \\
&=  \left(G_{a} - \left. H \right|_{t=a} \right) \mrmd a + G_{\bx(a)} \mrmd \bx(a) + G_{\bxi} \mrmd \bxi 
+ \left(G_{b} + \left. H \right|_{t=b} \right) \mrmd b + G_{\bx(b)} \mrmd \bx(b) + G_{\bnu} \mrmd \bnu \\
&\hphantom{=}-\left[ \blam^\mathsf{T} \left( \de \bx +   \dot {\bx} \dt \right) \right]_a^b+ \int_a^b \left[ \left(H_{\bx} + {\dot \blam}^\mathsf{T} \right) \de \bx+ \left(H_{\blam} - {\dot \bx}^\mathsf{T} \right) \de \blam + H_{\bu} \de \bu   \right] \dt \\
&=  \left(G_{a} - \left. H \right|_{t=a} \right) \mrmd a + G_{\bx(a)} \mrmd \bx(a) + G_{\bxi} \mrmd \bxi 
+ \left(G_{b} + \left. H \right|_{t=b} \right) \mrmd b + G_{\bx(b)} \mrmd \bx(b) + G_{\bnu} \mrmd \bnu \\
&\hphantom{=}-\left[ \blam^\mathsf{T} \mrmd \bx \right]_a^b+ \int_a^b \left[ \left(H_{\bx} + {\dot \blam}^\mathsf{T} \right) \de \bx+ \left(H_{\blam} - {\dot \bx}^\mathsf{T} \right) \de \blam + H_{\bu} \de \bu   \right] \dt \\
&=  \left(G_{a} - \left. H \right|_{t=a} \right) \mrmd a + \left( G_{\bx(a)} + \left. \blam^\mathsf{T} \right|_{t=a} \right) \mrmd \bx(a) + G_{\bxi} \mrmd \bxi \\
&\hphantom{=}+ \left(G_{b} + \left. H \right|_{t=b} \right) \mrmd b + \left( G_{\bx(b)} - \left. \blam^\mathsf{T} \right|_{t=b} \right) \mrmd \bx(b) + G_{\bnu} \mrmd \bnu \\
&\hphantom{=}+ \int_a^b \left[ \left(H_{\bx} + {\dot \blam}^\mathsf{T} \right) \de \bx+ \left(H_{\blam} - {\dot \bx}^\mathsf{T} \right) \de \blam + H_{\bu} \de \bu   \right] \dt.
\end{split}
\end{equation}
In the first equality, Leibnitz's rule is used to compute the differential of the integral. In the second equality, integration by parts is used. In the third equality, the formula $\mrmd \bx(t) = \de \bx(t)+ \dot \bx(t) \dt$, or more concisely $\mrmd \bx = \de \bx+ \dot \bx \dt$, is used.

The necessary conditions on $\bx$, $\blam$, and $\bu$ which make $\mrmd \tilde{J}=0$ are the DAEs defined for $a \le t \le b$
\begin{equation} \label{eq_pmp_dae}
\begin{split}
\dot {\bx} &= H_{\blam}^\mathsf{T} \left(t,\bx,\blam,\bu,\mu\right) = \mathbf{f}\left(t,\bx,\bu,\mu\right) \\
\dot {\blam} &= - H_{\bx}^\mathsf{T} \left(t,\bx,\blam,\bu,\mu\right) \\
\mathbf{0}_{m \times 1} &= H_{\bu}^\mathsf{T} \left(t,\bx,\blam,\bu,\mu\right),
\end{split}
\end{equation}
the left boundary conditions defined at time $t=a$
\begin{equation} \label{eq_pmp_lbc}
\left. H \right|_{t=a} = G_{a}, \quad \left. \blam \right|_{t=a} = -G_{\bx(a)}^\mathsf{T}, \quad G_{\bxi}^\mathsf{T} = \bsigma\left(a,\bx(a),\mu\right) = \mathbf{0}_{k_1 \times 1},
\end{equation}
and the right boundary conditions defined at time $t=b$
\begin{equation} \label{eq_pmp_rbc}
\left. H \right|_{t=b} = -G_{b}, \quad \left. \blam \right|_{t=b} = G_{\bx(b)}^\mathsf{T}, \quad G_{\bnu}^\mathsf{T} = \bpsi\left(b,\bx(b),\mu\right) = \mathbf{0}_{k_2 \times 1}.
\end{equation}
If the initial time $a$ is prescribed, then the left boundary condition $\left. H \right|_{t=a} = G_{a}$ is dropped. If the final time $b$ is prescribed, then the right boundary condition $\left. H \right|_{t=b} = -G_{b}$ is dropped. The necessary conditions \eqref{eq_pmp_dae}, \eqref{eq_pmp_lbc}, and \eqref{eq_pmp_rbc} constitute a differential-algebraic equation two-point boundary value problem (DAE TPBVP).}

If $H_{\bu \bu}$ is nonsingular, then the optimal control problem is said to be regular or nonsingular; otherwise if $H_{\bu \bu}$ is singular, then the optimal control problem is said to be singular. If $H_{\bu \bu}$ is nonsingular, then by the implicit function theorem, the algebraic equation $H_{\bu}=\mathbf{0}_{1 \times m}$ in \eqref{eq_pmp_dae} guarantees the existence of a unique function, say $\bpi$, for which 
\begin{equation} \label{eq_pi_def}
\bu = \bpi\left(t,\bx,\blam,\mu\right).
\end{equation}
If $H_{\bu \bu}$ is nonsingular, it may be possible to solve the algebraic equation $H_{\bu}=\mathbf{0}_{1 \times m}$ in \eqref{eq_pmp_dae} analytically for $\bu$ in terms of $t$, $\bx$, $\blam$, and $\mu$ to construct $\bpi$ explicitly in \eqref{eq_pi_def}; otherwise, the value $\bu$ of $\bpi$ in \eqref{eq_pi_def} may be constructed numerically in an efficient and accurate manner (with quadratic convergence) via a few iterations of Newton's method applied to $H_{\bu}=\mathbf{0}_{1 \times m}$ starting from an initial guess $\bu_0$ of $\bu$: 
\begin{equation} \label{eq_pi_Newton}
\bu_{i+1} = \bu_i-H_{\bu \bu}^{-1}\left(t,\bx,\blam,\bu_i,\mu\right)  H_{\bu}^\mathsf{T} \left(t,\bx,\blam,\bu_i,\mu\right).
\end{equation} 
Using \eqref{eq_pi_def}, the Hamiltonian may be re-expressed as a function of $t$, $\bx$, $\blam$, and $\mu$ via the regular or reduced Hamiltonian
\begin{equation} \label{eq_Hhat_def}
\hat{H}\left(t,\bx,\blam,\mu\right) \equiv H\left(t,\bx,\blam,\bpi\left(t,\bx,\blam,\mu\right),\mu\right).
\end{equation} 
Note that by construction of $\bpi$, 
\begin{equation} \label{eq_H_u_pi}
H_{\bu} \left(t,\bx,\blam,\bpi\left(t,\bx,\blam,\mu\right),\mu\right) = \mathbf{0}_{1 \times m}.
\end{equation}
By using the definition \eqref{eq_Hhat_def} of the regular Hamiltonian $\hat{H}$, invoking \eqref{eq_H_u_pi}, and defining
\begin{equation} \label{eq_fhat_def}
\hat{\mathbf{f}}\left(t,\bx,\blam,\mu\right) \equiv \mathbf{f}\left(t,\bx,\bpi\left(t,\bx,\blam,\mu\right),\mu\right),
\end{equation} 
it follows from the chain rule that
\begin{equation} \label{eq_Hhat_t}
\begin{split}
\hat{H}_{t} \left(t,\bx,\blam,\mu\right)  &=  H_{t}\left(t,\bx,\blam,\bpi\left(t,\bx,\blam,\mu\right),\mu\right)+H_{\bu}\left(t,\bx,\blam,\bpi\left(t,\bx,\blam,\mu\right),\mu\right) \bpi_{t}\left(t,\bx,\blam,\mu\right) \\
&= H_{t}\left(t,\bx,\blam,\bpi\left(t,\bx,\blam,\mu\right),\mu\right), 
\end{split}
\end{equation}
\begin{equation} \label{eq_Hhat_x}
\begin{split}
\hat{H}_{\bx} \left(t,\bx,\blam,\mu\right)  &=  H_{\bx}\left(t,\bx,\blam,\bpi\left(t,\bx,\blam,\mu\right),\mu\right)+H_{\bu}\left(t,\bx,\blam,\bpi\left(t,\bx,\blam,\mu\right),\mu\right) \bpi_{\bx}\left(t,\bx,\blam,\mu\right) \\
&= H_{\bx}\left(t,\bx,\blam,\bpi\left(t,\bx,\blam,\mu\right),\mu\right), 
\end{split}
\end{equation}
\begin{equation} \label{eq_Hhat_lambda}
\begin{split}
\hat{H}_{\blam} \left(t,\bx,\blam,\mu\right)  &=  H_{\blam}\left(t,\bx,\blam,\bpi\left(t,\bx,\blam,\mu\right),\mu\right)+H_{\bu}\left(t,\bx,\blam,\bpi\left(t,\bx,\blam,\mu\right),\mu\right) \bpi_{\blam}\left(t,\bx,\blam,\mu\right) \\
&= H_{\blam}\left(t,\bx,\blam,\bpi\left(t,\bx,\blam,\mu\right),\mu\right) \\
&= \mathbf{f}^\mathsf{T}\left(t,\bx,\bpi\left(t,\bx,\blam,\mu\right),\mu\right) \\
&= \hat{\mathbf{f}}^\mathsf{T}\left(t,\bx,\blam,\mu\right),
\end{split}
\end{equation}
and
\begin{equation} \label{eq_Hhat_mu}
\begin{split}
\hat{H}_{\mu} \left(t,\bx,\blam,\mu\right)  &=  H_{\mu}\left(t,\bx,\blam,\bpi\left(t,\bx,\blam,\mu\right),\mu\right)+H_{\bu}\left(t,\bx,\blam,\bpi\left(t,\bx,\blam,\mu\right),\mu\right) \bpi_{\mu}\left(t,\bx,\blam,\mu\right) \\
&= H_{\mu}\left(t,\bx,\blam,\bpi\left(t,\bx,\blam,\mu\right),\mu\right). 
\end{split}
\end{equation}
By using \eqref{eq_pi_def} to eliminate the algebraic equation $H_{\bu}=\mathbf{0}_{1 \times m}$ from \eqref{eq_pmp_dae}, by plugging \eqref{eq_Hhat_lambda} and \eqref{eq_Hhat_x} into the right hand sides of the ODEs in \eqref{eq_pmp_dae}, and by plugging the definition \eqref{eq_Hhat_def} into the left and right boundary conditions \eqref{eq_pmp_lbc} and \eqref{eq_pmp_rbc}, the necessary  conditions on $\bx$ and $\blam$ which a solution of \eqref{dyn_opt_problem_pmp} must satisfy are the ODEs defined for $a \le t \le b$
\begin{equation} \label{eq_pmp_bvp}
\begin{split}
\dot {\bx} &= \hat{H}_{\blam}^\mathsf{T} \left(t,\bx,\blam,\mu\right) = \hat{\mathbf{f}}\left(t,\bx,\blam,\mu\right), \\
\dot {\blam} &= - \hat{H}_{\bx}^\mathsf{T} \left(t,\bx,\blam,\mu\right),
\end{split}
\end{equation}
the left boundary conditions defined at time $t=a$
\begin{equation} \label{eq_pmp_lbc2}
\left. \hat{H} \right|_{t=a} = G_{a}, \quad \left. \blam \right|_{t=a} = -G_{\bx(a)}^\mathsf{T}, \quad G_{\bxi}^\mathsf{T} = \bsigma\left(a,\bx(a),\mu\right) = \mathbf{0}_{k_1 \times 1},
\end{equation}
and the right boundary conditions defined at time $t=b$
\begin{equation} \label{eq_pmp_rbc2}
\left. \hat{H} \right|_{t=b} = -G_{b}, \quad \left. \blam \right|_{t=b} = G_{\bx(b)}^\mathsf{T}, \quad G_{\bnu}^\mathsf{T} = \bpsi\left(b,\bx(b),\mu\right) = \mathbf{0}_{k_2 \times 1}.
\end{equation}
If the initial time $a$ is prescribed, then the left boundary condition $\left. \hat{H} \right|_{t=a} = G_{a}$ is dropped. If the final time $b$ is prescribed, then the right boundary condition $\left. \hat{H} \right|_{t=b} = -G_{b}$ is dropped. The necessary conditions \eqref{eq_pmp_bvp}, \eqref{eq_pmp_lbc2}, and \eqref{eq_pmp_rbc2} constitute an ODE TPBVP. Appendix~\ref{app_imp_details} provides implementation details for numerically solving the ODE TPBVP \eqref{eq_pmp_bvp}, \eqref{eq_pmp_lbc2}, and \eqref{eq_pmp_rbc2}.

A solution of the DAE TPBVP \eqref{eq_pmp_dae}, \eqref{eq_pmp_lbc}, and \eqref{eq_pmp_rbc} or of the ODE TPBVP \eqref{eq_pmp_bvp}, \eqref{eq_pmp_lbc2}, and \eqref{eq_pmp_rbc2} is said to be an extremal solution of the optimal control problem \eqref{dyn_opt_problem_pmp}. Note that an extremal solution only satisfies necessary conditions for a minimum of the optimal control problem \eqref{dyn_opt_problem_pmp}, so that an extremal solution is not guaranteed to be a local minimum of \eqref{dyn_opt_problem_pmp}. 

Since the DAE TPBVP \eqref{eq_pmp_dae}, \eqref{eq_pmp_lbc}, and \eqref{eq_pmp_rbc} and the ODE TPBVP \eqref{eq_pmp_bvp}, \eqref{eq_pmp_lbc2}, and \eqref{eq_pmp_rbc2} have small convergence radii \cite{BrHo1975applied,bryson1999dynamic,betts2010practical}, a continuation method (performing continuation in the parameter $\mu$) is often required to numerically solve them starting from a solution to a simpler optimal control problem \cite{ascher1994numerical}. The solution to the simpler optimal control problem might be obtained via analytics, the gradient method \cite{BrHo1975applied,bryson1999dynamic}, the method of successive approximations \cite{krylov1963method,krylov1972algorithm,lyubushin1982modifications,chernousko1982method}, or the direct method \cite{betts2010practical,patterson2014gpops}. For example in \cite{putkaradze2018constraint}, the continuation parameter $\mu$ is used to vary integrand cost function coefficients in $L$ in order to numerically solve the optimal control ODE TPBVPs for Suslov's problem 
via monotonic continuation, starting from an analytical solution to a singular optimal control problem. In Sections~\ref{sec_disk_sim} and \ref{sec_ball_sim}, the continuation parameter $\mu$ is used to vary integrand cost function coefficients in $L$ in order to numerically solve the optimal control ODE TPBVPs for the rolling disk and ball via predictor-corrector continuation, starting from a direct method solution to a simpler optimal control problem. Appendices~\ref{app_predictor_corrector} and \ref{app_sweep_predictor_corrector} describe predictor-corrector continuation methods for solving ODE TPBVPs and which are used to solve the optimal control ODE TPBVPs for the rolling disk and ball in Sections~\ref{sec_disk_sim} and \ref{sec_ball_sim}. 

\rem{Because path inequality constraints have been omitted from the optimal control problem \eqref{dyn_opt_problem_pmp}, the control is not restricted to lie in a compact set. Hence, an extremal solution solving \eqref{eq_pmp_dae}, \eqref{eq_pmp_lbc}, and \eqref{eq_pmp_rbc} or \eqref{eq_pmp_bvp}, \eqref{eq_pmp_lbc2}, and \eqref{eq_pmp_rbc2} does not lie on any boundary. If the control is restricted to lie in a compact set (due to path inequality constraints), then the control of an extremal solution may lie on the boundary of this compact set; if the control is discontinuous, hopping abruptly between points on the boundary, then it is said to be bang-bang.}

\section{Implementation Details for Solving the ODE TPBVP for a Regular Optimal Control Problem} \label{app_imp_details}

Details for numerically solving the ODE TPBVP \eqref{eq_pmp_bvp}, \eqref{eq_pmp_lbc2}, and \eqref{eq_pmp_rbc2} associated with the indirect method solution of a regular optimal control problem are presented here. There are two general methods, initial value and global, for numerically solving an ODE TPBVP. An initial value method, such as single or multiple shooting, subdivides the integration interval $[a,b]$ into a fixed, finite mesh and integrates the ODE on each mesh subinterval using a guess of the unknown initial conditions at one endpoint in each mesh subinterval. A root-finder is used to iteratively adjust the guesses of the unknown initial conditions until the solution segments are continuous at the internal mesh points and until the boundary conditions at the endpoints $a$ and $b$ are satisfied. A global method, such as a Runge-Kutta, collocation, or finite-difference scheme, subdivides the integration interval $[a,b]$ into a finite, adaptive mesh and solves a large nonlinear system of algebraic equations obtained by imposing the ODE constraints at a finite set of points in each mesh subinterval, by imposing continuity of the solution at internal mesh points, and by imposing the boundary conditions at the endpoints $a$ and $b$. By estimating the error in each mesh subinterval, a global method iteratively refines or adapts the mesh until a prescribed error tolerance is satisfied. Because initial value methods cannot integrate unstable ODEs, global methods are preferred \cite{ascher1994numerical,muir1999optimal,boisvert2011problem}.  

\subsection{Normalization and ODE Velocity Function}
There are many solvers available to numerically solve the ODE TPBVP \eqref{eq_pmp_bvp}, \eqref{eq_pmp_lbc2}, and \eqref{eq_pmp_rbc2}. For example, \mcode{bvp4c} \cite{shampine2000solving}, \mcode{bvp5c} \cite{kierzenka2008bvp}, \mcode{sbvp} \cite{auzinger2003collocation}, and \mcode{bvptwp} \cite{cash2013algorithm} (which encapsulates \mcode{twpbvp_m}, \mcode{twpbvpc_m}, \mcode{twpbvp_l}, \mcode{twpbvpc_l}, \mcode{acdc}, and \mcode{acdcc}) are \mcode{MATLAB} Runge-Kutta or collocation ODE TPBVP solvers, while COLSYS \cite{ascher1981algorithm}, COLNEW \cite{bader1987new}, COLMOD \cite{cash2001automatic}, COLCON \cite{bader1989continuation}, BVP\_M-2 \cite{boisvert2013runge}, TWPBVP \cite{cash1991deferred}, TWPBVPC \cite{cash2005new}, TWPBVPL \cite{bashir1998lobatto}, TWPBVPLC \cite{cash2006hybrid}, ACDC \cite{cash2001automatic}, and ACDCC \cite{cash2013algorithm} are Fortran Runge-Kutta or collocation ODE TPBVP solvers. The reader is referred to the Appendix in \cite{putkaradze2018constraint} for a comprehensive list of ODE TPBVP solvers. In order to numerically solve the ODE TPBVP \eqref{eq_pmp_bvp}, \eqref{eq_pmp_lbc2}, and \eqref{eq_pmp_rbc2}, many solvers (such as the global method \mcode{MATLAB} and Fortran solvers just listed) require that the ODE TPBVP be defined on a fixed time interval and any unknown parameters, such as $\bxi$, $\bnu$, $a$, and $b$, must often be modeled as dummy constant 
dependent variables with zero derivatives. In addition, to aid convergence, many solvers can exploit Jacobians of the ODE velocity function and of the two-point boundary condition function. Thus, \eqref{eq_pmp_bvp} is redefined on the normalized time interval $[0,1]$ through the change of independent variable $s \equiv \frac{t-a}{T}$, where $T \equiv  b-a$. Note that $t(s) = Ts+a$. Define the normalized state $\tilde \bx(s) \equiv \bx(t(s))$ and normalized costate $\tilde \blam(s) \equiv \blam(t(s))$. Define the expanded un-normalized ODE TPBVP dependent variable vector
\begin{equation} \label{eq_z_def}
\bz(t) \equiv \begin{bmatrix} \bx(t) \\ \blam(t) \\ \bxi \\ \bnu \\ a \\ b \end{bmatrix}.
\end{equation}
Defining $\tilde \bz(s) \equiv \bz(t(s))$, the expanded normalized ODE TPBVP dependent variable vector is
\begin{equation} \label{eq_tz_def}
\tilde \bz(s) \equiv \bz(t(s)) = \begin{bmatrix} \bx(t(s)) \\ \blam(t(s)) \\ \bxi \\ \bnu \\ a \\ b \end{bmatrix} = \begin{bmatrix} \tilde \bx(s) \\ \tilde \blam(s) \\ \bxi \\ \bnu \\ a \\ b \end{bmatrix}.
\end{equation}
By the chain rule, \eqref{eq_pmp_bvp}, and since $\dd{t(s)}{s} = T$,
\begin{equation} \label{eq_dtz}
\begin{split}
\dot {\tilde \bz}(s) = \dd{\tilde \bz(s)}{s} = \begin{bmatrix} \dot {\tilde \bx}(s) \\ \dot {\tilde \blam}(s) \\ \mathbf{0}_{(k_1+k_2+2) \times 1} \end{bmatrix}  = \dd{\bz(t(s))}{t} \dd{t(s)}{s} &= \begin{bmatrix} \dd {\bx(t(s))}{t} \\ \dd {\blam(t(s))}{t} \\ \mathbf{0}_{(k_1+k_2+2) \times 1} \end{bmatrix} \dd{t(s)}{s} \\ &= \begin{bmatrix} \hat{\mathbf{f}} \left(t(s), \bx(t(s)), \blam(t(s)),\mu\right)  \\ - \hat{H}_{\bx}^\mathsf{T} \left(t(s),\bx(t(s)),\blam(t(s)),\mu\right) \\ \mathbf{0}_{(k_1+k_2+2) \times 1}  \end{bmatrix} T \\ &= \begin{bmatrix} \hat{\mathbf{f}} \left(t(s),\tilde \bx(s),\tilde \blam(s),\mu\right)  \\ - \hat{H}_{\bx}^\mathsf{T} \left(t(s),\tilde \bx(s),\tilde \blam(s),\mu\right) \\ \mathbf{0}_{(k_1+k_2+2) \times 1}  \end{bmatrix} T.
\end{split}
\end{equation}
Define $\tilde \bPhi \left(s,\tilde \bz(s),\mu \right)$ to be the right-hand side of \eqref{eq_dtz}, i.e. the normalized ODE velocity function, so that
\begin{equation} \label{eq_Phi}
\tilde \bPhi \left(s,\tilde \bz(s),\mu \right) \equiv \begin{bmatrix} \hat{\mathbf{f}} \left(t(s),\tilde \bx(s),\tilde \blam(s),\mu\right)  \\ - \hat{H}_{\bx}^\mathsf{T} \left(t(s),\tilde \bx(s),\tilde \blam(s),\mu\right) \\ \mathbf{0}_{(k_1+k_2+2) \times 1}  \end{bmatrix} T.
\end{equation}
The Jacobian of $\tilde \bPhi$ with respect to $\tilde \bz(s)$ is
\begin{equation} \label{eq_DPhi_z}
\begin{split}
\tilde \bPhi_{\tilde \bz(s)}&\left(s,\tilde \bz(s),\mu\right) =  \\
&\begin{bmatrix} \hat{\mathbf{f}}_{\bx}T & \hat{\mathbf{f}}_{\blam}T & \mathbf{0}_{n \times (k_1+k_2)} & -\hat{\mathbf{f}}+\hat{\mathbf{f}}_{t}(1-s)T &  \hat{\mathbf{f}}+\hat{\mathbf{f}}_{t}sT   \\ - \hat{H}_{\bx \bx}T & -\hat{\mathbf{f}}_{\bx}^\mathsf{T}T & \mathbf{0}_{n \times (k_1+k_2)} & \hat{H}_{\bx}^\mathsf{T}-\hat{H}_{\bx t}(1-s)T &  -\hat{H}_{\bx}^\mathsf{T}-\hat{H}_{\bx t}sT \\ \mathbf{0}_{(k_1+k_2+2) \times n} & \mathbf{0}_{(k_1+k_2+2) \times n} & \mathbf{0}_{(k_1+k_2+2) \times (k_1+k_2)} & \mathbf{0}_{(k_1+k_2+2) \times 1} & \mathbf{0}_{(k_1+k_2+2) \times 1}  \end{bmatrix}
\end{split}
\end{equation} 
and the Jacobian of $\tilde \bPhi$ with respect to $\mu$ is
\begin{equation} \label{eq_DPhi_mu}
\tilde \bPhi_{\mu}\left(s,\tilde \bz(s),\mu\right) = \begin{bmatrix} \hat{\mathbf{f}}_{\mu}T  \\ - \hat{H}_{\bx \mu}T \\ \mathbf{0}_{(k_1+k_2+2) \times 1}   \end{bmatrix}.
\end{equation}
In \eqref{eq_DPhi_z} and \eqref{eq_DPhi_mu}, shorthand notation is used for conciseness and all zeroth and first derivatives of $\hat{\mathbf{f}}$ and all first and second derivatives of $\hat{H}$ are evaluated at $\left( s,\tilde \bz(s),\mu \right)$. An explanation of the meaning of the shorthand notation used to express all zeroth and first derivatives of $\hat{\mathbf{f}}$ and all first and second derivatives of $\hat{H}$ is given in Table~\ref{table_ODE_BVP_ODE_Jacobians}. In rows $n+1$ through $2n$ and columns $n+1$ through $2n$ of \eqref{eq_DPhi_z}, Clairaut's Theorem was used to obtain $\hat{H}_{\bx \blam} = \hat{H}_{\blam \bx}^\mathsf{T} = \hat{\mathbf{f}}_{\bx}^\mathsf{T}$, recalling from \eqref{eq_Hhat_lambda} that $\hat{H}_{\blam} = \hat{\mathbf{f}}^\mathsf{T}$.

\begin{table}[h!]
	\centering 
	{ 
		\setlength{\extrarowheight}{1.5pt}
		\begin{tabular}{| c c c c c c c |} 
			\hline
			\textbf{Shorthand} & $\mathbf{\vert}$ & \textbf{Extended Shorthand} & $\mathbf{\vert}$ &
			\textbf{Normalized} & $\mathbf{\vert}$ & \textbf{Un-Normalized} \\ 
			\hline\hline 
			
			$\hat{\mathbf{f}}$ &=& $\left. \hat{\mathbf{f}} \right|_{\left( s,\tilde \bz(s),\mu \right)}$ &=& $\hat{\mathbf{f}} \left(t(s),\tilde \bx(s),\tilde \blam(s),\mu\right) $ &=& $ \hat{\mathbf{f}} \left(t(s),\bx(t(s)),\blam(t(s)),\mu\right)$ \\ \hline
			
			$\hat{\mathbf{f}}_{\blam}$ &=& $\left. \hat{\mathbf{f}}_{\blam} \right|_{\left( s,\tilde \bz(s),\mu \right)} $ &=& $\hat{\mathbf{f}}_{\blam} \left(t(s),\tilde \bx(s),\tilde \blam(s),\mu\right) $ &=& $ \hat{\mathbf{f}}_{\blam} \left(t(s),\bx(t(s)),\blam(t(s)),\mu\right)$ \\ \hline
			
			$\hat{\mathbf{f}}_{\bx}$ &=& $\left. \hat{\mathbf{f}}_{\bx} \right|_{\left( s,\tilde \bz(s),\mu \right)} $ &=& $\hat{\mathbf{f}}_{\bx} \left(t(s),\tilde \bx(s),\tilde \blam(s),\mu\right) $ &=& $ \hat{\mathbf{f}}_{\bx} \left(t(s),\bx(t(s)),\blam(t(s)),\mu\right)$ \\ \hline
			
			$\hat{\mathbf{f}}_{t}$ &=& $\left. \hat{\mathbf{f}}_{t} \right|_{\left( s,\tilde \bz(s),\mu \right)} $ &=& $\hat{\mathbf{f}}_{t} \left(t(s),\tilde \bx(s),\tilde \blam(s),\mu\right) $ &=& $ \hat{\mathbf{f}}_{t} \left(t(s),\bx(t(s)),\blam(t(s)),\mu\right)$ \\ \hline
			
			$\hat{\mathbf{f}}_{\mu}$ &=& $\left. \hat{\mathbf{f}}_{\mu} \right|_{\left( s,\tilde \bz(s),\mu \right)} $ &=& $\hat{\mathbf{f}}_{\mu} \left(t(s),\tilde \bx(s),\tilde \blam(s),\mu\right) $ &=& $ \hat{\mathbf{f}}_{\mu} \left(t(s),\bx(t(s)),\blam(t(s)),\mu\right)$ \\ \hline
			
			$\hat{H}_{\bx}^\mathsf{T}$ &=& $\left. \hat{H}_{\bx}^\mathsf{T} \right|_{\left( s,\tilde \bz(s),\mu \right)} $&=&$\hat{H}_{\bx}^\mathsf{T} \left(t(s),\tilde \bx(s),\tilde \blam(s),\mu\right)$&=&$ \hat{H}_{\bx}^\mathsf{T} \left(t(s),\bx(t(s)),\blam(t(s)),\mu\right)$ \\ \hline
			
			$\hat{H}_{\bx \bx}$ &=& $\left. \hat{H}_{\bx \bx} \right|_{\left( s,\tilde \bz(s) ,\mu\right)} $ &=& $\hat{H}_{\bx \bx} \left(t(s),\tilde \bx(s),\tilde \blam(s),\mu\right) $ &=& $ \hat{H}_{\bx \bx} \left(t(s),\bx(t(s)),\blam(t(s)),\mu\right)$ \\ \hline
			
			$\hat{H}_{\bx t}$ &=& $\left. \hat{H}_{\bx t} \right|_{\left( s,\tilde \bz(s),\mu \right)} $ &=& $\hat{H}_{\bx t} \left(t(s),\tilde \bx(s),\tilde \blam(s),\mu\right) $ &=& $ \hat{H}_{\bx t} \left(t(s),\bx(t(s)),\blam(t(s)),\mu\right)$ \\ \hline
			
			$\hat{H}_{\bx \mu}$ &=& $\left. \hat{H}_{\bx \mu} \right|_{\left( s,\tilde \bz(s),\mu \right)} $ &=& $\hat{H}_{\bx \mu} \left(t(s),\tilde \bx(s),\tilde \blam(s),\mu\right) $ &=& $ \hat{H}_{\bx \mu} \left(t(s),\bx(t(s)),\blam(t(s)),\mu\right)$ \\ \hline
			
		\end{tabular} 
	}
	\caption{Explanation of shorthand notation for zeroth and first derivatives of $\hat{\mathbf{f}}$ and first and second derivatives of $\hat{H}$ used in \eqref{eq_DPhi_z} and \eqref{eq_DPhi_mu}.}
	\label{table_ODE_BVP_ODE_Jacobians}
\end{table}

Recall that $\hat{\mathbf{f}}$ is defined in \eqref{eq_fhat_def} in terms of $\mathbf{f}$ and $\bpi$. By using the chain rule, the first derivatives of $\hat{\mathbf{f}}$ that appear in \eqref{eq_DPhi_z}, \eqref{eq_DPhi_mu}, and Table~\ref{table_ODE_BVP_ODE_Jacobians} may be computed from first derivatives of $\mathbf{f}$ and $\bpi$ as follows:
\begin{equation} \label{eq_fhat_lam}
\hat{\mathbf{f}}_{\blam} \left(t,\bx,\blam,\mu\right)=\mathbf{f}_{\bu} \left(t,\bx,\bpi\left(t,\bx,\blam,\mu\right),\mu\right) \bpi_{\blam} \left(t,\bx,\blam,\mu\right),
\end{equation}
\begin{equation} \label{eq_fhat_x}
\begin{split}
\hat{\mathbf{f}}_{\bx} \left(t,\bx,\blam,\mu\right)=\mathbf{f}_{\bx} \left(t,\bx,\bpi\left(t,\bx,\blam,\mu\right),\mu\right)+\mathbf{f}_{\bu} \left(t,\bx,\bpi\left(t,\bx,\blam,\mu\right),\mu\right) \bpi_{\bx} \left(t,\bx,\blam,\mu\right),
\end{split}
\end{equation}
\begin{equation} \label{eq_fhat_t}
\begin{split}
\hat{\mathbf{f}}_{t} \left(t,\bx,\blam,\mu\right)=\mathbf{f}_{t} \left(t,\bx,\bpi\left(t,\bx,\blam,\mu\right),\mu\right)+\mathbf{f}_{\bu} \left(t,\bx,\bpi\left(t,\bx,\blam,\mu\right),\mu\right) \bpi_{t} \left(t,\bx,\blam,\mu\right),
\end{split}
\end{equation}
and
\begin{equation} \label{eq_fhat_mu}
\begin{split}
\hat{\mathbf{f}}_{\mu} \left(t,\bx,\blam,\mu\right)=\mathbf{f}_{\mu} \left(t,\bx,\bpi\left(t,\bx,\blam,\mu\right),\mu\right)+\mathbf{f}_{\bu} \left(t,\bx,\bpi\left(t,\bx,\blam,\mu\right),\mu\right) \bpi_{\mu} \left(t,\bx,\blam,\mu\right).
\end{split}
\end{equation}

Recall that by construction of $\bpi$, $H_{\bu} \left(t,\bx,\blam,\bpi\left(t,\bx,\blam,\mu\right),\mu\right) = \mathbf{0}_{1 \times m}$, as stated previously in \eqref{eq_H_u_pi}. Differentiating $H_{\bu} \left(t,\bx,\blam,\bpi\left(t,\bx,\blam,\mu\right),\mu\right) = \mathbf{0}_{1 \times m}$ with respect to $\blam$, $\bx$, $t$, and $\mu$, in turn, and using the chain rule gives
\begin{equation} \label{eq_H_ulam_pi}
H_{\bu \blam} \left(t,\bx,\blam,\bpi\left(t,\bx,\blam,\mu\right),\mu\right)+H_{\bu \bu} \left(t,\bx,\blam,\bpi\left(t,\bx,\blam,\mu\right),\mu\right) \bpi_{\blam} \left(t,\bx,\blam,\mu\right)= \mathbf{0}_{m \times n},
\end{equation}
\begin{equation} \label{eq_H_ux_pi}
H_{\bu \bx} \left(t,\bx,\blam,\bpi\left(t,\bx,\blam,\mu\right),\mu\right)+H_{\bu \bu} \left(t,\bx,\blam,\bpi\left(t,\bx,\blam,\mu\right),\mu\right) \bpi_{\bx} \left(t,\bx,\blam,\mu\right)= \mathbf{0}_{m \times n},
\end{equation}
\begin{equation} \label{eq_H_ut_pi}
H_{\bu t} \left(t,\bx,\blam,\bpi\left(t,\bx,\blam,\mu\right),\mu\right)+H_{\bu \bu} \left(t,\bx,\blam,\bpi\left(t,\bx,\blam,\mu\right),\mu\right) \bpi_{t} \left(t,\bx,\blam,\mu\right)= \mathbf{0}_{m \times 1},
\end{equation}
and
\begin{equation} \label{eq_H_umu_pi}
H_{\bu \mu} \left(t,\bx,\blam,\bpi\left(t,\bx,\blam,\mu\right),\mu\right)+H_{\bu \bu} \left(t,\bx,\blam,\bpi\left(t,\bx,\blam,\mu\right),\mu\right) \bpi_{\mu} \left(t,\bx,\blam,\mu\right)= \mathbf{0}_{m \times 1}.
\end{equation}
\eqref{eq_H_ulam_pi}, \eqref{eq_H_ux_pi}, \eqref{eq_H_ut_pi}, and \eqref{eq_H_umu_pi} may be solved for $\bpi_{\blam}$, $\bpi_{\bx}$, $\bpi_{t}$, and $\bpi_{\mu}$, respectively:
\begin{equation} \label{eq_pi_lam}
\begin{split}
\bpi_{\blam} \left(t,\bx,\blam,\mu\right) &= -H_{\bu \bu}^{-1} \left(t,\bx,\blam,\bpi\left(t,\bx,\blam,\mu\right),\mu\right)H_{\bu \blam} \left(t,\bx,\blam,\bpi\left(t,\bx,\blam,\mu\right),\mu\right) \\
&= -H_{\bu \bu}^{-1} \left(t,\bx,\blam,\bpi\left(t,\bx,\blam,\mu\right),\mu\right)\mathbf{f}_{\bu}^\mathsf{T} \left(t,\bx,\bpi\left(t,\bx,\blam,\mu\right),\mu\right),
\end{split}
\end{equation}
\begin{equation} \label{eq_pi_x}
\bpi_{\bx} \left(t,\bx,\blam,\mu\right)= -H_{\bu \bu}^{-1} \left(t,\bx,\blam,\bpi\left(t,\bx,\blam,\mu\right),\mu\right) H_{\bu \bx} \left(t,\bx,\blam,\bpi\left(t,\bx,\blam,\mu\right),\mu\right),
\end{equation}
\begin{equation} \label{eq_pi_t}
\bpi_{t} \left(t,\bx,\blam,\mu\right)= -H_{\bu \bu}^{-1} \left(t,\bx,\blam,\bpi\left(t,\bx,\blam,\mu\right),\mu\right) H_{\bu t} \left(t,\bx,\blam,\bpi\left(t,\bx,\blam,\mu\right),\mu\right),
\end{equation}
and
\begin{equation} \label{eq_pi_mu}
\bpi_{\mu} \left(t,\bx,\blam,\mu\right)= -H_{\bu \bu}^{-1} \left(t,\bx,\blam,\bpi\left(t,\bx,\blam,\mu\right),\mu\right) H_{\bu \mu} \left(t,\bx,\blam,\bpi\left(t,\bx,\blam,\mu\right),\mu\right).
\end{equation}
In \eqref{eq_pi_lam},  Clairaut's Theorem was used to obtain $H_{\bu \blam} = H_{\blam \bu}^\mathsf{T} = {\mathbf{f}}_{\bu}^\mathsf{T}$, since $H_{\blam} = \mathbf{f}^\mathsf{T}$. As will be stated again later, \eqref{eq_pi_lam}, \eqref{eq_pi_x}, \eqref{eq_pi_t}, and \eqref{eq_pi_mu} are especially useful if the value of $\bpi$ is constructed numerically via Newton's method as in \eqref{eq_pi_Newton}; if $\bpi$ is given analytically, then it should be possible to construct $\bpi_{\blam}$, $\bpi_{\bx}$, $\bpi_{t}$, and $\bpi_{\mu}$ via manual, symbolic, or automatic differentiation of the analytical formula for $\bpi$.

Since Clairaut's Theorem guarantees that 
\begin{equation}
H_{\bx \bu} \left(t,\bx,\blam,\bu,\mu\right) =
H_{\bu \bx}^\mathsf{T} \left(t,\bx,\blam,\bu,\mu\right)
\end{equation}
and
\begin{equation}
H_{\bu \bu} \left(t,\bx,\blam,\bu,\mu\right) = 
H_{\bu \bu}^\mathsf{T} \left(t,\bx,\blam,\bu,\mu\right),
\end{equation}
\eqref{eq_H_ux_pi} may be solved for $H_{\bx \bu} \left(t,\bx,\blam,\bpi\left(t,\bx,\blam,\mu\right),\mu\right)$:
\begin{equation} \label{eq_H_xu_pi}
\begin{split}
H_{\bx \bu} \left(t,\bx,\blam,\bpi\left(t,\bx,\blam,\mu\right),\mu\right) &=
H_{\bu \bx}^\mathsf{T} \left(t,\bx,\blam,\bpi\left(t,\bx,\blam,\mu\right),\mu\right) \\
&=-\bpi_{\bx}^\mathsf{T} \left(t,\bx,\blam,\mu\right) H_{\bu \bu}^\mathsf{T} \left(t,\bx,\blam,\bpi\left(t,\bx,\blam,\mu\right),\mu\right) \\ &=-\bpi_{\bx}^\mathsf{T} \left(t,\bx,\blam,\mu\right) H_{\bu \bu} \left(t,\bx,\blam,\bpi\left(t,\bx,\blam,\mu\right),\mu\right). 
\end{split}
\end{equation}
By differentiating \eqref{eq_Hhat_x} with respect to $\bx$, $t$, and $\mu$, using the chain rule, and exploiting \eqref{eq_H_xu_pi}, the second derivatives of $\hat{H}$ that appear in \eqref{eq_DPhi_z}, \eqref{eq_DPhi_mu}, and Table~\ref{table_ODE_BVP_ODE_Jacobians} may be computed from first derivatives of $\mathbf{\bpi}$ and second derivatives of $H$ as follows: 
\begin{equation} \label{eq_Hhat_xx}
\begin{split}
\hat{H}_{\bx \bx} \left(t,\bx,\blam,\mu\right)&=H_{\bx \bx} \left(t,\bx,\blam,\bpi\left(t,\bx,\blam,\mu\right),\mu\right)+H_{\bx \bu} \left(t,\bx,\blam,\bpi\left(t,\bx,\blam,\mu\right),\mu\right) \bpi_{\bx} \left(t,\bx,\blam,\mu\right) \\
& = H_{\bx \bx} \left(t,\bx,\blam,\bpi\left(t,\bx,\blam,\mu\right),\mu\right)\\&\hphantom{=}-\bpi_{\bx}^\mathsf{T} \left(t,\bx,\blam,\mu\right) H_{\bu \bu} \left(t,\bx,\blam,\bpi\left(t,\bx,\blam,\mu\right),\mu\right) \bpi_{\bx} \left(t,\bx,\blam,\mu\right),
\end{split}
\end{equation}
\begin{equation} \label{eq_Hhat_xt}
\begin{split}
\hat{H}_{\bx t} \left(t,\bx,\blam,\mu\right)&=H_{\bx t} \left(t,\bx,\blam,\bpi\left(t,\bx,\blam,\mu\right),\mu\right)+H_{\bx \bu} \left(t,\bx,\blam,\bpi\left(t,\bx,\blam,\mu\right),\mu\right) \bpi_{t} \left(t,\bx,\blam,\mu\right) \\
& = H_{\bx t} \left(t,\bx,\blam,\bpi\left(t,\bx,\blam,\mu\right),\mu\right)\\&\hphantom{=}-\bpi_{\bx}^\mathsf{T} \left(t,\bx,\blam,\mu\right) H_{\bu \bu} \left(t,\bx,\blam,\bpi\left(t,\bx,\blam,\mu\right),\mu\right) \bpi_{t} \left(t,\bx,\blam,\mu\right),
\end{split}
\end{equation}
and
\begin{equation} \label{eq_Hhat_xmu}
\begin{split}
\hat{H}_{\bx \mu} \left(t,\bx,\blam,\mu\right)&=H_{\bx \mu} \left(t,\bx,\blam,\bpi\left(t,\bx,\blam,\mu\right),\mu\right)+H_{\bx \bu} \left(t,\bx,\blam,\bpi\left(t,\bx,\blam,\mu\right),\mu\right) \bpi_{\mu} \left(t,\bx,\blam,\mu\right) \\
& = H_{\bx \mu} \left(t,\bx,\blam,\bpi\left(t,\bx,\blam,\mu\right),\mu\right)\\&\hphantom{=}-\bpi_{\bx}^\mathsf{T} \left(t,\bx,\blam,\mu\right) H_{\bu \bu} \left(t,\bx,\blam,\bpi\left(t,\bx,\blam,\mu\right),\mu\right) \bpi_{\mu} \left(t,\bx,\blam,\mu\right).
\end{split}
\end{equation}

If the value of $\bpi$ is constructed numerically via Newton's method as in \eqref{eq_pi_Newton} rather than analytically, then \eqref{eq_pi_lam}, \eqref{eq_pi_x}, \eqref{eq_pi_t}, and \eqref{eq_pi_mu} should be used to evaluate $\bpi_{\blam}$, $\bpi_{\bx}$, $\bpi_{t}$, and $\bpi_{\mu}$, which appear in the formulas for $\hat{\mathbf{f}}_{\blam}$ given in \eqref{eq_fhat_lam}, $\hat{\mathbf{f}}_{\bx}$ given in \eqref{eq_fhat_x}, $\hat{\mathbf{f}}_{t}$ given in \eqref{eq_fhat_t}, $\hat{\mathbf{f}}_{\mu}$ given in \eqref{eq_fhat_mu}, $\hat{H}_{\bx \bx}$ given in \eqref{eq_Hhat_xx}, $\hat{H}_{\bx t}$ given in \eqref{eq_Hhat_xt}, and $\hat{H}_{\bx \mu}$ given in \eqref{eq_Hhat_xmu}. The second equation in \eqref{eq_Hhat_xx}, \eqref{eq_Hhat_xt}, and \eqref{eq_Hhat_xmu} is given because it may be more computationally efficient than the first equation if $\bpi$ is given analytically, so that \eqref{eq_pi_lam}, \eqref{eq_pi_x}, \eqref{eq_pi_t}, and \eqref{eq_pi_mu} need not be used to evaluate $\bpi_{\blam}$, $\bpi_{\bx}$, $\bpi_{t}$, and $\bpi_{\mu}$. 

\subsection{Two-Point Boundary Condition Function}
Now the boundary conditions \eqref{eq_pmp_lbc2}-\eqref{eq_pmp_rbc2} are considered. Letting
\begin{equation} \label{eq_Ups_1}
\bUpsilon_1\left( \bz(a),\bz(b),\mu\right) \equiv \begin{bmatrix} \hat{H} \left(a,\bx(a), \blam(a),\mu\right) \\  \blam(a) \\ \bsigma\left(a,\bx(a),\mu\right) \\ \hat{H} \left(b,\bx(b), \blam(b),\mu\right) \\ \blam(b) \\ \bpsi\left(b,\bx(b),\mu\right) \end{bmatrix},
\end{equation}
\begin{equation}  \label{eq_Ups_2}
\bUpsilon_2\left(\bz(a),\bz(b),\mu\right) \equiv \begin{bmatrix} G_a\left(a,\bx(a), \bxi, b, \bx(b),\bnu,\mu\right) \\ -G_{\bx(a)}^\mathsf{T}\left(a,\bx(a), \bxi, b, \bx(b),\bnu,\mu\right) \\ \mathbf{0}_{k_1\times 1} \\ -G_b\left(a,\bx(a), \bxi, b, \bx(b),\bnu,\mu\right) \\ G_{\bx(b)}^\mathsf{T}\left(a,\bx(a), \bxi, b, \bx(b),\bnu,\mu\right) \\ \mathbf{0}_{k_2 \times 1} \end{bmatrix},
\end{equation}
and
\begin{equation} \label{eq_Ups}
\begin{split}
\bUpsilon\left(\bz(a),\bz(b),\mu\right) &\equiv \bUpsilon_1\left( \bz(a),\bz(b),\mu\right) -  \bUpsilon_2\left(\bz(a),\bz(b),\mu\right) \\ &=\begin{bmatrix} \hat{H} \left(a,\bx(a), \blam(a),\mu\right) \\  \blam(a) \\ \bsigma\left(a,\bx(a),\mu\right) \\ \hat{H} \left(b,\bx(b), \blam(b),\mu\right) \\ \blam(b) \\ \bpsi\left(b,\bx(b),\mu\right) \end{bmatrix} - \begin{bmatrix} G_a\left(a,\bx(a), \bxi, b, \bx(b),\bnu,\mu\right) \\ -G_{\bx(a)}^\mathsf{T}\left(a,\bx(a), \bxi, b, \bx(b),\bnu,\mu\right) \\ \mathbf{0}_{k_1 \times 1} \\ -G_b\left(a,\bx(a), \bxi, b, \bx(b),\bnu,\mu\right) \\ G_{\bx(b)}^\mathsf{T}\left(a,\bx(a), \bxi, b, \bx(b),\bnu,\mu\right) \\ \mathbf{0}_{k_2 \times 1} \end{bmatrix},
\end{split}
\end{equation}
the boundary conditions \eqref{eq_pmp_lbc2}-\eqref{eq_pmp_rbc2} in un-normalized dependent variables are given by the two-point boundary condition function
\begin{equation} \label{eq_Ups_BC}
\bUpsilon\left(\bz(a),\bz(b),\mu\right) = \mathbf{0}_{(2n+k_1+k_2+2) \times 1}.
\end{equation}
The Jacobians of $\bUpsilon$ with respect to $\bz(a)$, $\bz(b)$, and $\mu$ are
\begin{equation} \label{eq_Jac_Ups1_za}
\bUpsilon_{\bz(a)}\left(\bz(a),\bz(b),\mu\right) = \bUpsilon_{1,\bz(a)}\left( \bz(a),\bz(b),\mu\right) -  \bUpsilon_{2,\bz(a)}\left(\bz(a),\bz(b),\mu\right),
\end{equation}
\begin{equation} \label{eq_Jac_Ups1_zb}
\bUpsilon_{\bz(b)}\left(\bz(a),\bz(b),\mu\right) = \bUpsilon_{1,\bz(b)}\left( \bz(a),\bz(b),\mu\right) -  \bUpsilon_{2,\bz(b)}\left(\bz(a),\bz(b),\mu\right),
\end{equation}
and
\begin{equation} \label{eq_Jac_Ups1_mu}
\bUpsilon_{\mu}\left(\bz(a),\bz(b),\mu\right) = \bUpsilon_{1,\mu}\left( \bz(a),\bz(b),\mu\right) -  \bUpsilon_{2,\mu}\left(\bz(a),\bz(b),\mu\right),
\end{equation}
where
\begin{equation} \label{eq_J_Ups_1za}
\bUpsilon_{1,\bz(a)}\left( \bz(a),\bz(b),\mu\right) = \begin{bmatrix} \left. \hat{H}_{\bx} \right|_a & \left. \hat{\mathbf{f}}^\mathsf{T} \right|_a & \mathbf{0}_{1 \times k_1} & \mathbf{0}_{1 \times k_2} & \left. \hat{H}_t \right|_a & 0 \\ \mathbf{0}_{n \times n} & I_{n \times n} & \mathbf{0}_{n \times k_1} & \mathbf{0}_{n \times k_2} & \mathbf{0}_{n \times 1} & \mathbf{0}_{n \times 1} \\ \bsigma_{\bx(a)} & \mathbf{0}_{k_1 \times n}  & \mathbf{0}_{k_1 \times k_1} & \mathbf{0}_{k_1 \times k_2} & \bsigma_{a} & \mathbf{0}_{k_1 \times 1} \\ \mathbf{0}_{1 \times n} & \mathbf{0}_{1 \times n} & \mathbf{0}_{1 \times k_1} & \mathbf{0}_{1 \times k_2} & 0 & \left. \hat{H}_t \right|_b \\ \mathbf{0}_{n \times n} & \mathbf{0}_{n \times n} & \mathbf{0}_{n \times k_1} & \mathbf{0}_{n \times k_2} & \mathbf{0}_{n \times 1} & \mathbf{0}_{n \times 1} \\ \mathbf{0}_{k_2 \times n} & \mathbf{0}_{k_2 \times n}  & \mathbf{0}_{k_2 \times k_1} & \mathbf{0}_{k_2 \times k_2} & \mathbf{0}_{k_2 \times 1} & \bpsi_{b}   \end{bmatrix},
\end{equation}
\begin{equation} \label{eq_J_Ups_1zb}
\bUpsilon_{1,\bz(b)}\left( \bz(a),\bz(b),\mu\right) = \begin{bmatrix} \mathbf{0}_{1 \times n} & \mathbf{0}_{1 \times n} & \mathbf{0}_{1 \times k_1} & \mathbf{0}_{1 \times k_2} & \left. \hat{H}_t \right|_a & 0 \\ \mathbf{0}_{n \times n} & \mathbf{0}_{n \times n} & \mathbf{0}_{n \times k_1} & \mathbf{0}_{n \times k_2} & \mathbf{0}_{n \times 1} & \mathbf{0}_{n \times 1} \\ \mathbf{0}_{k_1 \times n} & \mathbf{0}_{k_1 \times n}  & \mathbf{0}_{k_1 \times k_1} & \mathbf{0}_{k_1 \times k_2} & \bsigma_{a} & \mathbf{0}_{k_1 \times 1} \\ \left. \hat{H}_{\bx} \right|_b & \left. \hat{\mathbf{f}}^\mathsf{T} \right|_b & \mathbf{0}_{1 \times k_1} & \mathbf{0}_{1 \times k_2} & 0 & \left. \hat{H}_t \right|_b \\ \mathbf{0}_{n \times n} & I_{n \times n} & \mathbf{0}_{n \times k_1} & \mathbf{0}_{n \times k_2} & \mathbf{0}_{n \times 1} & \mathbf{0}_{n \times 1} \\ \bpsi_{\bx(b)} & \mathbf{0}_{k_2 \times n}  & \mathbf{0}_{k_2 \times k_1} & \mathbf{0}_{k_2 \times k_2} & \mathbf{0}_{k_2 \times 1} & \bpsi_{b}   \end{bmatrix},
\end{equation}
\begin{equation} \label{eq_J_Ups_1mu}
\bUpsilon_{1,\mu}\left( \bz(a),\bz(b),\mu\right) = \begin{bmatrix} \left. \hat{H}_{\mu} \right|_a \\ \mathbf{0}_{n \times 1} \\  \bsigma_{\mu} \\ \left. \hat{H}_{\mu} \right|_b  \\ \mathbf{0}_{n \times 1} \\ \bpsi_{\mu} \end{bmatrix},
\end{equation}
\begin{equation}  \label{eq_J_Ups_2za}
\begin{split}
\bUpsilon_{2,\bz(a)}\left( \bz(a),\bz(b),\mu\right) &= \begin{bmatrix} G_{a\bx(a)} & \mathbf{0}_{1 \times n}  & G_{a\bxi} & G_{a\bnu} & G_{aa} & G_{ab} \\ -G_{\bx(a)\bx(a)} & \mathbf{0}_{n \times n} & -G_{\bx(a)\bxi} & -G_{\bx(a)\bnu} & -G_{\bx(a)a} & -G_{\bx(a)b} \\ \mathbf{0}_{k_1 \times n} & \mathbf{0}_{k_1 \times n}  & \mathbf{0}_{k_1 \times k_1} & \mathbf{0}_{k_1 \times k_2} & \mathbf{0}_{k_1 \times 1} & \mathbf{0}_{k_1 \times 1} \\ -G_{b\bx(a)} & \mathbf{0}_{1 \times n}  & -G_{b\bxi} & -G_{b\bnu} & -G_{ba} & -G_{bb} \\ G_{\bx(b)\bx(a)} & \mathbf{0}_{n \times n} & G_{\bx(b)\bxi} & G_{\bx(b)\bnu} & G_{\bx(b)a} & G_{\bx(b)b} \\ \mathbf{0}_{k_2 \times n} & \mathbf{0}_{k_2 \times n}  & \mathbf{0}_{k_2 \times k_1} & \mathbf{0}_{k_2 \times k_2} & \mathbf{0}_{k_2 \times 1} & \mathbf{0}_{k_2 \times 1}   \end{bmatrix} \\
&= \begin{bmatrix} G_{a\bx(a)} & \mathbf{0}_{1 \times n}  & \bsigma_a^\mathsf{T} & \mathbf{0}_{1 \times k_2} & G_{aa} & G_{ab} \\ -G_{\bx(a)\bx(a)} & \mathbf{0}_{n \times n} & -\bsigma_{\bx(a)}^\mathsf{T} & \mathbf{0}_{n \times k_2} & -G_{\bx(a)a} & -G_{\bx(a)b} \\ \mathbf{0}_{k_1 \times n} & \mathbf{0}_{k_1 \times n}  & \mathbf{0}_{k_1 \times k_1} & \mathbf{0}_{k_1 \times k_2} & \mathbf{0}_{k_1 \times 1} & \mathbf{0}_{k_1 \times 1} \\ -G_{b\bx(a)} & \mathbf{0}_{1 \times n}  & \mathbf{0}_{1 \times k_1} & -\bpsi_b^\mathsf{T} & -G_{ba} & -G_{bb} \\ G_{\bx(b)\bx(a)} & \mathbf{0}_{n \times n} & \mathbf{0}_{n \times k_1} & \bpsi_{\bx(b)}^\mathsf{T} & G_{\bx(b)a} & G_{\bx(b)b} \\ \mathbf{0}_{k_2 \times n} & \mathbf{0}_{k_2 \times n}  & \mathbf{0}_{k_2 \times k_1} & \mathbf{0}_{k_2 \times k_2} & \mathbf{0}_{k_2 \times 1} & \mathbf{0}_{k_2 \times 1}   \end{bmatrix},
\end{split}
\end{equation}
\begin{equation} \label{eq_J_Ups_2zb}
\begin{split}
\bUpsilon_{2,\bz(b)}\left( \bz(a),\bz(b),\mu\right) &= \begin{bmatrix} G_{a\bx(b)} & \mathbf{0}_{1 \times n}  & G_{a\bxi} & G_{a\bnu} & G_{aa} & G_{ab} \\ -G_{\bx(a)\bx(b)} & \mathbf{0}_{n \times n} & -G_{\bx(a)\bxi} & -G_{\bx(a)\bnu} & -G_{\bx(a)a} & -G_{\bx(a)b} \\ \mathbf{0}_{k_1 \times n} & \mathbf{0}_{k_1 \times n}  & \mathbf{0}_{k_1 \times k_1} & \mathbf{0}_{k_1 \times k_2} & \mathbf{0}_{k_1 \times 1} & \mathbf{0}_{k_1 \times 1} \\ -G_{b\bx(b)} & \mathbf{0}_{1 \times n}  & -G_{b\bxi} & -G_{b\bnu} & -G_{ba} & -G_{bb} \\ G_{\bx(b)\bx(b)} & \mathbf{0}_{n \times n} & G_{\bx(b)\bxi} & G_{\bx(b)\bnu} & G_{\bx(b)a} & G_{\bx(b)b} \\ \mathbf{0}_{k_2 \times n} & \mathbf{0}_{k_2 \times n}  & \mathbf{0}_{k_2 \times k_1} & \mathbf{0}_{k_2 \times k_2} & \mathbf{0}_{k_2 \times 1} & \mathbf{0}_{k_2 \times 1}   \end{bmatrix} \\
&= \begin{bmatrix} G_{a\bx(b)} & \mathbf{0}_{1 \times n}  & \bsigma_a^\mathsf{T} & \mathbf{0}_{1 \times k_2} & G_{aa} & G_{ab} \\ -G_{\bx(a)\bx(b)} & \mathbf{0}_{n \times n} & -\bsigma_{\bx(a)}^\mathsf{T} & \mathbf{0}_{n \times k_2} & -G_{\bx(a)a} & -G_{\bx(a)b} \\ \mathbf{0}_{k_1 \times n} & \mathbf{0}_{k_1 \times n}  & \mathbf{0}_{k_1 \times k_1} & \mathbf{0}_{k_1 \times k_2} & \mathbf{0}_{k_1 \times 1} & \mathbf{0}_{k_1 \times 1} \\ -G_{b\bx(b)} & \mathbf{0}_{1 \times n}  & \mathbf{0}_{1 \times k_1} & -\bpsi_b^\mathsf{T} & -G_{ba} & -G_{bb} \\ G_{\bx(b)\bx(b)} & \mathbf{0}_{n \times n} & \mathbf{0}_{n \times k_1} & \bpsi_{\bx(b)}^\mathsf{T} & G_{\bx(b)a} & G_{\bx(b)b} \\ \mathbf{0}_{k_2 \times n} & \mathbf{0}_{k_2 \times n}  & \mathbf{0}_{k_2 \times k_1} & \mathbf{0}_{k_2 \times k_2} & \mathbf{0}_{k_2 \times 1} & \mathbf{0}_{k_2 \times 1}   \end{bmatrix},
\end{split}
\end{equation}
and
\begin{equation} \label{eq_J_Ups_2mu}
\bUpsilon_{2,\mu}\left( \bz(a),\bz(b),\mu\right) = \begin{bmatrix} G_{a\mu} \\ -G_{\bx(a)\mu} \\  \mathbf{0}_{k_1 \times 1} \\  -G_{b\mu} \\  G_{\bx(b)\mu} \\ \mathbf{0}_{k_2 \times 1} \end{bmatrix}.
\end{equation}
In equations \eqref{eq_J_Ups_1za}, \eqref{eq_J_Ups_1zb}, and \eqref{eq_J_Ups_1mu}, $\hat{\mathbf{f}}$ and all first derivatives of $\hat{H}$ in row $1$ are evaluated at $\left(a,\bx(a), \blam(a),\mu\right)$ as shown in Table~\ref{table_ODE_BVP_BC_H_a_Jacobians}, all first derivatives of $\bsigma$ in rows  $n+2$ through $n+1+k_1$ are evaluated at $\left(a,\bx(a),\mu\right)$ as shown in Table~\ref{table_ODE_BVP_BC_sigma_a_Jacobians}, $\hat{\mathbf{f}}$ and all first derivatives of $\hat{H}$ in row $n+2+k_1$ are evaluated at $\left(b,\bx(b), \blam(b),\mu\right)$ as shown in Table~\ref{table_ODE_BVP_BC_H_b_Jacobians}, and all first derivatives of $\bpsi$ in rows  $2n+3+k_1$ through $2n+2+k_1+k_2$ are evaluated at $\left(b,\bx(b),\mu\right)$ as shown in Table~\ref{table_ODE_BVP_BC_psi_b_Jacobians}. Since $\hat{H}_{\blam} = \hat{\mathbf{f}}^\mathsf{T}$, $\left. \hat{H}_{\blam} \right|_a = \left. \hat{\mathbf{f}}^\mathsf{T} \right|_a$ in row $1$ and columns $n+1$ through $2n$ of \eqref{eq_J_Ups_1za} and $\left. \hat{H}_{\blam} \right|_b = \left. \hat{\mathbf{f}}^\mathsf{T} \right|_b$ in row $n+2+k_1$ and columns $n+1$ through $2n$ of \eqref{eq_J_Ups_1zb}. In equations \eqref{eq_J_Ups_2za}, \eqref{eq_J_Ups_2zb}, and \eqref{eq_J_Ups_2mu}, all second derivatives of $G$ are evaluated at $\left(a,\bx(a), \bxi, b, \bx(b),\bnu,\mu\right)$, while all first derivatives of $\bsigma$ and $\bpsi$ are evaluated at $\left(a,\bx(a),\mu\right)$ and $\left(b,\bx(b),\mu\right)$, respectively, as shown in Tables~\ref{table_ODE_BVP_BC_sigma_a_Jacobians} and \ref{table_ODE_BVP_BC_psi_b_Jacobians}. To simplify \eqref{eq_J_Ups_2za} and \eqref{eq_J_Ups_2zb}, Clairaut's Theorem, $G_{\bxi}^\mathsf{T} = \bsigma$, and $G_{\bnu}^\mathsf{T} = \bpsi$ are used to get $G_{a\bxi} = G_{\bxi a}^\mathsf{T}= \bsigma_{a}^\mathsf{T}$, $G_{a\bnu} = G_{\bnu a}^\mathsf{T}= \bpsi_{a}^\mathsf{T}=\mathbf{0}_{1 \times k_2}$, $G_{\bx(a)\bxi} = G_{\bxi \bx(a)}^\mathsf{T} = \bsigma_{\bx(a)}^\mathsf{T}$, $G_{\bx(a)\bnu} = G_{\bnu \bx(a)}^\mathsf{T} = \bpsi_{\bx(a)}^\mathsf{T}=\mathbf{0}_{n \times k_2}$, $G_{b\bxi}=G_{\bxi b}^\mathsf{T}=\bsigma_{b}^\mathsf{T}=\mathbf{0}_{1 \times k_1}$, $G_{b\bnu}=G_{\bnu b}^\mathsf{T}=\bpsi_{b}^\mathsf{T}$, $G_{\bx(b)\bxi}=G_{\bxi \bx(b)}^\mathsf{T} = \bsigma_{\bx(b)}^\mathsf{T}=\mathbf{0}_{n \times k_1}$, and $G_{\bx(b)\bnu}=G_{\bnu\bx(b)}^\mathsf{T}=\bpsi_{\bx(b)}^\mathsf{T}$.

\begin{table}[h!]
	\centering 
	{ 
		\setlength{\extrarowheight}{1.5pt}
		\begin{tabular}{| c c c c c |} 
			\hline
			\textbf{Shorthand} & $\mathbf{\vert}$ & \textbf{Meaning} & $\mathbf{\vert}$ &
			\textbf{Simplification} \\ 
			\hline\hline 
			
			$\left. \hat{H}_{\bx} \right|_a$ &=& $\hat{H}_{\bx} \left( a,\bx(a),\blam(a),\mu \right) $&=& $H_{\bx} \left( a,\bx(a),\blam(a),\bpi\left(a,\bx(a),\blam(a),\mu\right),\mu \right) $ \\ \hline
			
			$\left. \hat{\mathbf{f}}^\mathsf{T} \right|_a$ &=& $\hat{\mathbf{f}}^\mathsf{T} \left( a,\bx(a),\blam(a),\mu \right) $&=& $\mathbf{f}^\mathsf{T} \left( a,\bx(a),\blam(a),\bpi\left(a,\bx(a),\blam(a),\mu\right),\mu \right) $ \\ \hline
			
			$\left. \hat{H}_t \right|_a$ &=& $\hat{H}_t \left( a,\bx(a),\blam(a),\mu \right) $&=& $H_t \left( a,\bx(a),\blam(a),\bpi\left(a,\bx(a),\blam(a),\mu\right),\mu \right) $ \\ \hline
			
			$\left. \hat{H}_\mu \right|_a$ &=& $\hat{H}_\mu \left( a,\bx(a),\blam(a),\mu \right) $&=& $H_\mu \left( a,\bx(a),\blam(a),\bpi\left(a,\bx(a),\blam(a),\mu\right),\mu \right) $ \\ \hline
			
		\end{tabular} 
	}
	\caption{Explanation of shorthand notation for $\hat{\mathbf{f}}$ and first derivatives of $\hat{H}$ evaluated at $a$ used in \eqref{eq_J_Ups_1za}, \eqref{eq_J_Ups_1zb}, and \eqref{eq_J_Ups_1mu}. Note that $\left. \hat{H}_{\blam} \right|_a = \hat{H}_{\blam} \left( a,\bx(a),\blam(a),\mu \right) =\left. \hat{\mathbf{f}}^\mathsf{T} \right|_a$.}
	\label{table_ODE_BVP_BC_H_a_Jacobians}
\end{table}

\begin{table}[h!]
	\centering 
	{ 
		\setlength{\extrarowheight}{1.5pt}
		\begin{tabular}{| c c c |} 
			\hline
			\textbf{Shorthand} & $\mathbf{\vert}$ & \textbf{Meaning} \\ 
			\hline\hline 
			
			$\bsigma_{\bx(a)}$ &=& $\bsigma_{\bx(a)} \left( a,\bx(a),\mu \right) $ \\ \hline
			
			$\bsigma_{a}$ &=& $\bsigma_{a} \left( a,\bx(a),\mu \right) $ \\ \hline
			
			$\bsigma_{\mu}$ &=& $\bsigma_{\mu} \left( a,\bx(a),\mu \right) $ \\ \hline
			
		\end{tabular} 
	}
	\caption{Explanation of shorthand notation for first derivatives of $\bsigma$ used in \eqref{eq_J_Ups_1za}, \eqref{eq_J_Ups_1zb}, and \eqref{eq_J_Ups_1mu}.}
	\label{table_ODE_BVP_BC_sigma_a_Jacobians}
\end{table}

\begin{table}[h!]
	\centering 
	{ 
		\setlength{\extrarowheight}{1.5pt}
		\begin{tabular}{| c c c c c |} 
			\hline
			\textbf{Shorthand} & $\mathbf{\vert}$ & \textbf{Meaning} & $\mathbf{\vert}$ &
			\textbf{Simplification} \\ 
			\hline\hline 
			
			$\left. \hat{H}_{\bx} \right|_b$ &=& $\hat{H}_{\bx} \left( b,\bx(b),\blam(b),\mu \right) $&=& $H_{\bx} \left( b,\bx(b),\blam(b),\bpi\left(b,\bx(b),\blam(b),\mu\right),\mu \right) $ \\ \hline
			
			$\left. \hat{\mathbf{f}}^\mathsf{T} \right|_b$ &=& $\hat{\mathbf{f}}^\mathsf{T} \left( b,\bx(b),\blam(b),\mu \right) $&=& $\mathbf{f}^\mathsf{T} \left( b,\bx(b),\blam(b),\bpi\left(b,\bx(b),\blam(b),\mu\right),\mu \right) $ \\ \hline
			
			$\left. \hat{H}_t \right|_b$ &=& $\hat{H}_t \left( b,\bx(b),\blam(b),\mu \right) $&=& $H_t \left( b,\bx(b),\blam(b),\bpi\left(b,\bx(b),\blam(b),\mu\right),\mu \right) $ \\ \hline
			
			$\left. \hat{H}_\mu \right|_b$ &=& $\hat{H}_\mu \left( b,\bx(b),\blam(b),\mu \right) $&=& $H_\mu \left( b,\bx(b),\blam(b),\bpi\left(b,\bx(b),\blam(b),\mu\right),\mu \right) $ \\ \hline
			
		\end{tabular} 
	}
	\caption{Explanation of shorthand notation for $\hat{\mathbf{f}}$ and first derivatives of $\hat{H}$ evaluated at $b$ used in \eqref{eq_J_Ups_1za}, \eqref{eq_J_Ups_1zb}, and \eqref{eq_J_Ups_1mu}. Note that $\left. \hat{H}_{\blam} \right|_b= \hat{H}_{\blam} \left( b,\bx(b),\blam(b),\mu \right) =\left. \hat{\mathbf{f}}^\mathsf{T} \right|_b$.}
	\label{table_ODE_BVP_BC_H_b_Jacobians}
\end{table}

\begin{table}[h!]
	\centering 
	{ 
		\setlength{\extrarowheight}{1.5pt}
		\begin{tabular}{| c c c |} 
			\hline
			\textbf{Shorthand} & $\mathbf{\vert}$ & \textbf{Meaning} \\ 
			\hline\hline 
			
			$\bpsi_{\bx(b)}$ &=& $\bpsi_{\bx(b)} \left( b,\bx(b),\mu \right) $ \\ \hline
			
			$\bpsi_{b}$ &=& $\bpsi_{b} \left( b,\bx(b),\mu \right) $ \\ \hline
			
			$\bpsi_{\mu}$ &=& $\bpsi_{\mu} \left( b,\bx(b),\mu \right) $ \\ \hline
			
		\end{tabular} 
	}
	\caption{Explanation of shorthand notation for first derivatives of $\bpsi$ used in \eqref{eq_J_Ups_1za}, \eqref{eq_J_Ups_1zb}, and \eqref{eq_J_Ups_1mu}.}
	\label{table_ODE_BVP_BC_psi_b_Jacobians}
\end{table} 

To express the boundary conditions \eqref{eq_pmp_lbc2}-\eqref{eq_pmp_rbc2} in terms of normalized dependent variables, let $\tilde \bUpsilon_1\left(\tilde \bz(0),\tilde \bz(1),\mu\right) \equiv \bUpsilon_1\left(\bz(a),\bz(b),\mu\right)$, $\tilde \bUpsilon_2\left(\tilde \bz(0),\tilde \bz(1),\mu\right) \equiv \bUpsilon_2\left(\bz(a),\bz(b),\mu\right)$, and $\tilde \bUpsilon\left(\tilde \bz(0),\tilde \bz(1),\mu\right) \equiv \bUpsilon\left(\bz(a),\bz(b),\mu\right)$. Thus
\begin{equation} \label{eq_tUps_1}
\tilde \bUpsilon_1\left(\tilde \bz(0),\tilde \bz(1),\mu\right) =\begin{bmatrix} \hat{H} \left(a,\tilde \bx(0), \tilde \blam(0),\mu\right) \\ \tilde \blam(0) \\ \bsigma\left(a,\tilde \bx(0),\mu\right) \\ \hat{H} \left(b,\tilde \bx(1), \tilde \blam(1),\mu\right) \\ \tilde \blam(1) \\ \bpsi\left(b,\tilde \bx(1),\mu\right) \end{bmatrix},
\end{equation}
\begin{equation} \label{eq_tUps_2}
\tilde \bUpsilon_2\left(\tilde \bz(0),\tilde \bz(1),\mu\right) = \begin{bmatrix} G_a\left(a,\tilde \bx(0), \bxi, b, \tilde \bx(1),\bnu,\mu\right) \\ -G_{\bx(a)}^\mathsf{T}\left(a,\tilde \bx(0), \bxi, b, \tilde \bx(1),\bnu,\mu\right) \\ \mathbf{0}_{k_1 \times 1} \\ -G_b\left(a,\tilde \bx(0), \bxi, b, \tilde \bx(1),\bnu,\mu\right) \\ G_{\bx(b)}^\mathsf{T}\left(a,\tilde \bx(0), \bxi, b, \tilde \bx(1),\bnu,\mu\right) \\ \mathbf{0}_{k_2 \times 1} \end{bmatrix},
\end{equation}
and
\begin{equation} \label{eq_tUps}
\begin{split}
\tilde \bUpsilon\left(\tilde \bz(0),\tilde \bz(1),\mu\right) &=\tilde \bUpsilon_1\left(\tilde \bz(0),\tilde \bz(1),\mu\right) - \tilde \bUpsilon_2\left(\tilde \bz(0),\tilde \bz(1),\mu\right) \\ &=\begin{bmatrix} \hat{H} \left(a,\tilde \bx(0), \tilde \blam(0),\mu\right) \\ \tilde \blam(0) \\ \bsigma\left(a,\tilde \bx(0),\mu\right) \\ \hat{H} \left(b,\tilde \bx(1), \tilde \blam(1),\mu\right) \\ \tilde \blam(1) \\ \bpsi\left(b,\tilde \bx(1),\mu\right) \end{bmatrix} - \begin{bmatrix} G_a\left(a,\tilde \bx(0), \bxi, b, \tilde \bx(1),\bnu,\mu\right) \\ -G_{\bx(a)}^\mathsf{T}\left(a,\tilde \bx(0), \bxi, b, \tilde \bx(1),\bnu,\mu\right) \\ \mathbf{0}_{k_1 \times 1} \\ -G_b\left(a,\tilde \bx(0), \bxi, b, \tilde \bx(1),\bnu,\mu\right) \\ G_{\bx(b)}^\mathsf{T}\left(a,\tilde \bx(0), \bxi, b, \tilde \bx(1),\bnu,\mu\right) \\ \mathbf{0}_{k_2 \times 1} \end{bmatrix},
\end{split}
\end{equation}
and the boundary conditions \eqref{eq_pmp_lbc2}-\eqref{eq_pmp_rbc2} in normalized dependent variables are given by the normalized two-point boundary condition function
\begin{equation} \label{eq_tUps_BC}
\tilde \bUpsilon\left(\tilde \bz(0),\tilde \bz(1),\mu\right) = \mathbf{0}_{(2n+k_1+k_2+2) \times 1}.
\end{equation}
The Jacobians of $\tilde \bUpsilon$ with respect to $\tilde \bz(0)$, $\tilde \bz(1)$, and $\mu$ are
\begin{equation} \label{eq_Jac_tUps_tz0}
\tilde \bUpsilon_{\tilde \bz(0)}\left(\tilde \bz(0),\tilde \bz(1),\mu\right) = \tilde \bUpsilon_{1,\tilde \bz(0)}\left( \tilde \bz(0),\tilde \bz(1),\mu\right) -  \bUpsilon_{2,\tilde \bz(0)}\left(\tilde \bz(0),\tilde \bz(1),\mu\right),
\end{equation}
\begin{equation} \label{eq_Jac_tUps_tz1}
\tilde \bUpsilon_{\tilde \bz(1)}\left(\tilde \bz(0),\tilde \bz(1),\mu\right) = \tilde \bUpsilon_{1,\tilde \bz(1)}\left(\tilde \bz(0),\tilde \bz(1),\mu\right) -  \tilde \bUpsilon_{2,\tilde \bz(1)}\left(\tilde \bz(0),\tilde \bz(1),\mu\right),
\end{equation}
and
\begin{equation} \label{eq_Jac_tUps_mu}
\tilde \bUpsilon_{\mu}\left(\tilde \bz(0),\tilde \bz(1),\mu\right) = \tilde \bUpsilon_{1,\mu}\left(\tilde \bz(0),\tilde \bz(1),\mu\right) -  \tilde \bUpsilon_{2,\mu}\left(\tilde \bz(0),\tilde \bz(1),\mu\right),
\end{equation}
where the equality between the Jacobians of $\tilde \bUpsilon$, $\tilde \bUpsilon_1$, and $\tilde \bUpsilon_2$ with respect to $\tilde \bz(0)$, $\tilde \bz(1)$, and $\mu$ and the Jacobians of $\bUpsilon$, $\bUpsilon_1$, and $\bUpsilon_2$ with respect to $\bz(0)$, $\bz(1)$, and $\mu$ is given in Table~\ref{table_ODE_BVP_BC_Jacobians}.

\begin{table}[h!]
	\centering 
	{ 
		\setlength{\extrarowheight}{1.5pt}
		\begin{tabular}{| c c c |} 
			\hline
			\textbf{Normalized} & $\mathbf{\vert}$ & \textbf{Un-Normalized} \\ 
			\hline\hline 
			$\tilde \bUpsilon_{\tilde \bz(0)}\left(\tilde \bz(0),\tilde \bz(1),\mu\right)$ & = & $\bUpsilon_{\bz(a)}\left(\bz(a),\bz(b),\mu\right)$ \\ \hline
			$\tilde \bUpsilon_{1,\tilde \bz(0)}\left(\tilde \bz(0),\tilde \bz(1),\mu\right)$&= &$\bUpsilon_{1,\bz(a)}\left(\bz(a),\bz(b),\mu\right)$ \\ \hline
			$\tilde \bUpsilon_{2,\tilde \bz(0)}\left(\tilde \bz(0),\tilde \bz(1),\mu\right)$&= &$\bUpsilon_{2,\bz(a)}\left(\bz(a),\bz(b),\mu\right)$ \\ \hline
			$\tilde \bUpsilon_{\tilde \bz(1)}\left(\tilde \bz(0),\tilde \bz(1),\mu\right) $&= &$ \bUpsilon_{\bz(b)}\left(\bz(a),\bz(b),\mu\right)$\\ \hline $\tilde \bUpsilon_{1,\tilde \bz(1)}\left(\tilde \bz(0),\tilde \bz(1),\mu\right)$&= &$\bUpsilon_{1,\bz(b)}\left(\bz(a),\bz(b),\mu\right)$\\ \hline$\tilde \bUpsilon_{2,\tilde \bz(1)}\left(\tilde \bz(0),\tilde \bz(1),\mu\right)$&= &$\bUpsilon_{2,\bz(b)}\left(\bz(a),\bz(b),\mu\right)$\\ \hline $\tilde \bUpsilon_{\mu}\left(\tilde \bz(0),\tilde \bz(1),\mu\right) $&= &$ \bUpsilon_{\mu}\left(\bz(a),\bz(b),\mu\right)$ \\ \hline $\tilde \bUpsilon_{1,\mu}\left(\tilde \bz(0),\tilde \bz(1),\mu\right)$&= &$\bUpsilon_{1,\mu}\left(\bz(a),\bz(b),\mu\right)$ \\ \hline $\tilde \bUpsilon_{2,\mu}\left(\tilde \bz(0),\tilde \bz(1),\mu\right)$&= &$\bUpsilon_{2,\mu}\left(\bz(a),\bz(b),\mu\right)$ \\ \hline
		\end{tabular} 
	}
	\caption{Equality between Jacobians of two-point boundary condition functions in normalized and un-normalized coordinates.}
	\label{table_ODE_BVP_BC_Jacobians}
\end{table}

Special care must be taken when implementing the Jacobians \eqref{eq_Jac_tUps_tz0} and \eqref{eq_Jac_tUps_tz1}. Since the unknown constants $\bxi$, $\bnu$, $a$, and $b$ appear at the end of both $\tilde \bz(0)$ and $\tilde \bz(1)$, the unknown constants from only one of $\tilde \bz(0)$ and $\tilde \bz(1)$ are actually used to construct each term in $\tilde \bUpsilon$ involving $\bxi$, $\bnu$, $a$, and $b$. The trailing columns in \eqref{eq_Jac_tUps_tz0} are actually the Jacobian of $\tilde \bUpsilon$ with respect to $\bxi$, $\bnu$, $a$, and $b$ in $\tilde \bz(0)$, while the trailing columns in \eqref{eq_Jac_tUps_tz1} are actually the Jacobian of $\tilde \bUpsilon$ with respect to $\bxi$, $\bnu$, $a$, and $b$ in $\tilde \bz(1)$. Thus, the trailing columns in \eqref{eq_Jac_tUps_tz0} and \eqref{eq_Jac_tUps_tz1} corresponding to the Jacobian of $\tilde \bUpsilon$ with respect to $\bxi$, $\bnu$, $a$, and $b$ should not coincide in a software implementation. For example, if the unknown constants are extracted from $\tilde \bz(0)$ to construct $\tilde \bUpsilon$, $\tilde \bUpsilon_{\tilde \bz(0)}$ is as shown in \eqref{eq_Jac_tUps_tz0} while the trailing columns in \eqref{eq_Jac_tUps_tz1} corresponding to the Jacobian of $\tilde \bUpsilon$ with respect to the unknown constants in $\tilde \bz(1)$ should be all zeros. Alternatively, if the unknown constants are extracted from $\tilde \bz(1)$ to construct $\tilde \bUpsilon$, $\tilde \bUpsilon_{\tilde \bz(1)}$ is as shown in \eqref{eq_Jac_tUps_tz1} while the trailing columns in \eqref{eq_Jac_tUps_tz0} corresponding to the Jacobian of $\tilde \bUpsilon$ with respect to the unknown constants in $\tilde \bz(0)$ should be all zeros.

\subsection{Final Details}
In equations \eqref{eq_z_def}, \eqref{eq_tz_def}, \eqref{eq_dtz}, \eqref{eq_Phi}, and \eqref{eq_DPhi_mu}, the second to last row is needed only if the initial time $a$ is free and the last row is needed only if the final time $b$ is free. In equation \eqref{eq_DPhi_z}, the second to last row and column are needed only if the initial time $a$ is free and the last row and column are needed only if the final time $b$ is free. 

In equations \eqref{eq_Ups_1}, \eqref{eq_Ups_2}, \eqref{eq_Ups}, \eqref{eq_Ups_BC}, \eqref{eq_Jac_Ups1_mu}, \eqref{eq_J_Ups_1mu}, \eqref{eq_J_Ups_2mu}, \eqref{eq_tUps_1}, \eqref{eq_tUps_2}, \eqref{eq_tUps}, \eqref{eq_tUps_BC}, and \eqref{eq_Jac_tUps_mu} the first row is needed only if the initial time $a$ is free 
and row $n+k_1+2$  is needed only if the final time $b$ is free. In equations \eqref{eq_Jac_Ups1_za}, \eqref{eq_Jac_Ups1_zb}, \eqref{eq_J_Ups_1za}, \eqref{eq_J_Ups_1zb}, \eqref{eq_J_Ups_2za}, \eqref{eq_J_Ups_2zb}, \eqref{eq_Jac_tUps_tz0}, and \eqref{eq_Jac_tUps_tz1} the first row and second to last column are needed only if the initial time $a$ is free and row $n+k_1+2$ and the last column are needed only if the final time $b$ is free.

In order to numerically solve the ODE TPBVP \eqref{eq_pmp_bvp}, \eqref{eq_pmp_lbc2}, and \eqref{eq_pmp_rbc2} without continuation or with a monotonic continuation solver (such as \mcode{acdc} or \mcode{acdcc}), the solver should be provided \eqref{eq_Phi}, \eqref{eq_DPhi_z}, \eqref{eq_tUps}, \eqref{eq_Jac_tUps_tz0}, and \eqref{eq_Jac_tUps_tz1}. In order to numerically solve the ODE TPBVP \eqref{eq_pmp_bvp}, \eqref{eq_pmp_lbc2}, and \eqref{eq_pmp_rbc2} with a non-monotonic continuation solver (such as the predictor-corrector methods discussed in Appendices~\ref{app_predictor_corrector} and \ref{app_sweep_predictor_corrector}), the solver should be provided \eqref{eq_Phi}, \eqref{eq_DPhi_z}, \eqref{eq_DPhi_mu}, \eqref{eq_tUps}, \eqref{eq_Jac_tUps_tz0}, \eqref{eq_Jac_tUps_tz1}, and \eqref{eq_Jac_tUps_mu}. 

The first and second derivatives required to construct \eqref{eq_Phi}, \eqref{eq_DPhi_z}, \eqref{eq_DPhi_mu}, \eqref{eq_tUps}, \eqref{eq_Jac_tUps_tz0}, \eqref{eq_Jac_tUps_tz1}, and \eqref{eq_Jac_tUps_mu} are generally quite tedious to derive manually. Instead, symbolic differentiation \cite{hardy2008computer}, complex/bicomplex step differentiation \cite{squire1998using,martins2001connection,martins2003complex,lantoine2012using}, dual/hyper-dual numbers \cite{fike2011development,fike2011optimization,fike2012automatic,neuenhofen2018review}, and automatic differentiation \cite{corliss2002automatic,naumann2012art} are computational alternatives. If fact, it may be shown that the use of dual/hyper-dual numbers to compute first and second derivatives is equivalent to automatic differentiation \cite{neuenhofen2018review}. While symbolic differentiation suffers from expression explosion and complex/bicomplex step differentiation only applies to real analytic functions, dual/hyper-dual numbers and automatic differentiation are more robust and broadly-applicable. 

Therefore, while \eqref{eq_Phi}, \eqref{eq_DPhi_z}, \eqref{eq_DPhi_mu}, \eqref{eq_tUps}, \eqref{eq_Jac_tUps_tz0}, \eqref{eq_Jac_tUps_tz1}, and \eqref{eq_Jac_tUps_mu} are complicated, they may be readily constructed numerically through automatic differentiation of $H$, $\bpi$, $\hat{\mathbf{f}}$, $G$, $\bsigma$, and $\bpsi$ if $\bpi$ is given analytically and of $H$, $\mathbf{f}$, $G$, $\bsigma$, and $\bpsi$ if the value of $\bpi$ is constructed numerically via Newton's method as in \eqref{eq_pi_Newton}. There are many free automatic differentiation toolboxes available \cite{autodiff}, such as the \mcode{MATLAB} automatic differentiation toolbox \mcode{ADiGator} \cite{weinstein2017algorithm, weinstein2015utilizing}. Moreover, \mcode{ADiGator} is able to construct vectorized automatic derivatives, which is extremely useful for realizing the vectorized version of \eqref{eq_Phi}, \eqref{eq_DPhi_z}, and \eqref{eq_DPhi_mu}, as the non-vectorized version of these equations execute too slowly in \mcode{MATLAB} to solve the ODE TPBVP \eqref{eq_pmp_bvp}, \eqref{eq_pmp_lbc2}, and \eqref{eq_pmp_rbc2} in a timely manner. \rem{The use of vectorized automatic differentiation to realize \eqref{eq_Phi}, \eqref{eq_DPhi_z}, and \eqref{eq_DPhi_mu} in \mcode{MATLAB} is noteworthy, because not only is it tedious to manually derive the non-vectorized version of these equations, but it is terribly difficult to manually derive the vectorized version of these equations.}
	
\section{Predictor-Corrector Continuation Method for Solving an ODE TPBVP} \label{app_predictor_corrector}

\subsection{Introduction}
Suppose it is desired to solve the ODE TPBVP:
\begin{equation} \label{eqn_ode_tpbvp}
\begin{split}
\dd{}{s} \mathbf{y}(s) &= \mathbf{F} \left(s,\mathbf{y}(s),\ACP \right), \\
\mathbf{G}\left( \mathbf{y} (a),\mathbf{y}(b),\ACP \right) &= \mathbf{0}_{n \times 1},
\end{split}
\end{equation}
where $a,b \in \mathbb{R}$ are prescribed with $a < b$, $s \in \left[a,b\right] \subset \mathbb{R}$ is the independent variable, $n \in \mathbb{N}$ is the prescribed number of dependent variables in $\mathbf{y}$, $\mathbf{y} \colon \left[a,b\right] \to \mathbb{R}^n$ is an unknown function which must be solved for, $\ACP \in \mathbb{R}$ is a prescribed scalar parameter, $\mathbf{F} \colon \left[a,b\right] \times \mathbb{R}^n \times \mathbb{R} \to \mathbb{R}^n$ is a prescribed ODE velocity function defining the velocity of $\mathbf{y}$, and $\mathbf{G} \colon \mathbb{R}^n \times \mathbb{R}^n \times \mathbb{R} \to \mathbb{R}^{n}$ is a prescribed two-point boundary condition function. Observe that if $n=1$,  $\mathbf{y}$, $\mathbf{F}$, and $\mathbf{G}$ are scalar-valued functions, while if $n>1$, $\mathbf{y}$, $\mathbf{F}$, and $\mathbf{G}$ are vector-valued functions. The Jacobian of $\mathbf{F}$ with respect to $\mathbf{y}$ is  $\mathbf{F}_{\mathbf{y}} \colon \left[a,b\right] \times \mathbb{R}^n \times \mathbb{R} \to \mathbb{R}^{n \times n}$ and the Jacobian of $\mathbf{F}$ with respect to $\ACP$ is $\mathbf{F}_{\ACP} \colon \left[a,b\right] \times \mathbb{R}^n \times \mathbb{R} \to \mathbb{R}^{n \times 1}$. The Jacobian of $\mathbf{G}$ with respect to $\mathbf{y}(a)$ is $\mathbf{G}_{\mathbf{y}(a)} \colon \mathbb{R}^n \times \mathbb{R}^n \times \mathbb{R} \to \mathbb{R}^{n \times n}$, the Jacobian of $\mathbf{G}$ with respect to $\mathbf{y}(b)$ is $\mathbf{G}_{\mathbf{y}(b)} \colon \mathbb{R}^n \times \mathbb{R}^n \times \mathbb{R} \to \mathbb{R}^{n \times n}$, and the Jacobian of $\mathbf{G}$ with respect to $\ACP$ is $\mathbf{G}_{\ACP} \colon \mathbb{R}^n \times \mathbb{R}^n \times \mathbb{R} \to \mathbb{R}^{n \times 1}$. If $\mathbf{F}$ is linear in $\mathbf{y}$ and $\mathbf{G}$ is linear in $\mathbf{y} (a)$ and $\mathbf{y}(b)$, then \eqref{eqn_ode_tpbvp} is said to be a linear ODE TPBVP; otherwise, \eqref{eqn_ode_tpbvp} is said to be a nonlinear ODE TPBVP.

Note that a solution $\mathbf{y}$ to \eqref{eqn_ode_tpbvp} depends on the given value of the scalar parameter $\ACP$, so a solution to \eqref{eqn_ode_tpbvp} will be denoted by the pair $\left(\mathbf{y},\ACP \right)$. Usually it is not possible to solve \eqref{eqn_ode_tpbvp} analytically. Instead, a numerical method such as a shooting, finite-difference, or Runge-Kutta method (collocation is a special kind of Runge-Kutta method) must be utilized to construct an approximate solution to \eqref{eqn_ode_tpbvp}. All such numerical methods require an initial solution guess and convergence to a solution is guaranteed only if the initial solution guess is sufficiently near the solution. Thus, solving \eqref{eqn_ode_tpbvp} numerically requires construction of a good initial solution guess. 

One way to construct a good initial solution guess for \eqref{eqn_ode_tpbvp} is through continuation in the scalar parameter $\ACP$. If $\left(\mathbf{y}_{\mathrm I} ,\ACP_{\mathrm I} \right)$ solves \eqref{eqn_ode_tpbvp} and it is desired to solve \eqref{eqn_ode_tpbvp} for $\ACP = \ACP_{\mathrm F}$, it may be possible to construct a finite sequence of solutions $\left\{\left(\mathbf{y}_j ,\ACP_j \right)\right\}_{j=1}^{J+1}$ starting at the known solution $\left(\mathbf{y}_1 ,\ACP_1 \right) = \left(\mathbf{y}_{\mathrm I} ,\ACP_{\mathrm I} \right)$ and ending at the desired solution $\left(\mathbf{y}_{J+1} ,\ACP_{J+1} \right) = \left(\mathbf{y}_{\mathrm F} ,\ACP_{\mathrm F} \right)$, using the previous solution $\left(\mathbf{y}_j ,\ACP_j \right)$ as an initial solution guess for the numerical solver to obtain the next solution $\left(\mathbf{y}_{j+1} ,\ACP_{j+1} \right)$, $1 \le j \le J$, in the sequence. $J \in \mathbb{N}$ denotes the number of solutions in the sequence after the known solution. 

This appendix describes a particular such continuation method, called predictor-corrector continuation, for solving \eqref{eqn_ode_tpbvp}. The treatment given here follows \cite{trefethen2013numerical}. In the literature, predictor-corrector continuation is also called  path-following \cite{kitzhofer2009pathfollowing}, predictor-corrector path-following \cite{allgower1993continuation}, and differential path-following \cite{caillau2012differential}. AUTO \cite{doedel2007auto}, COLCON \cite{bader1989continuation}, and the algorithm presented in \cite{weinmuller1988pathfollowing} are Fortran predictor-corrrector continuation codes, while \mcode{bvpsuite1.1} \cite{kitzhofer2010new}, \mcode{Chebfun}'s \mcode{followpath} \cite{trefethen2013numerical}, and \mcode{COCO} \cite{dankowicz2013recipes} are \mcode{MATLAB} predictor-corrrector continuation codes. All these codes rely on global methods for solving ODE BVPs (e.g. Runge-Kutta, collocation, and finite-difference schemes), which are more robust than initial value methods for solving ODE BVPs (i.e. single and multiple shooting) because initial value methods cannot integrate unstable ODEs \cite{ascher1994numerical,muir1999optimal,boisvert2011problem}.

Before delving into the details, some functional analysis is reviewed which is necessary to understand how the predictor-corrector continuation method is applied to solve \eqref{eqn_ode_tpbvp}. 

\subsection{A Hilbert Space}
Let $\mathcal{H} = \left\{ \left(\mathbf{y} ,\ACP \right)  : \mathbf{y} \in L^2\left(\left[a,b\right], \mathbb{R}^n \right), \ACP \in \mathbb{R}  \right\}$. $\mathcal{H}$ is a Hilbert space over $\mathbb{R}$. If $\alpha, \beta \in \mathbb{R}$ and $\left(\mathbf{y} ,\ACP \right), \left(\tilde {\mathbf{y}} , \tilde \ACP \right) \in \mathcal{H}$, then
\begin{equation}
\alpha \left(\mathbf{y} ,\ACP \right) + \beta \left(\tilde {\mathbf{y}} , \tilde \ACP \right)  = \left(\alpha \mathbf{y}+  \beta \tilde {\mathbf{y}} ,\alpha \ACP +\beta \tilde \ACP \right),
\end{equation}
the inner product on $\mathcal{H}$ is
\begin{equation}
\left<\left(\mathbf{y} ,\ACP \right) , \left(\tilde {\mathbf{y}} , \tilde \ACP \right) \right> = \int_a^b \mathbf{y}^\mathsf{T}(s)  \tilde {\mathbf{y}}(s)  \mathrm{d}s + \ACP \tilde \ACP,
\end{equation}
and the norm on $\mathcal{H}$, induced by the inner product, is
\begin{equation}
\left \Vert \left(\mathbf{y} ,\ACP \right) \right \Vert = \left<\left(\mathbf{y} ,\ACP \right) , \left(\mathbf{y} , \ACP \right) \right>^{\frac{1}{2}} = \left[ \int_a^b \mathbf{y}^\mathsf{T}(s)  \mathbf{y}(s)  \mathrm{d}s + \ACP^2 \right]^{\frac{1}{2}}.
\end{equation}
\revisionNACO{R2Q3}{$\left(\mathbf{y} ,\ACP \right) \in \mathcal{H}$ and $\left(\tilde {\mathbf{y}} , \tilde \ACP \right) \in \mathcal{H}$ are said to be orthogonal iff $\left<\left(\mathbf{y} ,\ACP \right) , \left(\tilde {\mathbf{y}} , \tilde \ACP \right) \right> = 0$, and $\left(\mathbf{y} ,\ACP \right) \in \mathcal{H}$ is said to be of unit length iff $\left \Vert \left(\mathbf{y} ,\ACP \right) \right \Vert = 1$.}

\rem{
$\left(\mathbf{y} ,\ACP \right) \in \mathcal{H}$ and $\left(\tilde {\mathbf{y}} , \tilde \ACP \right) \in \mathcal{H}$ are said to be orthogonal if
\begin{equation}
\left<\left(\mathbf{y} ,\ACP \right) , \left(\tilde {\mathbf{y}} , \tilde \ACP \right) \right> = \int_a^b \mathbf{y}^\mathsf{T}(s)  \tilde {\mathbf{y}}(s)  \mathrm{d}s + \ACP \tilde \ACP = 0.
\end{equation}
$\left(\mathbf{y} ,\ACP \right) \in \mathcal{H}$ is said to be of unit length if
\begin{equation}
\left \Vert \left(\mathbf{y} ,\ACP \right) \right \Vert = \left<\left(\mathbf{y} ,\ACP \right) , \left(\mathbf{y} , \ACP \right) \right>^{\frac{1}{2}} = \left[ \int_a^b \mathbf{y}^\mathsf{T}(s)  \mathbf{y}(s)  \mathrm{d}s + \ACP^2 \right]^{\frac{1}{2}}=1.
\end{equation} }

\subsection{The Fr\'echet Derivative and Newton's Method}
Given a function $\mathbf{F}: \mathbb{R}^n \to \mathbb{R}^m$, recall that ordinary vector calculus defines the Jacobian of $\mathbf{F}$ as the function $\mathbf{F}': \mathbb{R}^n \to \mathbb{R}^{m \times n}$ such that $\mathbf{F}'\left(\bx\right)$ is the linearization of $\mathbf{F}$ at $\bx \in \mathbb{R}^n$. Given normed vector spaces $V$ and $W$ and an open subset $U$ of $V$, the Fr\'echet derivative is an extension of the Jacobian to an operator $\mathcal{F}: U \to W$. Before giving the definition of the Fr\'echet derivative, recall that $L\left(V,W\right)$ denotes the space of continuous linear operators from $V$ to $W$.  Now for the definition of the Fr\'echet derivative, which comes from Definition 2.2.4 of \cite{trefethen2013numerical}.
\begin{definition}
	Suppose that $V$ and $W$ are normed vector spaces, and let $U$ be an open subset of $V$. Then the operator $\mathcal{F}: U \to W$ is said to
	be Fr\'echet differentiable at $u \in U$ if and only if there exists an operator $\mathcal{L} \in L\left(V,W\right)$ such that 
	\begin{equation}
	\lim_{\left \Vert h \right \Vert_V \to 0} \frac{\left \Vert \mathcal{F}\left(u+h\right)-\mathcal{F}\left(u\right)-\mathcal{L}h \right \Vert_W}{\left \Vert h \right \Vert_V} = 0.
	\end{equation}
	The operator $\mathcal{L}$ is then called the Fr\'echet derivative of $\mathcal{F}$ at $u$, often denoted by $\mathcal{F}'(u)$.
	If $\mathcal{F}$ is Fr\'echet differentiable at all points in $U$, $\mathcal{F}$ is said to be Fr\'echet differentiable in $U$.
\end{definition}

Given a function $\mathbf{H}: \mathbb{R}^m \to \mathbb{R}^m$, Newton's method is an algorithm to solve $\mathbf{H}\left(\bx\right)=\mathbf{0}$ for $\bx \in \mathbb{R}^m$ and $\mathbf{0} \in \mathbb{R}^m$ when $\mathbf{H}$ satisfies certain mild conditions. Starting from an initial solution guess $\bx_0 \in \mathbb{R}^m$ sufficiently close to a solution, Newton's method converges to a solution of $\mathbf{H}\left(\bx\right)=\mathbf{0}$ by iteratively solving the equations
\begin{equation} \label{eq_newton}
\mathbf{H}'\left(\bx_k\right) \de \bx_k = - \mathbf{H}\left(\bx_k\right), \quad \bx_{k+1} = \bx_k + \de \bx_k,
\end{equation}
starting at $k=0$, where $\mathbf{H}'$ denotes the Jacobian of $\mathbf{H}$ and $\bx_k,\de \bx_k \in \mathbb{R}^m$ for $k \ge 0$. The iteration in \eqref{eq_newton} continues until $\mathbf{H}\left(\bx_k\right) \approx \mathbf{0}$ (or $\de \bx_k \approx \mathbf{0}$) or until $k$ exceeds a maximum iteration threshold. 

Now consider an operator $\mathcal{H}: U \subset V \to W$, where $V$ and $W$ are Banach spaces and $U$ is an open subset of $V$. Kantorovich \cite{kantorovich1948newton} provided an extension of Newton's method to solve $\mathcal{H}(u)=0$ for $u \in U$ and $0 \in W$ when $\mathcal{H}$ satisfies certain mild conditions.  Starting from an initial solution guess $u_0 \in U$ sufficiently close to a solution, Kantorovich's extension of Newton's method converges to a solution of $\mathcal{H}(u)=0$ by iteratively solving the equations
\begin{equation} \label{eq_kantorovich}
\mathcal{H}'\left(u_k\right) \de u_k = - \mathcal{H}\left(u_k\right), \quad u_{k+1} = u_k + \de u_k,
\end{equation}
starting at $k=0$, where $\mathcal{H}'$ denotes the Fr\'echet derivative of $\mathcal{H}$ and $u_k,\de u_k \in U$ for $k \ge 0$. The iteration in \eqref{eq_kantorovich} continues until $\mathcal{H}\left(u_k\right) \approx 0$ (or $\de u_k \approx 0$) or until $k$ exceeds a maximum iteration threshold.

\subsection{The Davidenko ODE IVP}
To motivate the predictor-corrector continuation method, the Davidenko ODE IVP is first presented. Let $\mathcal{C} = \left\{\left(\mathbf{y} ,\ACP \right) : \left(\mathbf{y} ,\ACP \right) \mathrm{\, solves \,} \eqref{eqn_ode_tpbvp}   \right\}$ denote the solution manifold of \eqref{eqn_ode_tpbvp}. Suppose the solution manifold $\mathcal{C}$ is parameterized by arclength $\nu$, so that an element of $\mathcal{C}$ is $\left(\mathbf{y}(\nu) ,\ACP(\nu) \right)$, the tangent $\left(\mathbf{v}(\nu) ,\tau(\nu) \right)$ to $\mathcal{C}$ at $\left(\mathbf{y}(\nu) ,\ACP(\nu) \right)$ satisfies $\left \Vert \left(\mathbf{v}(\nu) ,\tau(\nu) \right) \right \Vert^2 = \int_a^b \mathbf{v}^\mathsf{T}(s,\nu)  \mathbf{v}(s,\nu)  \mathrm{d}s + \left[ \tau(\nu) \right]^2=1$ (i.e. $\left(\mathbf{v}(\nu) ,\tau(\nu) \right)$ is a unit tangent), and the solution manifold $\mathcal{C}$ can be described as a solution curve. With this arclength parameterization, $\mathbf{y} \colon \left[a,b\right] \times \mathbb{R}  \to \mathbb{R}^n$, $\ACP \colon \mathbb{R} \to \mathbb{R}$, $\mathbf{v} \colon \left[a,b\right] \times \mathbb{R}  \to \mathbb{R}^n$, $\tau \colon \mathbb{R} \to \mathbb{R}$, $\mathbf{y}(\nu)$ is shorthand for $\mathbf{y}(\cdot,\nu) \colon \left[a,b\right] \to \mathbb{R}^n$, and $\mathbf{v}(\nu)$ is shorthand for $\mathbf{v}(\cdot,\nu) \colon \left[a,b\right] \to \mathbb{R}^n$. Note that the components of the unit tangent $\left(\mathbf{v}(\nu) ,\tau(\nu) \right)$ to $\mathcal{C}$ at $\left(\mathbf{y}(\nu) ,\ACP(\nu) \right)$ are given explicitly by $\mathbf{v}(s,\nu) = \pp{\mathbf{y}(s,\nu)}{\nu}$ and $\tau(\nu)=\dd{\lambda(\nu)}{\nu}$. 

The Fr\'echet derivative of the ODE TPBVP \eqref{eqn_ode_tpbvp} with respect to $\nu$ about the solution $\left(\mathbf{y}(\nu) ,\ACP(\nu) \right)$, in conjunction with the arclength constraint and the initial condition $\left(\mathbf{y}_{\mathrm I} ,\ACP_{\mathrm I} \right)$, gives the nonlinear ODE IVP in the independent arclength variable $\nu$:
\begin{equation} \label{eqn_frechet_ode_tpbvp}
\begin{split}
\dd{}{s} \mathbf{v}(s,\nu) &= \mathbf{F}_\mathbf{y} \left(s,\mathbf{y}(s,\nu),\ACP(\nu) \right) \mathbf{v}(s,\nu) + \mathbf{F}_{\ACP} \left(s,\mathbf{y}(s,\nu),\ACP(\nu) \right) \tau(\nu), \\
\mathbf{0}_{n \times 1} &= \mathbf{G}_{\mathbf{y}(a)} \left(\mathbf{y}(a,\nu),\mathbf{y}(b,\nu),\ACP(\nu) \right) \mathbf{v}(a,\nu)+\mathbf{G}_{\mathbf{y}(b)} \left(\mathbf{y}(a,\nu),\mathbf{y}(b,\nu),\ACP(\nu) \right) \mathbf{v}(b,\nu)\\&\hphantom{=}+\mathbf{G}_{\ACP} \left(\mathbf{y}(a,\nu),\mathbf{y}(b,\nu),\ACP(\nu) \right) \tau(\nu), \\
\left \Vert \left(\mathbf{v}(\nu) ,\tau(\nu) \right) \right \Vert^2 &=  
\left<\left(\mathbf{v}(\nu) ,\tau(\nu) \right) , \left(\mathbf{v}(\nu) , \tau(\nu) \right) \right> =
\int_a^b \mathbf{v}^\mathsf{T}(s,\nu)  \mathbf{v}(s,\nu)  \mathrm{d}s + \left[\tau(\nu)\right]^2=1, \\
\left(\mathbf{y}(\nu_0) ,\ACP(\nu_0) \right) &= \left(\mathbf{y}_{\mathrm I} ,\ACP_{\mathrm I} \right),
\end{split}
\end{equation}
which must be solved for $\left(\mathbf{y}(\nu) ,\ACP(\nu) \right)$ starting at $\nu_0$ from an initial solution $\left(\mathbf{y}_{\mathrm I} ,\ACP_{\mathrm I} \right)$ of \eqref{eqn_ode_tpbvp}. \eqref{eqn_frechet_ode_tpbvp} is called the Davidenko ODE IVP and its solution is called the Davidenko flow \cite{davidenko1953new,davidenko1953approximate,rall1968davidenko,deuflhard2011newton,trefethen2013numerical,boyd2014solving}. The first two equations in \eqref{eqn_frechet_ode_tpbvp} constitute the Fr\'echet derivative of the ODE TPBVP \eqref{eqn_ode_tpbvp}, the third equation is the arclength constraint, and the final equation is the initial condition. By introducing a dummy scalar-valued function $w$ to represent the integrand of the arclength constraint, \eqref{eqn_frechet_ode_tpbvp} can be re-written:
\begin{equation} \label{eqn_frechet_ode_tpbvp_rewrite}
\begin{split}
\dd{}{s} \mathbf{v}(s,\nu) &= \mathbf{F}_\mathbf{y} \left(s,\mathbf{y}(s,\nu),\ACP(\nu) \right) \mathbf{v}(s,\nu) + \mathbf{F}_{\ACP} \left(s,\mathbf{y}(s,\nu),\ACP(\nu) \right) \tau(\nu), \\
\dd{}{s} w(s,\nu) &= \mathbf{v}^\mathsf{T}(s,\nu)  \mathbf{v}(s,\nu), \\
\mathbf{0}_{n \times 1} &=
\mathbf{G}_{\mathbf{y}(a)} \left(\mathbf{y}(a,\nu),\mathbf{y}(b,\nu),\ACP(\nu) \right) \mathbf{v}(a,\nu)+\mathbf{G}_{\mathbf{y}(b)} \left(\mathbf{y}(a,\nu),\mathbf{y}(b,\nu),\ACP(\nu) \right) \mathbf{v}(b,\nu)\\&\hphantom{=}+\mathbf{G}_{\ACP} \left(\mathbf{y}(a,\nu),\mathbf{y}(b,\nu),\ACP(\nu) \right) \tau(\nu), \\
w(a,\nu) &= 0, \\
w(b,\nu)+\left[\tau(\nu)\right]^2-1 &=0,\\
\left(\mathbf{y}(\nu_0) ,\ACP(\nu_0) \right) &= \left(\mathbf{y}_{\mathrm I} ,\ACP_{\mathrm I} \right).
\end{split}
\end{equation}
Again, letting $\nu$ vary, \eqref{eqn_frechet_ode_tpbvp_rewrite} is a nonlinear ODE IVP which must be solved for $\left(\mathbf{y}(\nu) ,\ACP(\nu) \right)$ (i.e. $\mathbf{y} \colon \left[a,b\right] \times \mathbb{R}  \to \mathbb{R}^n$ and $\ACP \colon \mathbb{R} \to \mathbb{R}$) starting at $\nu_0$ from an initial solution $\left(\mathbf{y}_{\mathrm I} ,\ACP_{\mathrm I} \right)$ of \eqref{eqn_ode_tpbvp}.
However, for a fixed $\nu$, \eqref{eqn_frechet_ode_tpbvp_rewrite} is a nonlinear ODE TPBVP which must be solved for $\mathbf{v}(\cdot,\nu) \colon \left[a,b\right] \to \mathbb{R}^n$, $\tau(\nu) \in \mathbb{R}$, and $w(\cdot,\nu) \colon \left[a,b\right] \to \mathbb{R}$ and where the independent variable is $s \in \left[a,b\right]$. 

\begin{figure}[h]
	\centering
	\includegraphics[width=0.8\linewidth]{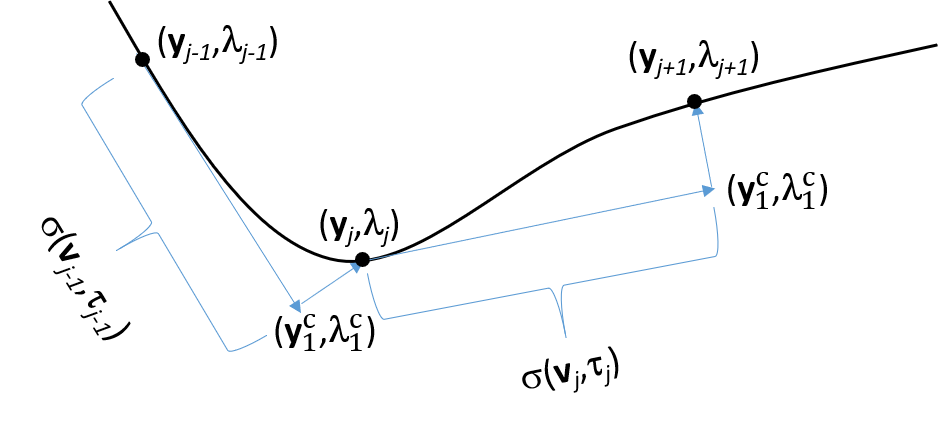}
	\caption{Predictor-corrector continuation.}
	\label{fig:pc_continuation}
\end{figure}

As explained in Chapter~5 of \cite{deuflhard2011newton}, it is inadvisable to integrate the Davidenko ODE IVP \eqref{eqn_frechet_ode_tpbvp}, or equivalently \eqref{eqn_frechet_ode_tpbvp_rewrite}. Instead, a predictor-corrector continuation method, depicted in Figure~\ref{fig:pc_continuation} and explained in detail in the following subappendices, is used to generate a solution sequence $\left\{\left(\mathbf{y}_j ,\ACP_j \right)\right\}_{j=1}^J$ which is a discrete subset of the Davidenko flow such that $\left(\mathbf{y}_1 ,\ACP_1 \right) = \left(\mathbf{y}_{\mathrm I} ,\ACP_{\mathrm I} \right)$. 

\subsection{Construct the Tangent}
Given a solution $\left(\mathbf{y}_j ,\ACP_j \right)$ to \eqref{eqn_ode_tpbvp} and a unit tangent $\left(\mathbf{v}_{j-1} ,\tau_{j-1} \right)$ to the previous solution $\left(\mathbf{y}_{j-1} ,\ACP_{j-1} \right)$ to \eqref{eqn_ode_tpbvp}, we seek to construct a tangent $\left(\mathbf{v}_j ,\tau_j \right)$ to the solution curve $\mathcal{C}$ at $\left(\mathbf{y}_j ,\ACP_j \right)$ which is roughly of unit length. The arclength constraint is
\begin{equation} \label{eqn_arclength_constraint}
\left \Vert \left(\mathbf{v}_j ,\tau_j \right) \right \Vert^2 = 
\left<\left(\mathbf{v}_j ,\tau_j \right) , \left(\mathbf{v}_j , \tau_j \right) \right>=
\int_a^b \mathbf{v}_j^\mathsf{T}(s)  \mathbf{v}_j(s)  \mathrm{d}s + \tau_j^2=1,
\end{equation}
which is nonlinear in the tangent $\left(\mathbf{v}_j ,\tau_j \right)$. An alternative constraint, the pseudo-arclength constraint, is
\begin{equation} \label{eqn_pseudo_constraint}
\left<\left({\mathbf{v}}_{j-1} ,\tau_{j-1} \right) , \left(\mathbf{v}_j , \tau_j \right) \right>=
\int_a^b {\mathbf{v}}_{j-1}^\mathsf{T}(s)  \mathbf{v}_j(s)  \mathrm{d}s + \tau_{j-1} \tau_j=1,
\end{equation}
which, in contrast to the arclength constraint \eqref{eqn_arclength_constraint}, is linear in the tangent $\left(\mathbf{v}_j ,\tau_j \right)$.
The linearization (i.e. Fr\'echet derivative) of the ODE TPBVP \eqref{eqn_ode_tpbvp} about the solution $\left(\mathbf{y}_j ,\ACP_j \right)$, in conjunction with the pseudo-arclength condition \eqref{eqn_pseudo_constraint}, gives the linear ODE TPBVP:
\begin{equation} \label{eqn_tan_ode_tpbvp}
\begin{split}
\dd{}{s} \mathbf{v}_j(s) &= \mathbf{F}_\mathbf{y} \left(s,\mathbf{y}_j(s),\ACP_j \right) \mathbf{v}_j(s) \\ &\hphantom{=} + \mathbf{F}_{\ACP} \left(s,\mathbf{y}_j(s),\ACP_j \right) \tau_j,  \\
\dd{}{s} \tau_j &= 0, \\
\dd{}{s} w(s) &= {\mathbf{v}}_{j-1}^\mathsf{T}(s) \mathbf{v}_j(s), \\
\mathbf{G}_{\mathbf{y}(a)} \left(\mathbf{y}_j(a),\mathbf{y}_j(b),\ACP_j \right) \mathbf{v}_j(a)+\mathbf{G}_{\mathbf{y}(b)} \left(\mathbf{y}_j(a),\mathbf{y}_j(b),\ACP_j \right) \mathbf{v}_j(b)\\+\mathbf{G}_{\ACP} \left(\mathbf{y}_j(a),\mathbf{y}_j(b),\ACP_j \right) \tau_j &= \mathbf{0}_{n \times 1}, \\
w(a) &= 0, \\
w(b) + \tau_{j-1} \tau_j - 1 &= 0,
\end{split}
\end{equation}
which must be solved for $\mathbf{v}_j \colon \left[a,b\right] \to \mathbb{R}^n$, $\tau_j \in \mathbb{R}$, and $w \colon \left[a,b\right] \to \mathbb{R}$ and where $\left(\mathbf{v}_j ,\tau_j \right)$ is a tangent to $\mathcal{C}$ at $\left(\mathbf{y}_j ,\ACP_j \right)$. Note that the first, second, and third equations in \eqref{eqn_tan_ode_tpbvp} are the ODEs, while the fourth, fifth, and sixth equations constitute the boundary conditions. The first, second, and fourth equations in \eqref{eqn_tan_ode_tpbvp} are the linearization (i.e. Fr\'echet derivative) of \eqref{eqn_ode_tpbvp} about the solution $\left(\mathbf{y}_j ,\ACP_j \right)$ and ensure that a tangent is produced, while the third, fifth, and sixth equations in \eqref{eqn_tan_ode_tpbvp} enforce the pseudo-arclength condition \eqref{eqn_pseudo_constraint}. The initial solution guess to solve \eqref{eqn_tan_ode_tpbvp} is $\left(\mathbf{v}_j ,\tau_j \right) = \left( \mathbf{v}_{j-1} ,\tau_{j-1} \right)$ and $w(s) = \int_a^s \mathbf{v}_{j-1}^\mathsf{T}(\tilde s)  \mathbf{v}_{j-1}(\tilde s)  \mathrm{d} \tilde s$, $s \in \left[a,b\right]$, for $j \ge 1$. For $j = 1$, define $\left(\mathbf{v}_0 ,\tau_0 \right) = \left(\mathbf{0}_{n \times 1} ,1 \right)$. Note that the construction of the initial guess for $w$ can be realized efficiently via the \mcode{MATLAB} routine \mcode{cumtrapz}. 

Note that the linear ODE TPBVP \eqref{eqn_tan_ode_tpbvp} can be solved numerically via the \mcode{MATLAB} routines \mcode{sbvp} or \mcode{bvptwp}, which offers 4 algorithms: \mcode{twpbvp_m}, \mcode{twpbvpc_m}, \mcode{twpbvp_l}, and \mcode{twpbvpc_l}; moreover, \mcode{sbvp} and \mcode{bvptwp} have special algorithms to solve linear ODE TPBVPs. Since $\mathbf{y}_j$ and $\mathbf{v}_{j-1}$ are usually only known at a discrete set of points in $\left[a,b\right]$, the values of these functions at the other points in $\left[a,b\right]$ must be obtained through interpolation in order to numerically solve \eqref{eqn_tan_ode_tpbvp}. The \mcode{MATLAB} routine \mcode{interp1} performs linear, cubic, pchip, makima, and spline interpolation and may be utilized to interpolate $\mathbf{y}_j$ and $\mathbf{v}_{j-1}$ while solving \eqref{eqn_tan_ode_tpbvp}.

Because the numerical solvers usually converge faster when provided Jacobians of the ODE velocity function and of the two-point boundary condition function, these are computed below.
Let 
\begin{equation}
\mathbf{x} = \begin{bmatrix} \mathbf{v}_{j}  \\ \tau_j  \\ w \end{bmatrix}.
\end{equation}
The ODE velocity function in \eqref{eqn_tan_ode_tpbvp} is
\begin{equation}
\mathbf{H}^{\mathrm t} \left(s, \mathbf{x}(s) \right) = \mathbf{H}^{\mathrm t} \left(s,\mathbf{v}_j(s),\tau_j,w(s) \right) = \begin{bmatrix} \mathbf{F}_\mathbf{y} \left(s,\mathbf{y}_j(s),\ACP_j \right) \mathbf{v}_j(s) + \mathbf{F}_{\ACP} \left(s,\mathbf{y}_j(s),\ACP_j \right) \tau_j \\ 0 \\ {\mathbf{v}}_{j-1}^\mathsf{T}(s) \mathbf{v}_j(s) \end{bmatrix}.
\end{equation}
The Jacobian of the ODE velocity function $\mathbf{H}^{\mathrm t}$ with respect to $\mathbf{x}$ is
\begin{equation}
\begin{split}
\mathbf{H}_\mathbf{x}^{\mathrm t} \left(s,\mathbf{x}(s) \right) &=
\mathbf{H}_\mathbf{x}^{\mathrm t} \left(s,\mathbf{v}_j(s),\tau_j,w(s) \right) \\ &= \begin{bmatrix} \mathbf{F}_{\mathbf{y}} \left(s,\mathbf{y}_j(s),\ACP_j \right) & \mathbf{F}_{\ACP} \left(s,\mathbf{y}_j(s),\ACP_j \right) & \mathbf{0}_{n \times 1} \\ \mathbf{0}_{1 \times n} & 0 & 0 \\ {\mathbf{v}}_{j-1}^\mathsf{T}(s) & 0 & 0 \end{bmatrix}.
\end{split}
\end{equation}
The two-point boundary condition in \eqref{eqn_tan_ode_tpbvp} is
\begin{equation}
\mathbf{K}^{\mathrm t} \left(\mathbf{x}(a),\mathbf{x}(b) \right) = \mathbf{0}_{\left(n+2 \right) \times 1},
\end{equation}
where $\mathbf{K}^{\mathrm t}$ is the two-point boundary condition function
\begin{equation}
\begin{split}
\mathbf{K}^{\mathrm t}& \left(\mathbf{x}(a),\mathbf{x}(b) \right) =\\& \begin{bmatrix} \mathbf{G}_{\mathbf{y}(a)} \left(\mathbf{y}_j(a),\mathbf{y}_j(b),\ACP_j \right) \mathbf{v}_j(a)+\mathbf{G}_{\mathbf{y}(b)} \left(\mathbf{y}_j(a),\mathbf{y}_j(b),\ACP_j \right) \mathbf{v}_j(b) +\mathbf{G}_{\ACP} \left(\mathbf{y}_j(a),\mathbf{y}_j(b),\ACP_j \right) \tau_j \\ w(a) \\  w(b) + \tau_{j-1} \tau_j - 1 \end{bmatrix}.
\end{split}
\end{equation}
The Jacobians of the two-point boundary condition function $\mathbf{K}^{\mathrm t}$ with respect to $\mathbf{x}(a)$ and $\mathbf{x}(b)$ are
\begin{equation} \label{eq_K_t_xa}
\mathbf{K}_{\mathbf{x}(a)}^{\mathrm t} \left(\mathbf{x}(a),\mathbf{x}(b) \right) = \begin{bmatrix} \mathbf{G}_{\mathbf{y}(a)} \left(\mathbf{y}_j(a),\mathbf{y}_j(b),\ACP_j \right) & \mathbf{G}_{\ACP} \left(\mathbf{y}_j(a),\mathbf{y}_j(b),\ACP_j \right) & \mathbf{0}_{n \times 1} \\ \mathbf{0}_{1 \times n} & 0 & 1 \\ \mathbf{0}_{1 \times n} & \tau_{j-1}  & 0 \end{bmatrix}
\end{equation}
and
\begin{equation} \label{eq_K_t_xb}
\mathbf{K}_{\mathbf{x}(b)}^{\mathrm t} \left(\mathbf{x}(a),\mathbf{x}(b) \right) = \begin{bmatrix} \mathbf{G}_{\mathbf{y}(b)} \left(\mathbf{y}_j(a),\mathbf{y}_j(b),\ACP_j \right) & \mathbf{G}_{\ACP} \left(\mathbf{y}_j(a),\mathbf{y}_j(b),\ACP_j \right) & \mathbf{0}_{n \times 1} \\ \mathbf{0}_{1 \times n} & 0 & 0 \\ \mathbf{0}_{1 \times n} & \tau_{j-1}  & 1\end{bmatrix}.
\end{equation}
Special care must be taken when implementing the Jacobians \eqref{eq_K_t_xa} and \eqref{eq_K_t_xb}. Since the unknown constant $\tau_j$ appears as the second to last element of both $\mathbf{x}(a)$ and $\mathbf{x}(b)$, $\tau_j$ from only one of $\mathbf{x}(a)$ and $\mathbf{x}(b)$ is actually used to construct each term in $\mathbf{K}^{\mathrm t}$ involving $\tau_j$. The middle column of \eqref{eq_K_t_xa} is actually the derivative of $\mathbf{K}^{\mathrm t}$ with respect to the $\tau_j$ in  $\mathbf{x}(a)$, while the middle column of \eqref{eq_K_t_xb} is actually the derivative of $\mathbf{K}^{\mathrm t}$ with respect to the $\tau_j$ in  $\mathbf{x}(b)$. Thus, the middle columns in \eqref{eq_K_t_xa} and \eqref{eq_K_t_xb} corresponding to the derivative of $\mathbf{K}^{\mathrm t}$ with respect to $\tau_j$ should not coincide in a software implementation. For example, if $\mathbf{K}^{\mathrm t}$ is constructed from the $\tau_j$ in $\mathbf{x}(a)$, $\mathbf{K}_{\mathbf{x}(a)}^{\mathrm t}$ is as shown in \eqref{eq_K_t_xa} while the middle column of \eqref{eq_K_t_xb} corresponding to the derivative of $\mathbf{K}^{\mathrm t}$ with respect to the $\tau_j$ in $\mathbf{x}(b)$ is all zeros. Alternatively, if $\mathbf{K}^{\mathrm t}$ is constructed from the $\tau_j$ in $\mathbf{x}(b)$, $\mathbf{K}_{\mathbf{x}(b)}^{\mathrm t}$ is as shown in \eqref{eq_K_t_xb} while the middle column of \eqref{eq_K_t_xa} corresponding to the derivative of $\mathbf{K}^{\mathrm t}$ with respect to the $\tau_j$ appearing in  $\mathbf{x}(a)$ is all zeros.

\subsection{Normalize the Tangent}
The tangent $\left(\mathbf{v}_j ,\tau_j \right)$ at $\left(\mathbf{y}_j ,\ACP_j \right)$ obtained by solving \eqref{eqn_tan_ode_tpbvp} in the previous step is only roughly of unit length. A unit tangent at $\left(\mathbf{y}_j ,\ACP_j \right)$ is obtained from $\left(\mathbf{v}_j ,\tau_j \right)$ through normalization:
\begin{equation} \label{eqn_norm_tan_pred} 
\left(\mathbf{v}_j ,\tau_j \right) \gets \frac{1}{\kappa} \left(\mathbf{v}_j ,\tau_j \right),
\end{equation}
where 
\begin{equation} \label{eqn_norm_tan_pred_kappa} 
\kappa = \left \Vert \left(\mathbf{v}_j ,\tau_j \right) \right \Vert = 
\left<\left(\mathbf{v}_j ,\tau_j \right) , \left(\mathbf{v}_j , \tau_j \right) \right>^\frac{1}{2}=
\left[ \int_a^b \mathbf{v}_j^\mathsf{T}(s)  \mathbf{v}_j(s)  \mathrm{d}s + \tau_j^2\right]^\frac{1}{2}.
\end{equation}
The integration operator to construct the normalization scalar $\kappa$ in \eqref{eqn_norm_tan_pred_kappa} can be realized via the \mcode{MATLAB} routine \mcode{trapz}.

\subsection{Construct the Tangent Predictor}
The unit tangent $\left(\mathbf{v}_j ,\tau_j \right)$ constructed in \eqref{eqn_norm_tan_pred} is used to obtain a guess (the so-called ``tangent predictor") $\left(\mathbf{y}_1^{\mathrm c} ,\ACP_1^{\mathrm c} \right)$ for the next solution $\left(\mathbf{y}_{j+1} ,\ACP_{j+1} \right)$ as follows:
\begin{equation} \label{eqn_intermediate_value}
\left(\mathbf{y}_1^{\mathrm c} ,\ACP_1^{\mathrm c} \right) =  \left(\mathbf{y}_j ,\ACP_j \right) + \sigma \left(\mathbf{v}_j ,\tau_j \right),
\end{equation}
where $\sigma \in \left[ \sigma_{\mathrm{min}} , \sigma_{\mathrm{max}} \right]$ is a steplength and where $0 < \sigma_{\mathrm{min}} \le \sigma_{\mathrm{max}}$. Concretely, $\sigma_{\mathrm{min}}$ might be $.0001$ and $\sigma_{\mathrm{max}}$ might be $\frac{1}{2}$. $\sigma$ is adapted during the predictor-corrector continuation method based on the corrector step, discussed in the next subappendix. Initially, the value of $\sigma$ is set to $\sigma_{\mathrm{init}} \in \left[ \sigma_{\mathrm{min}} , \sigma_{\mathrm{max}} \right]$. The notation $\left(\mathbf{y}_1^{\mathrm c} ,\ACP_1^{\mathrm c} \right)$ is used to denote the tangent predictor in \eqref{eqn_intermediate_value} because, as discussed in the next subappendix, the tangent predictor is used as the initial corrector in an iterative Newton's method that projects the tangent predictor onto $\mathcal{C}$.

\subsection{Construct the Corrector}
Since the tangent predictor $\left(\mathbf{y}_1^{\mathrm c} ,\ACP_1^{\mathrm c} \right)$ constructed in \eqref{eqn_intermediate_value} does not necessarily lie on $\mathcal{C}$, $\left(\mathbf{y}_1^{\mathrm c} ,\ACP_1^{\mathrm c} \right)$ must be projected onto $\mathcal{C}$ to obtain the next solution (the so-called ``corrector") $\left(\mathbf{y}_{j+1} ,\ACP_{j+1} \right)$. This projection process is the corrector step. In order to perform the projection efficiently, the difference between the next solution and the tangent predictor, $\left(\mathbf{y}_{j+1} ,\ACP_{j+1} \right)-\left(\mathbf{y}_1^{\mathrm c} ,\ACP_1^{\mathrm c} \right)$, should be orthogonal to the unit tangent $\left( \mathbf{v}_j ,\tau_j \right)$. That is, the orthogonality constraint is
\begin{equation} \label{eqn_ortho_constraint}
\begin{split}
\left<\left({\mathbf{v}}_j ,  \tau_j \right) , \left(\mathbf{y}_{j+1} , \ACP_{j+1} \right) - \left(\mathbf{y}_1^{\mathrm c} ,\ACP_1^{\mathrm c} \right) \right> &=
\left<\left({\mathbf{v}}_j ,  \tau_j \right) , \left(\mathbf{y}_{j+1} - \mathbf{y}_1^{\mathrm c} , \ACP_{j+1}-\ACP_1^{\mathrm c} \right) \right> \\  &=
\int_a^b {\mathbf{v}}_j^\mathsf{T}(s) \left[\mathbf{y}_{j+1}(s) - \mathbf{y}_1^{\mathrm c}(s)  \right]  \mathrm{d}s +  \tau_{j} \left[\ACP_{j+1}-\ACP_1^{\mathrm c}  \right]=0.
\end{split}
\end{equation}
The tangent predictor $\left(\mathbf{y}_1^{\mathrm c} ,\ACP_1^{\mathrm c} \right)$ can be iteratively corrected by applying Newton's method to \eqref{eqn_ode_tpbvp}, while enforcing the orthogonality constraint \eqref{eqn_ortho_constraint}, to generate a sequence of correctors $\left\{ \left(\mathbf{y}_k^{\mathrm c} ,\ACP_k^{\mathrm c} \right) \right\}_{k=1}^ {K+1} $. Applying Newton's method to the ODE TPBVP \eqref{eqn_ode_tpbvp} about the current corrector $\left(\mathbf{y}_k^{\mathrm c} ,\ACP_k^{\mathrm c} \right)$, in conjunction with the orthogonality constraint \eqref{eqn_ortho_constraint}, gives the linear ODE TPBVP:
\begin{equation} \label{eqn_corrn_ode_tpbvp}
\begin{split}
\dd{}{s} \delta \mathbf{y}_k^{\mathrm c}(s) &= \mathbf{F}_\mathbf{y} \left(s,\mathbf{y}_k^{\mathrm c}(s),\ACP_k^{\mathrm c} \right) \delta \mathbf{y}_k^{\mathrm c}(s) \\&\hphantom{=} + \mathbf{F}_{\ACP} \left(s,\mathbf{y}_k^{\mathrm c}(s),\ACP_k^{\mathrm c} \right)  \delta \ACP_k^{\mathrm c} \\
&\hphantom{=} -\dd{}{s} \mathbf{y}_k^{\mathrm c}(s) + \mathbf{F} \left(s,\mathbf{y}_k^{\mathrm c}(s),\ACP_k^{\mathrm c} \right), \\
\dd{}{s} \delta \ACP_k^{\mathrm c} &= 0, \\
\dd{}{s} w(s) &= {\mathbf{v}}_j^\mathsf{T}(s) \delta \mathbf{y}_k^{\mathrm c}(s),  \\
\mathbf{G}_{\mathbf{y}(a)}  \left( \mathbf{y}_k^{\mathrm c} (a),\mathbf{y}_k^{\mathrm c}(b),\ACP_k^{\mathrm c} \right)  \delta \mathbf{y}_k^{\mathrm c}(a)+\mathbf{G}_{\mathbf{y}(b)}  \left( \mathbf{y}_k^{\mathrm c} (a),\mathbf{y}_k^{\mathrm c}(b),\ACP_k^{\mathrm c} \right)  \delta \mathbf{y}_k^{\mathrm c}(b)\\+\mathbf{G}_{\ACP} \left( \mathbf{y}_k^{\mathrm c} (a),\mathbf{y}_k^{\mathrm c}(b),\ACP_k^{\mathrm c} \right)  \delta \ACP_k^{\mathrm c}  +
\mathbf{G} \left( \mathbf{y}_k^{\mathrm c} (a),\mathbf{y}_k^{\mathrm c}(b),\ACP_k^{\mathrm c} \right) &= \mathbf{0}_{n \times 1}, \\
w(a) &= 0, \\
w(b) + \tau_{j} \delta \ACP_k^{\mathrm c} &= 0,
\end{split}
\end{equation}
which must be solved for $\delta \mathbf{y}_k^{\mathrm c} \colon \left[a,b\right] \to \mathbb{R}^n$,  $\delta \ACP_k^{\mathrm c} \in \mathbb{R}$, and $w \colon \left[a,b\right] \to \mathbb{R}$ and where $\left( \delta \mathbf{y}_k^{\mathrm c}, \delta \ACP_k^{\mathrm c} \right)$ represents a correction to the current corrector $\left( \mathbf{y}_k^{\mathrm c},\ACP_k^{\mathrm c} \right)$. Note that the first, second, and third equations in \eqref{eqn_corrn_ode_tpbvp} are the ODEs, while the fourth, fifth, and sixth equations constitute the boundary conditions. The first, second, and fourth equations in \eqref{eqn_corrn_ode_tpbvp} are the result of applying Newton's method to \eqref{eqn_ode_tpbvp} about the current corrector $\left( \mathbf{y}_k^{\mathrm c},\ACP_k^{\mathrm c} \right)$, while the third, fifth, and sixth equations in \eqref{eqn_corrn_ode_tpbvp} enforce the orthogonality constraint \eqref{eqn_ortho_constraint}. \eqref{eqn_corrn_ode_tpbvp} must be solved iteratively for at most $K$ iterations, so that $1 \le k \le K$. The initial solution guess to solve \eqref{eqn_corrn_ode_tpbvp} at the beginning of each iteration is $\left( \delta \mathbf{y}_k^{\mathrm c}, \delta \ACP_k^{\mathrm c} \right)=\left(\mathbf{0}_{n \times 1} ,0 \right)$ and $w(s) = 0$, $s \in \left[a,b\right]$. The initial corrector about which Newton's method is applied in the first iteration is the tangent predictor $\left(\mathbf{y}_1^{\mathrm c} ,\ACP_1^{\mathrm c} \right)$. At the end of each iteration, the corrector about which Newton's method is applied for the next iteration is updated via $\left( \mathbf{y}_{k+1}^{\mathrm c},\ACP_{k+1}^{\mathrm c} \right)=\left( \mathbf{y}_k^{\mathrm c},\ACP_k^{\mathrm c} \right)+\left( \delta \mathbf{y}_k^{\mathrm c}, \delta \ACP_k^{\mathrm c} \right)$. At the end of each iteration, convergence to $\mathcal{C}$ should be tested via:
\begin{equation} \label{eqn_conv_test}
\frac { \left \Vert \left( \delta \mathbf{y}_k^{\mathrm c} , \delta \ACP_k^{\mathrm c} \right) \right \Vert } { \left \Vert \left(\mathbf{y}_1^{\mathrm c} ,\ACP_1^{\mathrm c} \right) \right \Vert } =
\frac{\left[ \int_a^b \left[\delta \mathbf{y}_k^{\mathrm c}(s)\right]^\mathsf{T}  \delta \mathbf{y}_k^{\mathrm c}(s)  \mathrm{d}s + \left[  \delta \ACP_k^{\mathrm c} \right]^2 \right]^\frac{1}{2}}{\left[ \int_a^b \left[\mathbf{y}_1^{\mathrm c}(s) \right]^\mathsf{T}  \mathbf{y}_1^{\mathrm c}(s)  \mathrm{d}s + \left[\ACP_1^{\mathrm c} \right]^2 \right]^\frac{1}{2}}< \gamma,
\end{equation}
where $\gamma$ is a small threshold such as $.001$. Since Newton's method enjoys quadratic convergence near a solution, only a few (say  $K=5$) iterative solves of \eqref{eqn_corrn_ode_tpbvp} should be attempted. If convergence has not been attained in $K$ iterations, the steplength $\sigma$ should be reduced: 
\begin{equation} \label{eqn_sigma_reduce}
\sigma \gets \sigma_{\mathrm r} \sigma,
\end{equation}
where $\sigma_{\mathrm r}$ is a reduction scale factor such as $\frac{1}{4}$ and the corrector step should be restarted at the new tangent predictor $\left(\mathbf{y}_1^{\mathrm c} ,\ACP_1^{\mathrm c} \right)= \left(\mathbf{y}_j ,\ACP_j \right) + \sigma \left(\mathbf{v}_j ,\tau_j \right)$, based on the updated value of $\sigma$ realized in \eqref{eqn_sigma_reduce}. If, as a result of the reduction realized in \eqref{eqn_sigma_reduce}, $\sigma < \sigma_{\mathrm{min}}$, the algorithm should halt and predictor-corrector continuation failed. However, if convergence has been achieved in $k \le K$ iterations, the next solution can be taken to be $\left(\mathbf{y}_{j+1} ,\ACP_{j+1} \right) = \left( \mathbf{y}_{k+1}^{\mathrm c},\ACP_{k+1}^{\mathrm c} \right)$ or the corrector can be further polished as explained in the next subappendix. Moreover, if convergence has been achieved rapidly in no more than $k_{\mathrm{fast}}$ iterations, where $1 \le  k_{\mathrm{fast}} \le K$ and, concretely, $k_{\mathrm{fast}}$ might be 3, then the steplength $\sigma$ may be increased:
\begin{equation}
\sigma \gets \min \left\{ \sigma_{\mathrm i} \sigma, \sigma_{\mathrm{max}} \right\},
\end{equation}
where $\sigma_{\mathrm i}$ is an increase scale factor such as $2$.

Note that the linear ODE TPBVP \eqref{eqn_corrn_ode_tpbvp} can be solved numerically via the \mcode{MATLAB} routines \mcode{sbvp} or \mcode{bvptwp}, which offers 4 algorithms: \mcode{twpbvp_m}, \mcode{twpbvpc_m}, \mcode{twpbvp_l}, and \mcode{twpbvpc_l}; moreover, \mcode{sbvp} and \mcode{bvptwp} have special algorithms to solve linear ODE TPBVPs. Since $\mathbf{y}_k^{\mathrm c}$, $\dd{}{s} \mathbf{y}_k^{\mathrm c}$, and ${\mathbf{v}}_j$ are usually only known at a discrete set of points in $\left[a,b\right]$, the values of these functions at the other points in $\left[a,b\right]$ must be obtained through interpolation in order to numerically solve \eqref{eqn_corrn_ode_tpbvp}. The \mcode{MATLAB} routine \mcode{interp1} performs linear, cubic, pchip, makima, and spline interpolation and may be utilized to interpolate $\mathbf{y}_k^{\mathrm c}$, $\dd{}{s} \mathbf{y}_k^{\mathrm c}$, and ${\mathbf{v}}_j$ while solving \eqref{eqn_corrn_ode_tpbvp}.

Because the numerical solvers usually converge faster when provided Jacobians of the ODE velocity function and of the two-point boundary condition function, these are computed below.
Let 
\begin{equation}
\mathbf{x} = \begin{bmatrix} \delta \mathbf{y}_k^{\mathrm c}  \\ \delta \ACP_k^{\mathrm c}  \\ w \end{bmatrix}.
\end{equation}
The ODE velocity function in \eqref{eqn_corrn_ode_tpbvp} is
\begin{equation}
\begin{split}
\mathbf{H}^{\mathrm c} \left(s, \mathbf{x}(s) \right) &= \mathbf{H}^{\mathrm c} \left(s,\delta \mathbf{y}_k^{\mathrm c}(s),\delta \ACP_k^{\mathrm c},w(s) \right) \\
&= \begin{bmatrix}  \mathbf{F}_\mathbf{y} \left(s,\mathbf{y}_k^{\mathrm c}(s),\ACP_k^{\mathrm c} \right) \delta \mathbf{y}_k^{\mathrm c}(s) + \mathbf{F}_{\ACP} \left(s,\mathbf{y}_k^{\mathrm c}(s),\ACP_k^{\mathrm c} \right)  \delta \ACP_k^{\mathrm c} -\dd{}{s} \mathbf{y}_k^{\mathrm c}(s) + \mathbf{F} \left(s,\mathbf{y}_k^{\mathrm c}(s),\ACP_k^{\mathrm c} \right) \\ 0 \\ {\mathbf{v}}_j^\mathsf{T}(s) \delta \mathbf{y}_k^{\mathrm c}(s)  \end{bmatrix}.
\end{split}
\end{equation}
The Jacobian of the ODE velocity function $\mathbf{H}^{\mathrm c}$ with respect to $\mathbf{x}$ is
\begin{equation}
\begin{split}
\mathbf{H}_\mathbf{x}^{\mathrm c} \left(s,\mathbf{x}(s) \right) &=
\mathbf{H}_\mathbf{x}^{\mathrm c} \left(s,\delta \mathbf{y}_k^{\mathrm c}(s),\delta \ACP_k^{\mathrm c},w(s) \right) \\ &= \begin{bmatrix} \mathbf{F}_{\mathbf{y}} \left(s,\mathbf{y}_k^{\mathrm c}(s),\ACP_k^{\mathrm c} \right) & \mathbf{F}_{\ACP} \left(s,\mathbf{y}_k^{\mathrm c}(s),\ACP_k^{\mathrm c} \right) & \mathbf{0}_{n \times 1} \\ \mathbf{0}_{1 \times n} & 0 & 0 \\ { \mathbf{v}}_j^\mathsf{T}(s) & 0 & 0 \end{bmatrix}.
\end{split}
\end{equation}
The two-point boundary condition in \eqref{eqn_corrn_ode_tpbvp} is
\begin{equation}
\mathbf{K}^{\mathrm c} \left(\mathbf{x}(a),\mathbf{x}(b) \right) = \mathbf{0}_{\left(n+2 \right) \times 1},
\end{equation}
where $\mathbf{K}^{\mathrm c}$ is the two-point boundary condition function
\begin{equation}
\mathbf{K}^{\mathrm c} \left(\mathbf{x}(a),\mathbf{x}(b) \right) = \begin{bmatrix} \mathbf{G}_{\mathbf{y}(a)}  \left( \mathbf{y}_k^{\mathrm c} (a),\mathbf{y}_k^{\mathrm c}(b),\ACP_k^{\mathrm c} \right)  \delta \mathbf{y}_k^{\mathrm c}(a)+\mathbf{G}_{\mathbf{y}(b)}  \left( \mathbf{y}_k^{\mathrm c} (a),\mathbf{y}_k^{\mathrm c}(b),\ACP_k^{\mathrm c} \right)  \delta \mathbf{y}_k^{\mathrm c}(b)\\+\mathbf{G}_{\ACP} \left( \mathbf{y}_k^{\mathrm c} (a),\mathbf{y}_k^{\mathrm c}(b),\ACP_k^{\mathrm c} \right)  \delta \ACP_k^{\mathrm c}  +
\mathbf{G} \left( \mathbf{y}_k^{\mathrm c} (a),\mathbf{y}_k^{\mathrm c}(b),\ACP_k^{\mathrm c} \right) \\ w(a) \\  w(b) + \tau_{j} \delta \ACP_k^{\mathrm c} \end{bmatrix}.
\end{equation}
The Jacobians of the two-point boundary condition function $\mathbf{K}^{\mathrm c}$ with respect to $\mathbf{x}(a)$ and $\mathbf{x}(b)$ are
\begin{equation} \label{eq_K_c_xa}
\mathbf{K}_{\mathbf{x}(a)}^{\mathrm c} \left(\mathbf{x}(a),\mathbf{x}(b) \right) = \begin{bmatrix} \mathbf{G}_{\mathbf{y}(a)}\left( \mathbf{y}_k^{\mathrm c} (a),\mathbf{y}_k^{\mathrm c}(b),\ACP_k^{\mathrm c} \right) & \mathbf{G}_{\ACP}\left( \mathbf{y}_k^{\mathrm c} (a),\mathbf{y}_k^{\mathrm c}(b),\ACP_k^{\mathrm c} \right) & \mathbf{0}_{n \times 1} \\ \mathbf{0}_{1 \times n} & 0 & 1 \\ \mathbf{0}_{1 \times n} & \tau_{j}  & 0 \end{bmatrix}
\end{equation}
and
\begin{equation} \label{eq_K_c_xb}
\mathbf{K}_{\mathbf{x}(b)}^{\mathrm c} \left(\mathbf{x}(a),\mathbf{x}(b) \right) = \begin{bmatrix} \mathbf{G}_{\mathbf{y}(b)}\left( \mathbf{y}_k^{\mathrm c} (a),\mathbf{y}_k^{\mathrm c}(b),\ACP_k^{\mathrm c} \right) & \mathbf{G}_{\ACP}\left( \mathbf{y}_k^{\mathrm c} (a),\mathbf{y}_k^{\mathrm c}(b),\ACP_k^{\mathrm c} \right) & \mathbf{0}_{n \times 1} \\ \mathbf{0}_{1 \times n} & 0 & 0 \\ \mathbf{0}_{1 \times n} & \tau_{j}  & 1\end{bmatrix}.
\end{equation}
Special care must be taken when implementing the Jacobians \eqref{eq_K_c_xa} and \eqref{eq_K_c_xb}. Since the unknown constant $\delta \ACP_k^{\mathrm c}$ appears as the second to last element of both $\mathbf{x}(a)$ and $\mathbf{x}(b)$, $\delta \ACP_k^{\mathrm c}$ from only one of $\mathbf{x}(a)$ and $\mathbf{x}(b)$ is actually used to construct each term in $\mathbf{K}^{\mathrm c}$ involving $\delta \ACP_k^{\mathrm c}$. The middle column of \eqref{eq_K_c_xa} is actually the derivative of $\mathbf{K}^{\mathrm c}$ with respect to the $\delta \ACP_k^{\mathrm c}$ in  $\mathbf{x}(a)$, while the middle column of \eqref{eq_K_c_xb} is actually the derivative of $\mathbf{K}^{\mathrm c}$ with respect to the $\delta \ACP_k^{\mathrm c}$ in  $\mathbf{x}(b)$. Thus, the middle columns in \eqref{eq_K_c_xa} and \eqref{eq_K_c_xb} corresponding to the derivative of $\mathbf{K}^{\mathrm c}$ with respect to $\delta \ACP_k^{\mathrm c}$ should not coincide in a software implementation. For example, if $\mathbf{K}^{\mathrm c}$ is constructed from the $\delta \ACP_k^{\mathrm c}$ in $\mathbf{x}(a)$, $\mathbf{K}_{\mathbf{x}(a)}^{\mathrm c}$ is as shown in \eqref{eq_K_c_xa} while the middle column of \eqref{eq_K_c_xb} corresponding to the derivative of $\mathbf{K}^{\mathrm c}$ with respect to the $\delta \ACP_k^{\mathrm c}$ in $\mathbf{x}(b)$ is all zeros. Alternatively, if $\mathbf{K}^{\mathrm c}$ is constructed from the $\delta \ACP_k^{\mathrm c}$ in $\mathbf{x}(b)$, $\mathbf{K}_{\mathbf{x}(b)}^{\mathrm c}$ is as shown in \eqref{eq_K_c_xb} while the middle column of \eqref{eq_K_c_xa} corresponding to the derivative of $\mathbf{K}^{\mathrm c}$ with respect to the $\delta \ACP_k^{\mathrm c}$ appearing in  $\mathbf{x}(a)$ is all zeros.

\subsection{Polish the Corrector}
The final corrector $\left( \mathbf{y}_{k+1}^{\mathrm c},\ACP_{k+1}^{\mathrm c} \right)$ from the previous step can be further polished by finding $\left(\mathbf{y}_{j+1} ,\ACP_{j+1} \right)$ that solves \eqref{eqn_ode_tpbvp} while satisfying the orthogonality constraint \eqref{eqn_ortho_constraint}. This yields the ODE TPBVP:
\begin{equation} \label{eqn_corrp_ode_tpbvp}
\begin{split}
\dd{}{s} \mathbf{y}_{j+1}(s) &= \mathbf{F} \left(s,\mathbf{y}_{j+1}(s),\ACP_{j+1} \right), \\
\dd{}{s} \ACP_{j+1} &= 0, \\
\dd{}{s} w(s) &= {\mathbf{v}}_j^\mathsf{T}(s) \left[\mathbf{y}_{j+1}(s) - \mathbf{y}_1^{\mathrm c}(s)  \right],  \\
\mathbf{G}\left( \mathbf{y}_{j+1} (a),\mathbf{y}_{j+1}(b),\ACP_{j+1} \right) &= \mathbf{0}_{n \times 1}, \\
w(a) &= 0, \\
w(b) + \tau_j \left[\ACP_{j+1}-\ACP_1^{\mathrm c}  \right] &= 0,
\end{split}
\end{equation}
which must be solved for $\mathbf{y}_{j+1} \colon \left[a,b\right] \to \mathbb{R}^n$,  $\ACP_{j+1} \in \mathbb{R}$, and $w \colon \left[a,b\right] \to \mathbb{R}$. Note that the first, second, and third equations in \eqref{eqn_corrp_ode_tpbvp} are the ODEs, while the fourth, fifth, and sixth equations constitute the boundary conditions. The first, second, and fourth equations in \eqref{eqn_corrp_ode_tpbvp} ensure that the solution lies on $\mathcal{C}$ (i.e. satisfies \eqref{eqn_ode_tpbvp}), while the third, fifth, and sixth equations in \eqref{eqn_corrp_ode_tpbvp} enforce the orthogonality constraint \eqref{eqn_ortho_constraint}. The initial solution guess to solve \eqref{eqn_corrp_ode_tpbvp} is the final corrector $\left( \mathbf{y}_{k+1}^{\mathrm c},\ACP_{k+1}^{\mathrm c} \right)$ from the previous step and $w(s) = 0$, $s \in \left[a,b\right]$. 

Note that the ODE TPBVP \eqref{eqn_corrp_ode_tpbvp} can be solved numerically via the \mcode{MATLAB} routines \mcode{sbvp} or \mcode{bvptwp}, which offers 4 algorithms: \mcode{twpbvp_m}, \mcode{twpbvpc_m}, \mcode{twpbvp_l}, and \mcode{twpbvpc_l}. Since $\mathbf{y}_1^{\mathrm c}$ and ${ \mathbf{v}}_j$ are usually only known at a discrete set of points in $\left[a,b\right]$, the values of these functions at the other points in $\left[a,b\right]$ must be obtained through interpolation in order to numerically solve \eqref{eqn_corrp_ode_tpbvp}. The \mcode{MATLAB} routine \mcode{interp1} performs linear, cubic, pchip, makima, and spline interpolation and may be utilized to interpolate $\mathbf{y}_1^{\mathrm c}$ and ${\mathbf{v}}_j$ while solving \eqref{eqn_corrp_ode_tpbvp}.

Because the numerical solvers usually converge faster when provided Jacobians of the ODE velocity function and of the two-point boundary condition function, these are computed below.
Let 
\begin{equation}
\mathbf{x} = \begin{bmatrix} \mathbf{y}_{j+1}  \\ \ACP_{j+1}  \\ w \end{bmatrix}.
\end{equation}
The ODE velocity function in \eqref{eqn_corrp_ode_tpbvp} is
\begin{equation}
\mathbf{H}^{\mathrm p} \left(s, \mathbf{x}(s) \right) = \mathbf{H}^{\mathrm p} \left(s,\mathbf{y}_{j+1}(s),\ACP_{j+1},w(s) \right) = \begin{bmatrix} \mathbf{F} \left(s,\mathbf{y}_{j+1}(s),\ACP_{j+1} \right) \\ 0 \\ {\mathbf{v}}_j^\mathsf{T}(s) \left[\mathbf{y}_{j+1}(s) - \mathbf{y}_1^{\mathrm c}(s)  \right]  \end{bmatrix}.
\end{equation}
The Jacobian of the ODE velocity function $\mathbf{H}^{\mathrm p}$ with respect to $\mathbf{x}$ is
\begin{equation}
\begin{split}
\mathbf{H}_\mathbf{x}^{\mathrm p} \left(s,\mathbf{x}(s) \right) &=
\mathbf{H}_\mathbf{x}^{\mathrm p} \left(s,\mathbf{y}_{j+1}(s),\ACP_{j+1},w(s) \right) \\ &= \begin{bmatrix} \mathbf{F}_{\mathbf{y}} \left(s,\mathbf{y}_{j+1}(s),\ACP_{j+1} \right) & \mathbf{F}_{\ACP} \left(s,\mathbf{y}_{j+1}(s),\ACP_{j+1} \right) & \mathbf{0}_{n \times 1} \\ \mathbf{0}_{1 \times n} & 0 & 0 \\ {\mathbf{v}}_j^\mathsf{T}(s) & 0 & 0 \end{bmatrix}.
\end{split}
\end{equation}
The two-point boundary condition in \eqref{eqn_corrp_ode_tpbvp} is
\begin{equation}
\mathbf{K}^{\mathrm p} \left(\mathbf{x}(a),\mathbf{x}(b) \right) = \mathbf{0}_{\left(n+2 \right) \times 1},
\end{equation}
where $\mathbf{K}^{\mathrm p}$ is the two-point boundary condition function
\begin{equation}
\mathbf{K}^{\mathrm p} \left(\mathbf{x}(a),\mathbf{x}(b) \right) = \begin{bmatrix} \mathbf{G}\left( \mathbf{y}_{j+1} (a),\mathbf{y}_{j+1}(b),\ACP_{j+1} \right) \\ w(a) \\  w(b) + \tau_j \left[\ACP_{j+1}-\ACP_1^{\mathrm c}  \right] \end{bmatrix}.
\end{equation}
The Jacobians of the two-point boundary condition function $\mathbf{K}^{\mathrm p}$ with respect to $\mathbf{x}(a)$ and $\mathbf{x}(b)$ are
\begin{equation} \label{eq_K_p_xa}
\mathbf{K}_{\mathbf{x}(a)}^{\mathrm p} \left(\mathbf{x}(a),\mathbf{x}(b) \right) = \begin{bmatrix} \mathbf{G}_{\mathbf{y}(a)}\left( \mathbf{y}_{j+1} (a),\mathbf{y}_{j+1}(b),\ACP_{j+1} \right) & \mathbf{G}_{\ACP}\left( \mathbf{y}_{j+1} (a),\mathbf{y}_{j+1}(b),\ACP_{j+1} \right) & \mathbf{0}_{n \times 1} \\ \mathbf{0}_{1 \times n} & 0 & 1 \\ \mathbf{0}_{1 \times n} & \tau_j  & 0 \end{bmatrix}
\end{equation}
and
\begin{equation} \label{eq_K_p_xb}
\mathbf{K}_{\mathbf{x}(b)}^{\mathrm p} \left(\mathbf{x}(a),\mathbf{x}(b) \right) = \begin{bmatrix} \mathbf{G}_{\mathbf{y}(b)}\left( \mathbf{y}_{j+1} (a),\mathbf{y}_{j+1}(b),\ACP_{j+1} \right) & \mathbf{G}_{\ACP}\left( \mathbf{y}_{j+1} (a),\mathbf{y}_{j+1}(b),\ACP_{j+1} \right) & \mathbf{0}_{n \times 1} \\ \mathbf{0}_{1 \times n} & 0 & 0 \\ \mathbf{0}_{1 \times n} & \tau_j  & 1\end{bmatrix}.
\end{equation}
Special care must be taken when implementing the Jacobians \eqref{eq_K_p_xa} and \eqref{eq_K_p_xb}. Since the unknown constant $\ACP_{j+1}$ appears as the second to last element of both $\mathbf{x}(a)$ and $\mathbf{x}(b)$, $\ACP_{j+1}$ from only one of $\mathbf{x}(a)$ and $\mathbf{x}(b)$ is actually used to construct each term in $\mathbf{K}^{\mathrm p}$ involving $\ACP_{j+1}$. The middle column of \eqref{eq_K_p_xa} is actually the derivative of $\mathbf{K}^{\mathrm p}$ with respect to the $\ACP_{j+1}$ in  $\mathbf{x}(a)$, while the middle column of \eqref{eq_K_p_xb} is actually the derivative of $\mathbf{K}^{\mathrm p}$ with respect to the $\ACP_{j+1}$ in  $\mathbf{x}(b)$. Thus, the middle columns in \eqref{eq_K_p_xa} and \eqref{eq_K_p_xb} corresponding to the derivative of $\mathbf{K}^{\mathrm p}$ with respect to $\ACP_{j+1}$ should not coincide in a software implementation. For example, if $\mathbf{K}^{\mathrm p}$ is constructed from the $\ACP_{j+1}$ in $\mathbf{x}(a)$, $\mathbf{K}_{\mathbf{x}(a)}^{\mathrm p}$ is as shown in \eqref{eq_K_p_xa} while the middle column of \eqref{eq_K_p_xb} corresponding to the derivative of $\mathbf{K}^{\mathrm p}$ with respect to the $\ACP_{j+1}$ in $\mathbf{x}(b)$ is all zeros. Alternatively, if $\mathbf{K}^{\mathrm p}$ is constructed from the $\ACP_{j+1}$ in $\mathbf{x}(b)$, $\mathbf{K}_{\mathbf{x}(b)}^{\mathrm p}$ is as shown in \eqref{eq_K_p_xb} while the middle column of \eqref{eq_K_p_xa} corresponding to the derivative of $\mathbf{K}^{\mathrm p}$ with respect to the $\ACP_{j+1}$ appearing in  $\mathbf{x}(a)$ is all zeros.

\subsection{Pseudocode for Predictor-Corrector Continuation}
Below is pseudocode that realizes the predictor-corrector continuation method.
\begin{algorithm}
	\caption{Predictor-Corrector Continuation for Nonlinear ODE TPBVPs.}
	\textbf{Input:} ODE velocity function $\mathbf{F} \colon \left[a,b\right] \times \mathbb{R}^n \times \mathbb{R} \to \mathbb{R}^n$, two-point boundary condition function $\mathbf{G} \colon \mathbb{R}^n \times \mathbb{R}^n \times \mathbb{R} \to \mathbb{R}^{n}$, and their Jacobians $\mathbf{F}_{\mathbf{y}} \colon \left[a,b\right] \times \mathbb{R}^n \times \mathbb{R} \to \mathbb{R}^{n \times n}$, $\mathbf{F}_{\ACP} \colon \left[a,b\right] \times \mathbb{R}^n \times \mathbb{R} \to \mathbb{R}^{n \times 1}$, $\mathbf{G}_{\mathbf{y}(a)} \colon \mathbb{R}^n \times \mathbb{R}^n \times \mathbb{R} \to \mathbb{R}^{n \times n}$, $\mathbf{G}_{\mathbf{y}(b)} \colon \mathbb{R}^n \times \mathbb{R}^n \times \mathbb{R} \to \mathbb{R}^{n \times n}$, and $\mathbf{G}_{\ACP} \colon \mathbb{R}^n \times \mathbb{R}^n \times \mathbb{R} \to \mathbb{R}^{n \times 1}$. Initial point on the solution curve $\mathcal{C}$, $\left(\mathbf{y}_1,\ACP_1 \right)$. Maximum number of points not including the initial point to be computed on $\mathcal{C}$, $J$. Initial tangent steplength, $\sigma_{\mathrm{init}}$. Minimum and maximum tangent steplengths permitted, $\sigma_{\mathrm{min}}$ and $\sigma_{\mathrm{max}}$. Tangent steplength reduction and increase scale factors, $\sigma_{\mathrm r}$ and $\sigma_{\mathrm i}$. Maximum number of Newton correction steps permitted, $K$. Maximum number of Newton correction steps for which a tangent steplength increase may occur if convergence is obtained, $k_{\mathrm{fast}}$. Newton correction convergence threshold, $\gamma$. Tangent direction at the first solution, $d$. $d$ may be $-1$ or $1$. If $d$ is $-1$ or $1$, the first tangent is scaled by $d$. \texttt{polish} is a Boolean that determines whether the Newton corrector solution is polished by solving \eqref{eqn_corrp_ode_tpbvp}.  \\
	\textbf{Output:} A solution curve $\boldsymbol{c}$ or a $\texttt{flag}$ indicating that the curve could not be traced. 
	\begin{algorithmic}[1]
		\Function{PAC\_BVP}{$\mathbf{F},\mathbf{G},\mathbf{F}_{\mathbf{y}},\mathbf{F}_{\ACP},\mathbf{G}_{\mathbf{y}(a)},\mathbf{G}_{\mathbf{y}(b)},\mathbf{G}_{\ACP},\left(\mathbf{y}_1,\ACP_1 \right),J,\sigma_{\mathrm{init}},\sigma_{\mathrm{min}},\sigma_{\mathrm{max}},\sigma_{\mathrm r},\sigma_{\mathrm i},K,k_{\mathrm{fast}},\gamma,d,\texttt{polish}$}
		
		\State $\sigma \gets \sigma_{\mathrm{init}}$ \Comment{Set the initial tangent steplength.}
		\State $\boldsymbol{c}(1) \gets \left(\mathbf{y}_{1} ,\ACP_{1} \right)$ \Comment{Store the initial solution on $\mathcal{C}$.}
		\State $\left(\mathbf{v}_0 ,\tau_0 \right) \gets \left(\mathbf{0}_{n \times 1} ,1 \right)$ \Comment{Select an initial unit tangent. This choice forces $\tau_1=1$.}
		\For {$j = 1$ to $J$}  \Comment{Trace the solution curve $\mathcal{C}$.}
		\State Obtain a tangent $\left(\mathbf{v}_j ,\tau_j \right)$ to $\mathcal{C}$ at $\left(\mathbf{y}_j ,\ACP_j \right)$ by solving \eqref{eqn_tan_ode_tpbvp} starting from $\left(\mathbf{v}_{j-1} ,\tau_{j-1} \right)$.
		\State $\kappa \gets \left \Vert \left(\mathbf{v}_j ,\tau_j \right) \right \Vert$
		\If {$j==1 $} \Comment{Choose the direction of the tangent at the initial solution, based on $d$.}
		\State $\kappa \gets \sgn\left(d\right) \kappa$
		\EndIf
		\State $\left(\mathbf{v}_j ,\tau_j \right) \gets \frac{1}{\kappa} \left(\mathbf{v}_j ,\tau_j \right)$ \Comment{Normalize the tangent.}
		\State \texttt{reject} $\gets$ \textbf{TRUE}
		\While{\texttt{reject}}
		\State $\left( \mathbf{y}_1^{\mathrm c},\ACP_1^{\mathrm c} \right) \gets  \left(\mathbf{y}_j ,\ACP_j \right) + \sigma \left(\mathbf{v}_j ,\tau_j \right)$ \Comment{Take a tangent step of length $\sigma$.}
		\For{$k=1$ to $K$} \Comment{Newton correction counter.}
		\State Obtain a Newton correction $\left( \delta \mathbf{y}_k^{\mathrm c}, \delta \ACP_k^{\mathrm c} \right)$ to $\left( \mathbf{y}_k^{\mathrm c},\ACP_k^{\mathrm c} \right)$ by solving \eqref{eqn_corrn_ode_tpbvp}.
		\State $\left( \mathbf{y}_{k+1}^{\mathrm c},\ACP_{k+1}^{\mathrm c} \right) \gets \left( \mathbf{y}_k^{\mathrm c},\ACP_k^{\mathrm c} \right)+\left( \delta \mathbf{y}_k^{\mathrm c}, \delta \ACP_k^{\mathrm c} \right)$ \Comment{Construct the Newton corrector}.
		\If {$\frac { \left \Vert \left( \delta \mathbf{y}_k^{\mathrm c} , \delta \ACP_k^{\mathrm c} \right) \right \Vert } { \left \Vert \left( \mathbf{y}_1^{\mathrm c},\ACP_1^{\mathrm c} \right) \right \Vert }< \gamma$} \Comment{Test for convergence to $\mathcal{C}$.}
		\State \texttt{reject} $\gets$ \textbf{FALSE}	
		\If {\texttt{polish}}
		\State \parbox[t]{\dimexpr\linewidth-\algorithmicindent}{Obtain the next solution $\left(\mathbf{y}_{j+1} ,\ACP_{j+1} \right)$ on $\mathcal{C}$ by solving \eqref{eqn_corrp_ode_tpbvp} starting from $\left( \mathbf{y}_{k+1}^{\mathrm c},\ACP_{k+1}^{\mathrm c} \right)$.\strut}
		\Else
		\State $\left(\mathbf{y}_{j+1} ,\ACP_{j+1} \right) \gets \left( \mathbf{y}_{k+1}^{\mathrm c},\ACP_{k+1}^{\mathrm c} \right)$ \Comment{Accept the Newton corrector solution.}
		\EndIf 	
		\State $\boldsymbol{c}(j+1) \gets \left(\mathbf{y}_{j+1} ,\ACP_{j+1} \right)$ \Comment{Store the new solution on $\mathcal{C}$.}
		\If{$k \le k_{\mathrm{fast}}$} \Comment{Test for rapid Newton convergence.}
		\State $\sigma \gets \min \left\{ \sigma_{\mathrm i} \sigma, \sigma_{\mathrm{max}} \right\}$ \Comment{Rapid Newton convergence, so increase the tangent \hphantom{blah blah blah blah blah blah blah blah blah blah blah } steplength.}
		\EndIf
		\State \textbf{break} \Comment{Break out of the for loop since convergence to $\mathcal{C}$ has been achieved.}
		\EndIf
		\EndFor
		\If{\texttt{reject}}
		\State $\sigma \gets \sigma_{\mathrm r} \sigma$ \Comment{Too many Newton steps taken, so reduce the tangent steplength.}
		\If{$\sigma < \sigma_{\mathrm{min}}$}
		\State \textbf{print} ``Unable to trace $\mathcal{C}$ because the tangent steplength is too small: $\sigma < \sigma_{\mathrm{min}}$.''
		\State \Return \texttt{flag}
		\EndIf
		\EndIf
		\EndWhile
		\EndFor
		\State \Return $\boldsymbol{c}$
		\EndFunction
	\end{algorithmic}
\end{algorithm}

\rem{
\begin{algorithm} \addtocounter{algorithm}{-1} 
	\caption{Predictor-Corrector Continuation for Nonlinear ODE TPBVPs. Part 2.}
	\begin{algorithmic}[1]
		\algrestore{bkbreak}
		\If{\texttt{reject}}
		\State $\sigma \gets \sigma_{\mathrm r} \sigma$ \Comment{Too many Newton steps taken, so reduce the tangent steplength.}
		\If{$\sigma < \sigma_{\mathrm{min}}$}
		\State \textbf{print} ``Unable to trace $\mathcal{C}$ because the tangent steplength is too small: $\sigma < \sigma_{\mathrm{min}}$.''
		\State \Return \texttt{flag}
		\EndIf
		\EndIf
		\EndWhile
		\EndFor
		\State \Return $\boldsymbol{c}$
		\EndFunction
	\end{algorithmic}
\end{algorithm}
}

\section{Sweep Predictor-Corrector Continuation Method for Solving an ODE TPBVP} \label{app_sweep_predictor_corrector}

\subsection{Introduction}
In this appendix, an alternative predictor-corrector continuation method is presented that exploits a monotonic continuation ODE TPBVP solver, such as \mcode{bvptwp}'s \mcode{acdc} or \mcode{acdcc}, to monotonically increase (i.e. sweep) the tangent steplength $\sigma$ from $0$ up until a maximum threshold $\sigma_{\mathrm{max}}$ is reached or until the next turning point is reached.

\subsection{Construct the Tangent}
Given a solution $\left(\mathbf{y}_j ,\ACP_j \right)$ to \eqref{eqn_ode_tpbvp}, we seek to construct a unit tangent $\left(\mathbf{v}_j ,\tau_j \right)$ to the solution curve $\mathcal{C}$ at $\left(\mathbf{y}_j ,\ACP_j \right)$. Recall the arclength constraint
\begin{equation} \label{eqn_arclength_constraint2}
\left \Vert \left(\mathbf{v}_j ,\tau_j \right) \right \Vert^2 = 
\left<\left(\mathbf{v}_j ,\tau_j \right) , \left(\mathbf{v}_j , \tau_j \right) \right>=
\int_a^b \mathbf{v}_j^\mathsf{T}(s)  \mathbf{v}_j(s)  \mathrm{d}s + \tau_j^2=1.
\end{equation}
The linearization (i.e. Fr\'echet derivative) of the ODE TPBVP \eqref{eqn_ode_tpbvp} about the solution $\left(\mathbf{y}_j ,\ACP_j \right)$, in conjunction with the arclength constraint \eqref{eqn_arclength_constraint2}, gives the nonlinear ODE TPBVP:
\begin{equation} \label{eqn_tan_nl_ode_tpbvp}
\begin{split}
\dd{}{s} \mathbf{v}_j(s) &= \mathbf{F}_\mathbf{y} \left(s,\mathbf{y}_j(s),\ACP_j \right) \mathbf{v}_j(s) \\
&\hphantom{=} + \mathbf{F}_{\ACP} \left(s,\mathbf{y}_j(s),\ACP_j \right) \tau_j,  \\
\dd{}{s} \tau_j &= 0, \\
\dd{}{s} w(s) &=  {\mathbf{v}}_j^\mathsf{T}(s) \mathbf{v}_j(s),  \\
\mathbf{G}_{\mathbf{y}(a)} \left(\mathbf{y}_j(a),\mathbf{y}_j(b),\ACP_j \right) \mathbf{v}_j(a)+\mathbf{G}_{\mathbf{y}(b)} \left(\mathbf{y}_j(a),\mathbf{y}_j(b),\ACP_j \right) \mathbf{v}_j(b)\\+\mathbf{G}_{\ACP} \left(\mathbf{y}_j(a),\mathbf{y}_j(b),\ACP_j \right) \tau_j &= \mathbf{0}_{n \times 1}, \\
w(a) &= 0, \\
w(b) +  \tau_j^2 - 1 &= 0,
\end{split}
\end{equation}
which must be solved for $\mathbf{v}_j \colon \left[a,b\right] \to \mathbb{R}^n$, $\tau_j \in \mathbb{R}$, and $w \colon \left[a,b\right] \to \mathbb{R}$ and where $\left(\mathbf{v}_j ,\tau_j \right)$ is a unit tangent to $\mathcal{C}$ at $\left(\mathbf{y}_j ,\ACP_j \right)$. Note that the first, second, and third equations in \eqref{eqn_tan_nl_ode_tpbvp} are the ODEs, while the fourth, fifth, and sixth equations constitute the boundary conditions. The first, second, and fourth equations in \eqref{eqn_tan_nl_ode_tpbvp} are the linearization (i.e. Fr\'echet derivative) of \eqref{eqn_ode_tpbvp} about the solution $\left(\mathbf{y}_j ,\ACP_j \right)$ and ensure that a tangent is produced, while the third, fifth, and sixth equations in \eqref{eqn_tan_nl_ode_tpbvp} enforce the arclength constraint \eqref{eqn_arclength_constraint2} ensuring that the tangent is of unit length. The initial solution guess to solve \eqref{eqn_tan_nl_ode_tpbvp} is $\left(\mathbf{v}_j ,\tau_j \right) = \left(\mathbf{0}_{n \times 1} ,1 \right)$ and $w(s) = 0$, $s \in \left[a,b\right]$.

Note that the ODE TPBVP \eqref{eqn_tan_nl_ode_tpbvp} can be solved numerically via the \mcode{MATLAB} routines \mcode{sbvp} or \mcode{bvptwp}, which offers 4 algorithms: \mcode{twpbvp_m}, \mcode{twpbvpc_m}, \mcode{twpbvp_l}, and \mcode{twpbvpc_l}. Since $\mathbf{y}_j$ is usually only known at a discrete set of points in $\left[a,b\right]$, the values of this function at the other points in $\left[a,b\right]$ must be obtained through interpolation in order to numerically solve \eqref{eqn_tan_nl_ode_tpbvp}. The \mcode{MATLAB} routine \mcode{interp1} performs linear, cubic, pchip, makima, and spline interpolation and may be utilized to interpolate $\mathbf{y}_j$ while solving \eqref{eqn_tan_nl_ode_tpbvp}.

Because the numerical solvers usually converge faster when provided Jacobians of the ODE velocity function and of the two-point boundary condition function, these are computed below.
Let 
\begin{equation}
\mathbf{x} = \begin{bmatrix} \mathbf{v}_{j}  \\ \tau_j  \\ w \end{bmatrix}.
\end{equation}
The ODE velocity function in \eqref{eqn_tan_nl_ode_tpbvp} is
\begin{equation}
\mathbf{H}^{\mathrm t} \left(s, \mathbf{x}(s) \right) = \mathbf{H}^{\mathrm t} \left(s,\mathbf{v}_j(s),\tau_j,w(s) \right) = \begin{bmatrix} \mathbf{F}_\mathbf{y} \left(s,\mathbf{y}_j(s),\ACP_j \right) \mathbf{v}_j(s) + \mathbf{F}_{\ACP} \left(s,\mathbf{y}_j(s),\ACP_j \right) \tau_j \\ 0 \\ {\mathbf{v}}_j^\mathsf{T}(s) \mathbf{v}_j(s) \end{bmatrix}.
\end{equation}
The Jacobian of the ODE velocity function $\mathbf{H}^{\mathrm t}$ with respect to $\mathbf{x}$ is
\begin{equation}
\begin{split}
\mathbf{H}_\mathbf{x}^{\mathrm t} \left(s,\mathbf{x}(s) \right) &=
\mathbf{H}_\mathbf{x}^{\mathrm t} \left(s,\mathbf{v}_j(s),\tau_j,w(s) \right) \\ &= \begin{bmatrix} \mathbf{F}_{\mathbf{y}} \left(s,\mathbf{y}_j(s),\ACP_j \right) & \mathbf{F}_{\ACP} \left(s,\mathbf{y}_j(s),\ACP_j \right) & \mathbf{0}_{n \times 1} \\ \mathbf{0}_{1 \times n} & 0 & 0 \\ 2 {\mathbf{v}}_j^\mathsf{T}(s) & 0 & 0 \end{bmatrix}.
\end{split}
\end{equation}
The two-point boundary condition in \eqref{eqn_tan_nl_ode_tpbvp} is
\begin{equation}
\mathbf{K}^{\mathrm t} \left(\mathbf{x}(a),\mathbf{x}(b) \right) = \mathbf{0}_{\left(n+2 \right) \times 1},
\end{equation}
where $\mathbf{K}^{\mathrm t}$ is the two-point boundary condition function
\begin{equation}
\begin{split}
\mathbf{K}^{\mathrm t}& \left(\mathbf{x}(a),\mathbf{x}(b) \right) = \\
& \begin{bmatrix} \mathbf{G}_{\mathbf{y}(a)} \left(\mathbf{y}_j(a),\mathbf{y}_j(b),\ACP_j \right) \mathbf{v}_j(a)+\mathbf{G}_{\mathbf{y}(b)} \left(\mathbf{y}_j(a),\mathbf{y}_j(b),\ACP_j \right) \mathbf{v}_j(b) +\mathbf{G}_{\ACP} \left(\mathbf{y}_j(a),\mathbf{y}_j(b),\ACP_j \right) \tau_j \\ w(a) \\  w(b) + \tau_j^2 - 1 \end{bmatrix}.
\end{split}
\end{equation}
The Jacobians of the two-point boundary condition function $\mathbf{K}^{\mathrm t}$ with respect to $\mathbf{x}(a)$ and $\mathbf{x}(b)$ are
\begin{equation} \label{eqs_K_t_xa}
\mathbf{K}_{\mathbf{x}(a)}^{\mathrm t} \left(\mathbf{x}(a),\mathbf{x}(b) \right) = \begin{bmatrix} \mathbf{G}_{\mathbf{y}(a)} \left(\mathbf{y}_j(a),\mathbf{y}_j(b),\ACP_j \right) & \mathbf{G}_{\ACP} \left(\mathbf{y}_j(a),\mathbf{y}_j(b),\ACP_j \right) & \mathbf{0}_{n \times 1} \\ \mathbf{0}_{1 \times n} & 0 & 1 \\ \mathbf{0}_{1 \times n} & 2 \tau_j  & 0 \end{bmatrix}
\end{equation}
and
\begin{equation} \label{eqs_K_t_xb}
\mathbf{K}_{\mathbf{x}(b)}^{\mathrm t} \left(\mathbf{x}(a),\mathbf{x}(b) \right) = \begin{bmatrix} \mathbf{G}_{\mathbf{y}(b)} \left(\mathbf{y}_j(a),\mathbf{y}_j(b),\ACP_j \right) & \mathbf{G}_{\ACP} \left(\mathbf{y}_j(a),\mathbf{y}_j(b),\ACP_j \right) & \mathbf{0}_{n \times 1} \\ \mathbf{0}_{1 \times n} & 0 & 0 \\ \mathbf{0}_{1 \times n} & 2 \tau_j  & 1\end{bmatrix}.
\end{equation}
Special care must be taken when implementing the Jacobians \eqref{eqs_K_t_xa} and \eqref{eqs_K_t_xb}. Since the unknown constant $\tau_j$ appears as the second to last element of both $\mathbf{x}(a)$ and $\mathbf{x}(b)$, $\tau_j$ from only one of $\mathbf{x}(a)$ and $\mathbf{x}(b)$ is actually used to construct each term in $\mathbf{K}^{\mathrm t}$ involving $\tau_j$. The middle column of \eqref{eqs_K_t_xa} is actually the derivative of $\mathbf{K}^{\mathrm t}$ with respect to the $\tau_j$ in  $\mathbf{x}(a)$, while the middle column of \eqref{eqs_K_t_xb} is actually the derivative of $\mathbf{K}^{\mathrm t}$ with respect to the $\tau_j$ in  $\mathbf{x}(b)$. Thus, the middle columns in \eqref{eqs_K_t_xa} and \eqref{eqs_K_t_xb} corresponding to the derivative of $\mathbf{K}^{\mathrm t}$ with respect to $\tau_j$ should not coincide in a software implementation. For example, if $\mathbf{K}^{\mathrm t}$ is constructed from the $\tau_j$ in $\mathbf{x}(a)$, $\mathbf{K}_{\mathbf{x}(a)}^{\mathrm t}$ is as shown in \eqref{eqs_K_t_xa} while the middle column of \eqref{eqs_K_t_xb} corresponding to the derivative of $\mathbf{K}^{\mathrm t}$ with respect to the $\tau_j$ in $\mathbf{x}(b)$ is all zeros. Alternatively, if $\mathbf{K}^{\mathrm t}$ is constructed from the $\tau_j$ in $\mathbf{x}(b)$, $\mathbf{K}_{\mathbf{x}(b)}^{\mathrm t}$ is as shown in \eqref{eqs_K_t_xb} while the middle column of \eqref{eqs_K_t_xa} corresponding to the derivative of $\mathbf{K}^{\mathrm t}$ with respect to the $\tau_j$ appearing in  $\mathbf{x}(a)$ is all zeros.

\subsection{Determine the Tangent Direction}
The unit tangent $\left(\mathbf{v}_j ,\tau_j \right)$ at $\left(\mathbf{y}_j ,\ACP_j \right)$ obtained by solving \eqref{eqn_tan_nl_ode_tpbvp} must be scaled so that the sweep predictor-corrector continuation method does not reverse direction. As shown in \cite{kitzhofer2009pathfollowing}, the correct direction for the unit tangent is obtained via:
\begin{equation} \label{eqn_dir_tan_pred} 
\left(\mathbf{v}_j ,\tau_j \right) \gets \sgn\left( \kappa \right) \left(\mathbf{v}_j ,\tau_j \right),
\end{equation}
where $\kappa$ is the inner product of the previous and current unit tangents:
\begin{equation} \label{eqn_ip_tan_pred_kappa} 
\kappa = 
\left<\left(\mathbf{v}_{j-1} ,\tau_{j-1} \right) , \left(\mathbf{v}_j , \tau_j \right) \right>  =
\int_a^b \mathbf{v}_{j-1}^\mathsf{T}(s)  \mathbf{v}_j(s)  \mathrm{d}s + \tau_{j-1} \tau_j.
\end{equation}
The integration operator to construct the inner product $\kappa$ in \eqref{eqn_ip_tan_pred_kappa} can be realized via the \mcode{MATLAB} routine \mcode{trapz}. With the sign direction selected by \eqref{eqn_dir_tan_pred}, the inner product of the previous and current unit tangents is positive:
\begin{equation} \label{eqn_ip2_tan_pred_kappa} 
\left<\left(\mathbf{v}_{j-1} ,\tau_{j-1} \right) , \left(\mathbf{v}_j , \tau_j \right) \right>  =
\int_a^b \mathbf{v}_{j-1}^\mathsf{T}(s)  \mathbf{v}_j(s)  \mathrm{d}s + \tau_{j-1} \tau_j > 0.
\end{equation}

\subsection{Sweep along the Tangent}
By monotonically increasing (or sweeping) the tangent steplength $\sigma$ from $0$, the current solution $\left(\mathbf{y}_j ,\ACP_j \right)$ and its unit tangent $\left(\mathbf{v}_j ,\tau_j \right)$ can be used to find the next solution $\left(\mathbf{y}_{j+1} ,\ACP_{j+1} \right)$ that solves \eqref{eqn_ode_tpbvp} while satisfying the orthogonality constraint:
\begin{equation} \label{eqn_ortho_constraint2}
\begin{split}
&\left<\left({\mathbf{v}}_j , \tau_j \right) , \left(\mathbf{y}_{j+1}, \ACP_{j+1} \right) - \left( \left( \mathbf{y}_j,\ACP_j \right) + \sigma \left( {\mathbf{v}}_j, \tau_j \right) \right) \right>  
=\left<\left({\mathbf{v}}_j , \tau_j \right) , \left(\mathbf{y}_{j+1} - \left( \mathbf{y}_j + \sigma {\mathbf{v}}_j  \right)  , \ACP_{j+1}-\left( \ACP_j + \sigma  \tau_j  \right) \right) \right> \\
&= \int_a^b {\mathbf{v}}_j^\mathsf{T}(s) \left[\mathbf{y}_{j+1}(s) - \left( \mathbf{y}_j(s) + \sigma {\mathbf{v}}_j(s)  \right)   \right]  \mathrm{d}s +   \tau_j \left[\ACP_{j+1}-\left( \ACP_j + \sigma  \tau_j  \right) \right]=0.
\end{split}
\end{equation}
This yields the ODE TPBVP:
\begin{equation} \label{eqn_sweep_ode_tpbvp}
\begin{split}
\dd{}{s} \mathbf{y}_{j+1}(s) &= \mathbf{F} \left(s,\mathbf{y}_{j+1}(s),\ACP_{j+1} \right), \\
\dd{}{s} \ACP_{j+1} &= 0, \\
\dd{}{s} w(s) &= {\mathbf{v}}_j^\mathsf{T}(s) \left[\mathbf{y}_{j+1}(s) - \left( \mathbf{y}_j(s) + \sigma {\mathbf{v}}_j(s)  \right)  \right],  \\
\mathbf{G}\left( \mathbf{y}_{j+1} (a),\mathbf{y}_{j+1}(b),\ACP_{j+1} \right) &= \mathbf{0}_{n \times 1}, \\
w(a) &= 0, \\
w(b) + \tau_j \left[\ACP_{j+1}-\left( \ACP_j + \sigma \tau_j  \right)  \right] &= 0,
\end{split}
\end{equation}
which must be solved for $\mathbf{y}_{j+1} \colon \left[a,b\right] \to \mathbb{R}^n$,  $\ACP_{j+1} \in \mathbb{R}$, and $w \colon \left[a,b\right] \to \mathbb{R}$ by monotonically increasing (or sweeping) $\sigma$. Note that the first, second, and third equations in \eqref{eqn_sweep_ode_tpbvp} are the ODEs, while the fourth, fifth, and sixth equations constitute the boundary conditions. The first, second, and fourth equations in \eqref{eqn_sweep_ode_tpbvp} ensure that the solution lies on $\mathcal{C}$ (i.e. satisfies \eqref{eqn_ode_tpbvp}), while the third, fifth, and sixth equations in \eqref{eqn_sweep_ode_tpbvp} enforce the orthogonality constraint \eqref{eqn_ortho_constraint2}. The initial solution guess to solve \eqref{eqn_sweep_ode_tpbvp} is the current solution $\left( \mathbf{y}_j,\ACP_j \right)$ and $w(s) = 0$, $s \in \left[a,b\right]$. $\sigma$ starts at $0$, since the initial solution guess for $\left(\mathbf{y}_{j+1} ,\ACP_{j+1} \right)$ is $\left( \mathbf{y}_j,\ACP_j \right)$, and increases monotonically until the maximum threshold $\sigma_{\mathrm{max}}$ is reached or until the ODE TPBVP solver halts (due to reaching a turning point). 

Note that the ODE TPBVP \eqref{eqn_sweep_ode_tpbvp} can be solved numerically via the \mcode{MATLAB} routine \mcode{bvptwp}, which offers 2 continuation algorithms: \mcode{acdc} and \mcode{acdcc}. The continuation algorithms \mcode{acdc} and \mcode{acdcc} assume that the continuation parameter (in this case $\sigma$) is monotonically increasing or decreasing, so that they will halt at a turning point in the continuation parameter. Since $\mathbf{y}_j$ and ${ \mathbf{v}}_j$ are usually only known at a discrete set of points in $\left[a,b\right]$, the values of these functions at the other points in $\left[a,b\right]$ must be obtained through interpolation in order to numerically solve \eqref{eqn_sweep_ode_tpbvp}. The \mcode{MATLAB} routine \mcode{interp1} performs linear, cubic, pchip, makima, and spline interpolation and may be utilized to interpolate $\mathbf{y}_j$ and ${\mathbf{v}}_j$ while solving \eqref{eqn_sweep_ode_tpbvp}.

Because the numerical solvers usually converge faster when provided Jacobians of the ODE velocity function and of the two-point boundary condition function, these are computed below.
Let 
\begin{equation}
\mathbf{x} = \begin{bmatrix} \mathbf{y}_{j+1}  \\ \ACP_{j+1}  \\ w \end{bmatrix}.
\end{equation}
The ODE velocity function in \eqref{eqn_sweep_ode_tpbvp} is
\begin{equation}
\mathbf{H}^{\mathrm q} \left(s, \mathbf{x}(s),\sigma \right) = \mathbf{H}^{\mathrm q} \left(s,\mathbf{y}_{j+1}(s),\ACP_{j+1},w(s),\sigma \right) = \begin{bmatrix} \mathbf{F} \left(s,\mathbf{y}_{j+1}(s),\ACP_{j+1} \right) \\ 0 \\ {\mathbf{v}}_j^\mathsf{T}(s) \left[\mathbf{y}_{j+1}(s) - \left( \mathbf{y}_j(s) + \sigma { \mathbf{v}}_j(s)  \right)   \right]  \end{bmatrix}.
\end{equation}
The Jacobian of the ODE velocity function $\mathbf{H}^{\mathrm q}$ with respect to $\mathbf{x}$ is
\begin{equation}
\begin{split}
\mathbf{H}_\mathbf{x}^{\mathrm q} \left(s,\mathbf{x}(s),\sigma \right) &=
\mathbf{H}_\mathbf{x}^{\mathrm q} \left(s,\mathbf{y}_{j+1}(s),\ACP_{j+1},w(s),\sigma \right) \\ &= \begin{bmatrix} \mathbf{F}_{\mathbf{y}} \left(s,\mathbf{y}_{j+1}(s),\ACP_{j+1} \right) & \mathbf{F}_{\ACP} \left(s,\mathbf{y}_{j+1}(s),\ACP_{j+1} \right) & \mathbf{0}_{n \times 1} \\ \mathbf{0}_{1 \times n} & 0 & 0 \\ {\mathbf{v}}_j^\mathsf{T}(s) & 0 & 0 \end{bmatrix}.
\end{split}
\end{equation}
The two-point boundary condition in \eqref{eqn_sweep_ode_tpbvp} is
\begin{equation}
\mathbf{K}^{\mathrm q} \left(\mathbf{x}(a),\mathbf{x}(b),\sigma \right) = \mathbf{0}_{\left(n+2 \right) \times 1},
\end{equation}
where $\mathbf{K}^{\mathrm q}$ is the two-point boundary condition function
\begin{equation}
\mathbf{K}^{\mathrm q} \left(\mathbf{x}(a),\mathbf{x}(b),\sigma \right) = \begin{bmatrix} \mathbf{G}\left( \mathbf{y}_{j+1} (a),\mathbf{y}_{j+1}(b),\ACP_{j+1} \right) \\ w(a) \\  w(b) + \tau_j \left[\ACP_{j+1}-\left( \ACP_j + \sigma \tau_j  \right)  \right] \end{bmatrix}.
\end{equation}
The Jacobians of the two-point boundary condition function $\mathbf{K}^{\mathrm q}$ with respect to $\mathbf{x}(a)$ and $\mathbf{x}(b)$ are
\begin{equation} \label{eq_K_q_xa}
\mathbf{K}_{\mathbf{x}(a)}^{\mathrm q} \left(\mathbf{x}(a),\mathbf{x}(b),\sigma \right) = \begin{bmatrix} \mathbf{G}_{\mathbf{y}(a)}\left( \mathbf{y}_{j+1} (a),\mathbf{y}_{j+1}(b),\ACP_{j+1} \right) & \mathbf{G}_{\ACP}\left( \mathbf{y}_{j+1} (a),\mathbf{y}_{j+1}(b),\ACP_{j+1} \right) & \mathbf{0}_{n \times 1} \\ \mathbf{0}_{1 \times n} & 0 & 1 \\ \mathbf{0}_{1 \times n} & \tau_j  & 0 \end{bmatrix}
\end{equation}
and
\begin{equation} \label{eq_K_q_xb}
\mathbf{K}_{\mathbf{x}(b)}^{\mathrm q} \left(\mathbf{x}(a),\mathbf{x}(b),\sigma \right) = \begin{bmatrix} \mathbf{G}_{\mathbf{y}(b)}\left( \mathbf{y}_{j+1} (a),\mathbf{y}_{j+1}(b),\ACP_{j+1} \right) & \mathbf{G}_{\ACP}\left( \mathbf{y}_{j+1} (a),\mathbf{y}_{j+1}(b),\ACP_{j+1} \right) & \mathbf{0}_{n \times 1} \\ \mathbf{0}_{1 \times n} & 0 & 0 \\ \mathbf{0}_{1 \times n} & \tau_j  & 1\end{bmatrix}.
\end{equation}
Special care must be taken when implementing the Jacobians \eqref{eq_K_q_xa} and \eqref{eq_K_q_xb}. Since the unknown constant $\ACP_{j+1}$ appears as the second to last element of both $\mathbf{x}(a)$ and $\mathbf{x}(b)$, $\ACP_{j+1}$ from only one of $\mathbf{x}(a)$ and $\mathbf{x}(b)$ is actually used to construct each term in $\mathbf{K}^{\mathrm q}$ involving $\ACP_{j+1}$. The middle column of \eqref{eq_K_q_xa} is actually the derivative of $\mathbf{K}^{\mathrm q}$ with respect to the $\ACP_{j+1}$ in $\mathbf{x}(a)$, while the middle column of \eqref{eq_K_q_xb} is actually the derivative of $\mathbf{K}^{\mathrm q}$ with respect to the $\ACP_{j+1}$ in $\mathbf{x}(b)$. Thus, the middle columns in \eqref{eq_K_q_xa} and \eqref{eq_K_q_xb} corresponding to the derivative of $\mathbf{K}^{\mathrm q}$ with respect to $\ACP_{j+1}$ should not coincide in a software implementation. For example, if $\mathbf{K}^{\mathrm q}$ is constructed from the $\ACP_{j+1}$ in $\mathbf{x}(a)$, $\mathbf{K}_{\mathbf{x}(a)}^{\mathrm q}$ is as shown in \eqref{eq_K_q_xa} while the middle column of \eqref{eq_K_q_xb} corresponding to the derivative of $\mathbf{K}^{\mathrm q}$ with respect to the $\ACP_{j+1}$ in $\mathbf{x}(b)$ is all zeros. Alternatively, if $\mathbf{K}^{\mathrm q}$ is constructed from the $\ACP_{j+1}$ in $\mathbf{x}(b)$, $\mathbf{K}_{\mathbf{x}(b)}^{\mathrm q}$ is as shown in \eqref{eq_K_q_xb} while the middle column of \eqref{eq_K_q_xa} corresponding to the derivative of $\mathbf{K}^{\mathrm q}$ with respect to the $\ACP_{j+1}$ appearing in  $\mathbf{x}(a)$ is all zeros.

\subsection{Pseudocode for Sweep Predictor-Corrector Continuation}
Below is pseudocode that realizes the sweep predictor-corrector continuation method.

\begin{algorithm}
	\caption{Sweep Predictor-Corrector Continuation for Nonlinear ODE TPBVPs.}
	\textbf{Input:} ODE velocity function $\mathbf{F} \colon \left[a,b\right] \times \mathbb{R}^n \times \mathbb{R} \to \mathbb{R}^n$, two-point boundary condition function $\mathbf{G} \colon \mathbb{R}^n \times \mathbb{R}^n \times \mathbb{R} \to \mathbb{R}^{n}$, and their Jacobians $\mathbf{F}_{\mathbf{y}} \colon \left[a,b\right] \times \mathbb{R}^n \times \mathbb{R} \to \mathbb{R}^{n \times n}$, $\mathbf{F}_{\ACP} \colon \left[a,b\right] \times \mathbb{R}^n \times \mathbb{R} \to \mathbb{R}^{n \times 1}$, $\mathbf{G}_{\mathbf{y}(a)} \colon \mathbb{R}^n \times \mathbb{R}^n \times \mathbb{R} \to \mathbb{R}^{n \times n}$, $\mathbf{G}_{\mathbf{y}(b)} \colon \mathbb{R}^n \times \mathbb{R}^n \times \mathbb{R} \to \mathbb{R}^{n \times n}$, and $\mathbf{G}_{\ACP} \colon \mathbb{R}^n \times \mathbb{R}^n \times \mathbb{R} \to \mathbb{R}^{n \times 1}$. Initial point on the solution curve $\mathcal{C}$, $\left(\mathbf{y}_1,\ACP_1 \right)$. Maximum number of points not including the initial point to be computed on $\mathcal{C}$, $J$. $\sigma_{\mathrm{max}}$ is a vector of length $J$ such that $\sigma_{\mathrm{max}}(j)$ is the maximum tangent steplength permitted to obtain solution $j+1$. Tangent direction at the first solution, $d$. $d$ may be $-2$, $-1$, $1$, or $2$. If $d$ is $-1$ or $1$, the first tangent is scaled by $d$. If $d$ is $-2$ $(2)$, the first tangent is scaled so that $\lambda$ decreases (increases) in the first step. \\
	\textbf{Output:} A solution curve $\boldsymbol{c}$. 
	\begin{algorithmic}[1]
		\Function{PAC\_s3\_BVP}{$\mathbf{F},\mathbf{G},\mathbf{F}_{\mathbf{y}},\mathbf{F}_{\ACP},\mathbf{G}_{\mathbf{y}(a)},\mathbf{G}_{\mathbf{y}(b)},\mathbf{G}_{\ACP},\left(\mathbf{y}_1,\ACP_1 \right),J,\sigma_{\mathrm{max}},d$}
		
		\State $\boldsymbol{c}(1) \gets \left(\mathbf{y}_{1} ,\ACP_{1} \right)$ \Comment{Store the initial solution on $\mathcal{C}$.}
		\For {$j = 1$ to $J$}  \Comment{Trace the solution curve $\mathcal{C}$.}
		\State Obtain a unit tangent $\left(\mathbf{v}_j ,\tau_j \right)$ to $\mathcal{C}$ at $\left(\mathbf{y}_j ,\ACP_j \right)$ by solving \eqref{eqn_tan_nl_ode_tpbvp} starting from $\left(\mathbf{0}_{n \times 1} ,1 \right)$.
		
		\If {$j==1 $} \Comment{Choose the direction of the tangent at the initial solution, based on $d$.}
		\If {$\left(d==-2 \right.$ \textbf{OR} $\left. d==2 \right)$  \textbf{AND} $\tau_1 < 0$}
		\State $d \gets -d$ \Comment{Flip the sign of $d$ to get the desired tangent direction.}
		\EndIf
		\State $\kappa \gets d$
		\Else
		\State $\kappa \gets \left<\left(\mathbf{v}_{j-1} ,\tau_{j-1} \right) , \left(\mathbf{v}_j , \tau_j \right) \right>$	 \Comment{Ensure that the traced solution does not reverse direction.}	
		\EndIf
		\State $\left(\mathbf{v}_j ,\tau_j \right) \gets \sgn\left( \kappa \right) \left(\mathbf{v}_j ,\tau_j \right)$ \Comment{Choose the correct tangent direction.}
		\State \parbox[t]{\dimexpr\linewidth-\algorithmicindent}{Obtain the next solution $\left(\mathbf{y}_{j+1} ,\ACP_{j+1} \right)$ on $\mathcal{C}$ by solving \eqref{eqn_sweep_ode_tpbvp} starting from $\left( \mathbf{y}_j,\ACP_j \right)$ and monotonically increasing $\sigma$ starting from 0 and without exceeding $\sigma_{\mathrm{max}}(j)$.\strut}
		\State $\boldsymbol{c}(j+1) \gets \left(\mathbf{y}_{j+1} ,\ACP_{j+1} \right)$ \Comment{Store the new solution on $\mathcal{C}$.}
		\EndFor
		\State \Return $\boldsymbol{c}$
		\EndFunction
		
	\end{algorithmic}
\end{algorithm}

\end{document}